**W. B. VASANTHA KANDASAMY**
**FLORENTIN SMARANDACHE**

# FUZZY COGNITIVE MAPS AND NEUTROSOPHIC COGNITIVE MAPS

# Fuzzy Cognitive Maps and Neutrosophic Cognitive Maps

**W. B. Vasantha Kandasamy**
Department of Mathematics
Indian Institute of Technology, Madras
Chennai – 600036, India
e-mail: vasantha@iitm.ac.in
web: http://mat.iitm.ac.in/~wbv

**Florentin Smarandache**
Department of Mathematics
University of Gallup, NM
New Mexico, USA
e-mail: smarand@gallup.unm.edu

# CONTENTS

SILENCE = DEATH

DEDICATED TO

ALL THE VOICES OF DISSENT

# PREFACE

In a world of chaotic alignments, traditional logic with its strict boundaries of truth and falsity has not imbued itself with the capability of reflecting the reality. Despite various attempts to reorient logic, there has remained an essential need for an alternative system that could infuse into itself a representation of the real world. Out of this need arose the system of Neutrosophy, and its connected logic, Neutrosophic Logic. Neutrosophy is a new branch of philosophy that studies the origin, nature and scope of neutralities, as well as their interactions with different ideational spectra. This was introduced by one of the authors, Florentin Smarandache. A few of the mentionable characteristics of this mode of thinking are [90-94]: It proposes new philosophical theses, principles, laws, methods, formulas and movements; it reveals that the world is full of indeterminacy; it interprets the uninterpretable; regards, from many different angles, old concepts, systems and proves that an idea which is true in a given referential system, may be false in another, and vice versa; attempts to make peace in the war of ideas, and to make war in the peaceful ideas! The main principle of neutrosophy is: Between an idea <A> and its opposite <Anti-A>, there is a continuum-power spectrum of Neutralities. This philosophy forms the basis of Neutrosophic logic.

Neutrosophic logic grew as an alternative to the existing logics and it represents a mathematical model of uncertainty, vagueness, ambiguity, imprecision, undefined, unknown, incompleteness, inconsistency, redundancy, contradiction. It can be defined as a logic in which each proposition is estimated to have the percentage of truth in a subset T, the percentage of indeterminacy in a subset I, and the percentage of falsity in a subset F, is called Neutrosophic Logic. We use a subset of truth (or indeterminacy, or falsity), instead of using a number, because in many cases, we are not able to exactly determine the percentages of truth and of falsity but to approximate them: for example a proposition is between 30-40% true. The subsets are not necessarily intervals, but any sets (discrete, continuous, open or closed or half-open/ half-closed interval, intersections or unions of the previous sets, etc.) in accordance with the given proposition. A subset may have one element only in special cases of this logic. It is imperative to mention here that the Neutrosophic logic is a further generalization of the theory of Fuzzy Logic.

In this book we study the concepts of Fuzzy Cognitive Maps (FCMs) and their Neutrosophic analogue, the Neutrosophic Cognitive Maps (NCMs).

Fuzzy Cognitive Maps are fuzzy structures that strongly resemble neural networks, and they have powerful and far-reaching consequences as a mathematical tool for modeling complex systems. Prof. Bart Kosko, the guru of fuzzy logic, introduced the Fuzzy Cognitive Maps [54] in the year 1986. It was a fuzzy extension of the cognitive map pioneered in 1976 by political scientist Robert Axelrod [5], who used it to represent knowledge as an interconnected, directed, bilevel-logic graph. Till today there are over a hundred research papers which deal with FCMs, and the tool has been used to study real-world situations as varied as stock-investment analysis to supervisory system control, and child labor to community mobilization against the AIDS epidemic.

This book has been written with two aims: First, we seek to consolidate the vast amount of research that has been done around the concepts of FCMs, and also try to give an inclusive view of the various real-world problems to which FCMs have been applied. Though there are over a hundred research papers relating to FCMs, there is no book that deals exclusively with them — and we hope this book possibly bridges that gap. Second, we introduce here (for the first time) the concept of Neutrosophic Cognitive Maps (NCMs), which are a generalization of Fuzzy Cognitive Maps. The special feature of NCMs is their ability to handle indeterminacy in relations between two concepts, which is denoted by 'I'. This new structure — the NCM is capable of giving results with greater sensitivity than the FCM. It also allows a larger liberty of intuition by allowing an expert to express not just the positive, negative and absence of impacts but also the indeterminacy of impacts.

Practically speaking, we must be aware that even in our day-to-day lives, the indeterminacy and unpredictability of life, affect us almost as much as the determined factors. It is a major handicap in mathematical modeling that we are only able to give weightages for known concepts; and most of the time we exhibit an unconcern for indeterminate relationships between concepts, thereby presenting onto ourselves a skewed view. Prof. Bart Kosko, in his book *Heaven in a Chip: Fuzzy Visions of Society and Science in the Digital Age* writes that fuzzy theory can offer more choices and blur the hard lines of power that define the politics of our age. We have, in our effort at introducing indeterminacy into Fuzzy theory, and by our construction of neutrosophic structures of modeling, only extended the liberty of choice to a greater level. We have written this book as a maiden effort to inculcate into real-world problems the concept of indeterminacy, uncertainty and inconclusiveness.

This book is divided into three chapters. In chapter one, we recall the definition of Fuzzy Cognitive Maps, suggest properties about FCM models and give illustrations. We give details about the multifarious applications of FCMs which has been studied by many authors. A new notion called Fuzzy Relational Maps (FRMs) [a particularization of the FCMs] are introduced by us — the FRMs are applicable when the nodes of the FCMs can be divided into two disjoint classes, and they are more beneficial owing to their capability of being economic, time-saving and sensitive. In the concluding sections of the first chapter we deal with illustrations of FRMs and put forth the concept of combined FRMs and linked FRMs. The second chapter introduces the notion of Neutrosophic Cognitive Maps (NCMs), and in order to introduce this concept we have introduced neutrosophic graphs, neutrosophic fields, neutrosophic matrices and neutrosophic vector spaces. We provide many illustrations and applications relating NCMs. In this chapter, we compare and contrast NCMs and FCMs. We also define NRMs and illustrate its applications to real-world problems and compare NRMs and FRMs. In the final chapter we suggest problems relating to these concepts. An almost exhaustive bibliography relating to the theory of FCMs completes the book.

Some of the varied applications of FCMs and NCMs (and alternately FRMs and NRMs) which has been explained by us, in this book, include: modeling of supervisory systems; design of hybrid models for complex systems; mobile robots and in intimate technology such as office plants; analysis of business performance assessment; formalism debate and legal rules; creating metabolic and regulatory network models; traffic and transportation problems; medical diagnostics; simulation of strategic planning process in intelligent systems; specific language impairment; web-mining inference application; child labor problem; industrial relations: between employer and employee, maximizing production and profit; decision support in intelligent intrusion detection system; hyper-knowledge representation in strategy formation; female infanticide; depression in terminally ill patients and finally, in the theory of community mobilization and women empowerment relative to the AIDS epidemic.

It is worth mentioning here, that in this book we have not used degrees of uncertainties or indeterminates, although one can easily use these concepts. We have dealt with only simple FCMs and NCMs. We also wish to mention that this book has been written only for readers who are well versed with basic graph theory and matrix theory.

The authors thank Dr. Minh Perez of the American Research Press for his constant support and encouragement towards the writing of this book. We also thank Dr.Kandasamy, who diligently proofread four rough drafts of this book, Kama who patiently drew all the cognitive maps and ensured that the formatting and page making of the book was intact and Meena who collected the existing literature.

**W.B. VASANTHA KANDASAMY**
**FLORENTIN SMARANDACHE**

Chapter One

# BASIC CONCEPTS ABOUT FUZZY COGNITIVE MAPS AND FUZZY RELATIONAL MAPS

This chapter has seven sections. In section one we recall the definition and basic properties of Fuzzy Cognitive Maps (FCMs). In section two we give properties and models of FCMs and present some of its applications to problems such as the maximum utility of a route, Socio-economic problems and Symptom-disease model. In section three we give some illustration of FCMs. Applications of FCMs is dealt with in section four and section five is concerned with the introduction of a new model [125] called the Fuzzy Relational Maps (FRMs). Fuzzy Relational Model happens to be better than that of the FCM model in several ways mainly when the analysis of the data can be treated as two disjoint entities. Thus its application to the Employee-Employer problem, the study of maximizing production in Cement Industries by giving maximum satisfaction to employees and the notion of Fuzzy relational models illustrating FRMs and combined FRMs are dealt with in section six. Section seven introduces linked FRMs.

## 1.1 Definition of Fuzzy Cognitive Maps

In this section we recall the notion of Fuzzy Cognitive Maps (FCMs), which was introduced by Bart Kosko [54] in the year 1986. We also give several of its interrelated definitions. FCMs have a major role to play mainly when the data concerned is an unsupervised one. Further this method is most simple and an effective one as it can analyse the data by directed graphs and connection matrices.

**DEFINITION 1.1.1:** *An FCM is a directed graph with concepts like policies, events etc. as nodes and causalities as edges. It represents causal relationship between concepts.*

***Example 1.1.1:*** In Tamil Nadu (a southern state in India) in the last decade several new engineering colleges have been approved and started. The resultant increase in the production of engineering graduates in these years is disproportionate with the need of engineering graduates. This has resulted in thousands of unemployed and underemployed graduate engineers. Using an expert's opinion we study the effect of such unemployed people on the society. An expert spells out the five major concepts relating to the unemployed graduated engineers as

| | | |
|---|---|---|
| $E_1$ | – | Frustration |
| $E_2$ | – | Unemployment |
| $E_3$ | – | Increase of educated criminals |
| $E_4$ | – | Under employment |
| $E_5$ | – | Taking up drugs etc. |

The directed graph where $E_1, \ldots, E_5$ are taken as the nodes and causalities as edges as given by an expert is given in the following Figure 1.1.1:

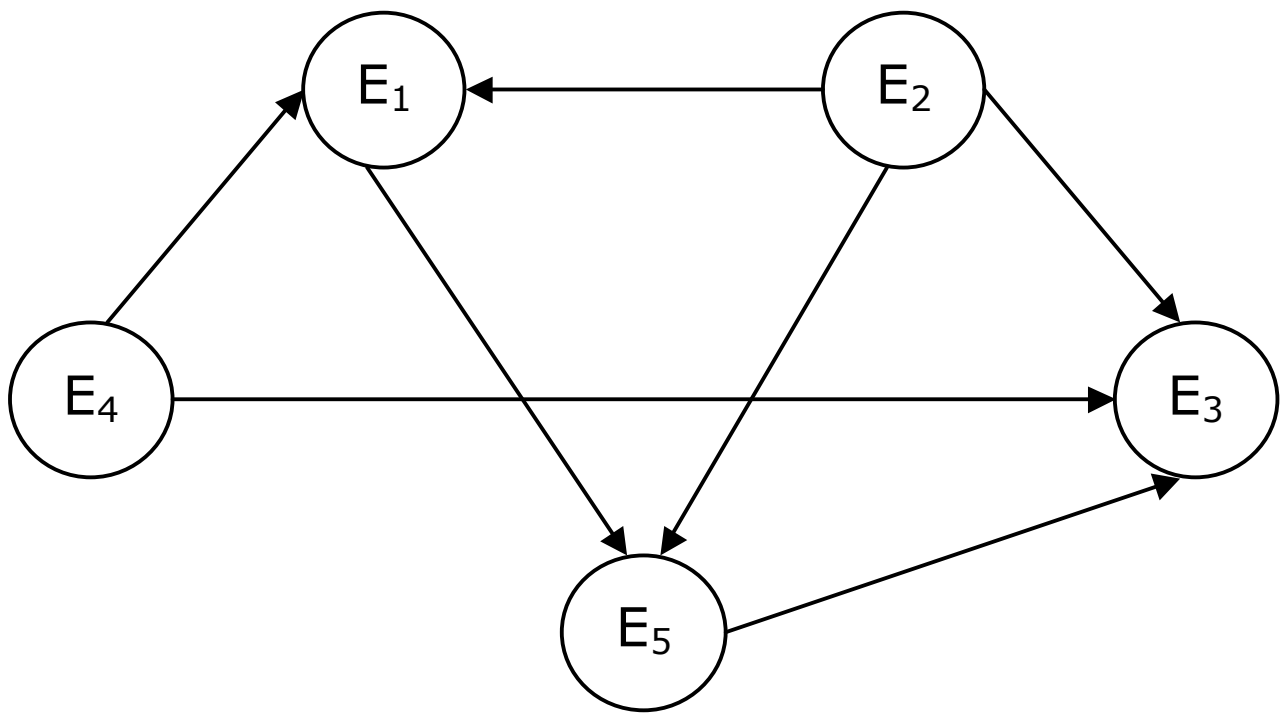

**FIGURE: 1.1.1**

According to this expert, increase in unemployment increases frustration. Increase in unemployment, increases the educated criminals. Frustration increases the graduates to take up to evils like drugs etc. Unemployment also leads to the increase in number of persons who take up to drugs, drinks etc. to forget their worries and unoccupied time. Under-employment forces then to do criminal acts like theft (leading to murder) for want of more money and so on. Thus one cannot actually get data for this but can use the expert's opinion for this unsupervised data to obtain some idea about the real plight of the situation. This is just an illustration to show how FCM is described by a directed graph.

{If increase (or decrease) in one concept leads to increase (or decrease) in another, then we give the value 1. If there exists no relation between two concepts the value 0 is given. If increase (or decrease) in one concept decreases (or increases) another, then we give the value –1. Thus FCMs are described in this way.}

**DEFINITION 1.1.2:** *When the nodes of the FCM are fuzzy sets then they are called as fuzzy nodes.*

**DEFINITION 1.1.3:** *FCMs with edge weights or causalities from the set {–1, 0, 1} are called simple FCMs.*

**DEFINITION 1.1.4:** *Consider the nodes / concepts $C_1, \ldots, C_n$ of the FCM. Suppose the directed graph is drawn using edge weight $e_{ij} \in \{0, 1, -1\}$. The matrix E be defined by $E = (e_{ij})$ where $e_{ij}$ is the weight of the directed edge $C_i C_j$. E is called the adjacency matrix of the FCM, also known as the connection matrix of the FCM.*

It is important to note that all matrices associated with an FCM are always square matrices with diagonal entries as zero.

**DEFINITION 1.1.5:** *Let $C_1, C_2, \ldots, C_n$ be the nodes of an FCM. $A = (a_1, a_2, \ldots, a_n)$ where $a_i \in \{0, 1\}$. A is called the instantaneous state vector and it denotes the on-off position of the node at an instant.*

$$a_i = 0 \text{ if } a_i \text{ is off and}$$
$$a_i = 1 \text{ if } a_i \text{ is on}$$

*for $i = 1, 2, \ldots, n$.*

**DEFINITION 1.1.6:** *Let $C_1$, $C_2$, ... , $C_n$ be the nodes of an FCM. Let $\overrightarrow{C_1C_2}$, $\overrightarrow{C_2C_3}$, $\overrightarrow{C_3C_4}$, … , $\overrightarrow{C_iC_j}$ be the edges of the FCM ($i \neq j$). Then the edges form a directed cycle. An FCM is said to be cyclic if it possesses a directed cycle. An FCM is said to be acyclic if it does not possess any directed cycle.*

**DEFINITION 1.1.7:** *An FCM with cycles is said to have a feedback.*

**DEFINITION 1.1.8:** *When there is a feedback in an FCM, i.e., when the causal relations flow through a cycle in a revolutionary way, the FCM is called a dynamical system.*

**DEFINITION 1.1.9:** *Let $\overline{C_1C_2}$, $\overline{C_2C_3}$, … ,$\overrightarrow{C_{n-1}C_n}$ be a cycle. When $C_i$ is switched on and if the causality flows through the edges of a cycle and if it again causes $C_i$ , we say that the dynamical system goes round and round. This is true for any node $C_i$, for $i = 1, 2, \ldots, n$. The equilibrium state for this dynamical system is called the hidden pattern.*

**DEFINITION 1.1.10:** *If the equilibrium state of a dynamical system is a unique state vector, then it is called a fixed point.*

***Example 1.1.2:*** Consider a FCM with $C_1$, $C_2$, …, $C_n$ as nodes. For example let us start the dynamical system by switching on $C_1$. Let us assume that the FCM settles down with $C_1$ and $C_n$ on i.e. the state vector remains as (1, 0, 0, …, 0, 1) this state vector (1, 0, 0, …, 0, 1) is called the fixed point.

**DEFINITION 1.1.11:** *If the FCM settles down with a state vector repeating in the form*

$$A_1 \to A_2 \to \ldots \to A_i \to A_1$$

*then this equilibrium is called a limit cycle.*

Methods of finding the hidden pattern are discussed in the following Section 1.2.

**DEFINITION 1.1.12:** *Finite number of FCMs can be combined together to produce the joint effect of all the FCMs. Let $E_1$, $E_2$, ... , $E_p$ be the adjacency matrices of the FCMs with nodes $C_1$, $C_2$, ..., $C_n$ then the combined FCM is got by adding all the adjacency matrices $E_1$, $E_2$, ..., $E_p$.*

*We denote the combined FCM adjacency matrix by $E = E_1 + E_2 + \ldots + E_p$.*

**NOTATION:** Suppose A = ($a_1$, … , $a_n$) is a vector which is passed into a dynamical system E. Then AE = ($a'_1$, … , $a'_n$) after thresholding and updating the vector suppose we get ($b_1$, … , $b_n$) we denote that by

$$(a'_1, a'_2, \ldots, a'_n) \to (b_1, b_2, \ldots, b_n).$$

Thus the symbol '→' means the resultant vector has been thresholded and updated.

FCMs have several advantages as well as some disadvantages. The main advantage of this method it is simple. It functions on expert's opinion. When the data happens to be an unsupervised one the FCM comes handy. This is the only known fuzzy technique that gives the hidden pattern of the situation. As we have a very well known theory, which states that the strength of the data depends on, the number of experts' opinion we can use combined FCMs with several experts' opinions.

At the same time the disadvantage of the combined FCM is when the weightages are 1 and –1 for the same $C_i$ $C_j$, we have the sum adding to zero thus at all times the connection matrices $E_1, \ldots, E_k$ may not be conformable for addition.

Combined conflicting opinions tend to cancel out and assisted by the strong law of large numbers, a consensus emerges as the sample opinion approximates the underlying population opinion. This problem will be easily overcome if the FCM entries are only 0 and 1.

We have just briefly recalled the definitions. For more about FCMs please refer Kosko [58].

## 1.2 Fuzzy Cognitive Maps – Properties and Models

Fuzzy cognitive maps (FCMs) are more applicable when the data in the first place is an unsupervised one. The FCMs work on the opinion of experts. FCMs model the world as a collection of classes and causal relations between classes.

FCMs are fuzzy signed directed graphs with feedback. The directed edge $e_{ij}$ from causal concept $C_i$ to concept $C_j$ measures how much $C_i$ causes $C_j$. The time varying concept function $C_i(t)$ measures the non negative occurrence of some fuzzy event, perhaps the strength of a political sentiment, historical trend or military objective.

FCMs are used to model several types of problems varying from gastric-appetite behavior, popular political developments etc. FCMs are also used to model in robotics like plant control.

The edges $e_{ij}$ take values in the fuzzy causal interval [–1, 1]. $e_{ij} = 0$ indicates no causality, $e_{ij} > 0$ indicates causal increase $C_j$ increases as $C_i$ increases (or $C_j$ decreases as $C_i$ decreases). $e_{ij} < 0$ indicates causal decrease or negative causality. $C_j$ decreases as $C_i$ increases (and or $C_j$ increases as $C_i$ decreases). Simple FCMs have edge values in {–1, 0, 1}. Then if causality occurs, it occurs to a maximal positive or negative degree. Simple FCMs provide a quick first approximation to an expert stand or printed causal knowledge.

We illustrate this by the following, which gives a simple FCM of a Socio-economic model. A Socio-economic model is constructed with Population, Crime, Economic condition, Poverty and Unemployment as nodes or concept. Here the simple trivalent directed graph is given by the following Figure 1.2.1 which is the experts opinion.

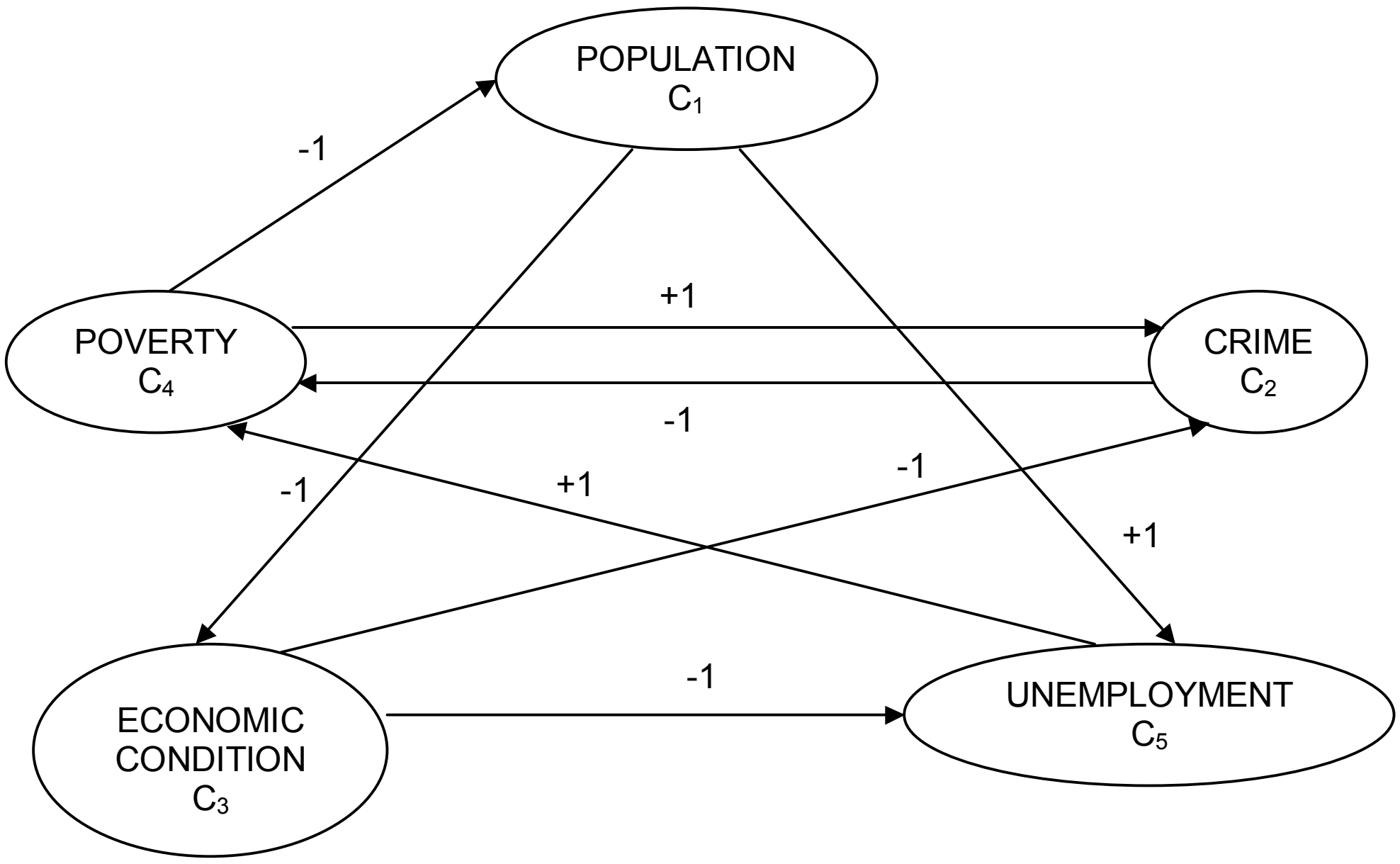

**FIGURE: 1.2.1**

Causal feedback loops abound in FCMs in thick tangles. Feedback precludes the graph-search techniques used in artificial-intelligence expert systems.

FCMs feedback allows experts to freely draw causal pictures of their problems and allows causal adaptation laws, infer causal links from simple data. FCM feedback forces us to abandon graph search, forward and especially backward chaining. Instead we view the FCM as a dynamical system and take its equilibrium behavior as a forward-evolved inference. Synchronous FCMs behave as Temporal Associative Memories (TAM). We can always, in case of a model, add two or more FCMs to produce a new FCM. The strong law of large numbers ensures in some sense that knowledge reliability increases with expert sample size.

We reason with FCMs. We pass state vectors C repeatedly through the FCM connection matrix E, thresholding or non-linearly transforming the result after each pass. Independent of the FCMs size, it quickly settles down to a temporal associative memory limit cycle or fixed point which is the hidden pattern of the system for that state vector C. The limit cycle or fixed-point inference summarizes the joint effects of all the interacting fuzzy knowledge.

***Example 1.2.1:*** Consider the 5 × 5 causal connection matrix E that represents the socio economic model using FCM given in figure in Figure 1.2.1.

$$E = \begin{bmatrix} 0 & 0 & -1 & 0 & 1 \\ 0 & 0 & 0 & -1 & 0 \\ 0 & -1 & 0 & 0 & -1 \\ -1 & 1 & 0 & 0 & 0 \\ 0 & 0 & 0 & 1 & 0 \end{bmatrix}$$

Concept nodes can represent processes, events, values or policies. Consider the first node $C_1 = 1$. We hold or clamp $C_1$ on the temporal associative memories recall process. Threshold signal functions synchronously update each concept after each pass, through the connection matrix E. We start with the population $C_1 = (1\ 0\ 0\ 0\ 0)$. The arrow indicates the threshold operation,

$$\begin{array}{lllllll} C_1 E & = & (0\ 0\ {-1}\ 0\ 1) & \rightarrow & (1\ 0\ 0\ 0\ 1) & = & C_2 \\ C_2 E & = & (0\ 0\ {-1}\ 1\ 1) & \rightarrow & (1\ 0\ 0\ 1\ 1) & = & C_3 \\ C_3 E & = & ({-1}\ 1\ {-1}\ 1\ 1) & \rightarrow & (1\ 1\ 0\ 1\ 1) & = & C_4 \\ C_4 E & = & ({-1}\ 1\ {-1}\ 0\ 1) & \rightarrow & (1\ 1\ 0\ 0\ 1) & = & C_5 \\ C_5 E & = & (0\ 0\ {-1}\ 0\ 1) & \rightarrow & (1\ 0\ 0\ 0\ 1) & = & C_6 \quad = \quad C_2. \end{array}$$

So the increase in population results in the unemployment problem, which is a limit cycle. For more about FCM refer Kosko [58] and for more about this socio economic model refer [112, 119].

This example illustrates the strengths and weaknesses of FCM analysis. FCM allows experts to represent factual and evaluative concepts in an interactive framework. Experts can quickly draw FCM pictures or respond to questionnaires. Experts can consent or dissent to the local causal structure and perhaps the global equilibrium. The FCM knowledge representation and inference-ing structure reduces to simple vector-matrix operations, favors integrated circuit implementation and allows extension to neural statistical or dynamical systems techniques. Yet an FCM equally encodes the experts' knowledge or ignorance, wisdom or prejudice. Worse, different experts differ in how they assign causal strengths to edges and in which concepts they deem causally relevant. The FCM seems merely to encode its designers' biases and may not even encode them accurately.

FCM combination provides a partial solution to this problem. We can additively superimpose each experts FCM in associative memory fashion, even though the FCM connection matrices $E_1, \ldots, E_K$ may not be conformable for addition. Combined conflicting opinions tend to cancel out and assisted the strong law of large numbers a consensus emerges as the sample opinion approximates the underlying population opinion. FCM combination allows knowledge researchers to construct FCMs with iterative interviews or questionnaire mailings.

The laws of large numbers require that the random samples be independent identically distributed random variables with finite variance. Independence models each experts individually. Identical distribution models a particular domain focus.

We combine arbitrary FCM connection matrices $F_1, F_2, \ldots, F_K$ by adding augmented FCM matrices. $F_1, \ldots, F_K$. Each augmented matrix $F_i$ has n-rows and n-columns n equals the total number of distinct concepts used by the experts. We permute the rows and columns of the augmented matrices to bring them into mutual coincidence. Then we add the $F_i$ point wise to yield the combined FCM matrix F.

$$F = \sum_i F_i$$

We can then use F to construct the combined FCM directed graph.

Even if each expert gives trivalent description in {–1, 0, 1}, the combined (and normalized) FCM entry $f_{ij}$ tends to be in {–1, 1}. The strong law of large numbers ensures that $f_{ij}$ provides a rational approximation to the underlying unknown population opinion of how much $C_i$ affects $C_j$. We can normalize $f_{ij}$ by the number K of experts. Experts tend to give trivalent evaluations more readily and more accurately than they give weighted evaluations. When transcribing interviews or documents, a knowledge engineer can more reliably determine an edge's sign than its magnitude.

Some experts may be more credible than others. We can weight each expert with non-negative credibility weight, weighing the augmented FCM matrix.

$$F = \Sigma\, w_i\, F_i.$$

The weights need not be in [0, 1]; the only condition is they should be non-negative. Different weights may produce different equilibrium limit cycles or fixed points as hidden patterns. We can also weigh separately any submatrix of each experts augmented FCM matrix.

Augmented FCM matrices imply that every expert causally discusses every concept $C_1, \ldots, C_n$. If an expert does not include $C_j$ in his FCM model the expert implicitly say that $C_j$ is not causally relevant. So the $j^{th}$ row and the $j^{th}$ column of his augmented connection matrix contains only zeros.

The only drawback which we felt while adopting FCM to several of the models is that we do not have a means to say or express if the relation between two causal concepts $C_i$ and $C_j$ is an indeterminate. So we in the next chapter will adopt in FCM the concept of indeterminacy and rename the Fuzzy Cognitive Maps as Neutrosophic Cognitive Maps as Neutrosophy enables one to accept the Truth, Falsehood and the Indeterminate. Such situations very often occur when we deal with unsupervised data that has more to do with feelings. Like political scenario, child labor, child's education, parent-children model, symptom-disease model, personality-medicine model in case of Homeopathy medicines, crime and punishment, in judicial problems (where evidences may be indeterminate), problem of the aged, female infanticide problem and so on and so forth.

Here we mention some examples in which certain factors are indeterminate but where we have used FCM to find solutions. So that in the next chapter we will use for these FCM model, Neutrosophic Cognitive Maps to study the systems.

***Example 1.2.2:*** In any nation the study of political situation i.e., the prediction of electoral winner or how people tend to prefer a particular politician and so on and so forth involves not only a lot of uncertainly for this, no data is available. They form an unsupervised data. Hence we study this model using FCM we see that while applying FCM to Indian politics the expert takes the following six nodes.

$x_1$ - Language
$x_2$ - Community

| | | |
|---|---|---|
| $x_3$ | - | Service to people, public figure configuration and Personality and nature |
| $x_4$ | - | Finance and Media |
| $x_5$ | - | Party's strength and opponent's strength |
| $x_6$ | - | Working members (Active volunteers) for the party. |

We using these six nodes and obtain several experts opinion. Now let us consider a case where in three or four parties join up to make a coalition party and stand in the election. Now each party has its own opinion and when they join up as an united front we will see the result. So, let us take the opinion of 4 experts on any arbitrary four concepts of the given six. We form the directed graph of the four experts opinion, we obtain the corresponding relational matrices. Finally we will get the combined effect. The Figures 1.2.2 to 1.2.5 correspond to directed graphs of the expert opinion.

First expert's opinion in the form of the directed graph and its relational matrix is given below.

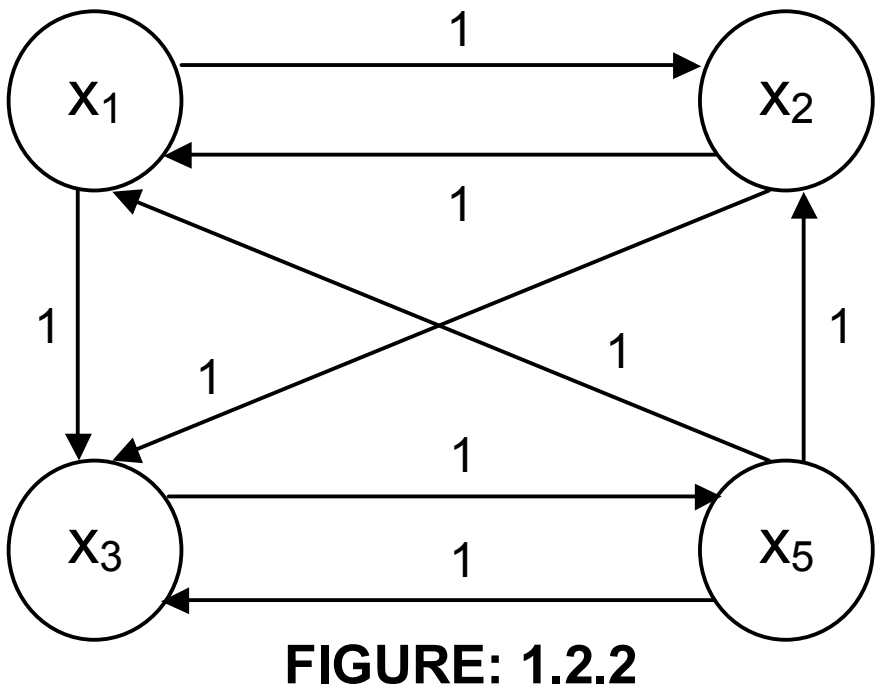

**FIGURE: 1.2.2**

$$E_1 = \begin{bmatrix} 0 & 1 & 1 & 0 & 0 & 0 \\ 1 & 0 & 1 & 0 & 0 & 0 \\ 0 & 0 & 0 & 0 & 1 & 0 \\ 0 & 0 & 0 & 0 & 0 & 0 \\ 1 & 1 & 1 & 0 & 0 & 0 \\ 0 & 0 & 0 & 0 & 0 & 0 \end{bmatrix}.$$

Second expert's opinion – directed graph is shown in Figure 1.2.3 and the relational matrix is given by $E_2$.

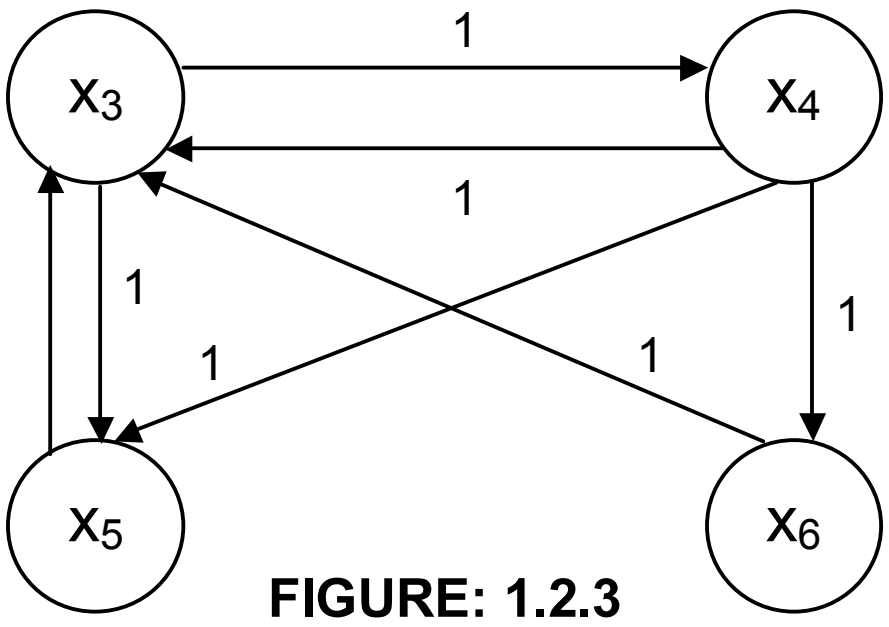

**FIGURE: 1.2.3**

$$E_2 = \begin{bmatrix} 0 & 0 & 0 & 0 & 0 & 0 \\ 0 & 0 & 0 & 0 & 0 & 0 \\ 0 & 0 & 0 & 1 & 1 & 0 \\ 0 & 0 & 1 & 0 & 1 & 1 \\ 0 & 0 & 1 & 0 & 0 & 0 \\ 0 & 0 & 1 & 0 & 0 & 0 \end{bmatrix}.$$

Third expert's opinion with the four conceptual nodes $x_2$, $x_3$, $x_4$ and $x_6$, its directed graph and the relational matrix are given in Figure 1.2.4 and its related matrix $E_3$.

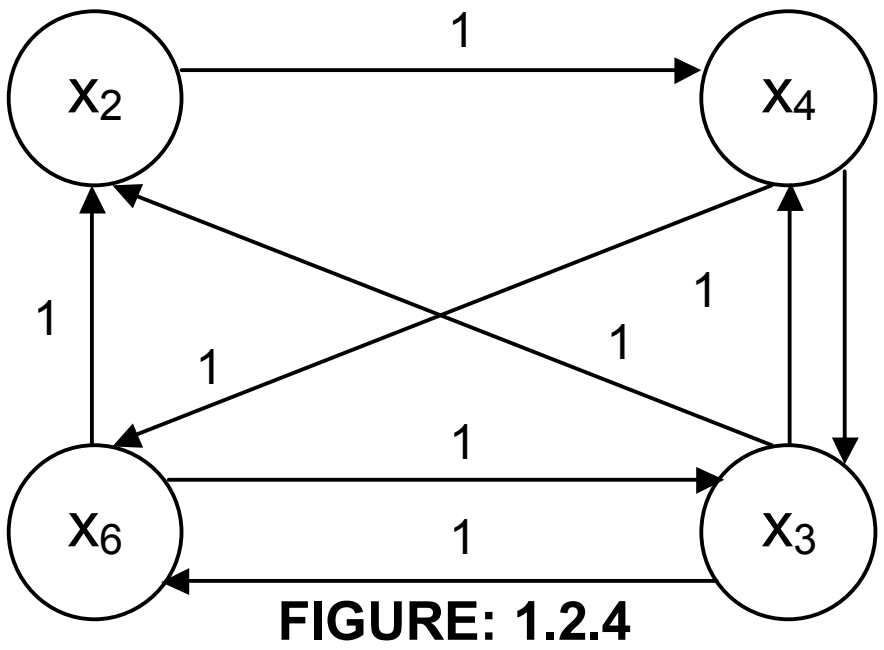

**FIGURE: 1.2.4**

$$E_3 = \begin{bmatrix} 0 & 0 & 0 & 0 & 0 & 0 \\ 0 & 0 & 0 & 1 & 0 & 0 \\ 0 & 1 & 0 & 1 & 0 & 1 \\ 0 & 0 & 1 & 0 & 0 & 1 \\ 0 & 0 & 0 & 0 & 0 & 0 \\ 0 & 1 & 1 & 0 & 0 & 0 \end{bmatrix}.$$

Directed graph and the relational matrix of a fourth expert using the concepts $x_2$, $x_5$, $x_4$ and $x_6$ is given in Figure 1.2.5 and the related matrix $E_4$ is given below.

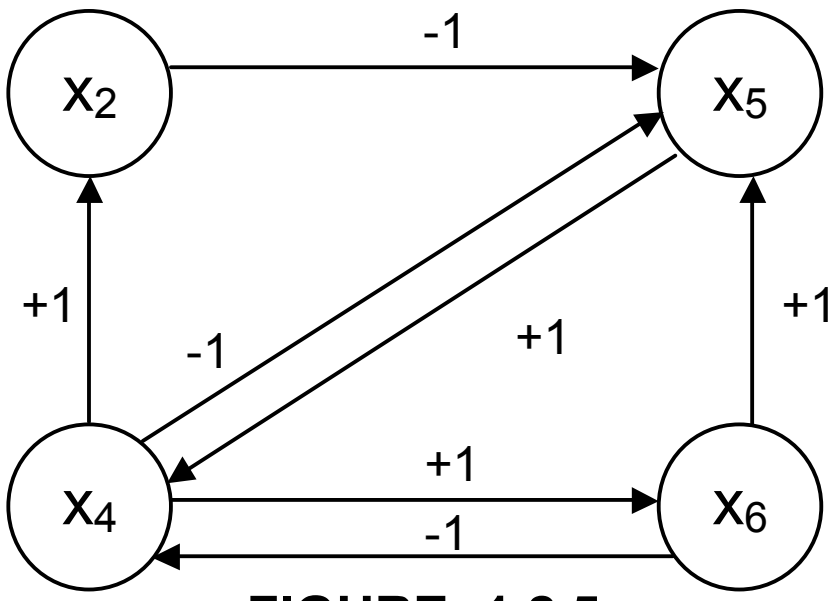

**FIGURE: 1.2.5**

$$E_4 = \begin{bmatrix} 0 & 0 & 0 & 0 & 0 & 0 \\ 0 & 0 & 0 & 0 & -1 & 0 \\ 0 & 0 & 0 & 0 & 0 & 0 \\ 0 & 1 & 0 & 0 & -1 & 1 \\ 0 & 0 & 0 & 1 & 0 & 0 \\ 0 & 0 & 0 & -1 & 1 & 0 \end{bmatrix}.$$

We note that each matrix contains two zero rows and two zero columns corresponding to the experts causally irrelevant concepts. We now combine the directed graph of the four experts and obtain the Figure 1.2.6.

The combined FCM matrix E, which is equal to the addition of the four matrices and its related directed graph, is as follows:

Combined directed graph of the combined FCM is given in Figure 1.2.6:

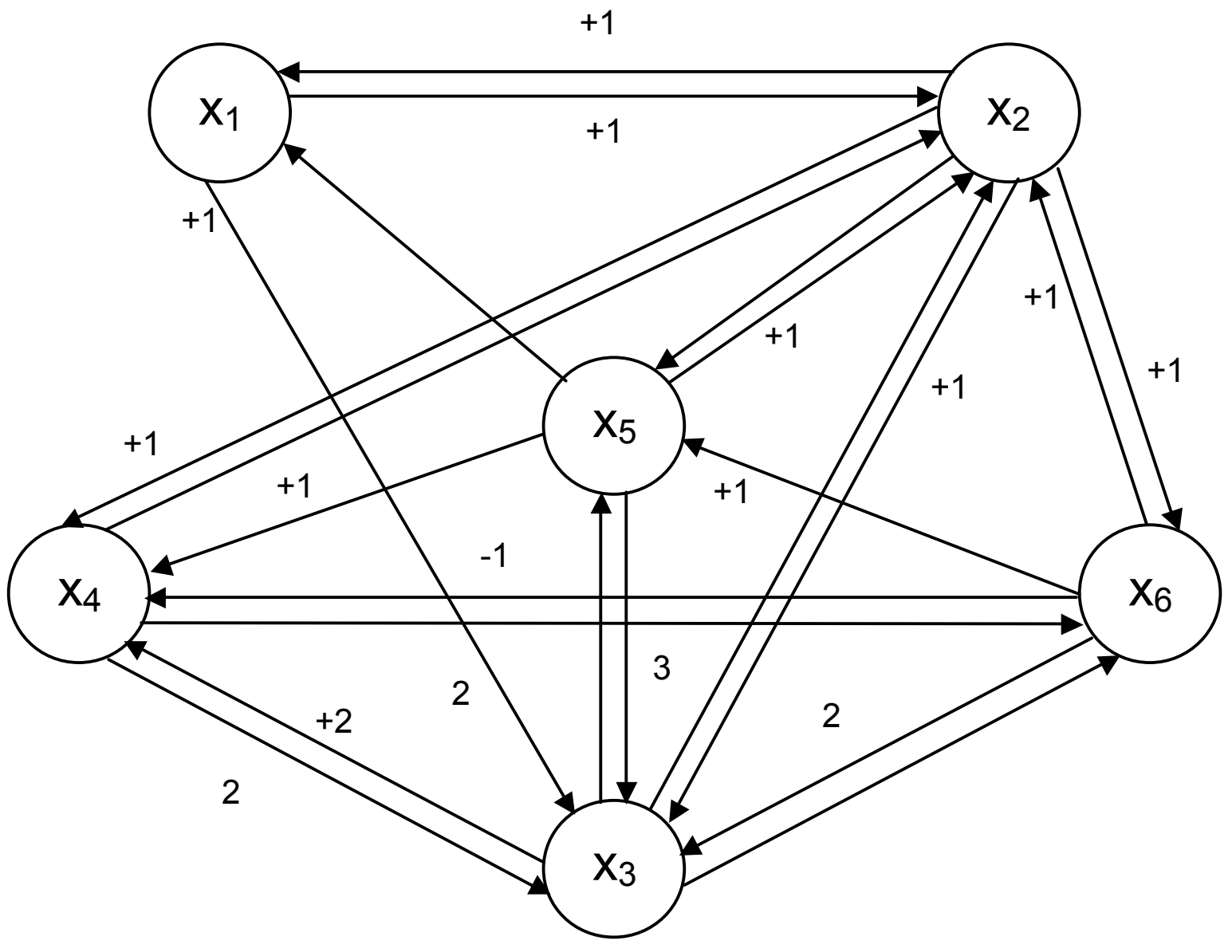

**FIGURE: 1.2.6**

$$E = \begin{bmatrix} 0 & 1 & 1 & 0 & 0 & 0 \\ 1 & 0 & 1 & 1 & -1 & 0 \\ 0 & 1 & 0 & 2 & 2 & 1 \\ 0 & 1 & 2 & 0 & 0 & 3 \\ 1 & 1 & 2 & 1 & 0 & 0 \\ 0 & 1 & 2 & -1 & 1 & 0 \end{bmatrix}.$$

We now consider a single state element through the matrix and multiply it with the above relational matrix and we see that nearly every state is active.

For example, let us consider the input with parties' strength and opponent's strength in the on state and rest of the coordinates as off state, vector matrix A = (0, 0, 0, 0, 1, 0). Now by passing on we get the matrix AE = (1, 1, 2, 1, 0, 0) $\rightarrow$ (1, 1, 2, 1, 1, 0).

Thus we see parties' strength and the opponent's strength is a very sensitive issue and it makes all other coordinate to be on except the state of the working members.

Likewise when the two concepts the community service to people and the public figure configuration is taken to be in the on state we see that on passing through the system all state becomes on.

Thus the FCMs gives us the hidden pattern. No other method can give these resultants that too with the unsupervised data. For more about this please refer [116].

***Example 1.2.3:*** The study of symptoms and its associations with disease in children happens to be very uncertain and difficult one.

For at times the doctor treats for the symptoms instead of treating for disease. So that when a ready-made model is made it may serve the better purpose for the doctor.

We have also adopted the FCM in case of Symptom-disease model in children [115]. To build the symptom-disease model for children we use the following 8 nodes of FCM, which are labeled as follows:

| | | |
|---|---|---|
| $C_1$ | - | Fever with cold |
| $C_2$ | - | Fever with vomiting (or) Loose Motions |
| $C_3$. | - | Fever with loss of appetite |
| $C_4$ | - | Fever with cough |
| $C_5$ | - | Respiratory diseases |
| $C_6$ | - | Gastroenteritis |
| $C_7$ | - | Jaundice |
| $C_8$ | - | Tuberculosis. |

The directed graph as given by the doctor who is taken as an expert is given in the Figure: 1.2.7, which is as follows:

The corresponding connection or adjacency matrix E is as follows:

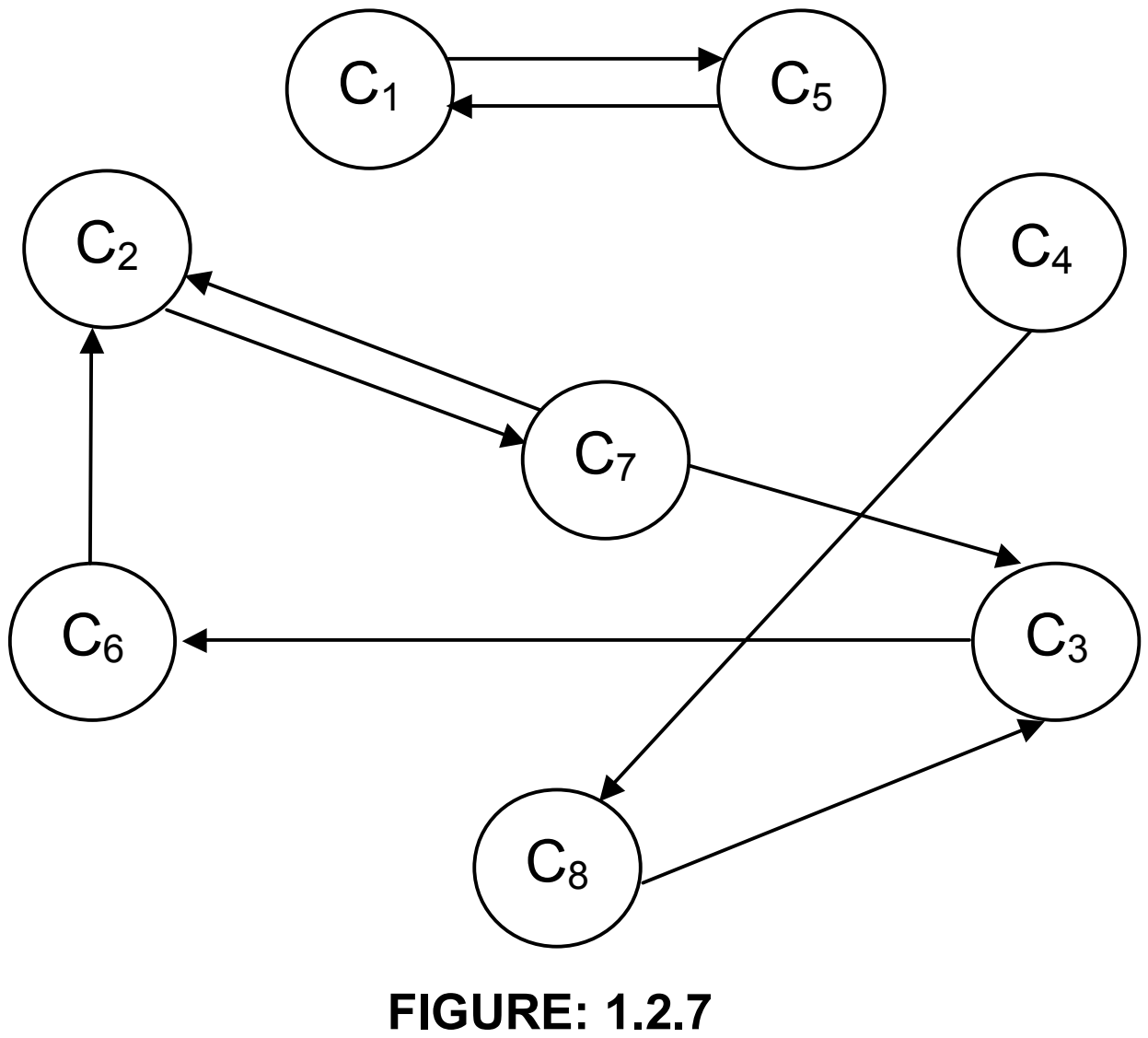

**FIGURE: 1.2.7**

$$E = \begin{bmatrix} 0 & 0 & 0 & 0 & 1 & 0 & 0 & 0 \\ 0 & 0 & 0 & 0 & 0 & 0 & 1 & 0 \\ 0 & 0 & 0 & 0 & 0 & 1 & 0 & 0 \\ 0 & 0 & 0 & 0 & 0 & 0 & 0 & 1 \\ 1 & 0 & 0 & 0 & 0 & 0 & 0 & 0 \\ 0 & 1 & 0 & 0 & 0 & 0 & 0 & 0 \\ 0 & 1 & 1 & 0 & 0 & 0 & 0 & 0 \\ 0 & 0 & 1 & 0 & 0 & 0 & 0 & 0 \end{bmatrix}.$$

Input the vector $A_1 = (1\ 0\ 0\ 0\ 0\ 0\ 0\ 0)$

$$A_1E \to (1\ 0\ 0\ 0\ 1\ 0\ 0\ 0) = A_2$$
$$A_2E \to (1\ 0\ 0\ 0\ 1\ 0\ 0\ 0) = A_3 = A_2.$$

According to this doctor's opinion, fever and cold induces the respiratory diseases. Suppose the vector $D_1 = (0\ 0\ 1\ 0\ 0\ 0\ 0\ 0)$ is taken as the input vector;

$$D_1E = (0\ 0\ 0\ 0\ 0\ 1\ 0\ 0) \to (0\ 0\ 1\ 0\ 0\ 1\ 0\ 0) = D_2$$
$$D_2E \to (0\ 1\ 1\ 0\ 0\ 1\ 0\ 0) = D_3$$
$$D_3E \to (0\ 1\ 1\ 0\ 0\ 1\ 1\ 0) = D_4$$
$$D_4E = D_4.$$

Thus when a child suffers with the symptom fever and loss of appetite then the doctor suspects the child may develop fever with vomiting or loose motion leading to the sickness of gastroenteritis and jaundice. For further results in this direction please refer [115].

Use of FCMs in the study of the maximum utility of a bus route in Madras city (in South India) happens to be a difficult one for the concept deals with the many aspects of modern metropolitan public transportation. Now we just illustrate how we have applied FCM to the problem of determining the maximum utility of a route [122]. Here we not only give the FCM model to find the maximum utility of a route but also we have developed a Java program to study implications of the model [123].

***Example 1.2.4:*** The application of FCM in this case is done in a very ingenious way. Here for the first time FCMs are used in identifying the maximum utilization of a time period in a day, in predicting the overall utility of the routes and in the stability analysis.

We represent the different time periods of a day, viz morning, noon, evening etc by the total number of individual hours i.e. $\{H_1, H_2, \ldots, H_{24}\}$, where $H_i$ represents the $i^{th}$ hour ending of the day. In order to estimate the utility rate of a route and to identity the peak period of a route, we consider the total number of time periods in a day and the various attributes acting on these time periods as the conceptual nodes viz. $\{C_1, \ldots, C_{10}\}$. It is in the hands of the designer of the problem to assign the time periods from $\{H_1, H_2, \ldots, H_{24}\}$ to each value of $C_i$. Here we take n = 10 say $\{C_1, C_2, \ldots, C_{10}\}$ where $C_1$ corresponds to the early hours $\{H_6, H_7\}$, $C_2$ corresponds to the morning hours $\{H_8, H_9, H_{10}\}$, $C_3$ corresponds to early noon hours $\{H_{11}, H_{12}, H_{13}\}$, $C_4$ refers to the evening $\{H_{14}, H_{15}, H_{16}\}$. $C_5$ corresponds to late evening hours $\{H_{17}, H_{18}, H_{19}\}$. $C_6$ corresponds to the night hours, $\{H_{20}, H_{21}, H_{22}\}$, $C_7$ indicates the number of passengers in each time period, $C_8$ corresponds to the total collection in each time period, $C_9$ denotes the number of trips made in each time period and $C_{10}$ corresponds to the hourly occupancy in each time period. We make an assumption i.e. an increase in the conceptual nodes depicting the time period say $C_i$ comprising of hours $\{H_i, H_{(i+1)}, \ldots, H_{(i+n)}\}$ will imply a graded rise in the hour of the day from hour $H_i$ to hour $H_{i+n}$.

To assess the interactions among the conceptual nodes, we collect expert opinion. The more number of experts the higher is the reliability in the knowledge base. The experts include the passengers traveling in the city transport service, the raw data obtained from the Pallavan Transport Corporation (currently renamed as Metropolitan Transport Corporation) (taken from the Madras city in India) denoting its route-wise loading pattern analysis. Though we have taken several experts opinion, here we give only the expert opinion of a regular passenger traveling along the route 18B, the adjacency matrix of it is given below:

$$E = \begin{bmatrix} 0 & 0 & 0 & 0 & 0 & 0 & 1 & 1 & -1 & 1 \\ 0 & 0 & 0 & 0 & 0 & 0 & 1 & 1 & 1 & 0 \\ 0 & 0 & 0 & 0 & 0 & 0 & -1 & -1 & -1 & -1 \\ 0 & 0 & 0 & 0 & 0 & 0 & -1 & -1 & -1 & 0 \\ 0 & 0 & 0 & 0 & 0 & 0 & 1 & -1 & 1 & 0 \\ 0 & 0 & 0 & 0 & 0 & 0 & -1 & -1 & -1 & -1 \\ 0 & 0 & 0 & 0 & 0 & 0 & 0 & 1 & 1 & 0 \\ 0 & 0 & 0 & 0 & 0 & 0 & 1 & 0 & 1 & 1 \\ 0 & 0 & 0 & 0 & 0 & 0 & 1 & 0 & 0 & -1 \\ 0 & 0 & 0 & 0 & 0 & 0 & 1 & 1 & -1 & 0 \end{bmatrix}.$$

To study the hidden pattern first we clamp the concept with the first node in the on state and rest of the nodes in the off state

$A_1$ = (1 0 0 0 0 0 0 0 0 0)
$A_1E$ = (0 0 0 0 0 0 1 1 –1 1) $\to$ (1 0 0 0 0 0 1 1 0 1) = $A_2$
$A_2E$ = (0 0 0 0 0 0 2 3 0 3) $\to$ (1 0 0 0 0 0 1 1 0 1) = $A_3$
$A_3E$ = (0 0 0 0 0 0 3 2 0 2) $\to$ (1 0 0 0 0 0 1 1 0 1) = $A_4 = A_3$.

We observe that increase in time period $C_1$ increases the number of passengers, the total collection and the number of trips. Next consider the increase of the time period $C_2$. The input vector is

$B_1$ = (0 1 0 0 0 0 0 0 0 0)
$B_1E$ = (0 0 0 0 0 0 1 1 1 0) $\to$ (0 1 0 0 0 0 1 1 1 0) = $B_2$
$B_2E$ = (0 0 0 0 0 0 3 2 3 0) $\to$ (0 1 0 0 0 0 1 1 1 0) = $B_2$.

Thus we observe that an increase in the time period $C_2$ leads to an increase in number of passengers, the total collection and the increased occupancy.

Consider the combined effect of increasing $C_2$ and $C_{10}$ i.e. increase of trips in the time period $\{H_8, H_{10}, H_{10}\}$.

| | | | | | | |
|---|---|---|---|---|---|---|
| $W_1$ | = | (0 1 0 0 0 0 0 0 0 1) | | | | |
| $W_1E$ | = | (0 0 0 0 0 0 2 2 0 0) | $\to$ | (0 1 0 0 0 0 1 1 0 1) | = | $W_2$ |
| $W_2E$ | = | (0 0 0 0 0 0 3 3 2 1) | $\to$ | (0 1 0 0 0 0 1 1 1 1) | = | $W_3$ |
| $W_3E$ | = | (0 0 0 0 0 0 4 3 2 0) | $\to$ | (0 1 0 0 0 0 1 1 1 1) | = | $W_4 = W_3$. |

Thus increase of time period $C_2$ and $C_{10}$ simultaneously increases the total number of passengers, the total collection and the occupancy rate.

### THE JAVA PSEUDO CODE TO EVALUATE THE FIXED POINTS FOR CLAMPED VECTORS

```
import java.applet.Applet;
import java.awt.*;
import java.lang.*;
import java.io.*;
import java.net.*;
import java.util.*;
public class cognitive extends java.applet.Applet
{
public void init()
   {
/* Through Text Label capture the Number of Concepts used
in this Study */
/* Define an Editable input window */
/* Define a non-Editable output window */
/* Add Button for Sample 1 to provide data for sample
matrix and the number of Concepts */
/* Add Evaluate Button: If clicked on this it provides
the values for Clamped vectors and the corresponding
Fixed point in the system */
/* Add a clear button to clear the values of concepts,
input window, output window */
/* Show all the values defined above */
}/* init() */

public boolean action(Event evt, Object arg)
{
String label = (String)arg;
if(evt.target instanceof Button)
{
if(label.equals("EVALUATE"))
{
/* Capture the number of rows in Matrix */
/* Call the function readandevaluate() to evaluate the
values */
/* show the calculated values */
```

```
}
else if (label.equals("Clear")
{
/*clear the rows, input window and output window */
}
else if (label.equals("Load Sample Matrix"))
{
/* set the number of values to 5 */
try
{
dataURL=new URL
("http://members.tripod.com/~mandalam/java1-
0/data3.txt");
try
{
/* get the sample data from the URL above and populate
the input window */
}
catch(IOException e)
{
System.err.println("Error:"+e);
}
}
catch(MalformedURLException e) {
return false;
}
return true;
}/* sample Matrix*/
  }
return false;
  }
void readandevaluate(String mystring)
  {
/*Tokenize the input window and get the initial Adjacency
matrix in A */
/* show the above read adjacency matrix A in output
window */
/* Determine the Hidden Pattern in B, the Clamped Vector
*/
/* Determine the new vector C */
/* apply Threshold to new vector C */
/* check if C is same as B */
/* Determine the new vector C and compare it with B */
/* if both are same this is the Fixed point in the
system, else continue by copying the values of C into B
*/
  }/* readandevaluate() */
}/* end of class Cognitive */
```

Several new results in this direction can be had from [122]. Another type of problem about transportation is as follows.

The problem studied in this case is for a fixed source S, a fixed destination D and a unique route from the source to the destination, with the assumption that all the passengers travel in the same route, we identify the preferences in the regular services at the peak hour of a day.

We have considered only the peak-hour since the passenger demand is very high only during this time period, where the transport sector caters to the demands of the different groups of people like the school children, the office goers, the vendors etc.

We have taken a total of eight characteristic of the transit system, which includes the level of service and the convenience factors. We have the following elements, Frequency of the service, in-vehicle travel time, the travel fare along the route, the speed of the vehicle, the number of intermediate points, the waiting time, the number of transfers and the crowd in the bus or equivalently the congestion in the service.

Before defining the cognitive structure of the relationship, we give notations to the concepts involved in the analysis as below.

| | | |
|---|---|---|
| $C_1$ | - | Frequency of the vehicles along the route |
| $C_2$ | - | In-vehicle travel time along the route |
| $C_3$ | - | Travel fare along the route |
| $C_4$ | - | Speed of the vehicles along the route |
| $C_5$ | - | Number of intermediate points in the route |
| $C_6$ | - | Waiting time |
| $C_7$ | - | Number of transfers in the route |
| $C_8$ | - | Congestion in the vehicle. |

The graphical representation of the inter-relationship between the nodes is given in the form of directed graph given in Figure: 1.2.8.

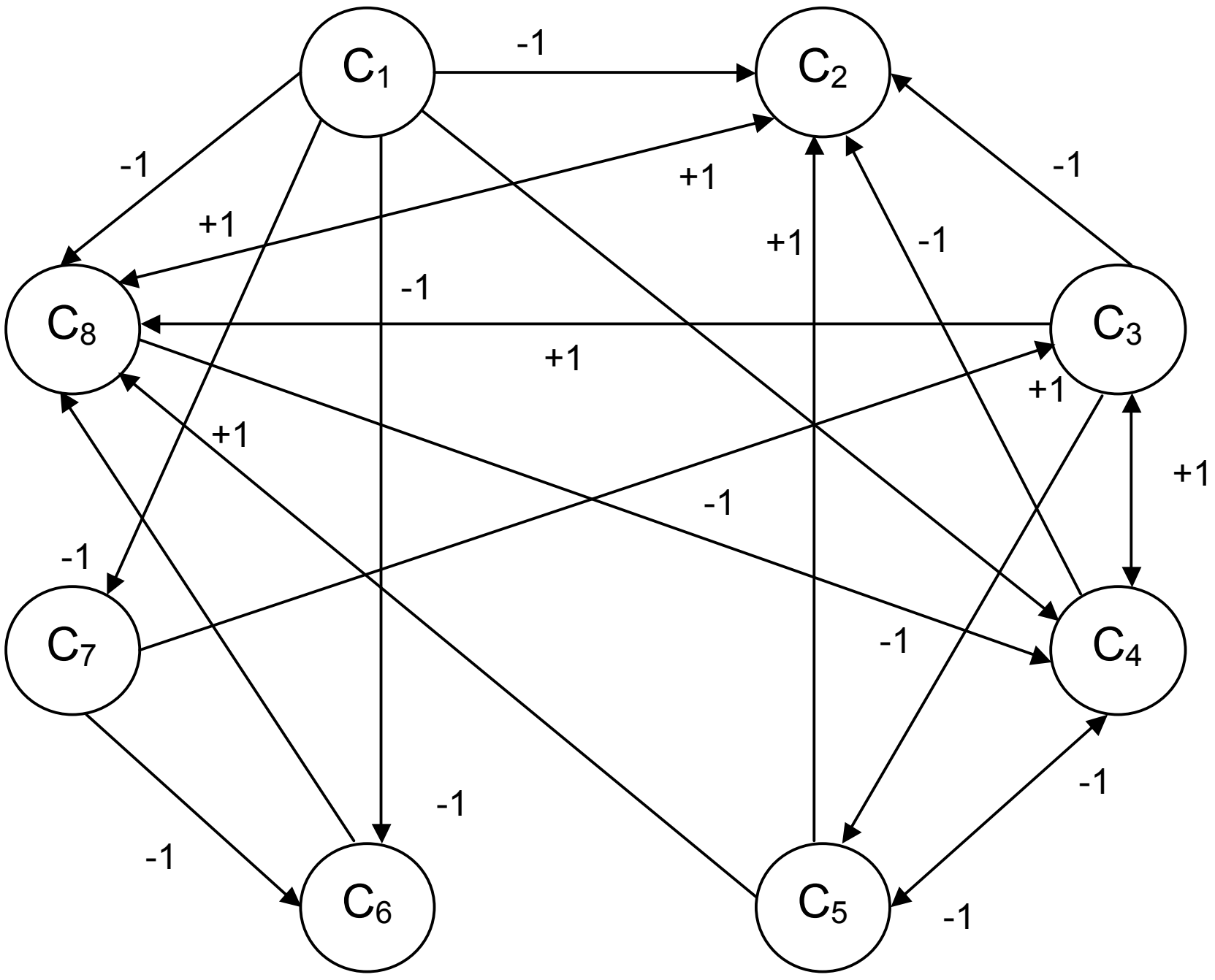

**FIGURE: 1.2.8**

From the above signed directed graph, we obtain a connection matrix E, since the number of concepts used here are eight, the connection matrix is a 8 × 8 matrix.
Thus we have $E = [A_y]_{8\times 8}$

$$E = \begin{bmatrix} 0 & -1 & 0 & 1 & 0 & -1 & -1 & -1 \\ 0 & 0 & 0 & 0 & 0 & 0 & 0 & 1 \\ 0 & -1 & 0 & 1 & -1 & 0 & 0 & -1 \\ 0 & -1 & 1 & 0 & -1 & 0 & 0 & 0 \\ 0 & 1 & 0 & -1 & 0 & 0 & 0 & 1 \\ 0 & 0 & 0 & 0 & 0 & 0 & 0 & 1 \\ 0 & 0 & 1 & 0 & 0 & -1 & 0 & 0 \\ 0 & 1 & 0 & -1 & 0 & 0 & 0 & 0 \end{bmatrix}.$$

The model implications i.e. the predictions regarding changes in the state behaviour of a model are determined by activating the individual involved elements. This is achieved by 'Clamping' concepts of interest and including iterative operations on the matrix. We use the equation.

$$I_{t+1} = O_{(t)} = I_{(t)} * E$$

where $I_{(t)}$ is the input vector at the $t^{th}$ iteration, E is the connection matrix $O_{(t)}$ is the output vector at the $t^{th}$ iteration, used as the input for the $(t+1)^{th}$ iteration.

Initially, we clamp the concept of interest. This is done by fixing the corresponding element in the input and the rest of the elements are given a value zero. Thus we have the input vector as a 1 × n vector (where n is the number of concepts involved). Now we do matrix multiplication operation of $I_{1\times n}$ on the connection matrix $E_{(n\times n)}$. The output O is again a vector of order 1 × n. Now we implement the threshold function – a binary conversion of the output vector.

Thus we have

$$O(x) = \begin{cases} 1 & \text{if } x > 1 \\ 0 & \text{otherwise} \end{cases}.$$

The output vector after the implementation of the threshold function is used as the input vector for the $(t+1)^{th}$ stage. This new input vector again is operated on the connection matrix. The working process described above can be expressed by this algorithm.

**ALGORITHM TO STUDY THE MODEL IMPLICATIONS**

Begin

Step 1: Read the input vector $I_{(t)}$

Step 2: Give the connection matrix, E.

Step 3: Calculate the output vector $O_{(t)} = I_{(t)} * E$

Step 4: Apply threshold to output vector:

$$O_{(t)} \cong I_{(t+1)}$$

Step 5: If $(I_{(t+1)} = I_t)$, stop

else go to step 1.

End

By studying the final state of the iterations, we infer from the model. An equilibrium in the system is attained when we have a set of repeated patterns. Repeating patterns can be fixed points or limiting cycles or a chaotic attractor. A 'fixed point' is a single recurring pattern such as $A_4 \Rightarrow A_4$, in the pattern $A_1 \Rightarrow A_2 \Rightarrow A_3 \Rightarrow A_4 \Rightarrow A_4 \Rightarrow A_1$ . A "limiting cycle" is a set of multiple repeating patterns such as $A_3 \Rightarrow A_4 \Rightarrow A_5$ in the pattern $A_1 \Rightarrow A_2 \Rightarrow A_3 \Rightarrow A_4 \Rightarrow A_5 \Rightarrow A_3 \Rightarrow A_4 \Rightarrow A_5$. Here, it is a limit cycle of length 3 or 3-L. A “Chaotic attractor” is a limit cycle consisting of repeating patterns of differing lengths.

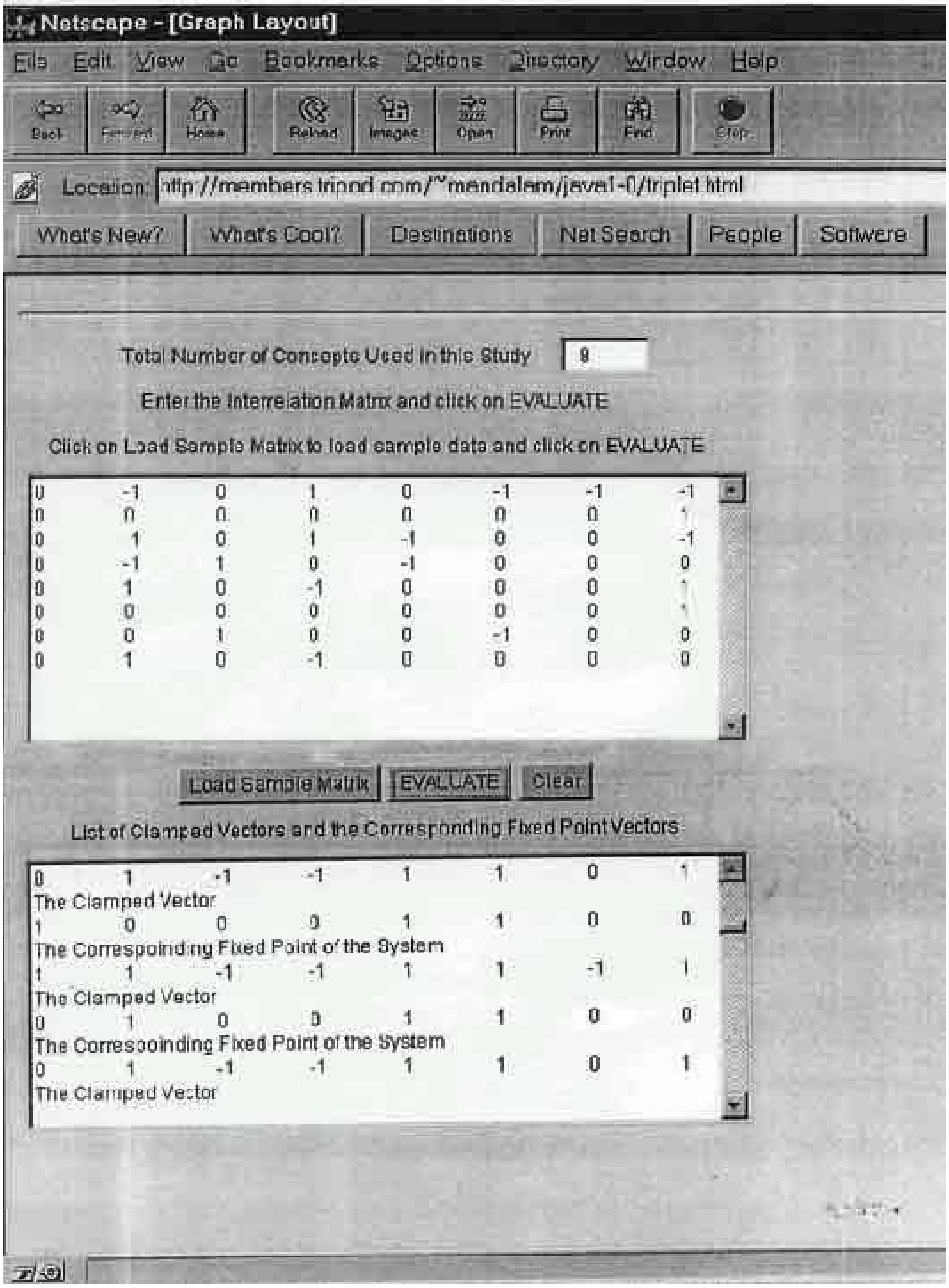

These equilibrium states are model inferences.

We give the Java screen depicting the loading of the simple matrix and the clamped vectors and their fixed points.

For a more detailed analysis of this refer [123].

***Example 1.2.5:*** Further, in the year 1999, Ozesmi U. [77] has used FCM to study the ecosystems of the Kizilirmar delta wetlands in Turkey using the expert opinions of the local people, NGOs (non-governmental organizations), government officials, stakeholder groups and vacation house owners. The data under study happens to be an unsupervised one and further his study was based on 31 FCM models, which were converted to adjacency matrices. Villages or the locals FCMS showed a large capacity to adapt to changing ecological and social conditions. They actively changed and challenged these conditions through the political process. Villagers were faced with many important forcing function they could not control. Most of the variables defined by villagers were related to agricultural and animal husbandry. Villagers or locals view must be given due weightage as mainly they are well versed about the situation and the probable problems. He has analyzed the FCMs using the following definitions:

**DEFINITION 1.2.1:** *The density of a cognitive map D is an index of connectivity.*

*$D = C / N(N-1)$ or alternatively $D = C / N^2$. This is known as the density equation where in C represents the number of connections possible between N variables if the number of connections possible between N variables can have a causal effect on themselves then the minimum number of connections is $N^2$.*

The structure of a FCM apart from the number of variables and connections can best be analyzed by finding the following variables:

**DEFINITION 1.2.2:** *For a given FCM the transmitter variables (are forcing, given, tails, independent) the receiver variables (ends, heads, dependent).*

*These variables are defined by their outdegree $[od(v_j)]$ and indegree $[id(v_j)]$. Outdegree is the row sum of absolute values of a variables in the adjacency matrix and shows the cumulative strengths of connections ($a_{ij}$) exiting the variable*

$$od(v_j) = \sum_{k=1}^{N} \bar{a}_{jk}$$

*indegree is the column sum of absolute values of a variable and shows the cumulative strength of variables entering the unit*

$$id(v_j) = \sum_{k=1}^{N} \bar{a}_{kj}\,.$$

*The immediate domain of a variable is the summation of its indegree (inarrows) and outdegree (outarrows) also called centrality. The contribution of a variable in a cognitive map can be understood by calculating its centrality (c) whether it is a*

*transmitter, receiver or ordinary variable. The centrality (c) of a variable is also called its total degree [td($v_i$)].*

$$c_i = td(v_i) = od(v_i) + id(v_i).$$

*Transmitter variables are units whose od($v_i$) is positive and their id($v_i$) is 0. Receiver variables are units whose od($v_i$) is 0 and their id($v_i$) is positive. Other variables which have both non-zero od($v_i$) and id($v_i$) are ordinary variables (mean).*

For more refer [10, 34, 35, 36].

The total number of receiver variables in a cognitive map can be considered an index of its complexity. Larger number of receiver indicate that the cognitive map considers many outcomes and implications that are a result of the system.

We just recall the notions of receiver to transmitter variable rates and the hierarchy index.

**DEFINITION 1.2.3:** *Many transmitter variables indicate thinking with top down influences, a formal hierarchical system. Many transmitter units shows the flatness of a cognitive map where causal arguments are not well elaborated. Then we can compare cognitive maps in terms of their complexity by number of receiver to transmitter variable ratios (R / T). Larger ratios will indicate more complex maps because they define more utility outcomes and are less controlling forcing functions. Another structural measure of a cognitive map is the hierarchy index (h)*

$$h = \frac{12}{(N-1)N(N+1)} \sum_i \left[ od(v_i) - \frac{\sum od(v_i)}{N} \right]^2$$

*where N is the total number of variables. When h is equal to one then the map is fully hierarchical and when h is equal to zero the system is fully democratic.*

This is an excellent piece of research where a nice analysis of "ecosystems" using FCM is made. Definitions from this paper are mainly recalled to enable the reader to apply to possible problems using these new definitions. For complete data refer [77].

## 1.3 Some more Illustrations of FCMs

In this section we give the illustrations of FCM models of eleven researchers in order to show to the reader the wide areas in which FCMs have found usage. We have only given an extract and a basic idea of how FCMs are used. For a detailed description one can refer the original research papers. These concepts are mainly given to show in the later chapter how Neutrosophic Cognitive Maps will be preferred in some cases. Also, we in this section recall the definitions of various types of FCMs which have been introduced by researchers including: Adaptive FCM, A new Balance Degree for FCM, Automatic Construction of FCM, Rule Based FCM etc.

### 1.3.1: FCM for Decision Support in an Intelligent Intrusion Detection System

As computer technology advances and the threats of computer crime increase, the apprehension and preemption of such infractions become more and more difficult and challenging. Over the years, intrusion detection has become a major area of research in computer science. Intrusion detection systems are often characterized based on two aspects [18].

a. the data source (host based/ multi-host based / network based), or
b. the model of intrusion detection (anomaly detection/misuse detection).

There are wide variations in the techniques used by intrusion detection systems.

The Intelligent Intrusion Detection System (IIDS) is a prototype developed as part of an intrusion detection research effort in the Department of Computer Science at Mississippi State University. The following unique features characterize IIDS:

1. Real time adaptive distributed and network based architecture.
2. Incorporation of both anomaly and misuse detection (i.e., misuse detection modules look for known patterns of attack while anomaly detection modules look for deviations from "normal" patterns of behaviour).
3. Integration of data mining algorithms with fuzzy logic and the use of genetic algorithms for optimization of membership functions and for feature selection.
4. Use of a decision engine to fuse information from different types of detection modules in order to make decisions about the overall health of the network.

We now recall from Siraj. A. et al [89] the description of a model of a decision engine incorporated in the intelligent intrusion detection system architecture that uses causal knowledge reasoning with FCMs. A primary task of the decision engine is to investigate the results generated by the misuse detection components that look at signatures of known attacks. For misuse detection the IIDS uses rule based detection mechanisms that work on each of the hosts of the network. Output from the misuse detection modules may be crisp or fuzzy. For certain types of attacks, the signature is either present or absent and the output of the module is binary. For other types of attacks like the number of failed logins, the output of the misuse detection module is a fuzzy measure of the degree of suspicion. The decision engine must assess results of the multiple misuse detection modules in order to compute the alert status for each machine and for each user account.

The following simple example illustrates how an FCM can be used to capture a suspicious event pattern. Consider the following scenario: the intruder is trying to access the network by breaking into a particular workstation. The intruder tries to login with several users' passwords and fails. One can identify such an attack scenario by observing the number of login failures, the user and machine involved, and the date and time of attack. This kind of attack should generate an alert for both- the machine and the user(s) concerned.

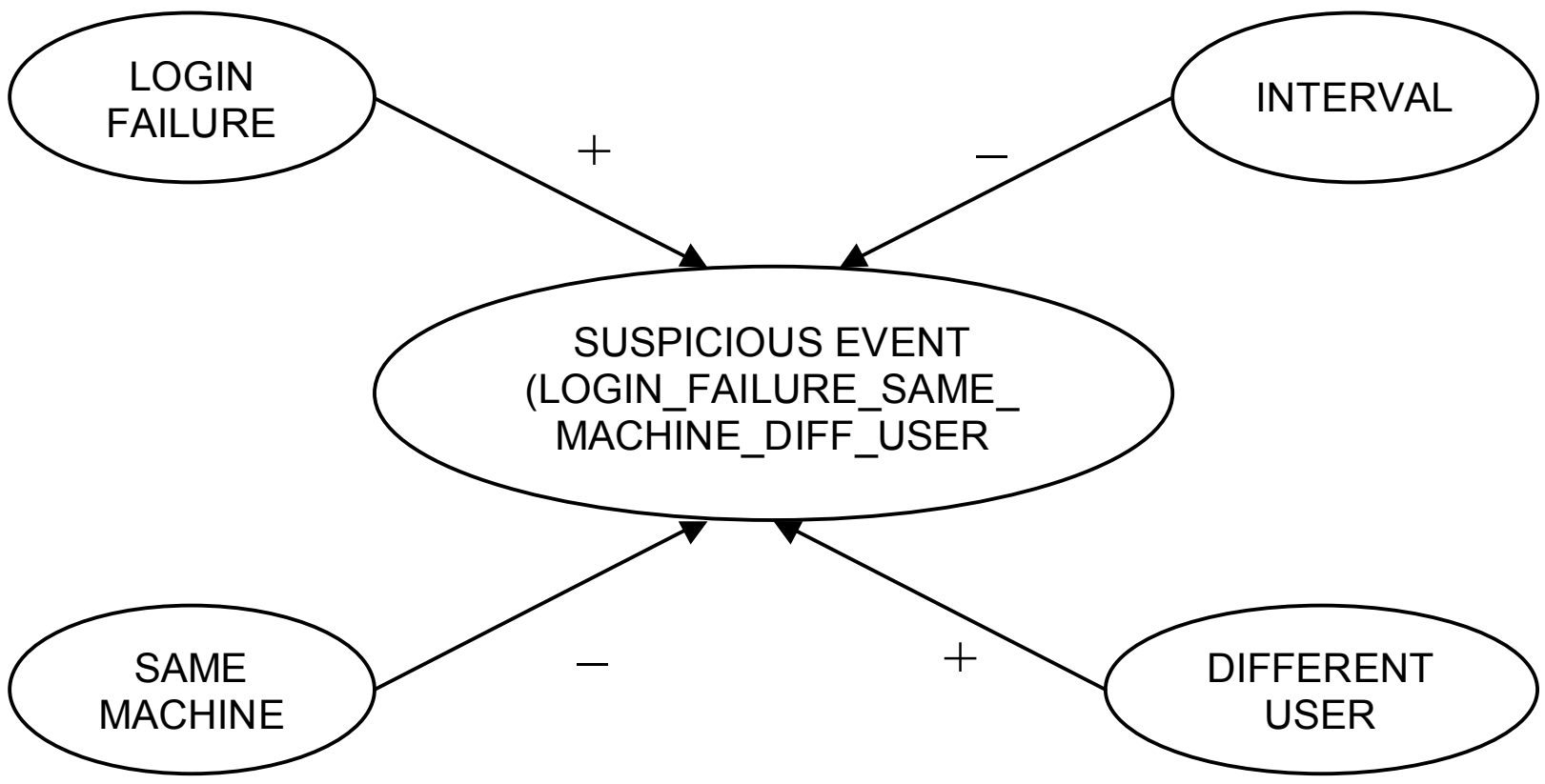

**FIGURE: 1.3.1**

Figure 1.3.1 shows an FCM that captures this scenario. The activation of the concept Login_Same_Machine_Diff_User by the other concepts can be implemented using a fuzzy rule-base. The fuzzy rules are used to map the multiple input concepts (the causes) to the output concept (the effect). Multiple fuzzy rules may be used to correspond to the knowledge described in an individual FCM. The FCM in Figure 1.3.1 can be implemented with a single fuzzy rule:

If number of login failures is moderate and interval is short and this happened for same machine and different users then, Login_Same_Machine_Diff_User is activated highly.

The kind of fuzzy cognitive modeling where fuzzy rules are used to support FCMs has been used for risk assessment in health care *see* Smith and Eloff [95]. Carvalho and Tome [19-23] report that rule-based FCMs are more effective than simple FCMs. Supporting fuzzy rules make FCMs fuzzy compatible and allows qualitative modeling [19, 20].

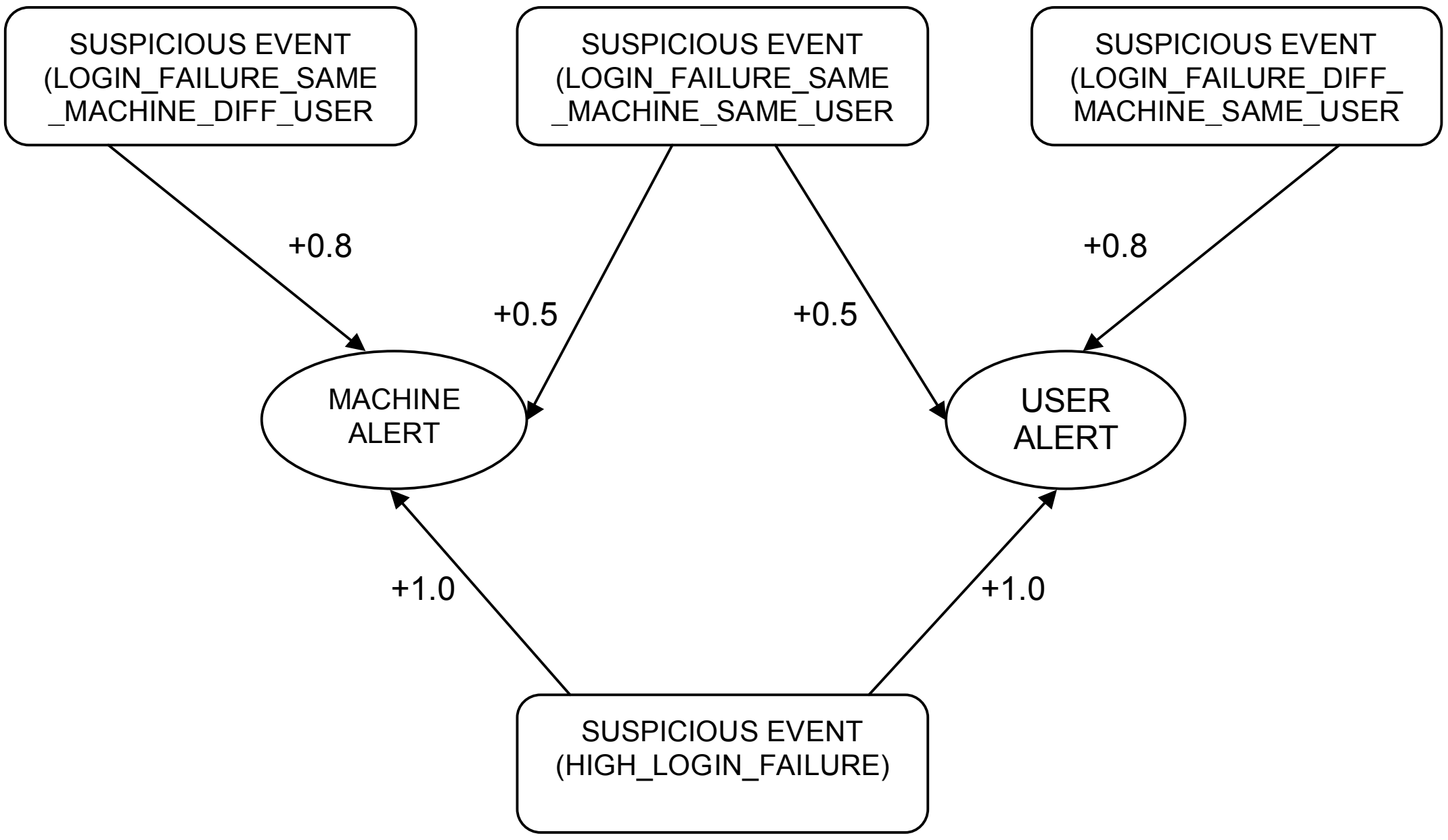

**FIGURE: 1.3.2**

An adjacency matrix is used to list the cause and effect relationships between the nodes. The runtime operation is simply observed by determining the next value of each concept from the current concept and connecting edge values [13, 14] and this can be represented as in [54], for each concept $C_i$ at $t_{n+1}$ time.

The anomaly detection modules measure the dissimilarities between normal and abnormal data and report anomalous behavior in terms of an alarm percentage and a similarity measurement. In order to assess the combined effect of outputs of the anomaly detection modules on the network's anomaly alert level, the decision engine uses causal analysis via FCMs. Figure 1.3.1 shows an example of such a FCM that can be used to infer the anomaly alert level of the network.

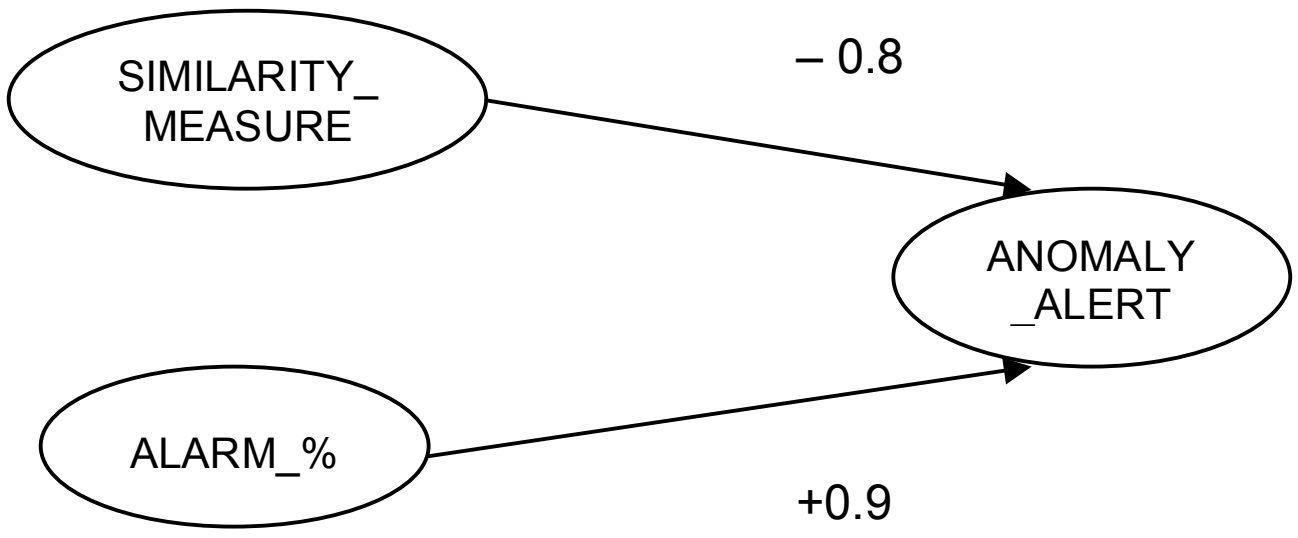

**FIGURE: 1.3.3**

Causal relations can be used to denote a set of fuzzy FCM rules. The FCM in Figure 1.3.3 describes two such causal relations. One denotes that as the similarity measure increases, the anomaly alert level decreases with 0.8 causality. Intuitively this can mean that when the similarity measure is very high, the anomaly alert is very low and, conversely, when the similarity measure is low, the anomaly alert is high.

So far Siraj et al [89] have used FCMs whose structure has been defined by human experts. Experiments are now ongoing that investigate the operational impact and effectiveness of employing this causal knowledge inferencing technique for decision support in the intelligent intrusion detection prototype. After reporting on the effectiveness of using FCMs in their problem domain, they [89] wish to explore how FCMs can be learned and trained like neural networks in this context. For a elaborate and interesting discussion of IIDS and FCMs, the reader is asked to refer [89].

### 1.3.2: A Strategic Planning Simulation Based on FCM Knowledge

K.C. Lee et al [65] in their paper on "*A Strategic Planning Simulation Based on FCM Knowledge and Differential Game*" use the mechanism of integrating fuzzy cognitive map knowledge with a strategic planning simulation where a FCM helps the decision makers understand the complex dynamics between a certain strategic goal and the related environmental factors.

They [65] argue that environments can be classified into three categories: uncontrollable, semi-controllable, and controllable. Environments such as economic conditions regardless of domestic or international interest rates, political stability, etc. are uncontrollable. Quality control activity, sales price, productivity, etc. are controllable because they might become favorable or unfavorable according to the efforts of the company. Semi-controllable environments possess two aspects of both uncontrollable and controllable environments. Examples are competitiveness, brand image, etc. Competitiveness represents a general competitive position of a firm in the target market, which may be determined by complicated interacting forces of various exogenous and /or endogenous factors described so far. Environments facing the new-entering company are assumed in this paper to encompass those three kinds above. According to this definition of environments, those environmental factors are listed as follows:

| | |
|---|---|
| Uncontrollable factors: | Economic conditions, Competitor's advertisements |
| Semi-controllable factors: | Competitiveness, Market demand |
| Controllable factors: | Quality control, Sales price, Productivity. |

Observations: To identify those factors relevant to the goal, a FCM should be built and then analyzed. Based on it K.C. Lee et al suppose that the strategic goal to be pursued is a market share. An appropriate FCM showing the cause-effect relationships between the strategic goal and those seven environmental factors described under 3 heads together with market share is depicted in the Figure 1.3.4.

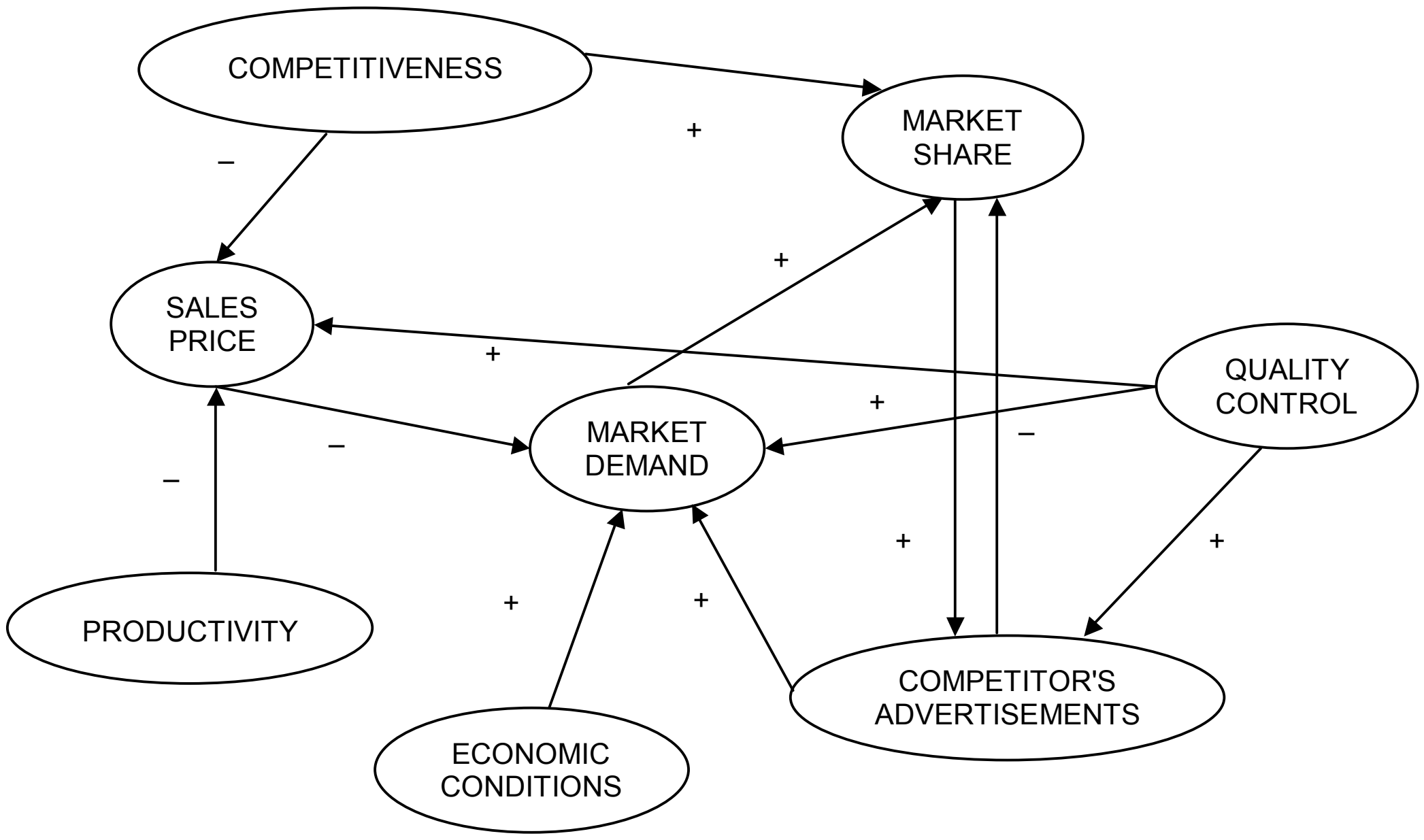

**FIGURE: 1.3.4**

Based on the causal relationships shown in Figure 1.3.4, an initial version of FCM matrix can be built as in Table 1. The description of notations is given below:

• • : Market Share
• • : Competitiveness
• • : Quality Control
• • : Sales Price
• • : Productivity
• • : Market Demand
• • : Economic Conditions
• • : Competitor's Advertisements

+1 : Positive cause-and-effect
0 : No cause-and-effect
–1 : Negative cause-and-effect

**Table 1. Initial version of FCM matrix**

| Cause | Effect | | | | | | | |
|---|---|---|---|---|---|---|---|---|
| | ● ● | ● ● | ● ● | ● ● | ● ● | ● ● | ● ● | ● ● |
| ● ● | | | | | | | | +1 |
| ● ● | +1 | | | | | | | |
| ● ● | | +1 | | | | +1 | | +1 |
| ● ● | –1 | –1 | | | | –1 | | |
| ● ● | | | | –1 | | | | |
| ● ● | +1 | | | | | +1 | | |
| ● ● | | | | | | | | |
| ● ● | –1 | | | | | +1 | | |

After considering the expert's collective opinion, the initial version of FCM matrix can be converted into a following refined version of FCM matrix in which causal values fall between –1 and 1.

**Table 2. Refined version of FCM matrix**

| ●●●●● | ●●●●●● | | | | | | | |
|---|---|---|---|---|---|---|---|---|
| | ● ● | ● ● | ● ● | ● ● | ● ● | ● ● | ● ● | ● ● |
| ● ● | | | | | | | | ● ●● |
| ● ● | ● ●● | | | | | | | |
| ● ● | | ● ●● | | | | ● ● | | ● ● |
| ● ● | ● ●● | ● ●● | | | | ● ●● | | |
| ● ● | | | | ● ●● | | | | |
| ● ● | ● ●● | | | | | | | |
| ● ● | | | | | | ● ●● | | |
| ● ● | ● ● | | | | | ● ●● | | |

From the first column in Table 2, it can easily be concluded that the competitor's advertisements affect the market share of the new-entering company most negatively. It is reasonable for the new-entering company to formulate appropriate strategies considering the competitor's advertisements sufficiently. In addition, the FCM matrix in Table 2 can be used usefully for clearly measuring the composite effects resulting from changes of multiple factors. For example, let us consider two cases:

[Case 1]: Suppose that three factors changed.

Then stimulus input vector may be obtained as follows:

Stimulus Input

V: Productivity = + 0.3
VII: Economic Conditions = – 0.3
VIII: Competitor's advertisements = +.5

This information can be organized into stimulus input vector: (0, 0, 0, 0, 0.3, 0, – 0.3, 0.5). Therefore multiplying this stimulus input vector with FCM matrix, K.C. Lee et al [65] obtain a consequence vector as follows:

(0, 0, 0, 0, 0.3, 0, – 0.3, 0.5) × FCM matrix = (– 0.5, 0, 0.24, 0, – 0.08, 0, 0).

The consequence vector may be interpreted such that changes in those three factors including productivity, economic conditions, and competitor's advertisements affect the market share negatively by .5, the sales price negatively by .24, and the market demand negatively by .08.

Therefore, [65] conclude that changes of three factors may affect the market share most unfavorably.

[Case 2]: If four factors change including sales price (– 0.1), productivity (+ 0.3), economic conditions (– 0.3), and competitor's advertisements (+0.2), stimulus input vector can be organized into (0, 0, 0, – 0.1, 0.3, 0, – 0.3, .2). Then a consequence vector may be obtained as (– 0.15, .02, 0, – 0.24, 0, –0.12, 0, 0). Composite effect of stimulus input vectors on the goal is – 0.15 which means that changes of those factors may affect the market share negatively by 0.15.

In addition to these two cases, if we change the state of competitor's advertisements into 0.7, then the effect is directly affecting the market share by – 0.7 and indirectly increasing the market demand by 0.14. In conclusion, K.C. Lee et al select the competitor's advertisements as a most relevant factor to the stated strategic goal, incorporating into one of variables of the following strategic simulation mechanism.

Since competitor's advertisements is selected as a factor most relevant to the strategic goal, advertising efforts are appropriate strategy for the new-entering company. In specific, strategy concerned is to appropriately allocate advertising behaviour.

K.C. Lee et al [65] found that competitive advertising model formulated in differential games is most suitable for this type of strategic simulation because two simulation entities including the new-entering company and the competitor must be equally considered into the simulation dynamics.

In the performance index given in K.C. Lee et al [65], $u_1$ is squared for analytical convenience to calculate the derivative of the optimal controls, $u^*$ and $v^*$. It is noteworthy that performance index of each company follows a goal hierarchy formulated by [65].

To solve this kind of differential game problems, an adaptive algorithm to take advantage of the competitor's non-optimal play. Let us suppose that the following sets of parameters shown in Table 3 are available to the new-entering company. M is given 500.0, $x_1(t_0) = 40.0$ and $x_2(t_0) = 100.0$, $u_1(t) = 1.00$ and $u_2(t) = 1.00$, and planning period is denoted as $t \in [0.0, 5.0]$, being descretized into integers for analytical simplicity.

After differential game-based strategy simulation with these parameters, [65] maintained $u_1^*$ for each case.

**Table 3. Sets of parameters for strategy formulation**

| Case | a1 | a2 | b1 | b2 | c1 | c2 | w1 | w2 |
|---|---|---|---|---|---|---|---|---|
| • • | .20 | .25 | 1.10 | 1.10 | .60 | .80 | 15.0 | 8.00 |
| • • | .25 | .25 | 2.00 | 1.50 | .50 | .60 | 22.5 | 8.00 |
| III | .25 | .25 | 2.00 | 1.50 | .50 | .60 | 25.0 | 5.00 |
| • • | .25 | .25 | 1.10 | 1.10 | .50 | .60 | 2.5 | 8.00 |
| • • | .20 | .25 | 1.10 | 1.10 | .50 | .60 | 15.0 | 8.00 |
| • • | .15 | .25 | 1.10 | 1.10 | .50 | .60 | 15.0 | 8.00 |
| • • | .25 | .40 | 1.10 | 1.10 | .50 | .60 | 15.0 | 8.00 |

Strategy evaluation is based on the criterion: "Will the proposed strategy be any good in the future?" For strategy evaluation, K.C. Lee et al [65] expected performance index JI for each case after strategic planning period, which is summarized in the following Table 4.

**Table 4. Performance value JI for each**

| Case | Performance Value (after strategic planning period) |
|---|---|
| I | 180.86 |
| II | 157.90 |
| III | 164.08 |
| IV | 170.56 |
| V | 152.77 |
| VI | 133.70 |
| VII | 174.67 |

The value of performance index in Table 4 are obtained from summing two terms: profits during strategic planning periods and relative market share at the terminal planning time. From the results shown in Table 4, [65] conclude that case I will result in the most promising strategy for the future. Then the following strategy (i.e., $u_1^*$) related to case I is implemented for the new-entering company.

**Table 5. Proposed strategy $u_1^*$ for the new-entering company**

| Planning Period | Strategy |
|---|---|
| 1 | .4118 |
| 2 | .1581 |
| 3 | .1912 |
| 4 | .1518 |
| 5 | .0808 |
| Sum | 1.0000 |

Values shown in Table 5 represent the normalized version of $u_1^*$ and then detailed strategy implementation may be accomplished such that if advertising budget is A, advertising expenditures to be invested at each planning period are obtained by multiplying A with the corresponding normalized value.

After strategy implementation, the expected sales amount ($x_1^*$) for the new-entering company is summarized in Table 6.

**Table 6. Expected sales amount**

| Time | Expected Sales Amount |
|---|---|
| 1 | 84.70 |
| 2 | 69.78 |
| 3 | 57.34 |
| 4 | 47.26 |
| 5 | 38.99 |

However, in a real situation, the expected sales amount shown in Table 6 cannot be fully accomplished due to a variety of unexpected changes in business environments. Then assume that deviations occur in the new-entering company's sales amount $x_1^*$ and also deviations are observed in the competitor's sales amount $x_2^*$, which are represented in Table 7.

**Table 7. Deviations observed in the sales amount**

| Planning Period | $x_1^*$ | $x_2^*$ |
|---|---|---|
| 1 | –1.0 | –2.0 |
| 2 | –0.5 | –1.5 |
| 3 | 0.0 | –1.0 |
| 4 | 0.5 | –0.5 |

Deviations shown in Table 7 will cause recycling, resulting in both adjusted strategies and sales amount for the new-entering company, which is summarized in Table 8.

**Table 8. Adjusted strategy and sales amount**

| Planning Period | Adjusted | |
|---|---|---|
| | Strategy | Sales Amount |
| 1 | .4095 | 84.7 |
| 2 | .1604 | 69.70 |
| 3 | .1918 | 103.00 |
| 4 | .1578 | 85.40 |
| 5 | .0805 | 69.70 |

Performance value adjusted after considering deviations in the competitor's sales amount are 238.72 which is greater than the original value 180.86. Also an increase in the sales amount is observed because [65] use the adaptive algorithm to take advantage of the competitor's non-optimal play. In addition, terminal market share increased from 7.80% to 13.95%.

Lee et al propose the FCM knowledge-based strategic planning simulation mechanism in which

1. The FCM-based causal knowledge provides overall interrelationships between strategic goals and environmental factors which are thought to be affecting the performance of strategic planning.

2. Strategic planning simulation system driven by differential game adapt to the changes in surrounding environments.

Main contributions of this paper to simulation are as follows:

1. Unlike the previous simulation approaches to strategic planning in which expert systems or decision support systems or neural networks were used to suggest business strategies, their approach applied FCM-based simulation analysis for environmental analysis. With this, [65] can identify and consider the most relevant environmental factor which seems to affect the expected target variable.

2. Those environmental factors identified through FCM-based simulation analysis can be incorporated into the strategic planning simulation framework with the aid of differential game mechanism.

3. Time variable, which has been neglected in the literature about conventional strategic planning simulation, can be effectively taken into consideration by the differential game-driven mechanism.

4. Synergistic framework for integrating FCM analysis and differential game is proposed for more effective strategic planning simulation under turbulent situation.

Experiments with a set of illustrative examples showed that approach of [65] can yield robust strategic planning simulation results which can vary according to the changes in environment. Future research topics are as follows: (1) how can we obtain a unified knowledge about environment when experts opinion differs, (2) what if we use more refined FCM in which edge values are defined more rigorously such as fuzzy partial relationship. Though we have tried to be very comprehensive and presented a lot of the research of K.C.Lee et al, the interested reader can refer [65].

### 1.3.3: Adaptive fuzzy cognitive maps for hyperknowledge representation in strategy formation process

Carlsson and Fuller [17] have shown that the effectiveness and usefulness of this hyperknowledge support system can be further advanced using Adaptive Fuzzy Cognitive Maps.

Strategic Management is defined as a system of action programs which form sustainable competitive advantages for a corporation, its divisions and its business units in a strategic planning period. A research team of the IAMSR institute has developed a support system for strategic management, called the Woodstrat, in two major Finnish forest industry corporations in 1992-96.

The system is modular and is built around the actual business logic strategic management in the two corporations, i.e. the main modules cover the market position (MP), the competitive position (CP), the productivity position (PROD), the profitability (PROF), the investments (INV) and the financing of investments (FIN). The innovation in Woodstrat is that these modules are linked together in a hyperknowledge fashion, i.e. when a strong market position is built in some market segment it will have an immediate impact on profitability through links running from key assumptions on expected developments to the projected income statement. There are similar links making the competitive position interact with the market position, and productivity position interact with both the market and the competitive positions, and with the profitability and financing positions.

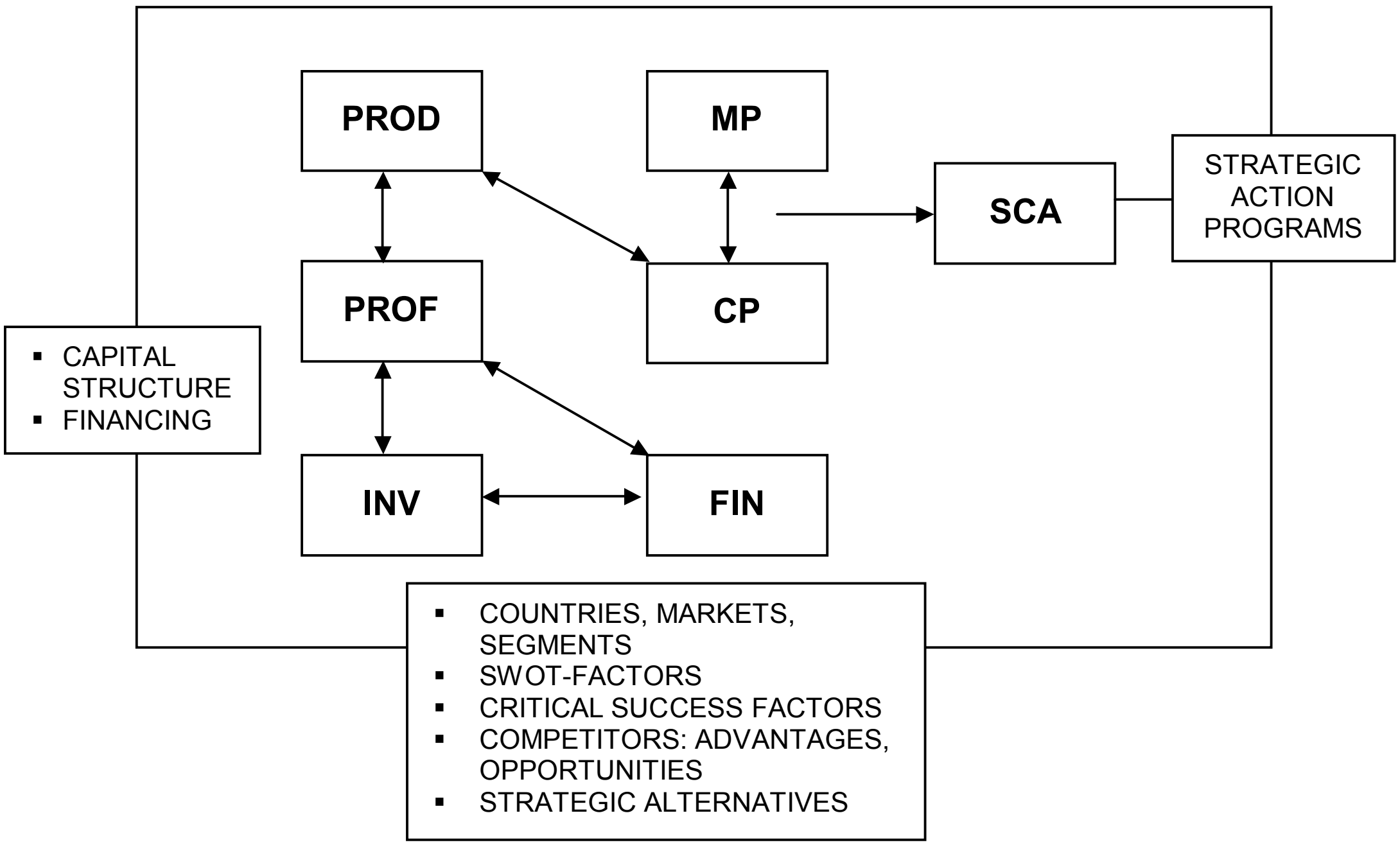

**FIGURE: 1.3.5**

The Woodstrat offers an intuitive and effective strategic planning support with object-oriented expert systems elements and a hyperknowledge user interface. In this paper Carlsson and Fuller [17] show that effectiveness and usefulness of a hyperknowledge support system can be further advanced using Adaptive Fuzzy Cognitive Maps.

When addressing strategic issues cognitive maps are used as action-oriented representations of the context the managers are discussing.

For the sake of simplicity, in their paper, Carlsson and Fuller [17] illustrate the strategy building process by the following fuzzy cognitive map with six states.

The causal connections between the states MP (Market position), CP (Competitive position), PROF (Profitability), FIN (Financing position), PROD (Productivity position) and INV (Investments) are derived from the opinions of managers' of different Strategic Business Units.

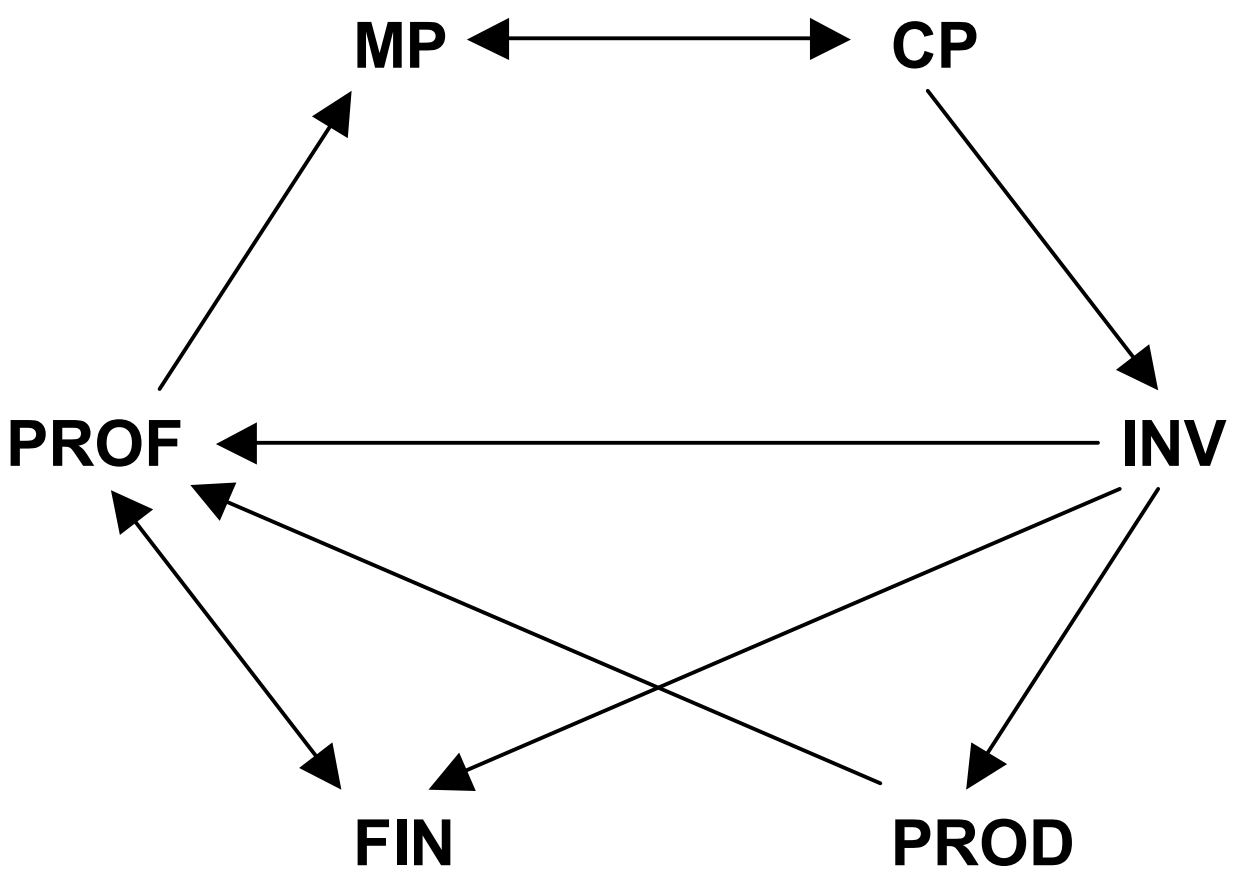

FIGURE: 1.3.6

Adaptive fuzzy cognitive maps can learn the weights from historical data. Once the FCM is trained it lets us play what-if games (e.g. What is demand goes up and prices remain stable? – i.e. we improve out MP) and can predict the future.

In the following Carlsson and Fuller describe a learning mechanism for the FCM of the strategy building process,

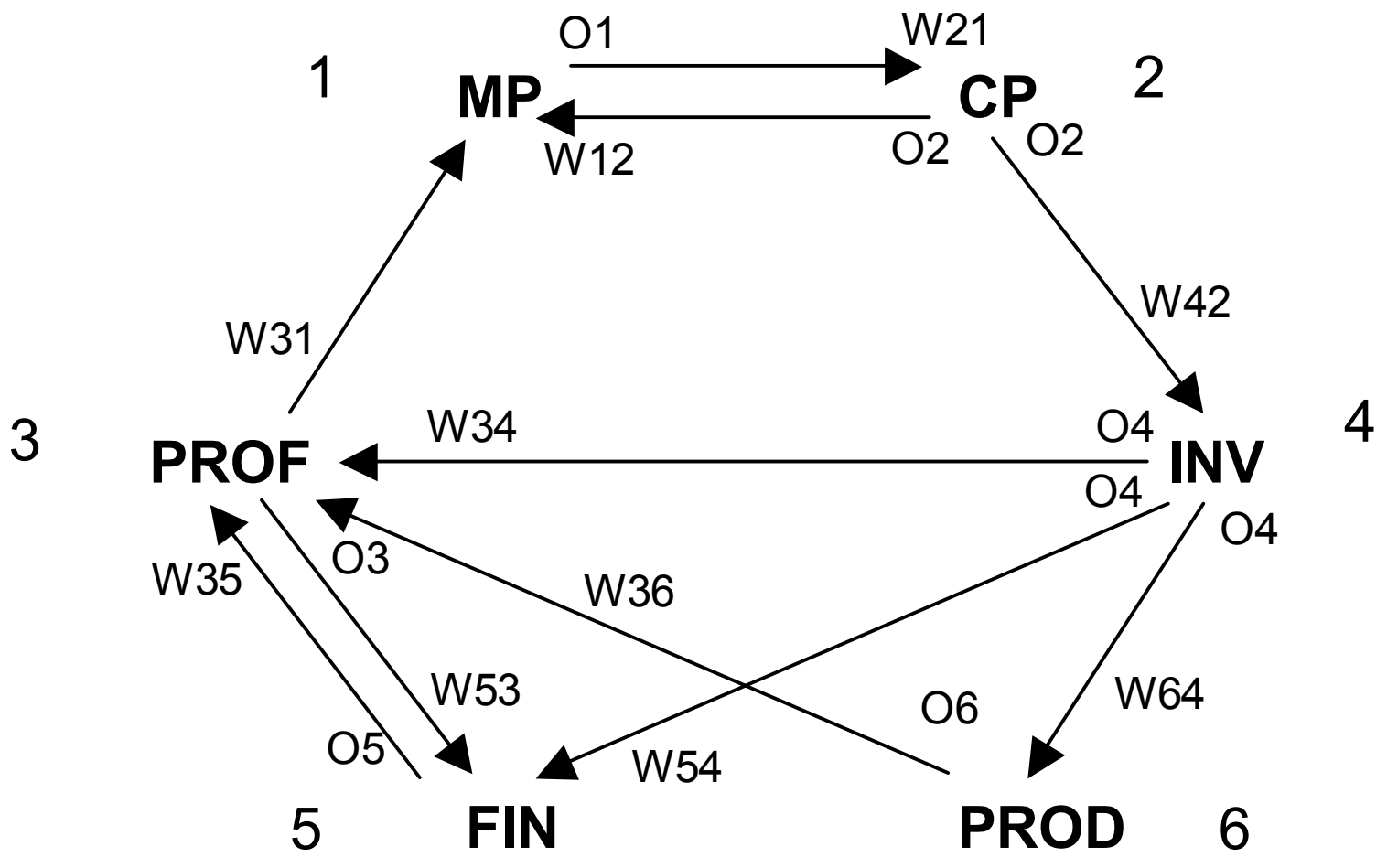

FIGURE: 1.3.7

In the following, Carlsson and Fuller [37] describe a learning mechanism for the FCM of the strategy building process, and illustrate the effectiveness of the map by a simple training set. Figure 1.3.6 shows the structure of the FCM of the strategy building process.

Inputs of states are computed as the weighted sum of the outputs of its causing states.

$$net = Wo$$

where W denotes the matrix of weights, is the vector of computed outputs, and net is the vector of inputs to the inputs to the states. In this case the weight matrix is given by

$$W = \begin{pmatrix} 0 & w_{12} & 0 & 0 & 0 & 0 \\ w_{21} & 0 & 0 & 0 & 0 & 0 \\ w_{31} & 0 & 0 & w_{34} & w_{35} & w_{36} \\ 0 & w_{42} & 0 & 0 & 0 & 0 \\ 0 & 0 & w_{53} & w_{54} & 0 & 0 \\ 0 & 0 & 0 & w_{64} & 0 & 0 \end{pmatrix}$$

where the zero elements denote no causal link between the states, and

$$net = \begin{bmatrix} net_1 \\ net_2 \\ net_3 \\ net_4 \\ net_5 \\ net_6 \end{bmatrix} = \begin{bmatrix} net(MP) \\ net(CP) \\ net(PROF) \\ net(INV) \\ net(FIN) \\ net(PROD) \end{bmatrix}, o = \begin{bmatrix} o_1 \\ o_2 \\ o_3 \\ o_4 \\ o_5 \\ o_6 \end{bmatrix} = \begin{bmatrix} o(MP) \\ o(CP) \\ o(PROF) \\ o(INV) \\ o(FIN) \\ o(PROD) \end{bmatrix}.$$

That is,

$$\begin{aligned} net_1 &= net\ (MP) &&= w_{12}o_2, \\ net_2 &= net\ (CP) &&= w_{21}\ o_1, \\ net_3 &= net\ (PROF) &&= w_{31}o_1 + w_{34}o_4 + w_{35}o_5 + w_{36}o_6, \\ net_4 &= net\ (INV) &&= w_{42}o_2, \\ net_5 &= net\ (FIN) &&= w_{54}o_4 + w_{53}o_3, \\ net_6 &= net(PROD) &&= w_{64}o_4 \end{aligned}$$

Carlsson and Fuller [17] have derived the extensions of Woodstrat to a fuzzy hyperknowledge support system that has main effect on the support approximate reasoning schemes in linking MP, CP etc. For an in depth analysis, please refer [17].

### 1.3.4: A New Balanced Degree for FCM

In this section we discuss the concept of balanced degrees for FCMs and the significant results that have been obtained by Tsadiras and Margaritis [107] in connection with it. An FCM concerning the public health with 7 nodes $C_1, C_2, \ldots, C_7$ is given, for example, concept $C_1$ (Number of people in the city) in Figure 1.3.8 is connected to concept $C_7$ (Bacteria) by two paths. The first is the $C_1$ (+)→ $C_4$ (+)→ $C_7$ and indirect causal relationship. The second path is $C_1$ (+)→ $C_3$ (+)→ $C_5$ (+)→ $C_7$ and indicates indirect causal relationship. In this case we cannot conclude if the total effect of concept $C_1$ to $C_7$ is positive or negative, or in other words if the increase of the population of the city will increase or decrease the bacteria.

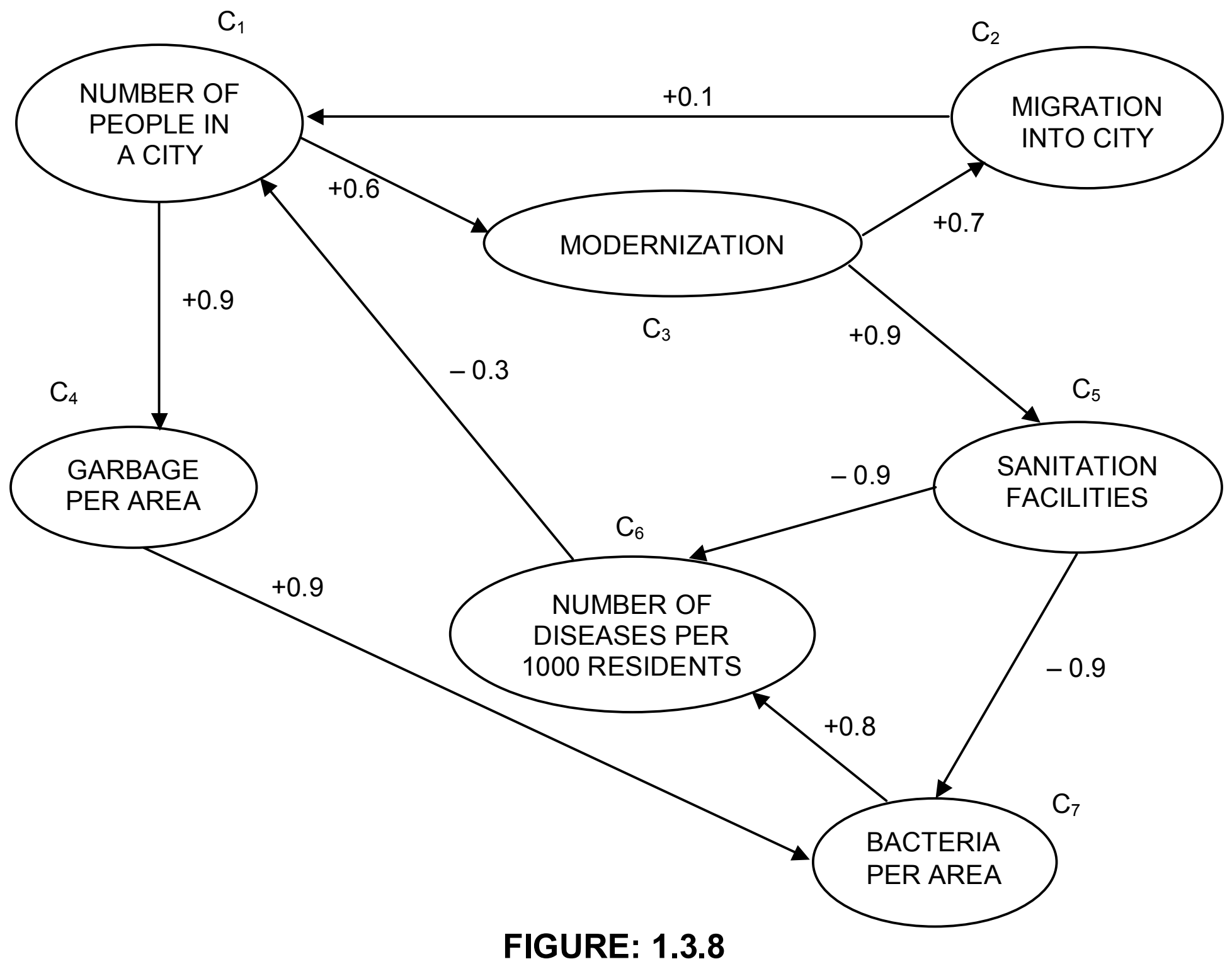

**FIGURE: 1.3.8**

An FCM is imbalanced if we can find two paths between the same two nodes that create causal relations of different sign. In the opposite case the FCM is balanced. The term "balanced" digraph is borrowed from Graph Theory. In an imbalanced FCM we cannot determine the sign of the total effect of a concept to another. Using the idea that as the length of the path increases, the indirect casual relation becomes weak, Eden [33-35] proposed that the total effect should have the sign of the shortest path between the two nodes. Tsadiras and Margaritis [107] on the other hand, proposed that the sign of the total effect should be the sign of the most important path where the most important path is the one that passes through the most important nodes.

The degree to which the digraph of the FCM is balanced or imbalanced is given by the Balance Degree of the digraph. Various types of Balance Degrees have been proposed, each of them being suitable for certain types of problems. Tsadiras and Margaritis [107] proposed the following balance degree $\bar{a}$.

$$\bar{a} = p / t \quad (1)$$

where p is the number of positive semicycles of the digraph and t is the total number of semicycles of the digraph (semicycle is a cycle in the digraph which can be created by not taking into consideration the direction of the arcs). The closer to 1 is the degree $\bar{a}$, the closer to balance is the digraph

$$\bar{a} = \frac{\sum_m p_m f(m)}{\sum_m t_m f(m)} \text{ or} \quad (2)$$

$$\bar{a} = \frac{\sum_m p_m f(m)}{\sum_m n_m f(m)} \tag{3}$$

where $p_m$ is the number of positive semicycles in the graph that have length m, $n_m$ is the number of negative semicycles in the graph that have length m, $t_m$ is the total number of semicycles in graph that have length m and f (m) a monotonously increasing function, for example f(m) = 1/ m, f(m) = 1 /$m^2$ or f(m) = 1/$2^m$.

[107] introduced the new Balance Degree r that checks strictly the balance of the digraph and is very useful for FCMs. This degree can take a value among the whole interval [0, 1]. The closer r is to 0 the closer to a completely balanced digraph is the graph. To calculate this new degree we check the signs of all the indirect causal relationships that exists among each pair of concepts $C_i$ and $C_j$. Let us assume that there exist $p_{ij}$ positive paths from concepts $C_i$ to concept $C_j$ and that $n_{ij}$ is the number of negative paths from concepts $C_i$ to concept $C_j$. Now there are two cases.

Case 1: $p_{ij} + n_{ij} = 2k$.

In this case the total number of paths between the two concepts is even. If i = j then we refer to cycles. The case of a completely balanced digraph (r = 0) occurs when we have either $p_{ij} = 0$ and $n_{ij} = 2k$ or $p_{ij} = 2k$ and $n_{ij} = 0$. The case of a completely imbalanced FCM (r = 1) occurs when $p_{ij} = n_{ij} = k$. The fraction

$$\frac{\min\{p_{ij} n_{ij}\}}{k}$$

gives a measurement of the balance according the above. If we have k + x, x > 0 positive (negative) paths and k – x negative (positive) paths then the above fraction gives (k – x) /k which is right measurement for the balance of the FCM.

Case 2: $p_{ij} + n_{ij} = 2k + 1$.

In this case the total number of paths between the two concepts is odd. The case of a completely balanced digraph (r = 0) occurs when we have either $p_{ij} = 0$ and $n_{ij} = 2k + 1$ or $p_{ij} = 2k + 1$ and $n_{ij} = 0$. The case of a completely imbalanced FCM (r = 1) occurs when $p_{ij} = k$ and $n_{ij} = k + 1$ or when $p_{ij} = k + 1$ and $n_{ij} = k$. The fraction

$$\frac{\min\{p_{ij} n_{ij}\}}{k+1}$$

gives a measurement of the balance according the above. If we have k – x, x > 0 positive (negative) paths and k + x + 1 negative (positive) paths then the above fraction gives (k – x) / (k + 1) which is right measurement for the balance of the FCM.

The degrees that the two cases above propose are different to the denominator. We can unify the two cases into one using the following fraction:

$$\frac{\min\{p_{ij} n_{ij}\}}{\operatorname{int}\left[\frac{p_{ij}+n_{ij}+1}{2}\right]} \quad (4)$$

where function int() returns the integer part of a decimal.

This new Balance Degree is useful because it provides us with a clear indication of the number of conflicts that exist in the digraph of the FCM. By that we can draw conclusions about the dynamical behaviour that we should expect from the FCM. This is important because the forecasting process of the FCMs requires the study of the dynamical behaviour of the FCM model. We have tried our level best to introduce all definitions and equations so that the reader can develop these concepts to Neutrosophic Cognitive Maps and apply them to several other problems.

### 1.3.5: Rule based FCM

FCMs are Causal Maps (a subset of Cognitive Maps that only allow causal relations) and in most applications, a FCM is indeed a man-trained Neural Network that is not fuzzy in a traditional sense and does not explore usual Fuzzy capacities; Carvalho and Tome [19-23]. They do not share the properties of other fuzzy systems and the causal maps end up being quantitative matrixes without any qualitative knowledge [21]. To avoid these limitations of existing approaches, Rule Based Fuzzy Cognitive Maps (RB-FCM) were introduced by Carvalho and Tome [19-23] and are being developed as a tool that models/ simulates real-world qualitative system dynamics.

One very important issue that has apparently been ignored (or avoided) when approaching Qualitative System Dynamics is <u>time</u>. Time is obviously essential in the study of System Dynamics (SD). When we are dealing with quantitative SD timing issues are naturally solved, since Time is <u>explicitly</u> expressed in the mathematical equations used to describe the relations between concepts. However, in qualitative SD, where natural language is more adequate to express the relations between concepts, time must be implicitly expressed in those relations (in the case of RBFCM, time must be somehow included in the fuzzy rule bases that define the relations). Without imposing this implicit timing knowledge it is impossible to guarantee a minimally acceptable simulation.

RB-FCM provide a representation of the dynamics of complex real-world qualitative systems with feedback and allow the simulation of the occurrence of events and their influence in the system. They are fuzzy directed graphs with feedback, which are composed of fuzzy nodes (Concepts), and fuzzy links (Relations). RB-FCM are true Cognitive Maps (CM) since are not limited to the representation of causal relations. Unlike FCM, Concepts are fuzzy variable defined by fuzzy linguistic variables, and fuzzy rule bases are used to express Relations. RB-FCM are essentially fuzzy rule based systems where we added fuzzy mechanisms to deal with feedback, and different kinds of relations in order to cope with the complexity and diversity of the qualitative systems we are trying to model.

RB-FCM are iterative: the current value of each concept is computed with its inputs previous values. The evolution of the system through time might reach equilibrium and converge to a single state or a cycle of states under certain conditions [19-23].

Introduction or removal of concepts and /or relations, or the change of state in one or more concepts affect the modeled systems in ways that are usually difficult or impossible to predict due to the complex feedback links. RB-FCM are a tool to predict the evolution through time caused by those changes. RB-FCM provide 2 kinds of Concepts (Levels and Variations), several kinds of Relations (Causal, Inference Common, similarity, ill-defined, Crisp, Level-Variation, and other common relations), and mechanisms that support invariant and time-variant probabilities, the possibility of subsystems including decision support systems to simulate the process of decision making by "actors" within the system [19-23].

With the introduction of mechanisms presented in this paper, which allow the modeling of delays and inhibition of certain relations when they have no influence on a given instant, RB-FCM can deal with timing issues and become a tool to represent and analyse the dynamics of qualitative systems.

RB-FCM Syntax is a language that was developed to describe real-world qualitative system in RB-FCM [21, 22]. It allows the definition of the concepts, relations, linguistic variables, membership functions, etc, that compose the system.

As time in RB-FCM must inevitably be represented implicitly in every relation; therefore, the responsibility of maintaining temporal coherence in the process of modeling a system relies heavily in the modeler. He or she must not only be responsible to find the nature and characteristics of the relation, but must also ensure that the magnitude of the relation is adequate to the time interval that it represents.

Due to the iterative nature of the RB-FCM simulation process, each system iteration must represent the flow of a given amount of time. In RB-FCM this period is called Base Time (B-Time). It represents the "resolution" of the simulation, i.e., the highest level of temporal detail that a simulation can provide in the modeled system.

B-Time must always be implicit while defining each rule in each rule base, especially in the case of causal relations, since Variations are always involved in those relations and there is an intrinsic relation between B-Time and Variations: the linguistic variable indicating the highest possible amount of change in each Variation (Increase_Very_Much, Decrease_Very_Much, Gets_Incredibly_Better, etc.), represents the largest amount of change that the physical entity modeled by the Variation is expected to suffer during the period of time that B-time models.

The choice of B-Time is highly dependent on the Real-world system being modeled. It could be one hour, one day, two days, one week, one year, one century or any other time period. It depends on the desired or advisable level of detail, complexity and intended long-term analysis of the system. Shorter B-Times usually need more detailed and complex rule bases and imply a much more careful approach to the precision and validity of the rules. Longer B-Times should provide more valid long-term simulations, but short-term detail, precision and validity will possibly be sacrificed.

The intrinsic nature of the system is also an important factor for B-Time selection. Chaotic or pseudo-chaotic systems obviously need shorter B-Times (1 hour, 1 day), since small deviations in the first steps of the simulation will cause huge errors in later steps. On the other hand, stable systems that tend to converge to a single state or single cycle of states usually handle very well long B-Times (1 month, 1 year), which should be fully utilized to allow the implementation of less detailed rule bases and provide faster simulations.

Delays are the most obvious timing issue, and the only one that has been addressed on FCM on some occasions. Sometimes the effect of a relation is not immediate (or at least can't felt during B-Time), and the delay must somehow be modeled on the RB-FCM. It is important to note that this is a different issue from the one mentioned in the previous section. We are not trying to model an event that produces a small discernible effect during a single iteration, but an event which has a real delay. A common example of this kind of timing events is the effect of oil-price variations on the cost of fuel or electricity: due to the fact that most countries have long term reserves, and the process of oil transport and refining is slow, then the people only feel those effects after a few months.

For more about RB-FCM refer [19-23]. We conclude that Rule Based Fuzzy Cognitive Maps are the best tools to carry out qualitative analysis.

### 1.3.6: Fuzzy Causal Relationship and Fuzzy Partially Causal Relationship

This section recalls the notions of fuzzy causal relationship and fuzzy partially causal relationship as given by Kim H.S. and K. C. Lee [50].

"A causally increases B" means that if A increases then B increases, and if A decreases then B also decreases. On the other hand, "A causally decreases B" means that if A increases then B decreases and if A decreases then B increases. So, in the concepts that constitute causal relationship, there must exist quantitative elements that can increase or decrease. Kosko [58] defined the concept $C_i$ that constitutes causal relationships in FCM as follows: $C_i = (Q_i \cup {}^{\sim}Q_i) \cap M_i$, where $Q_i$ is a quantity fuzzy set and ${}^{\sim}Q_i$ is a dis-quantity fuzzy set. ${}^{\sim}Q_i$ is the negation of $Q_i$. $M_i$ is a modifier fuzzy set that modifies $Q_i$ or ${}^{\sim}Q_i$ concretely. Each $Q_i$ and ${}^{\sim}Q_i$ partitions the whole set $C_i$. Double negation ${}^{\sim\sim}Q_i$ is equal to $Q_i$ implying that ${}^{\sim}Q_i$ is corresponding to $Q_i^c$ the complement of $Q_i$. However, negation does not mean antonym. For example, assume that $Q_i$ is "tall" and ${}^{\sim}Q_i$ is "short" in height. The complement of fuzzy set "tall" does not correspond to the fuzzy set "short". That is, in verbal representation, "not tall" does not necessarily mean "short". Therefore, if a dis-quantity fuzzy set ${}^{\sim}Q_i$ does not correspond to the complement of $Q_i$, we will call it as anti-quantity fuzzy set to clarify the subtle meaning in the dis-quantity fuzzy set.

There should be a modifier fuzzy set in the concept to modify the quantity fuzzy set and the dis-quantity fuzzy set. For example, a concept cannot be figured out if there are only quantitative sets such as "stability" or "instability".

**DEFINITION [50]:** *$C_i$ causes $C_j$ if and only if $(Q_i \cap M_i) \subset (Q_j \cap M_j)$ and $(\tilde{\ }Q_i \cap M_j) \subset (\tilde{\ }Q_j \cap M_j)$.*

**DEFINITION [50]:** *$C_i$ Causally decreases $C_j$ if and only if $(Q_i \cap M_i) \subset (\tilde{\ }Q_j \cap M_j)$ and $(\tilde{\ }Q_i \cap M_j) \subset (\tilde{\ }Q_j \cap M_j)$.*

Here "$\subset$" stands for fuzzy set inclusion (logical implication). For brevity, fuzzy causal relationship is abbreviated as FCR in the following.

Logical implication has an antecedent (or premise) and a consequent (or conclusion). A typical example of the logical implication is "If A, then B". In the logical implication, there are two important inference rules, namely, modus ponens and modus tollens. Modus ponens is defined as the following inference procedure:

| | |
|---|---|
| premise 1: | A |
| premise 2: | If A, Then B |
| consequence: | B. |

Modus tollens is defined as the following inference procedure:

| | |
|---|---|
| premise 1: | negation of B (or $\neg$B) |
| promise 2: | If A, Then B |
| consequence: | negation of A (or $\neg$A). |

The FCR is more complicated than the logical implication.

1. Modus tollens does not hold in the FCR. "When migration into city increases, then amount of garbage per area also increases". However, from this FCR we cannot conclude that if amount of garbage decreases then migration into city also decreases. (Here neutrosophic causality relationship is applicable for relations between several concepts are indeterminate).

2. In a logical implication "If A, Then B", we cannot draw inference in case when A is not true, i.e., when A is negated. However, in a FCR "A causes B", we can draw inference even if A is negated. For example, returning to the previous FCR, if migration into city decreases, then the amount of garbage also decreases.

Now, let us explore several characteristics of the FCR. At first, assume that there exist the following (FCR–1). We abbreviate the strength of FCR for the sake of explanatory convenience. Characteristics of the FCRs with real-valued causality strength will be discussed in the section:

Buying by institute investors in the stock market example [50].

$\xrightarrow[+]{}$ Increase of composite stock price, (FCR–1)

where the quantity fuzzy sets are "buying" and "increase". The modifier fuzzy sets are "institute in-vestors" and "composite stock price". By definition, the (FCR–1) is equivalent to the following:

Not buying by institute investors

$\xrightarrow[+]{}$ Not increase of composite stock price. (FCR–2)

Here "not buying" and "not increase" are the dis-quantity fuzzy sets. It is noteworthy that the negation term "not" is used in the (FCR-2). Let us represent (FCR-2) by anti-quanity fuzzy sets. Then the following holds:

Selling by institute investors

$\xrightarrow[+]{}$ Decrease of composite stock price. (FCR-3)

Here "selling and decrease" are the anti-quantity fuzzy sets. In general, it seems that (FCR-2) is equivalent to (FCR-3). However, rigorous distinction is necessary because the (FCR-2) in terms of dis-quantity fuzzy sets is not equivalent to (FCR-3) in terms of anti-quantity fuzzy sets. Now, it is necessary to define clearly what both quantity fuzzy sets and the dis-quantity fuzzy sets mean in case of formulating a FCR. From now on, we assume that "buying"-"selling" and "increase"-"decrease" are the corresponding pairs of quantity fuzzy sets ($Q_i$) and dis-quantity fuzzy sets( $\tilde{Q}_i$). In other words, we assume that (FCR-2) is equivalent to (FCR-3) in stock market given in example [50].

In FCM, it is not necessary that both equivalent (FCR–1) and (FCR-3) appear. Only the following form of FCR is sufficient:

Institute investors

$\xrightarrow[+]{}$ Composite stock price (FCR-4)

In (FCR-4), the quantity fuzzy sets or the dis-quantity fuzzy sets do not appear. Only the modifier is used. However, its meaning can be represented in the mathematical expression as $(Q_i \cup \tilde{Q}_i) \cap M_i$ which contains the quantity and dis-quantity fuzzy sets. (FCR-4) is, in fact, equivalent to both (FCR–1) and (FCR-3) . This is one of the important characteristics of FCM, which differs from the other knowledge representations.

(FCR–1) can be represented equivalently by using the negative causality form such as (FCR-5) where the quantity fuzzy set "increase" is replaced with the dis-quantity fuzzy set "decrease".

Buying by institute investors

$\xrightarrow[-]{}$ Decrease of composite stock price. (FCR-5).

(FCR-5) means that if buying by institute investors occurs, then the opposite situation of "decrease of composite stock price" takes place. The following also holds by the definition of negative causality;

Selling by institute investors

$\xrightarrow[-]{}$ Increase of composite stock price. (FCR-6)

This means that if selling by institute investors occurs, then the opposite situation of "increase of composite stock price" takes place, which is also equivalent to (FCR-3). Therefore, the FCRs such as (FCR–1), (FCR-3), (FCR-5), and (FCR-6) are the same with each other. From the discussion so far, the following theorem is suggested by Kim H.S. and Lee K.C. and for more elaborate information refer [50].

**THEOREM [50].** *When a concept $C_i$ is $(Q_i \cap M_i)$ and the negative concept $^{\sim}C_i$ is $(^{\sim}Q_i \cap M_i)$, the following FCRs are all equivalent:*

$$C_i \xrightarrow[+]{} C_j,\ {}^{\sim}C_i \xrightarrow[+]{} C_j,\ C_i \xrightarrow[-]{} {}^{\sim}C_j,\ {}^{\sim}C_i \xrightarrow[-]{} C_j$$

For a more detailed discussion please refer [50].

Now we proceed on to recall characteristics of the FCR with a real valued causality strength.

**DEFINITION [50]:** *$\{(a, b) \ / a \in A, b \in B\}^c = \{(c, d) \ / c \in A - \{a\}$ or $d \in B - \{b\}\}$.*

For example, if set A = $\{a_1, a_2\}$ and set B = $\{b_1, b_2\}$, then A × B = $\{(a_1, b_1), (a_1, b_2), (a_2, b_1), (a_2, b_2)\}$. Let a relation R = $\{(a_1, b_1)\}$. Then $R^c = \{(a_1, b_2), (a_2, b_1), (a_2, b_2)\}$. As defined, $C_i = (Q_i \cup {}^{\sim}Q_i) \cap M_i$. Then $(Q_i \cap M_i)^c = ({}^{\sim}Q_i \cap M_i)$ if $M_i$ is a whole set. Therefore, the complement of a (general) relation is defined as follows:

**DEFINITION [50]:** *When $C_i = (Q_i \cup {}^{\sim}Q_i) \cap M_i$. and $C_j = (Q_j \cup {}^{\sim}Q_j) \cap M_j$, the complement of a relation $\{((Q_i \cap M_i), (Q_j \cap M_j))\}$ is $\{(({}^{\sim}Q_i \cap M_i), (Q_j \cap M_j)), ((Q_i \cap M_i), ({}^{\sim}Q_j \cap M_j), ({}^{\sim}Q_i \cap M_i), ({}^{\sim}Q_j \cap M_j))\}$.*

However, if $((Q_i \cap M_i), (Q_j \cap M_j))$ is causal relation (i.e. not a general relation), then $({}^{\sim}Q_i \cap M_i), ({}^{\sim}Q_j \cap M_j))$ is equivalent to $((Q_i \cap M_i), (Q_j \cap M_j))$ by defintion of FCR. Likewise $(({}^{\sim}Q_i \cap M_i), (Q_j \cap M_j))$ is also equivalent to $((Q_i \cap M_i), ({}^{\sim}Q_j \cap M_j))$. Therefore, the complement of a causal relation is defined as follows.

**DEFINITION [50]:** *When $C_i = (Q_i \cup {}^{\sim}Q_i) \cap M_i$ and $C_j = (Q_j \cup {}^{\sim}Q_j) \cap M_j$, the complement of a causal relation $\{((Q_i \cap M_i), (Q_j \cap M_j))\}$ is $\{((Q_i \cap M_i), ({}^{\sim}Q_j \cap M_j))\}$.*

Let a causal relation R = $((Q_i \cap M_i), (Q_j \cap M_j))$. and x ∈ $\{(Q_i \cap M_i), (Q_j \cap M_j)\}$. Then the following equation holds;

$$\mu_R^C(x) = \begin{cases} 1 - \mu_R(x), & \mu_R(x) \ge 0, \\ -1 - \mu_R(x), & \mu_R(x) < 0. \end{cases}$$

**THEOREM [50]:** *When fuzzy causal concepts $C_i$, and $C_j$ are given, the following FCRs are all equivalent:*

$$C_i \xrightarrow[\alpha]{} C_j,\ {}^{\sim}C_i \xrightarrow[\alpha]{} {}^{\sim}C_j,\ C_i \xrightarrow[-\alpha]{} {}^{\sim}C_j,\ {}^{\sim}C_i \xrightarrow[-\alpha]{} C_j$$

*where $-1 \le \alpha \le 1$.*

**THEOREM [50]:** *When $\sim C_i$ is a negative concept of $C_i$ and the dis-quantity fuzzy set of $\sim C_i$ is equal to the complement of $C_i$'s quantity fuzzy set, then the following FCRs are all equivalent:*

$$C_i \xrightarrow[\alpha]{} C_j,\ C_i \xrightarrow[1-\alpha]{} {}^{\sim}C_j,\ C_i \xrightarrow[\alpha-1]{} C_j,$$

*where $0 < \alpha < 1$.*

**THEOREM [50]:** *When $\sim C_i$ is a negative concept of $C_i$ and the dis-quantity fuzzy set of $\sim C_i$ is equal to the complement of $C_i$'s quantity fuzzy set, then the following FCRs are all equivalent:*

$$C_i \xrightarrow[\alpha]{} C_j,\ C_i \xrightarrow[-1-\alpha]{} {}^{\sim}C_j,\ C_i \xrightarrow[\alpha+1]{} C_j,$$

*where $-1 < \alpha < 0$.*

**THEOREM [50]:** *When $\sim C_i$ is a negative concept of $C_i$ and the dis-quantity fuzzy set of $\sim C_i$ is equal to the complement of $C_i$'s quantity fuzzy set, then the following hold:*

$$C_i \xrightarrow[1]{} C_j,\quad implies\quad C_i \xrightarrow[+0]{} {}^{\sim}C_j,\ C_i \xrightarrow[-0]{} C_j$$
$$C_i \xrightarrow[1]{} C_j,\quad implies\quad C_i \xrightarrow[-0]{} {}^{\sim}C_j,\ C_i \xrightarrow[+0]{} C_j$$

There may be a case that even though $(Q_i \cap M_i) \subset (Q_j \cap M_j)$ is true, but $(\sim Q_i \cap M_i) \subset (\sim Q_j \cap M_j)$ is not true. Also such a case would happen that $(Q_i \cap M_i) \subset (\sim Q_j \cap M_j)$ is true, but $(\sim Q_i \cap M_i) \subset (Q_j \cap M_j)$ is not true. For example, there may be a stock market situation that institute investors' buying causes the increases of composite stock price but their selling cannot cause the decrease of composite stock price.

Buying by institute investors

$\xrightarrow[0.9]{}$ Increase of composite stock price,

Selling of institute investors

$\xrightarrow[0.6]{}$ Decrease of composite stock price,

where $\rightarrow$ means "partial causality". Both relations show another types of FCR, which will be termed as "fuzzy partially causal relation (FPCR)". To help readers understand more clearly, we will define FPCR.

**DEFINITION [50]:** *$C_i$ partially causes $C_j$ iff $(Q_i \cap M_i) \subset (Q_j \cap M_j)$*

**DEFINITION [50]:** *$C_i$ partially causally decrease $C_j$ iff $(Q_i \cap M_i) \subset (\sim Q_j \cap M_j)$.*

**DEFINITION [50]:** *A uni-quantitative concept $C_i$ is either $Q_i \cap M_i$ or $\sim Q_i \cap M_i$.*

**THEOREM [50]:** *Only FPCR exists between the uni-quantitative concepts.*

We have only recalled here the main definitions and an overall picture of the work carried out in the field of Fuzzy Causal Relationships and Fuzzy Partially Causal Relationships by Kim H.S. and Lee K.C. [50]. For more about this please refer [50].

### 1.3.7 Automatic Construction of FCMs

Here we just recall the theory of Automatic Construction of FCM based on the user provided data as introduced by Schneider M. et al [86]. It mainly deals with conversion of numerical vectors into fuzzy sets, closeness of relation between numerical vectors, determining whether relations between variables are direct or inverse and causality among variables. Since FCMs are dynamic systems, their resonant states are limit cycles. The limit cycle or hidden pattern is a FCM inference. The limit cycle or hidden pattern is a FCM inference. However, there are 2 important drawbacks in the conventional FCMs:

a. Lack of a concept of time: practically each causal has a different time delay.
b. They cannot deal with co-occurrence of multiple causes such as expressed by "and" conditions.

To overcome these drawbacks, the extended fuzzy cognitive maps (E-FCMs) were proposed by Hagiwara [41]. The features of the E-FCMs are:

a. Weights having non-linear membership functions.
b. Conditional weights.
c. Time-delay weights

The questions one might ask are:

a. Why equations are based on summation and not other operations such as max?
b. Why equation utilizes multiplication rather than min?
c. How to determine relative weights for experts?

These and other issues led us to explore the questions related to construction of the FCM. It is our belief that understanding how the FCM is constructed will contribute to resolving some of the above-raised questions with respect to inferencing.

Each demographic / economic concept for which numerical measurements (data) are available is represented by a numerical vector (V), whereas each element in the vector ($\nu$) represents a measurement from one country, where each vector consists of measurement from 106 countries. The method presented here requires that each numerical vector be treated as a fuzzy set [86].

One simplistic way to convert a numerical vector into a fuzzy set is:

1. Find a maximum value in V, and assign $\chi = 1$ to it, that is

$$MAX(\nu) \Rightarrow \chi v\ (\nu) = 1.$$

2. Find a minimum value in V, and assign $\chi = 0$ to it that is

$$MIN(\nu) \Rightarrow \chi v\ (\nu) = 0.$$

3. Project all other i vector elements ($v_i$) into the interval [0, 1], proportionally, that is

$$\chi v\ (\nu_I) = \frac{\nu_i - MIN(\nu)}{MAX(\nu) - MIN(\nu)}.$$

However, this procedure is legitimate only if the computed grades of membership truly reflect the membership of each element of the vector in the fuzzy set. In many cases, however, the straightforward procedure described above might result in assigning grades of membership in a way unsupportable by common sense reasoning.

In order to facilitate a more reasonable representation of the membership grades, the system should allow the user to provide a threshold (where relevant) for which every element above it will fully belong to the set with $\chi = 1$. In the same fashion the user may provide a threshold for which every element in the vector below it will be excluded from the set. That is, we assign the value 0 to its membership grade.

Let V be a numerical vector of n elements, and let $\nu$ be an element in V. Then the user can provide the upper threshold ($\alpha_u$) so that

$$\forall \upsilon\ (\upsilon \geq \alpha_u) \Rightarrow (\chi_\upsilon = 1)$$
$$\forall \upsilon\ (\upsilon \leq \alpha_u) \Rightarrow (\chi_\upsilon = 0).$$

When researching relation between two variables, the important question is: how closely are they related? The closeness of relation between two vectors is based on the concept of distance between vectors. The implementation of this concept in our system requires different computation for vectors, which are directly related, and vectors, which are inversely related. We shall restrict our examples to linear dependencies only.

If vectors $V_1$ and $V_2$ are directly related, then the closet relation between them is when for each i (i = 1, …, n), $\chi_1(v_i) = \chi_2(v_i)$. Let $d_i$ be the distance between the corresponding elements of $V_1$ and $V_2$ such that

$$d_i = |\chi_1 (v_i) - \chi_2 (v_I)|,$$

and let AD be the average distance between the vectors $V_1$ and $V_2$. Then utilizing the above equation

$$AD = \frac{\sum_{i=1}^{n} |d_i|}{n}$$

or

$$AD = \sqrt{\sum_{i=1}^{n} (d_i)^2}\ .$$

Once the average distance between the two vectors is computed, the closeness, or similarity S between the two vectors is computed using the following equation:

$$S = 1 - AD.$$

In the case of perfectly identical sets where each element $\chi_1(v_i) = \chi_2 (v_i)$, the average distance between the vectors (AD) is

$$AD = \frac{\sum_{i=1}^{n} |d_i|}{n} = 0$$

and the similarity between the two vectors is equal to

$$S = 1 - AD = 1 - 0 = 1$$

which indicates perfect similarity.

The greatest degree of dissimilarity is possible only in the case of binary sets. That is, for every element $V_1$ where $\chi_1 (v_i) = 1$, the corresponding value in $V_2$ has the grade of membership of 0, and vice versa. In this case

$$AD = \frac{\sum_{i=1}^{n} |d_i|}{n} = \frac{n}{n} = 1$$

and thus

$$S = 1 - AD = 1 - 1 = 0,$$

which indicates the extreme degree of dissimilarity.

In the case of fuzzy sets, where most of the elements have grades of membership between 0 and 1, no perfect similarity, or perfect dissimilarity is expected, but similarity to some degree is more reasonable outcome.

The method for computing closeness of relation (similarity) for inversely related vectors is in principle very similar to that of the directly related vectors, with the exception that the distance between the corresponding elements of the inversely related vectors $V_1$ and $V_2$ is

$$Id_i = \left| \chi_1 (v_i) - (1 - \chi_2 (v_I)) \right|$$

where $id_i$ is the distance between the corresponding elements of inversely related vectors.

The utilizing the above equation,

$$AD = \frac{\sum_{i=1}^{n} |id_i|}{n}$$

or

$$AD = \sqrt{\sum_{i=1}^{n} (id_i)^2} \; .$$

Once the average distance between the two vectors is computed, the inverse similarity between the two vectors is computed using the following equation: S = 1 – AD.

**TABLE 1: LIST OF VARIABLES OF THE MODEL**

| VARIABLE NAME | DESCRIPTION |
|---|---|
| Illiteracy rate | Adult illiteracy rate (% 1990) |
| Elementary education | Elementary education (% of age group, 1989) |
| Secondary education | Secondary education (% of age group, 1989) |
| Higher education | Higher education (% of age group, 1989) |
| Infant mortality | Infant mortality (per 1000, 1990) |
| Urban population | Urban population (as % of total population, 1990) |
| Death rate | Death rate (per 1000, 1990) |
| Birth rate | Birth rate (per 1000, 1990) |
| Population growth | Population growth rate (per 1000, 1990) |
| Average population growth | Average annual popl. growth rate (% 1980–1990) |
| GNP per capita | GNP per capita $, 1990) |
| ICP estimates of GDP/Capita | ICP estimates of GDP / capita ($, 1990) |
| Average growth of GNP per capita | GNP per capita (average annual rate of growth, 1965–1990) |
| Agriculture as % of GDP | Agriculture as % of GDP (1990) |
| Exports of machinery & transport equip. | Machinery & transport equipment (% of exports, 1990) |
| Exports of textiles and clothing | Textiles and clothing (% of exports, 1990) |
| Exports of primary commodities | Exports of primary commodities (% of exports) |
| Exports | Exports (millions $, 1990) |

Graphical illustration of the model is given in Figure 1.3.9.

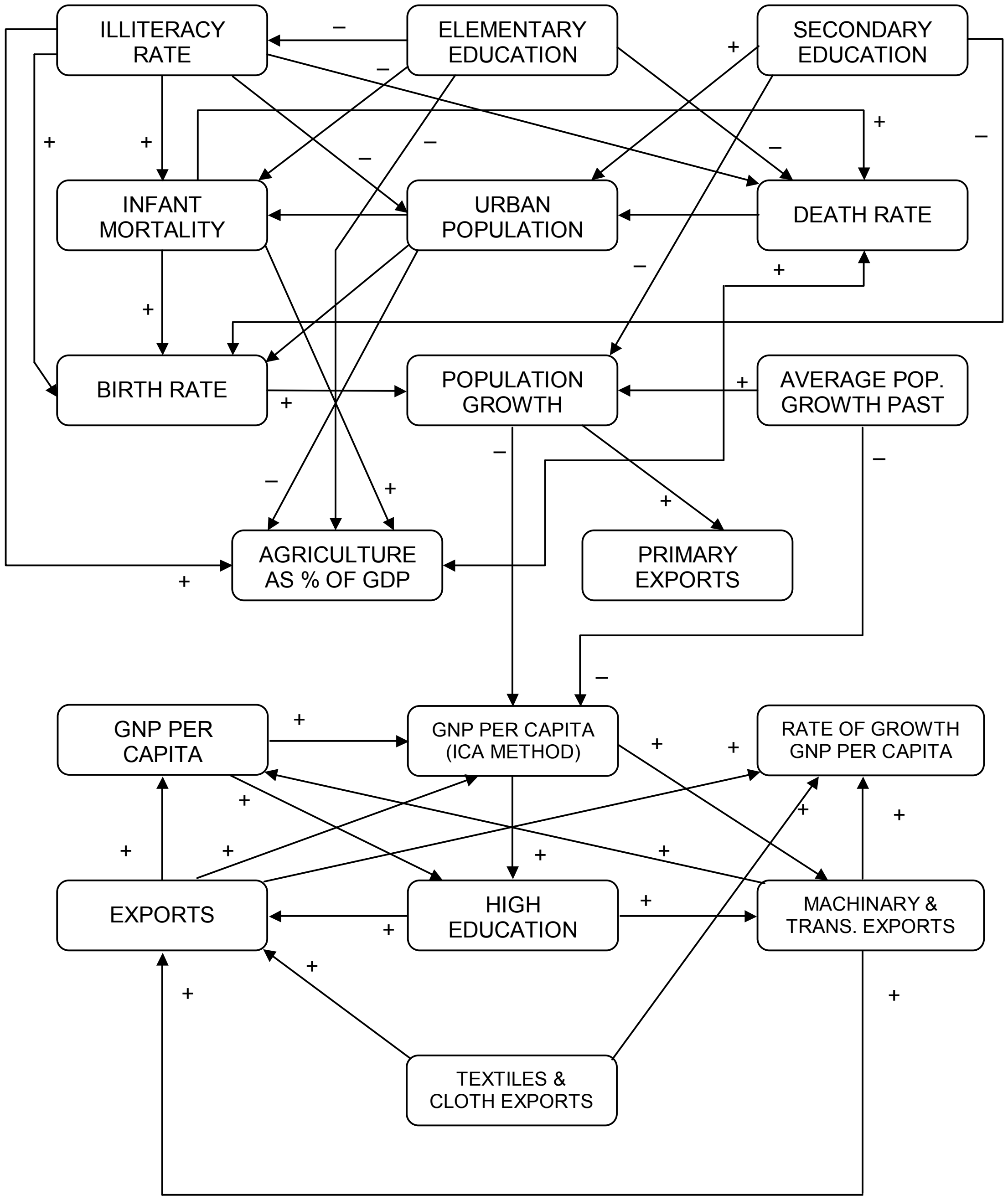

**FIGURE: 1.3.9**

Using the model Schneider et al [86] have described the automatic construction of FCM. We have only sketched the relevant facts, for more information, refer [86].

### 1.3.8: Fuzzy Cognitive State Map vs. Markovian model in the study of users web behaviour

Searching for information in general is complex, fuzzy and an uncertain process. The Internet is a system of storage and retrieval of information characterized by its enormous size, hypermedia, structure and distributed architecture. Much research has been done on the process of browsing for information.

Some seven concepts are taken as the nodes by Meghabghab [71] for studying which are given by $C_1, \ldots, C_7$.

A Markovian Modeling of Users Behaviour given by [71] is introduced and compared to a Fuzzy Cognitive Map (FCM) that represents the opinions of experts on how users surf the web. If we list all the concepts involved in the user's behaviour on the web: $C_1$ = Search, $C_2$ = Browse, $C_3$ = T.C. (time constraints), $C_4$ = I.C. (Information constraints), $C_5$ = Success, $C_6$ = Relevance, $C_7$ = Failure.

The directed graph of the FCM is given in Figure 1.3.10, using the connection matrix given by [71] one can determine the edge-weights.

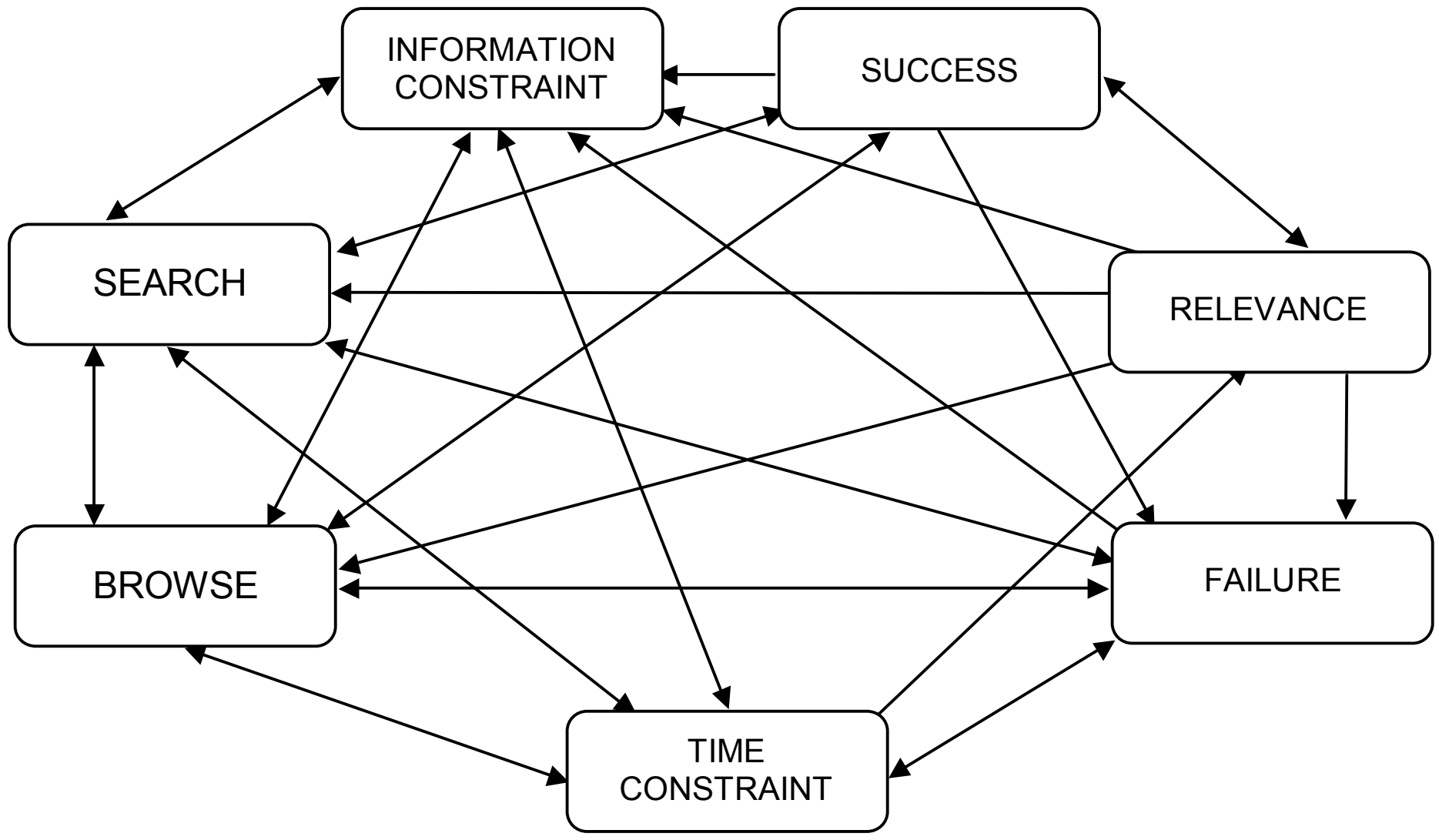

**FIGURE: 1.3.10**

Then E can be built by the expert's opinion and observation on user's behaviour on the web by the following matrix:

$$E = \begin{bmatrix} 0 & -1 & -1 & 1 & -1 & -1 & 1 \\ 1 & 0 & -1 & -1 & -1 & -1 & 1 \\ -1 & -1 & 0 & -1 & 0 & 1 & 1 \\ -1 & -1 & -1 & 0 & 0 & 0 & 0 \\ 1 & 1 & 1 & 1 & 0 & 1 & -1 \\ 1 & 1 & 0 & 1 & 1 & 0 & -1 \\ 1 & 1 & 1 & 1 & 0 & 0 & 0 \end{bmatrix}$$

Experts are divided on how all these concepts affect failure which is concept $C_7$. Some even give contradictory information on their answers why users fail to give an answer for a given query. Experts answered $C_7$ in the following way:

| | $C_1$ | $C_2$ | $C_3$ | $C_4$ | $C_5$ | $C_6$ | $C_7$ |
|---|---|---|---|---|---|---|---|
| $C_7$ | ±1 | ±1 | ±1 | ±1 | 0 | 0 | 0 |

The reason the positive values have been adopted in the last column of E as the final answer is that the negative values yield to minimum points that lead to chaotic behavior. If we are to look at "failure" as an input vector: $F_1$ = [0 0 0 0 0 0 1], then by applying in the matrix E we get: $F_1$* E = [1 1 1 1 1 0 –1 0], which translates into $F_2$ = [1 1 1 1 0 0 0]. By applying the same produce to $F_2$ again we get: $F_2$* E = [1 –3 –3 –1 –2 –1 +3], which translates into $F_3$ = [0 0 0 0 0 0 1]. Notice that and the network is back where it starts and converges to the limit cycle: $F_1 \rightarrow F_2 \rightarrow F_3 = F_1$.

If we look at success in an input vector $S_1$ = [0 0 0 0 1 0 0] and using in E we get: $S_1$ * E = [1 1 1 1 0 1 –1], which translated into $S_2$ = [ 1 1 1 1 0 1 0], $S_2$ * E = [0 –2 –3 0 –1 –1 2], which is translated into $S_3$ = [0 0 0 0 0 0 1], and $S_3$ * E = [ 1 1 1 1 1 0 –1 0], which translates into $S_4$ = [ 1 1 1 1 0 0 0]. $S_4$ * E = [–1 –3 –3 –1 –2 –1 3], which translates into $S_5$ = [0 0 0 0 0 0 1] = $S_3$. Notice how $S_5 = S_3$ and the network converges to the limit cycle $S_3 \rightarrow S_4 \rightarrow S_5 = S_3$.

The length of the limit cycle should be less than the number of concepts otherwise cross talk can occur.

How does user's behaviour change with time while searching the web? Do they learn new behavior or reinforce old ones. To answer such a question an adaptive FCM is built that uses correlation or Hebbian learning to encode some limit cycles in the FCMs or Temporal Associative Memories [58]. This method can store only few patterns. Differential Hebbian learning encodes change in a concept. To encode binary limit cycles in matrix E the TAM method sums the weighted correlation matrices between successive states. To encode a cycle a concept C has to be converted to a $X_i$ by changing the 0s into –1. Then E is a weighted sum:

$$E = X_1'{*}\, X_2 + X_2'{*}\, X_3 + X_3'{*}\, X_4 + \ldots + X_n'{*}X_1.$$

Then applying to our data of 7 concepts we have the following matrix E.

$$E = \begin{bmatrix} -1 & 1 & -1 & 3 & 1 & -1 & 1 \\ 1 & 3 & 1 & 5 & 3 & -3 & -1 \\ -1 & 1 & -1 & 2 & 5 & -1 & 1 \\ 3 & 1 & -1 & 3 & 1 & -5 & 1 \\ 1 & 3 & 5 & 1 & 3 & 1 & -1 \\ -1 & 1 & -1 & -1 & 5 & -1 & -3 \\ -3 & -5 & -3 & -7 & -1 & 1 & 3 \end{bmatrix}.$$

If we are to look at "failure" as an input vector: $F_1$ = [ 0 0 0 0 0 0 1], then $F_1$*E = [-3 –5 –3 –7 –1 1 3], which translates into $F_2$ = [0 0 0 0 0 1 1], $F_2$*E= [– 4 – 4 – 4 –8 4 0 0] which translates into $F_3$ = [0 0 0 0 1 0 0], $F_3$*E = [ 1 3 5 1 3 1 –1], which translated into $F_4$ = [ 1 1 1 1 1 1 0], $F_4$*E = [2 10 2 14 18 –10 –2], which translates into $F_5$= [ 1 1 1 1 1 0 0], $F_5$*E= [ 3 9 3 15 13 –9 1], which translates into $F_6$ = [ 1 1 1 1 1 0 1], $F_6$*E

= [ 0 4 0 8 12 –8 4], which translates into $F_7$ = [ 0 1 0 1 1 0 1], $F_7$*E = [ 2 2 2 2 6 –6 2], which translates into $F_8$ = [1 1 1 1 1 0 1].

Notice how $F_8 = F_6$ and the network converges to the limit cycle $F_6 \to F_7 \to F_8 = F_6$.

If we look at success in an input vector $S_1$ = [0 0 0 0 1 0 0], then $S_1$*E = [1 3 5 1 3 1 –1], which translates into $S_2$ = [1 1 1 1 1 1 0], $S_2$*E = [2 10 2 14 18 –10 –2], which translates into $S_3$ = [ 1 1 1 1 1 0 0], $S_3$*E = [3 9 3 15 13 –9 1], which translates into $S_4$ = [1 1 1 1 1 0 1], $S_4$*E = [ 0 4 0 8 12 –8 4], which translates into $S_5$ = [0 1 0 1 1 0 1], $S_5$*E = [ 2 2 2 2 6 –6 2], which translates into $S_6$ = [ 1 1 1 1 1 0 1]. Notice how $S_6 = S_4$ and the network converges to the limit cycle $S_4 \to S_5 \to S_6 = S_4$.

We can still work with different input vectors. The main conclusion made are:

They are in a redefined stage of thinking on the web at any given time. These stages are unique and cannot be violated at any time. For any time i different than j. state $S_i$, is different from state $S_j$, which can be stated as follows: $\forall\ i \neq j: S_i, \neq S_j$.

A markovian [71] modeling has been applied to look at the cognitive behavior of users on the Web in the graph of given in Figure 1.3.11. Markovian Modeling of user's cognitive map on the web ($S_1$ = Home, $S_2$ = Browse, $S_3$ = Search, $S_4$ = Select).

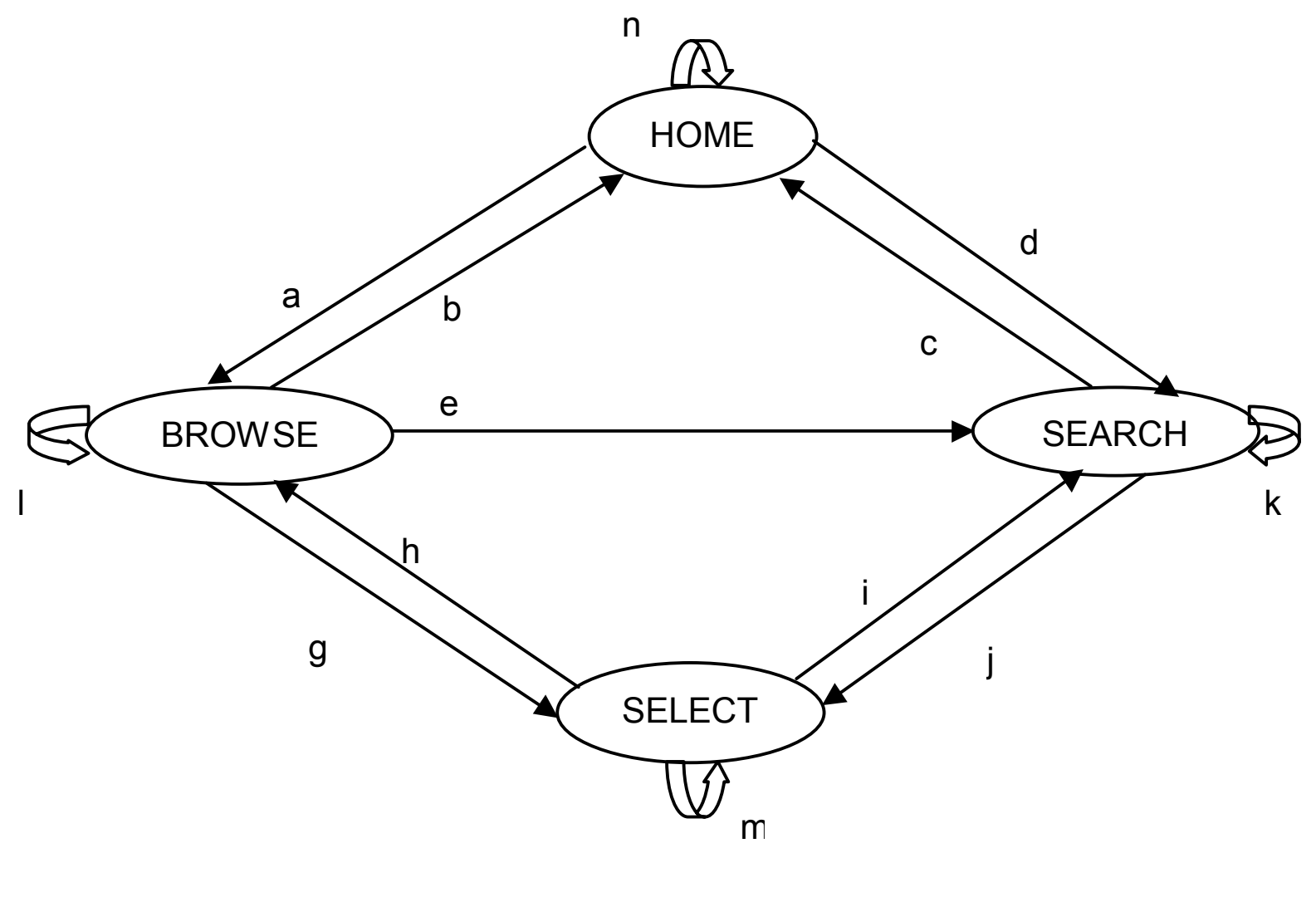

**FIGURE: 1.3.11**

20 students were selected to participate in the study. Their role is to elucidate the different values of the edges a, b, c, d, e, g, h. i, j, k, l, m, n that are needed in Figure 1.3.11. They will also verify the correctness of the cognitive model established in the paper. Any behavior that cannot be accounted for in Figure 1.3.11 will be taken into account and will help to modify or correct the model proposed. Also, the time spent in each task will help to explain the learning that takes place and the experience that builds up in the synapses of the user during a web browsing session. The difference between users in a web browsing session will help also to evaluate the failure or success rate of any web browsing session for a given query.

Fuzzy cognitive maps can model the causal graph of expert's perception of how users behave on the web. Although experts are still divided on the reasons users fail to answer a query, whether it is a fact based, research based or a task based query, Meghabghab [71] establishes a stable dynamical system that behaves well in a limit cycle. When asked what happens in a successful query the system responds well and uncovers a limit cycle. When asked what happens in the case of a failure the system responds well and uncovers a limit cycle.

For example, if a user loops while searching, loops browsing, stress out because of time constraints and information overload, he/she will fail. In the case of success, if a user loops while searching, loops while browsing, stresses out because of time constraints and information overload but uncovers relevant information, he/she will succeed.

This model shows that relevance is a key issue in separating success from failure. The same information presented to a user cannot be interpreted as being relevant to a user while it is being interpreted relevant to another. Relevance becomes a major issue in information interpretation and filtering. Relevance needs to be exploited as part of query expansion and relevance feedback to improve information retrieval. Popularity should not be confused with relevance.

Although Google and other engines use different ways to judge relevance, this is still a very tough issue and more research in this regard to provide a reliable automatic relevance algorithmic will help improve returned set results and reduce the high failing rate of users (40% in our case) in searching for information.

For more information, please refer [71].

### 1.3.9: Cognitive maps and Certainty Neuron Fuzzy Cognitive Maps

According to Tsadiras et al [108], each concept node can be imagined as a living entity (cell) that is positively or negatively activated and can be influenced by cells in its neighborhood. The natural behavior for that cell would be to lose some of its activation when there is no stimulation to maintain the activation. The same should apply to concepts of our structure.

Imagine a concept that is activated in a positive or negative manner to a certain activation level and receives no influence from other concepts. As time passes, its activation level decays towards the zero activation level. An example of such a concept is that of the population of a city. A positive activation implies an intention of the population of the city to increase. If this concept receives no influences from other concepts, this intention is not justified anymore, and because of that, it should gradually decrease to activation level zero.

To capture the above behavior the decay factor is introduced. The decay factor defines the fraction of the current activation that will be lost at each time step, according to the decay mechanism. It can take a value within [0, 1] , and the greater it is, the faster the cell becomes inactive when it receives no stimulation, and so the new activation level of a node is determined from its previous activation, the stimulation that it receives, and the decay term. The new activation level, after also taking into account the decay mechanism, is the following:

$$A_i^{t+1} = f_M(A_i^t, S_i^t) - d_i A_i^t$$

where $d_i$ is the decay factor of concept $C_i$. Clearly, the decay factor subtracts from the initial new activation a percentage of the current activation level. Using the activation level due to the decay can take a value less than –1 or greater than 1. This implies that the activations tend to be more negative or positive than they should be. In these extreme cases, the activation is set to –1 or 1, respectively.

Neurons that use as a transfer function the two-variable function were defined as certainty neurons [110]. Using this type of neuron in the FCM structure, the Certainty Neuron Fuzzy Cognitive Maps (CNFCMs) are built. Some preliminary results on the dynamical behavior of the CNFCMs [108,110] have shown that:

1. CNFCMs can reach equilibrium at fixed points in a direct or asymptotic way with the activation to be decimals in the interval [–1,1];

2. CNFCMs can exhibit a limit cycle behavior. The system falls in a loop of a specific period, and after a certain number of steps, it reaches the same state.

The equilibrium points of a CNFCM can be calculated. For the equilibrium state, $A_i^{t+1} = A_i^t$, for i = 1,…, n. This can only happen if $A_i^t$ and $S_i^t$ are of the same sign.

It is shown that if $A_i^t$ and $S_i^t$ are of opposite signs, then

$$A_i^t > 0 \Rightarrow f_M(A_i^t, S_i^t) < (A_i^t \Rightarrow f_M(A_i^t, S_i^t) - d_i A_i^t < A_i^t \Rightarrow A_i^{t+1} < A_i^t,$$

$$A_i^t < 0 \Rightarrow f_M(A_i^t, S_i^t) > A_i^t \Rightarrow f_M(A_i^t, S_i^t) - d_i A_i^t > A_i^t \Rightarrow A_i^{t+1} > A_i^t$$

and so equation $A_i^{t+1} = A_i^t$ can never be satisfied. The conclusion is that the equilibrium point is not reached. For a complete and detailed working refer [108].

We have that $A_i^t > 0$ $S_i^t > 0$. In that case,

$$\begin{aligned} A_i^{t+1} = A_i^t \quad &\Rightarrow \quad A_i^t + S_i^t(1 - A_i^t) - d_i A_i^t = A_i^t \\ &\Rightarrow \quad S_i^t - S_i^t A_i^t - d_i A_i^t = 0 \Rightarrow A_i^t(d_i + S_i^t) = S_i^t \\ &\Rightarrow \quad A_i^t = \frac{S_i^t}{(d_i + S_i^t)} \end{aligned}$$

Using the second rule, we have that $A_i^t < 0$, $S_i^t < 0$. In that case,

$$\begin{aligned} A_i^{t+1} = A_i^t \quad &\Rightarrow \quad A_i^t + S_i^t(1 + A_i^t) - d_i A_i^t = A_i^t \\ &\Rightarrow \quad S_i^t - S_i^t A_i^t - d_i A_i^t = 0 \Rightarrow A_i^t(d_i - S_i^t) = S_i^t \\ &\Rightarrow \quad A_i^t = \frac{S_i^t}{(d_i - S_i^t)} \end{aligned}$$

Now combining the results we get

$$A_i^{eq} = \frac{S_i^{eq}}{d_i + \left|S_i^{eq}\right|} = \frac{\sum_j w_{ji} A_j^{eq}}{d_i + \left|\sum_j w_{ji} A_j^{eq}\right|}, \quad i = 1, \cdots, n.$$

It is observed that if $d_i = 0$, $i = 1,\ldots, n$, then $A_i^{eq}$ equals either – 1 or 1. We can notice that the calculation of the equilibrium state can be seen as the solution of set of n nonlinear equations. The solutions of this set of equations produce the values of the activations $A_i^{eq}$ as far as such a state exists and the system does not enter into a limit cycle. Therefore, the CNFCM can be seen as a system that attempts to solve the set of n nonlinear equations, using a simple iterative method employing synchronous updating. It should be mentioned that if $A^{eq} = \left[A_1^{eq}, A_2^{eq}, \cdots, A_n^{eq}\right]$ is an equilibrium state, then $-A^{eq} = \left[-A_1^{eq}, -A_2^{eq}, \cdots, -A_n^{eq}\right]$ is also an equilibrium state as it can be easily verified. This is the reason that we should refer to pairs of equilibrium points.

It is shown that CNFCM increases the inference and representing capabilities of CM and FCM techniques, and allows the activation level of concepts to be decimals in the whole interval [–1, 1]. The wealthy inference capabilities and dynamical behavior of the CNFCM are presented in comparison to that of the classical FCMs.

For a complete work refer [108].

### 1.3.10: Fuzzy Cognitive Maps considering time relationships

As described previously, crisp cognitive maps assume that all effects are equally strong, by placing unit intensities (i.e. +1 or –1) on each arc. It might be more reasonable to place a different strength on each arc of a map, thus yielding a FCM. The strength is interpreted as the relative intensity of the effect, and can be a value in the interval [–1, 1].

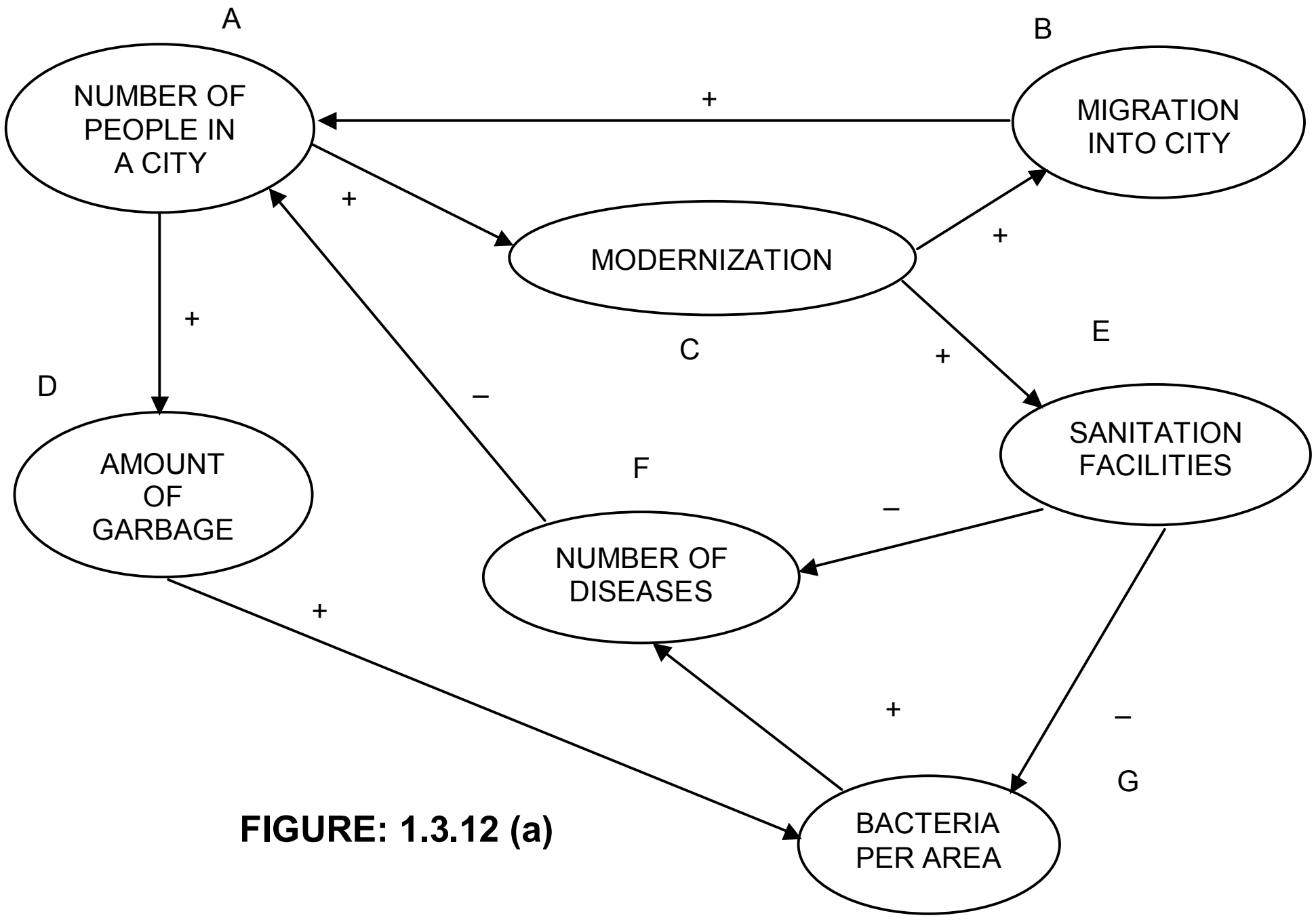

**FIGURE: 1.3.12 (a)**

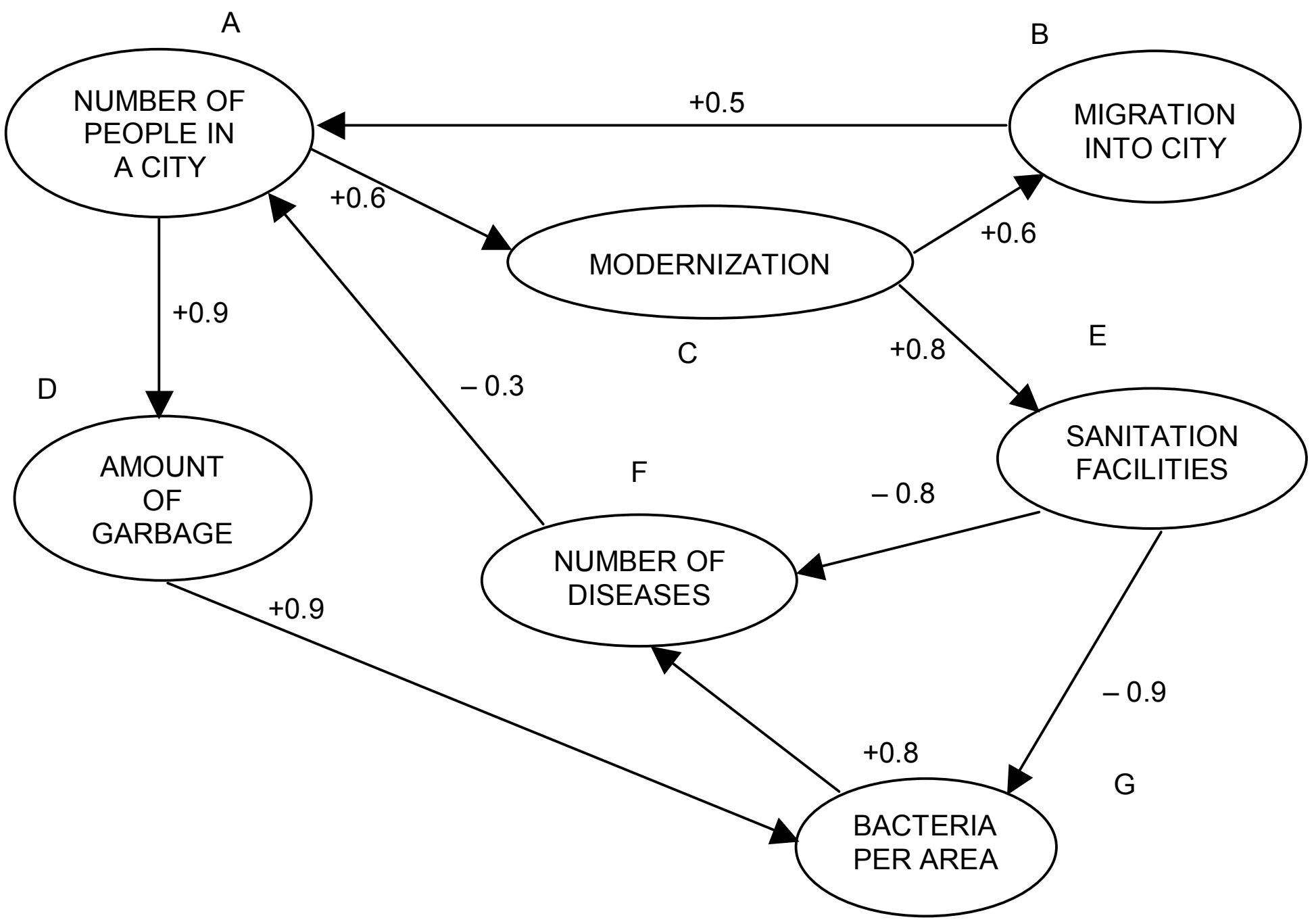

**FIGURE: 1.3.12 (b)**

Another omission in the FCM is the time relationship (or time lag) involved before a change in node $x_i$ has an effect on node $x_j$. For example, in Figure 1.3.12 an increase in "migration into city" will lead almost immediately to an increase in "number of people in a city", while an increase in "number of people in a city" might lead less immediately to an increase in "modernization". Up to the present, FCM assumes all effects take place in one unit time. Thus, a more realistic map would introduce a time lag corresponding to each effect.

However, the time lags are hard to estimate, and there is a trace-off between the generality of the map and the possibility of estimating its time lags in the realistic way. A method is that the time lags can be represented in fuzzy, relative terms by a domain expert: The time lags might be expressed as being in the set {I, N, L} or another appropriate set, where I, N and L, respectively, are the abbreviations of "Immediate", "Normal" and "Long". However the concept of indeterminacy I is left out which can be incorporated only in NCMs. The terms can be assigned to the relative fuzzy values, for example, {1, 2, 3} for each term. We can interpret that "3" is three delay units such 3 years, months, etc, "2" is two delay units, and "1", one unit. Note here that assessing precise values on time lags such as continuous values, in most conceptual FCMs, would be very difficult. Figure 1.3.13 shows a hypothetical FCM with time lags.

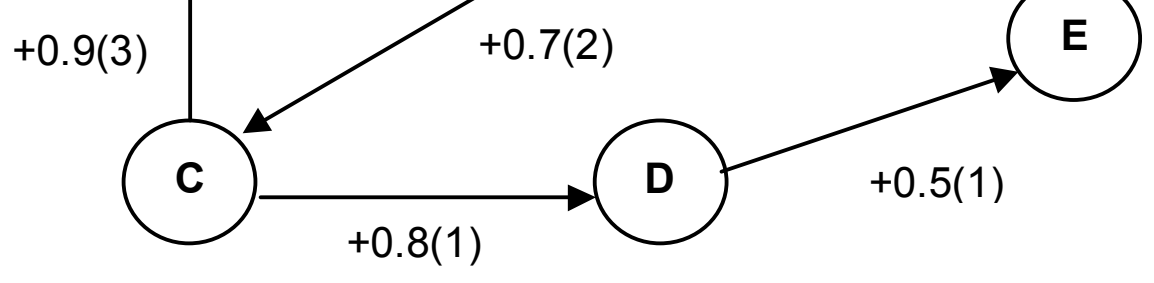

**FIGURE: 1.3.13**

Since considering time lags in the FCM, we can analyse system behaviors according to lapse of time, i.e. the change of causalities among factors for some discrete time t. Indirect effects at time t are represented by:

$$IE(t) = S^t \tag{1}$$

Then, total effects for some time t are given by:

$$TE(t) = OR^t_{k=1} (IE(k)). \tag{2}$$

Note in previous FCMs that we can only refer to the change of causalities for a stable state which is represented by total effect-matrix TE. Furthermore, the waiting time for a stable state would be meaningless, because of overlooking time relations. Park K.S. and Kim S.H. [78] formally define a FCM with time lags as Fuzzy Time Cognitive Map(FTCM), and illustrate a transformation that leads to convenient causal inference.

**DEFINITION [78]:** *A FTCM is a directed graph G = (X, E) consisting of a finite set X of N nodes,* $X = \{i\}^N_{i=1}$, *and a set E of edges (or arcs)* $e_{ij}$, $i, j \in X$, $E \subseteq X \times X$. *A note i* $\in X$ *is representative of a relevant factor of the decision domain. Each edge* $e_{ij}$ *of the FTCM has two kinds of relative causalities: the strength* $s_{ij} \in [-1, 1]$ *and the time lag* $t_{ij} \in [a, b]$, *where a and b are constants and* $0 < a \le b$.

As described in the previous section, Park and Kim [78] deal with time lags on not continuous but discrete space. The $t_{ij}$ is a value in a finite set Ω of M-many values $\Omega = \{t\}^M_{t=1}$. For example, if the time lag involved before a change in node i has an effect on node j, t (i → j) = 6 months, t(j → k) = 1 year, and t(k → l) = 2 years, then we can assign as $t_{ij} = 1$, $t_{jk} = 2$ and $t_{kl} = 4$.

In summary, the FTCM is identified by X, $S = \|s_{ij}\|_{NxN}$ and $T = \|t_{ij}\|_{NxN}$. Thus, it has a different time lag on each edge $e_{ij}$. In order to use equations (1) and (2) for causal inference, it is required that all edges must be translated into the edges with one unit-time (or the same time) lag. The next subsection contains a technique for translating into one unit-time lag.

For deriving the causal strength per one unit-time $s^*_{ij}$ on an edge $e_{ij}$, simply the arithmetic division might be considered (i.e. $s^*_{ij} = s_{ij} / t_{ij}$ for all i and j). But it would be misleading. The reason is that large time lag on an edge cannot lead to small causal strength on it. Rather, it should be interpreted that the t(i→j) is only large, without a change of the $s_{ij}$. Therefore, a value-preserving translation is required.

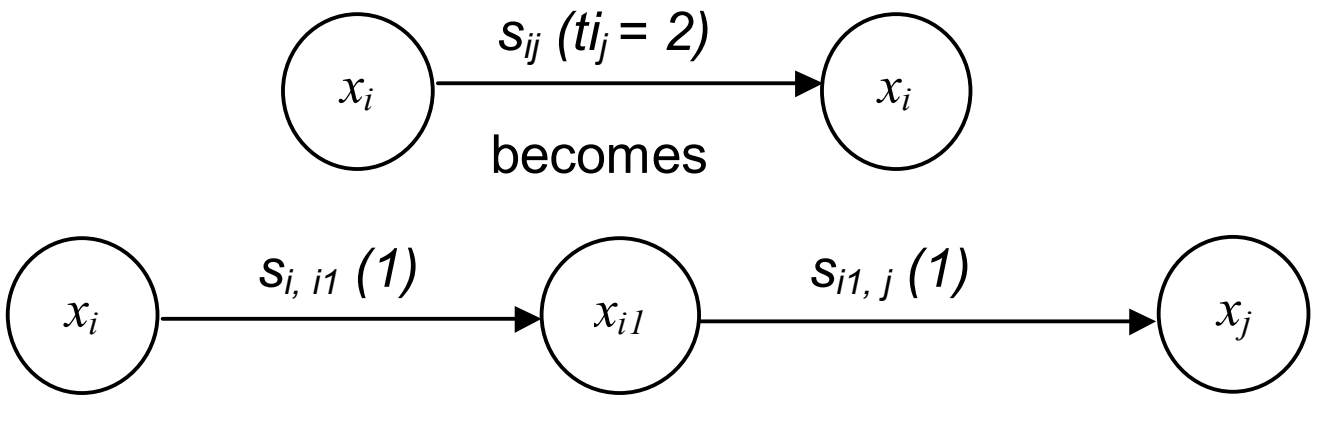

**FIGURE: 1.3.14**

For value-preserving translation, we introduce dummy nodes. For example, as shown in Figure 1.3.14, if a node i influences a node j through two delay units, then we introduce one dummy node iI so that the edge $e_{ij}$ of the initial FTCM becomes two edges, $e_{i,i1}$ and $e_{i1,j}$, of new FTCM. Thus, each time lag on the edges of new FTCM becomes one unit-time.

In Figure 1.3.14, the name or concept of node i1 equals to the name of the original mode i. It can be illustrated as follows: the causality of node i stays at node i1, since one unit time is elapsed. After two unit-time elapsed, the causality of i imparts to j. Then, values $s_{i,i1}$ and $s_{i1,j}$ can be assigned appropriately with the preservation of the original value: If the indirect effect amounts to specifying the weakest (i.e. minimum) causal link in a path, then $|s_{i,i1}| = |s_{i1,j}| = |s_{ij}|$, because

$$s_{ij} = \text{sign}\,(s_{i,i1})\,.\,\text{sign}\,(s_{i2,\,j}).\ \min\left(|s_{i,i1}|, |s_{i1,j}|\right).$$

The above observation is easily generalized as follows: If the time lag on $e_{ij}$ is m ($1 \leq m$) delay units, then we introduce m –1 dummy nodes which are the equal name with node i. Then, a set of new edges, $\{s_{i,i1}, s_{i1,i2},\ldots,s_{i(nt-1),j}\}$, comes out. Each time lag on the new edges is one unit-time. If the original strength on $e_{ij}$ is $s_{ij}$, then an appropriate strength on each new edge, $s_{p,p+1}$ for $(p, p+1) \in \{(I, i1), (i1, i2),\ldots, (i(m-1), j)\}$, is assigned according to what t-norm is used for propagating the indirect effects.

If product operator is used, then the following conditions have to be satisfied:

$$|s_{ij}| = |s| \text{ and sign } (s_{ij}) = \text{sign } (s),$$

where $s = s_{i,\,i1}.\ s_{i1,i2},\ldots,s_{i(m-1),j}$. Suppose $|s_{i,i1}| = |s_{i1,i2}| =,\ldots,= |s_{i(m-1),j}|$, then $|s_{p,p+1}|^{m} = |s_{ij}|$. Thus, $|s_{p,p+1}| = |s_{ij}|^{1/m}$. To identify the signs of the $s_{p,p+1}$, we consider the following four cases:

1. sign ($s_{ij}$) = +, and m is an odd number,
2. sign ($s_{ij}$) = +, and m is an even number,
3. sign ($s_{tj}$) = –, and m is an odd number,
4. sign ($s_{tj}$) = –, and m is an even number.

In the case of (1) through (3), $s_{p,p+1} = \text{sign}\,(s_{ij})\ |s_{ij}|^{1/m}$ for all (p,p+1). In the case of (4), the number of minus signs on the $s_{k,k-1} = -|s_{ij}|^{1/m}$ and $s_{q,q+1}=|s_{ij}|^{1/m}$, where $\{s_{q,q+1}\} = \{s_{p,p+1}\}\backslash S_{k,k+1}$ (The operator "\" denotes set subtraction, A\B = A-B). If minimum operator is used, then in the cases of (1) through (3), $\{s_{p,p+1}\}= s_{ij}$ for all (p,p+1), and in the case of (4), $s_{k,k+1} = -|s_{ij}|$, and $s_{q,q+1} = |s_{ij}|$.

Shown in Figure 1.3.15 is the translated FTCM of Figure 1.3.13 at the time that minimum operator is used for propagating indirect effects. Consequently, the use of dummy nodes translates the two matrices, S and T, into only one strength matrix, $S^*$, represents the causal strengths per one unit-time, we can use the equations of (1) and (2) for analyzing system behaviors according to lapse of time. Additionally, the size of the matrix $S^*$ would be greater than the size of the starting matrices, S and T. This entails an increased use of memory storage in the computer, but the computation for causal inferences is greatly facilitated.

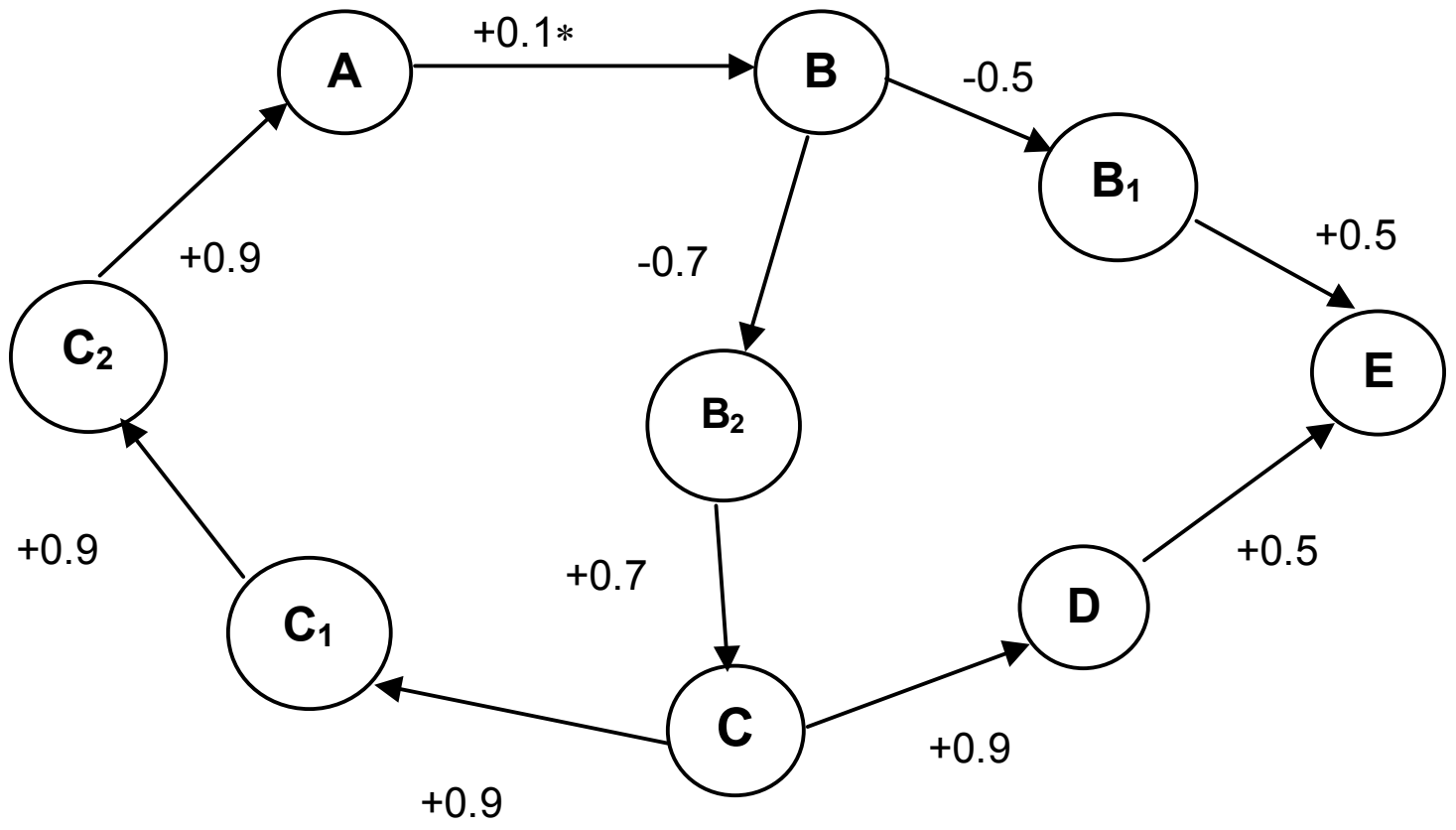

**FIGURE: 1.3.15**

Park and Kim [78] have proposed a FTCM that has different time relations on arrows among factors. We also presented a method of translating the FTCM into the FTCM with one unit-time lag. With the FTCM, we can analyse the change of causalities among factors for some discrete time t.

Readers who wish to learn more about FTCMs can refer [78] and other related papers.

### 1.3.11: Balanced Differential Learning Algorithm in FCM

Vazquiz [131] has studied about a balanced differential learning algorithm in FCM. Cognitive maps (CM) also called Causal Maps; a special kind of CM are the Fuzzy Cognitive Maps. The negative causal connections between A and B are described in the following Figure 1.3.16.

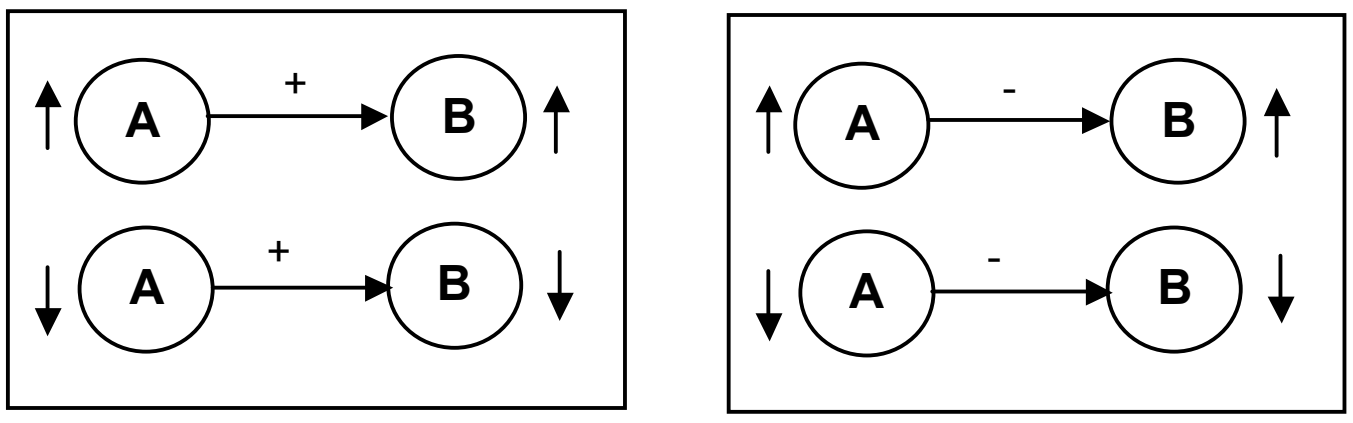

POSITIVE CAUSAL WEIGHT NEGATIVE CAUSAL WEIGHT

**FIGURE: 1.3.16**

The following table gives the linguistic labels used for the domain experts and a possible automatic translation of weights in the causal web of the FCM.

| Symbolic Values | Numeric Values |
|---|---|
| Affects a lot | 1.0 |
| Affects | 0.5 |
| Does not affect | 0.0 |
| Affects negatively | –0.5 |
| Affects negatively a lot | –1.0 |

A FCM using the above kind of table can be easily designed by the human domain experts and can be used to have a graphical model of the domain. Sometimes it is not possible to have a human expert to model the domain because the number and complexity of the variables involved are so large that the task to build by hand the FCM is very difficult. In this case the solution is to use adaptive fuzzy cognitive maps. An adaptive FCM learns its causal web from data. We have just mentioned a section the learning algorithm in FCM as given by [131]. As the main motivation is to apply to it in the neutrosophic setup we have not dealt it in this book elaborately but request the reader to refer [131] and other related papers.

## 1.4 Applications of FCM

Here we collect some of the applications of FCM as given by researchers. The notion of FCM is used in modeling of supervisory systems, in the design of hybrid models for complex systems, in the International Advanced Robotics Program, in business performance assessment, in analyzing legal problems etc. In fact the illustrations given in the earlier sections can also be viewed and studied as applications.

It is pertinent to mention here that FCM model has also been applied in case of transportation problems like passengers preference, peak-hour problem, stock investment analysis, etc.

### APPLICATION 1.4.1: Modeling Supervisory systems using FCM

In the application of FCM, Stylios and Groumpos [100] have given a new methodology, an FCM for modeling the supervisor of a complex control system. An FCM has been implemented in a simple process control problem that makes apparent the qualities and characteristic of the method. It has been observed how simply an FCM describes a system's behaviour and its flexibility in any change of system.

A more integrated approach i.e. a two-level structure where the FCM is the upper level, is used for more sophisticated supervisory control of manufacturing systems which is given below:
A two level structure approach

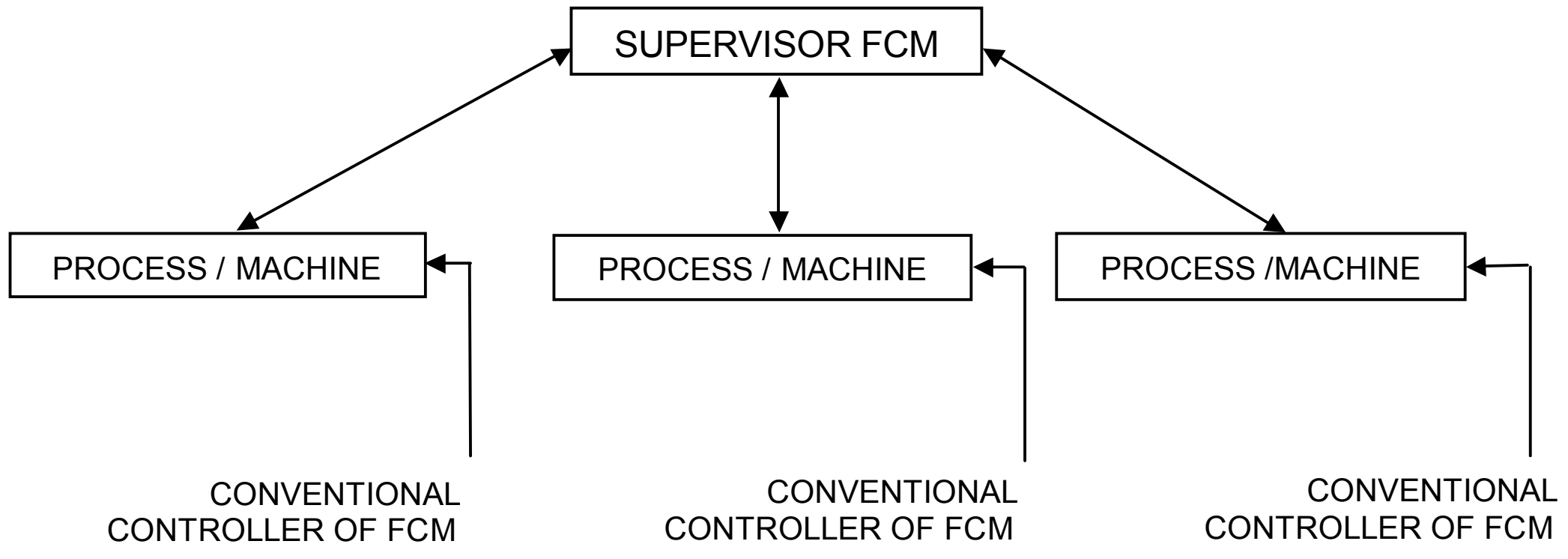

**FIGURE: 1.4.1**

The most important use of an FCM is for supervisory control of a conventional control element, thus complementing rather than replacing a conventional controller. Thus the FCM may replicate some of the knowledge and skills of the control engineer and it is built using a combination of knowledge representation techniques as causal models, production rules and object hierarchies and it is used to perform more demanding procedures such as failure detection, decision making and planning (tasks usually performed by a human supervisor of the controlled process).

In the upper level a supervisory FCM will include advanced features such as fault diagnosis, effects and cause analysis, prediction capabilities, decision analysis and strategic planning. The FCM will consist of concepts that stand for the irregular operation of some elements of the system, for failure mode variables for failure effects variables, for failure cause variable, severity of effects, design variables. Thus [100] have used the construction of a map mainly based on the operators heuristic knowledge about alarms, faults, what are their causes, and when they happen.

For more about the application of FCM in supervisory control please refer [100].

#### APPLICATION 1.4.2: FCM applied in the Design of Hybrid Models for Complex Systems

A Fuzzy Cognitive Map is used to aggregate separate models and to fit more precisely the plant behaviour at different operational conditions. One of the important advantages of fuzzy logic methodologies such as FCMs is that they are applicable to task oriented problems.

Here [39] FCMs are used to aggregate the set of different modeling technologies. FCMs best utilize existing experience in the operation of the system and are capable in modeling the behaviour of complex systems. FCMs seem to be a useful method in complex system modeling and control which will help the designer of a system in decision analysis. They are an appealing tool in the description of the modeling techniques, which teamed up with other methods will lead to be more sophisticated model and control design systems.

FCMs are used to aggregate the separate models and to perform a kind of maintenance of the system by integrating alternative modeling techniques. An augmented FCM can accomplish identification of the process models and cope with limited uncertainty situations. It may comprise different models, identification and estimation algorithms.

Using FCMs, which best utilize existing experience in the operation of the system and are capable in modeling the behaviour of the complex systems. FCMs seem to be a useful method in complex system modeling and control which will help the designer of a system in decision analysis. FCMs appear to be an appealing tool in the description of the modeling techniques, which teamed up with other method will lead to the more sophisticated model and control design systems.

For more about the application of FCM in the design of hybrid models for complex systems please refer Hadjiski et al [39].

**APPLICATION 1.4.3: Use of FCMs in Robotics: Mobile Robots performing Cutting Operations and Intimate Technologies like Office Plant**

In a bulletin released following the International Advanced Robotics Program (IARP) held in Roma, Italy between Nov 6-9 2002, it is said as follows: "FCMs are applicable for modeling the operation of mobile robots performing cutting operations. FCMs represent an approach to model the behaviour and operation of complex non-linear systems. The implementation of FCMs for modeling the mobile robot motion and cutting process are presented and the corresponding models developed." [85] Thus FCMs prove to be a very useful tool in mobile robots.

Bohlen M. and Mateas M. make use of FCMs in their research on Office Plant. Office Plant #1 (OP#1) is an exploration of a technological artifact, adapted to the office ecology, which fills the same social and emotional niche as an office plant. OP#1 monitors the ambient sound and light level, and, employing text classification techniques, also monitors its owner's email activity. Its robotic body, reminiscent of a plant in form, responds in slow, rhythmic movements to comment on the monitored activity. In addition, it makes its presence and present knowledge known through low, quiet, ambient sound. OP# is a new instantiation of our notion of intimate technology, that is, technologies which address human needs and desires as opposed to technologies which meet exclusively functional task specifications. OP#1 lives in a technological niche and interacts with users through their use of electronic mail. It acts as a companion and commentator on these activities.

In this section we describe the major artistic and technical concepts that underlie the design of OP#1. These concepts are: E-mail Space, Text Classification, Plant Behavior Architecture, and Sculptural Presence. In our practice we simultaneously explore both spaces; artistic and technical constraints and opportunities mutually inform each other.

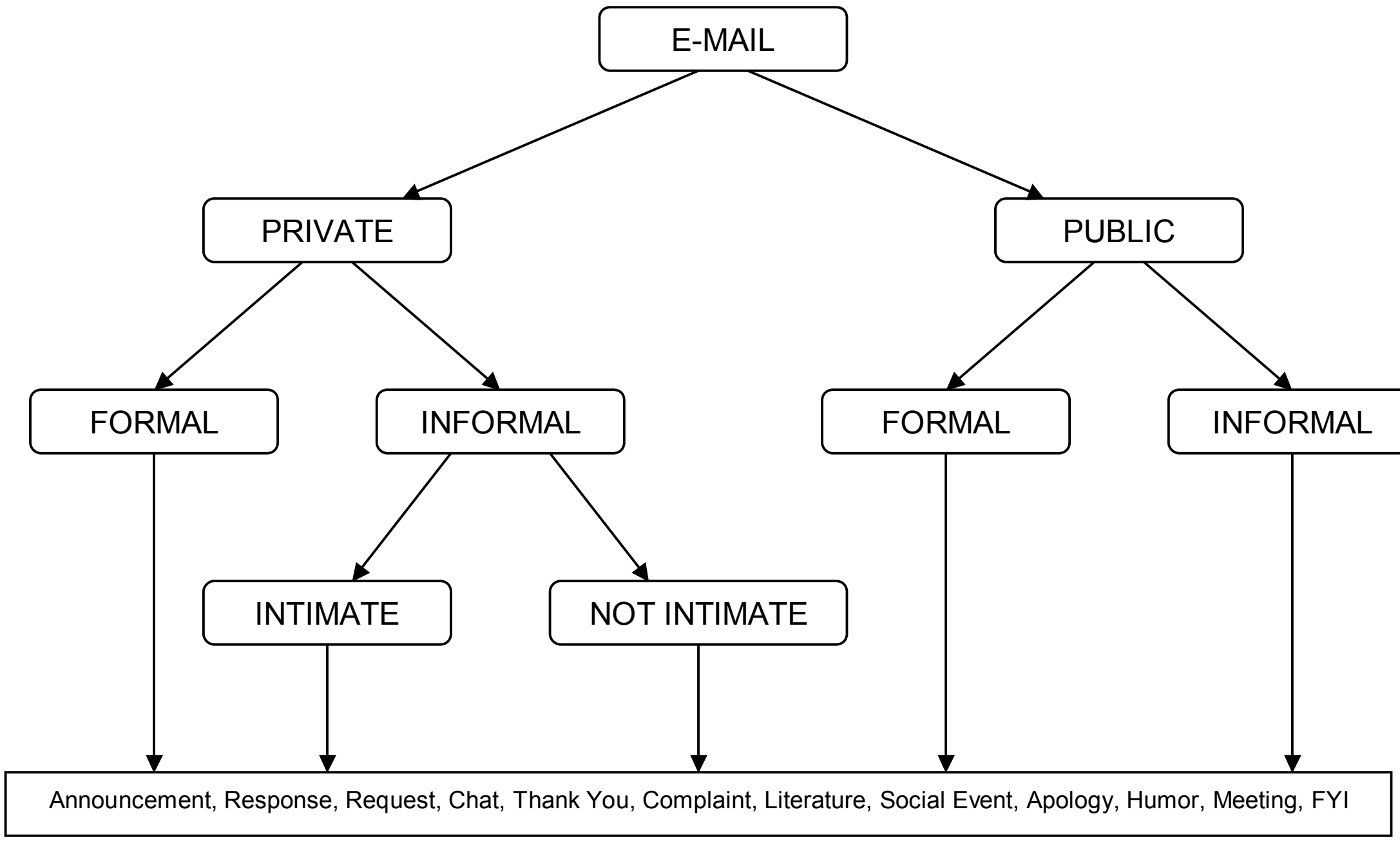

**FIGURE: 1.4.2**

Figure 1.4.2 depicts the category tree employed by OP#1. An email message is either private (addressed to a single person) or public (multiple address). The tone can be either formal or informal. Private, informal email can be intimate, that is, email addressed to close friends. After passing through this initial category tree, every message can be assigned one or more of the categories in the box at the bottom of the figure. In this categorization scheme, every message is assigned a set of labels. For example, a message may be a public, informal announcement, or a private, informal, humorous request

The state of the plant is dynamically modeled with a Fuzzy Cognitive Map (FCM) [58]. In a FCM, nodes representing actions and variables (states of the world) are connected in a network structure reminiscent of a neural network. FCMs are fuzzy signed digraphs with feedback. Nodes stand for fuzzy sets or events that occur to some degree. At any point in time, the total state of the system is defined by the vector of node values. In this implementation, the nodes represent actions. The action associated with the action node with the highest value is executed at each point in time. The values of nodes change over time as each nodes exerts positive and negative influence on the nodes it is connected to. The FCM approach is attractive because it can resolve contradictory inputs and maintains sufficient state to exhibit incremental effects. Figure 1.4.3 shows the FCM for OP#1.

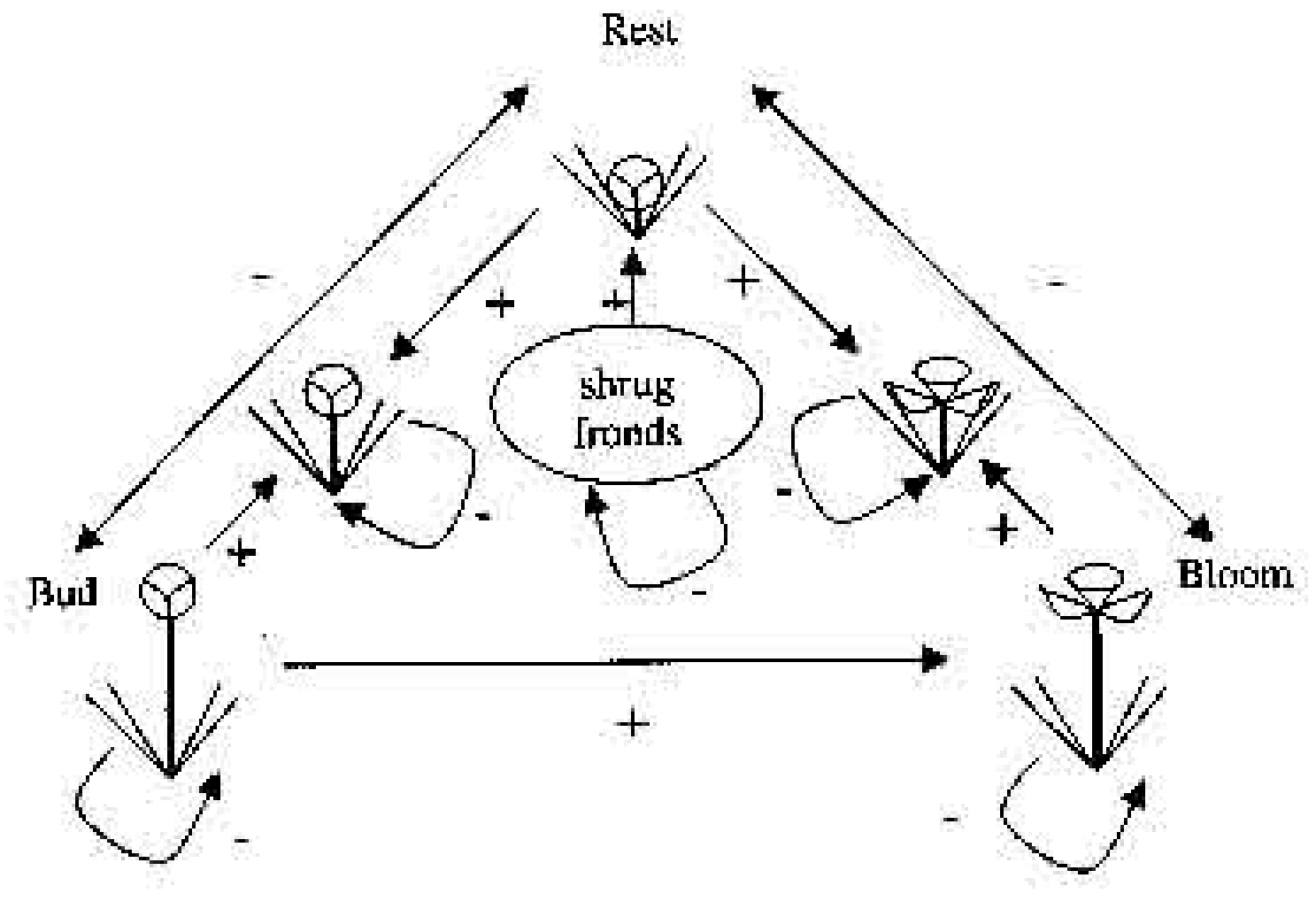

**FIGURE: 1.4.3**

Office Plant #1 is a desktop sculpture, an office machine that serves as a companion. In an environment of button pushing activity, OP#1, like a good piece of sculpture, is always on. OP#1 creates its own kind of variable presence in a user's email space, taking on various attitudes and falling into decided silence. In its reactions it is a commentator on the events it analyzes. It goes beyond mirroring these events and delivers reactions as if it understood the significance of the exchange. But effectively, OP#1 is mostly inactive. It has a well defined sense of doing nothing, yet. It is simply there and in that sense a tradition piece of sculpture. Its physicality is as important as its text classifying capabilities.

OP#1's activity cycle is given by a defined period of 24 hours. During the active office hours it is receptive to user presence. After hours it uncouples itself from the daily trivia as it moves into a contemplative space for regeneration.

Figure 1.4.3 shows the Fuzzy Cognitive Map relating physical plant states. The three primary physical postures of the plant are rest (protect), bud and bloom. In rest, the bulb is closed and fully lowered. The fronds occasionally move. In bud, the bulb is closed and fully extended. In bloom, the bulb is open and fully extended. When the node associated with a posture has the most activation energy, the plant performs the action of moving to this posture from its current posture frond is the action of tweaking one or more of the fronds. Activation energy from shrug flows towards rest. If the plant is shrugging too much, it moves into a protective posture. Activation energy from bud flows towards bloom; budding makes blooming more likely. Rest and bud, and rest and bloom, are mutually inhibitory. Rest and bud both spread their energy to an intermediate posture, and rest and bloom spread their energy to a second intermediate posture.

The combination of the mutual inhibition plus the intermediate posture will cause these pairs of states to compromise towards the intermediate posture will cause these pairs of states to compromise towards the intermediate posture. Finally, the self-inhibitory links tend to cause values in the system to decay; in the absence of input, the plant will not stay in a posture forever. When all of the nodes are zero, the plant will move towards the rest posture. As email is classified, energy is added to nodes, thus initiating the process of competition and cooperation between the nodes.

### APPLICATION 1.4.4: Use of FCMs to model & Analyse Business Performance Assessment

Kardaras and Mentzas [48] introduce Fuzzy Cognitive Maps (FCMs) as a computational modeling formalism which tackles the complexity and allows the analysis and simulation of the metrics.

The FCM proposed by Kardaras and Mentzas is called Impact Analysis Model (IAM) and depicts the strength (expressed in fuzzy terms such as weak, strong, etc.) and polarity (positive or negative) of the impact of selected business metrics from business areas such as competition, customers, products, suppliers, risk assessment, organizational effectiveness, etc., on the business performance. Reasoning mechanisms taken from the theory of FCMs, allow the deduction of additional causal relationships between variables, which can be used to assess the impact of Information Systems (IS) on a particular organizational area. IAM provides a comprehensive knowledge base of organizational and IS factors and their interrelationships which can be used as the foundation for the development and testing of a Business Metrics theory. Our approach, is highly customizable, enables users' participation and brainstorming, and allows the development of alternative scenarios with regard to current organizational situation assessment and alternative Business and IS strategies.

Kardaras and Mentzas [48] discuss a framework for developing business metrics. Their proposal reflects a broader view of traditional performance indicators which are mainly drawn from studies of time and motion or financial theories. The proposed business metrics are drawn from domains such as organizational effectiveness, business and Information Systems (IS) strategy, human aspects and organizational behaviour, IT issues, financial indicators, etc.

The business performance factors can be quite broad and therefore relevant to more than one business areas. There are factors which address the strategic organizational options and factors which represent choices and performance indicators at the operational business level. There are also factors which can be easily measured and others less tangible. In such cases our approach uses fuzzy terms (e.g. at the strategic level, where a threat from a competitor is "high") or surrogate measures (e.g. an indication of the level of the employees morale can be the absenteeism, accidents, etc.) which can be quantified. Figure 1.4.4 depicts the areas where metrics are drawn from.

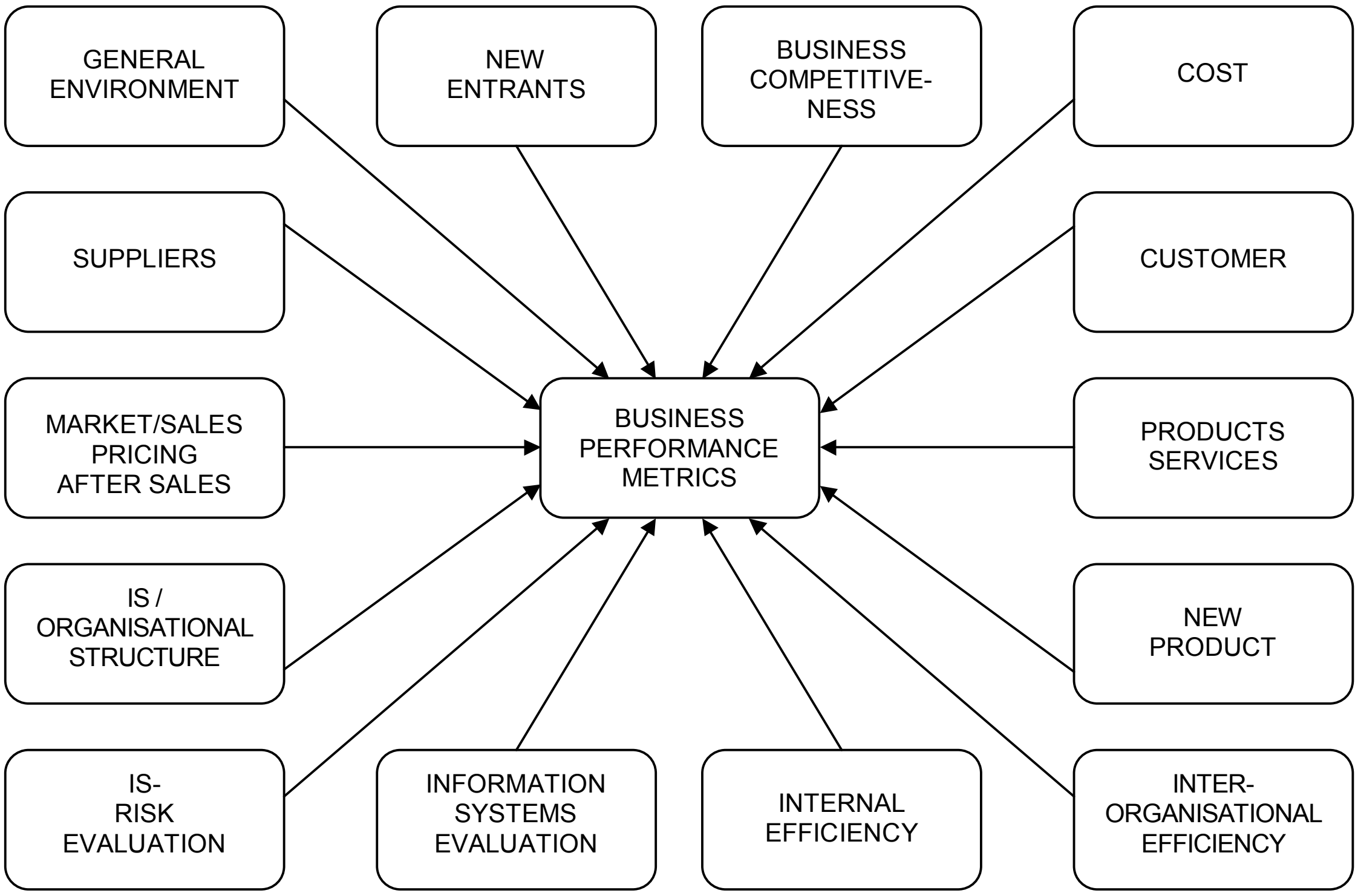

**FIGURE: 1.4.4**

For each of the business areas a set of Business Metrics (BM) is developed as a different sub model. Sub models are integrated however in order to be able to examine the impact from one area to all the others which it is related to. The decomposition of the (BM) into separate sub models according to business areas enables the close examination of specific organizational areas and to cope with the complexity.

The IAM framework distinguishes two levels of business metrics namely those at strategic level and those at process level. Three organizational perspectives namely: the Business; the Technological (i.e. IT), and the Social which are so modeled that they integrate the above three organizational domain perspectives.

Each of the organizational perspectives is considered from an Environmental Scanning View i.e. metrics reflect both an Internal and an External organizational point of view. The external view focuses on metrics which are drawn from domains such as competition, customers, etc. The internal one, focuses on the internal organizational effectiveness, e.g. flexibility of production, decision making efficiency, co-ordination and communication among processes or departments, etc.

Two metrics modes namely: the Options; and the Actions. Each one of the above views indicates the use of the metrics. Thus, the IAM metrics framework models the three organizational perspectives at both levels which are viewed from two different views and used for two different purposes.

More specifically, the business perspective considers issues such as business strategic objectives and business processes objectives; consistency with strategic goals; responsiveness to the external environment.

For an exhaustive description of the topic, the readers are requested to refer [48].

The IAM framework considers two different metrics modes with each one indicating the use of the metrics.
The Options Mode, indicates the use of the metrics for either,

- Competitive purposes i.e. metrics indicating organizational competitiveness, such as bargaining power, or
- Diagnostic purposes, e.g. metrics for identifying problems such as to much control, bureaucracy, etc.

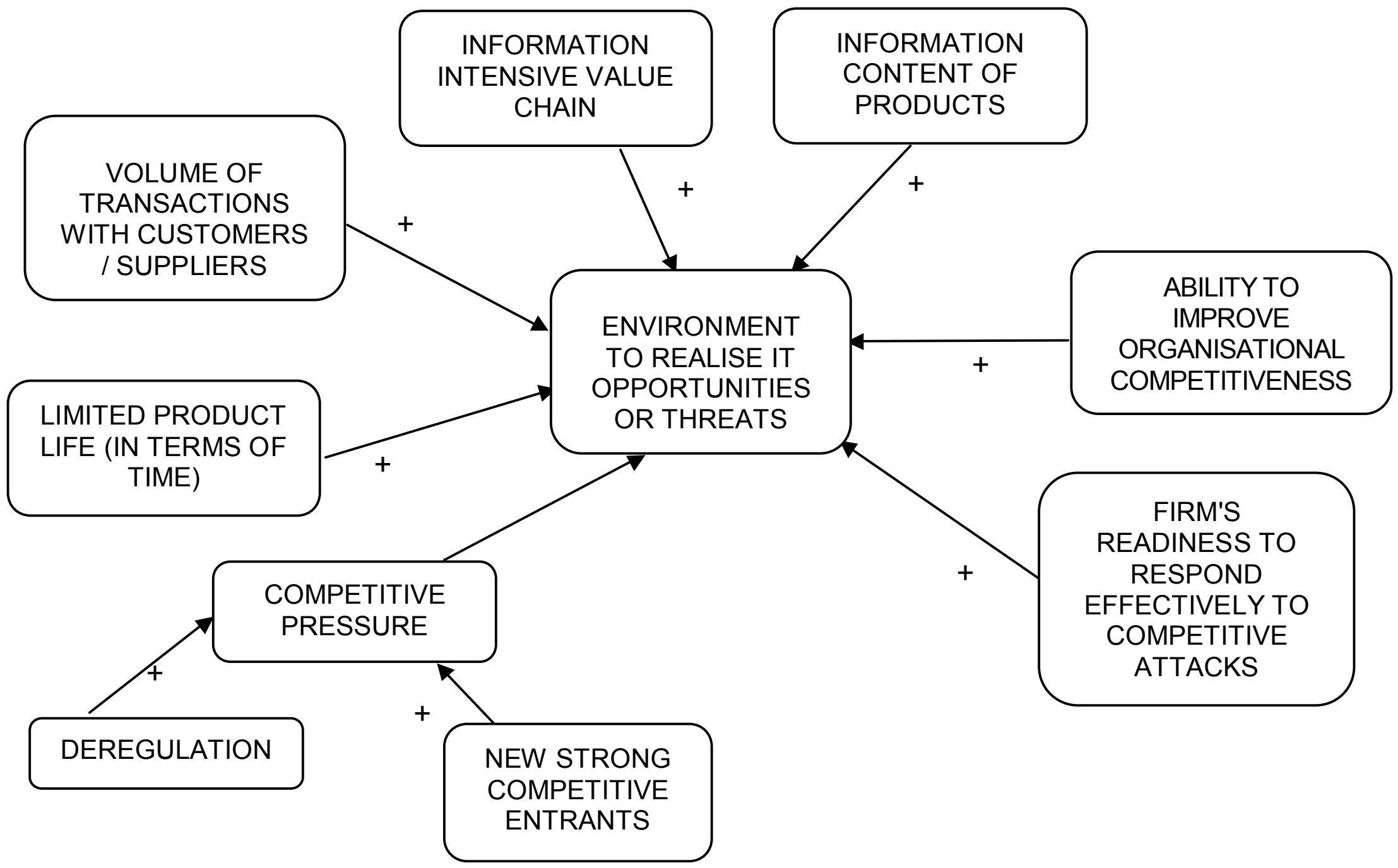

**FIGURE: 1.4.5**

Figure 1.4.5 shows an example of a FCM, which can be used to assess the business environment of an organization. It can be interpreted as follows: an increase in the cause variable "Volume of Transactions" with "Customers/Suppliers" affects positively, (since the sign of relationship is +), the effect variable "Environment to Realise IT Opportunities or Threats". The fuzzy weights are determined according to beliefs shared within a team responsible for business metrics development.

The business metrics approach developed by [48] reflects the following:

- Business strategy which is translated into business processes objectives i.e. time, cost, customers satisfaction, organizational and process effectiveness, etc.
- The participants' objectives and requirements, i.e. the human aspects and perspective of the business process, i.e. job satisfaction, individual and group effectiveness, etc.
- The technology. Our approach recognizes the significance of IT as an integral part of business process and important contributor to business success. The technology will be assessed in terms of: organizational fit, impact on business process, i.e. objectives, people, technology, overall organizational impact, etc.

The reader is expected to refer [48] for complete working of the paper which is innovative and interesting.

#### APPLICATION 1.4.5: Use of FCM in Legal Rules

Here we have provided a complete extract taken from Adams and Farber [1] for we find it thrilling to see how Fuzzy logic is used in Law.

"Judges tinker with the legal system all the time. Frequently forced to apply statutes in factual situations never dreamt of by the drafters of such legislations, judges labor to make the relevant statutory language "fit" the case at issue. The results of such tinkering or interpretation vary widely. Sometimes such tinkering produces out comes that are both unintended and unacceptable, leading to disastrous results. Other times judges, like the rest of us, get it right.

"Statutory intepretation is currently the subject of a lively scholarly debate. Should judges render a decision based merely on the words of the text before them or should they go one step further and attempt to discern the purpose of the legislation at issue? Is it legitimate for them to consult current social values? On one side of the debate are formalists, who eschew legislative history and current social values while being suspicious of the concept of statutory purpose. On the other side are antiformalists, who prefer to downplay textual arguments in favor of these other sources of guidance. Anti formalists endorse "practical reason"-meaning reliance on complex judgments regarding text, purpose, legal context, and societal norms.

"Although this is an important debate, the dispute between these viewpoints has obscured significant aspects of statutory interpretation while highlighting others. As a result, analysis of statutory interpretation has been skewed. First, the debaters have focused on the differences between formalism and antiformalism (or "practical reasoning"), ignoring considerable areas of overlap. Second, much of the discussion has involved stylized examples, often in the context of fairly simple statutes, rather than the complexities presented by actual judicial decisions under complicated modern statutory schemes. Third, neither side has a very clear model of the process that judges actually use to make decisions. […]

"This Article seeks to provide more functional models of judicial decision-making. Normally, formalists and antiformalists alike resort to “hand waving” when seeking to describe the judge’s cognitive process. This Article will attempt to provide some substance to its description of this process by referring to work by cognitive psychologists on expert judgments and to an emerging field known as fuzzy logic. The latter term-which sounds like a professor’s disparagement of an inept first-year student-requires a brief explanation.

"Two leading formalists, Justice Antonin Scalia and Judge Frank Easterbrook wrote the opinions. Specifically, Part I analyzes and critiques the Seventh Circuit decision in Levit v. Inger-soll R and Financial Corp. (Deprizio) and the Supreme Court’s decision in BEP v. Resolution Trust Corp., which utilize different varieties of formalism. Through an in-depth analysis of these cases, we probe the complex interactions between judicial interpretation, the business and legal communities, and the legislature."

For more about the facts of the case and the formalism debate, please refer [1].

Fuzzy logic provides a model of how to analyze such conflicts.

The example that follows develops an FCM virtual world based on the principal question in Deprizio: “whether the Trustee may recover from an outside creditor under section 550(a) (1) a transfer more than 90 days before the filing that is avoided under section 547 (b) because of a benefit for an inside creditor". The Augmented Deprizio FCM consists of Simple FCMs and Nested FCMs. Simple FCMs describe the judicial process: the trustee’s argument, the creditors’ argument, and Judge Easterbrook’s reasoning. Nested FCMs describe the degree of risk that the insiders and the creditors bear. Together, these simple FCMs describe the judicial process and the “Real World” of creditor, insider, and firm relationships. The purpose of this section is not to reach a definite conclusion about Deprizio but rather to introduce the FCM as a tool for understanding and perhaps predicting how judicial decisions affect business behavior.

The trustee’s argument is linear and based on a literal reading of the Bankruptcy Code. As a simple FCM, eight nodes arranged in four layers can describe the trustee’s argument. The first layer contains one node (T1) which simply represents the fact that the trustee has an argument. The second layer consists of facts specific to the case: Guarantor is an officer (T2); Guarantor has a claim (T3); and transfer from Firm to Creditor benefits Creditor (T4). The third layer of the trustee’s FCM incorporates Bankruptcy Code provisions to reach intermediate inferences: 101 (30) (B) (ii) transforms node T2 to node T5 (Guarantor is an insider) ; 101 (9) and 101 (4) (A) transform node T3 to node T6 (Guarantor is a creditor); and 547 (b) (4) transforms node T4 to T7 (Guarantor’s benefit is avoidable).

The fourth layer combines the intermediate inferences to reach the final conclusion that the trustee sought: "Lender may have to repay transfers received during the year before filing [for bankruptcy], even though Lender is not an insider." The trustee's simple FCM is shown in Figure 1.4.6.

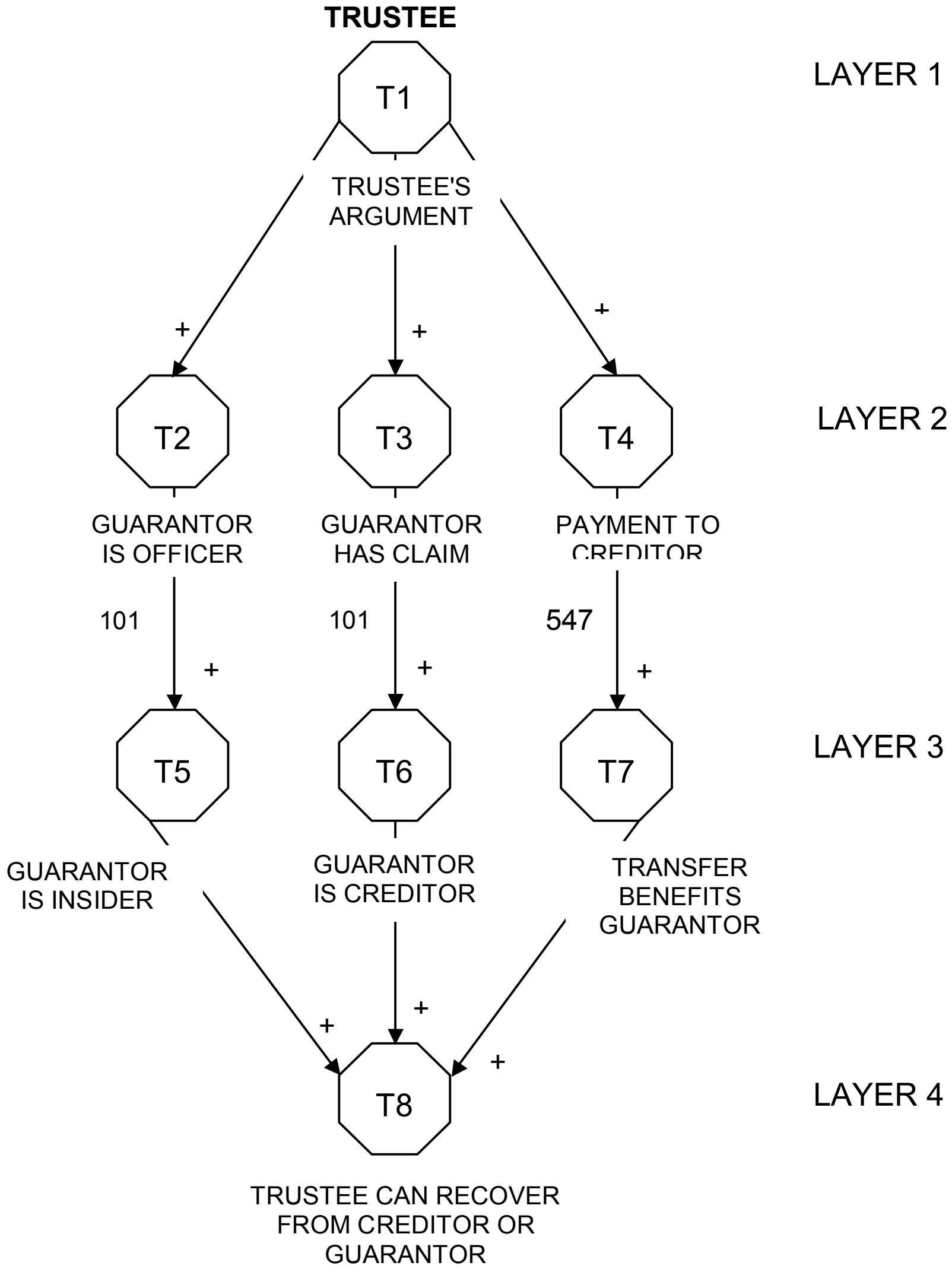

**FIGURE: 1.4.6**

The creditor's argument on the principal question relies on a two-benefit/two-transfer theory. The creditor's FCM, like that of the trustee, has four layers. However, only six nodes are used to model the creditor's argument. The first layer contains one node (C1) which imply represents the fact that the creditor has an argument. The second layer contains the basic theory: the payment benefits the Lender (C2) and the Guarantor (C3). The third layer suppose that a benefit equals a transfer and it uses the Code to determine if the resulting individual transfers are "avoidable". According to 547(b) (5), a transfer is avoidable only to the extent it gives the creditor more than the creditor would have received in a liquidation. Using this definition and the fact that the Guarantor's interest is junior to the Lender's, the Creditor's argument produces concept nodes C4 (transfer to Lender is unavoidable ) and C5 (transfer to Guaranotr is avoidable). Using 550 (a), the Creditor arrives at the fourth layer node C6. The node C6 represents the proposition that the trustee can recover only from the Guarantor for transfers made less than one year and beyond ninety days from filing. Figure 1.4.7 shows the Creditor's simple FCM.

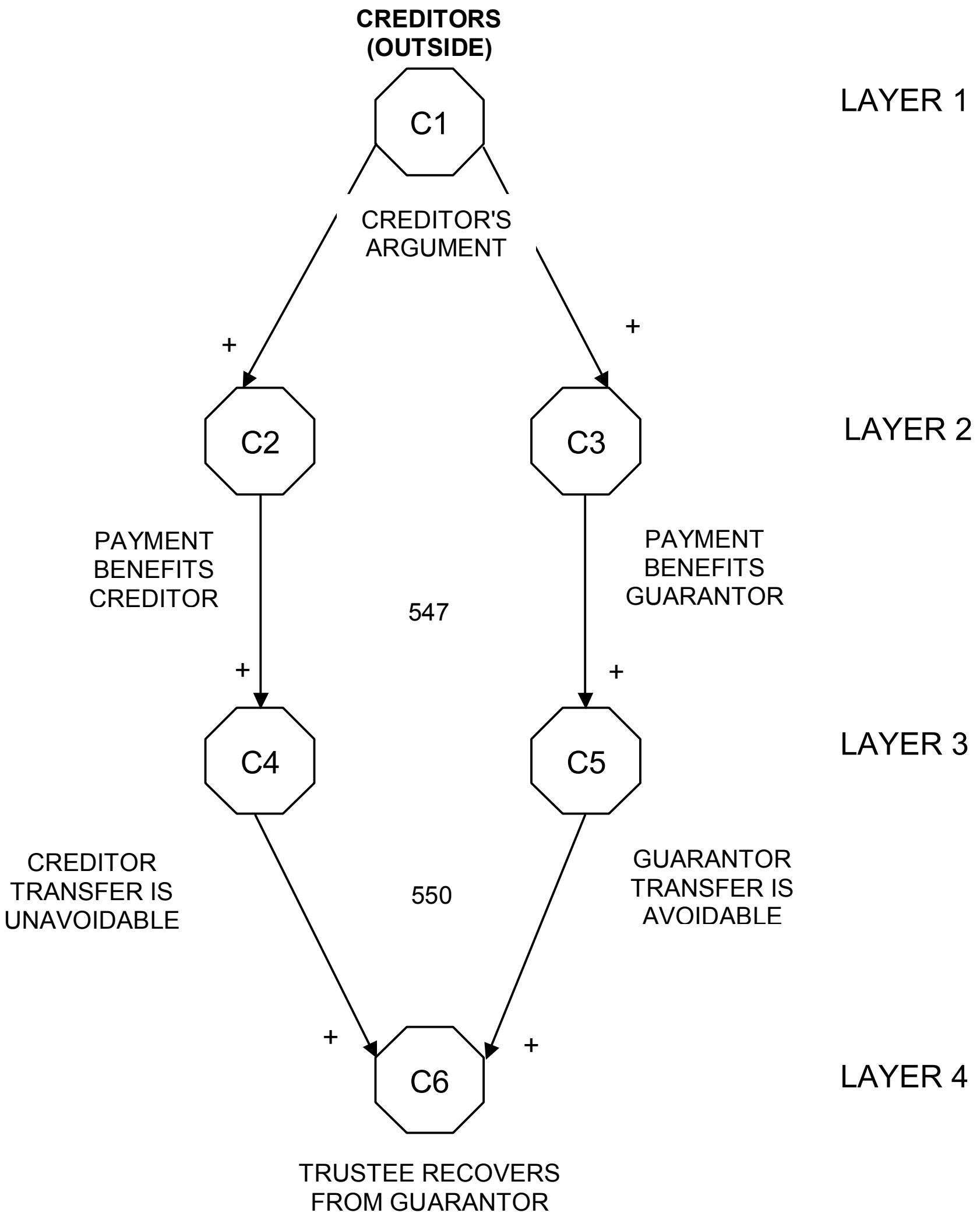

**FIGURE: 1.4.7**

A three layer simple FCM can likewise model Easterbrook's reasoning on both arguments. Easterbrook approaches the principal question using textualism and economic theory. These are represented in the first layer as nodes E1 (textualism) and E2 (economic theory). The second layer consists of three nodes. Node E3 represents Easterbrook's interpretation of 101 (9): a claim against the firm makes one a creditor. Node E4 represents his textualist reading of 547 (b) (1) and 101 (50): a payment is a transfer, a benefit is not a transfer. Node E5 contains Easterbrook's notion transforms economic node E5 to node E6 which asserts "all creditors gain from a rule of low that induces each to hold back. Figure 1.4.8 shows Easterbrook's simple FCM.

To achieve an augmented FCM representation of the entire case, the three simple FCMs are placed side by side. The next step is to link the FCMs with appropriate edge weights {-1, +1}. Once linked, the augmented FCM connection matrix can be formed from the individual FCM connection matrices. Converting the FCMs to connection matrix form will become apparent later as it facilitates the programming, or "bookkeeping", of more complex relationships in larger augmented FCMs. The matrices for this model can be found in [1].

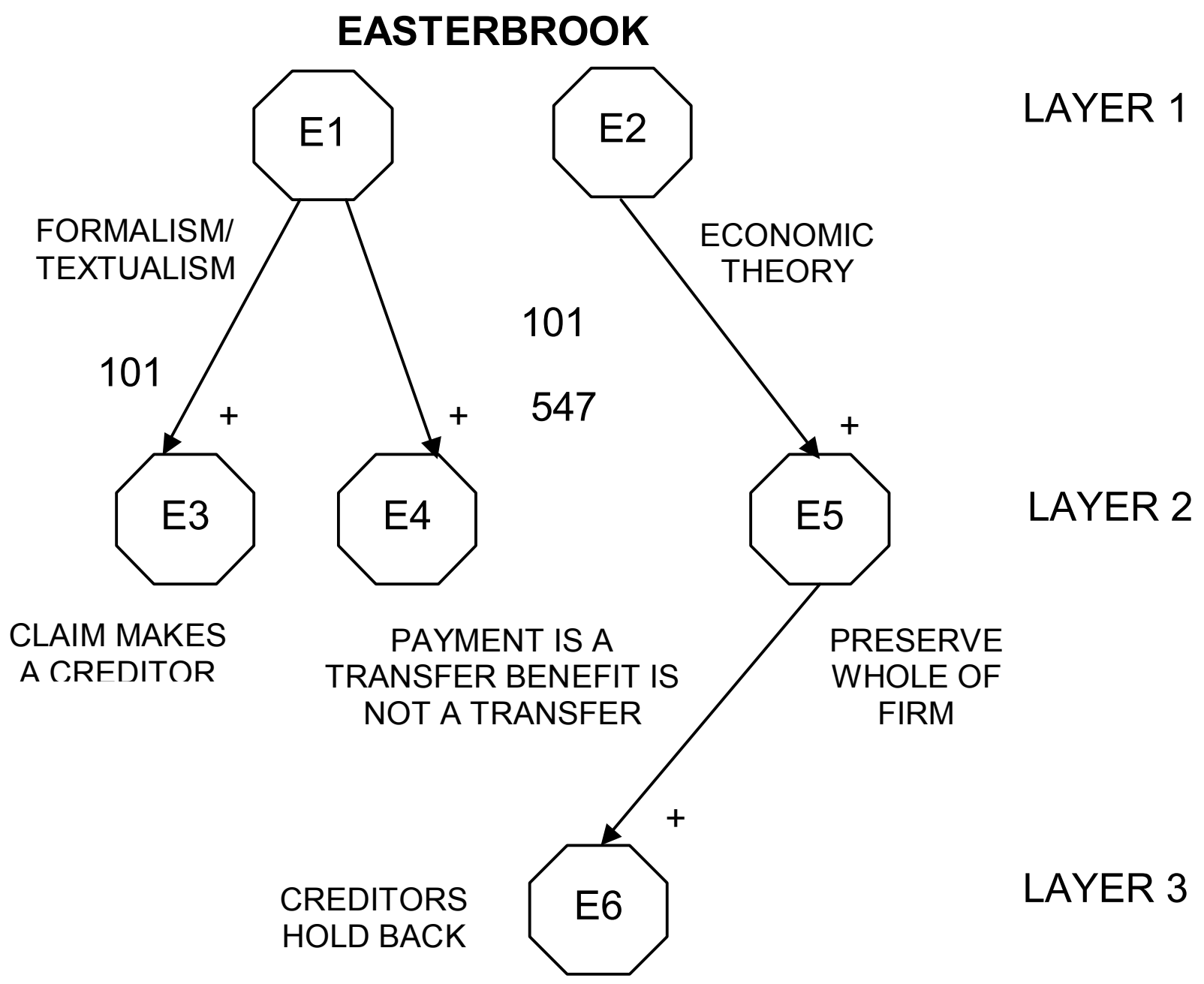

**FIGURE: 1.4.8**

When a firm borrows money, both the creditor and the firm engage in a transaction that reflects their respective levels of risk. The level of risk depends on many variables, including the prevailing rule of law. Deprizio acts to heighten the creditor's risk level, extending the preference window from 90 to 365 days. The heightened level to risk induces creditors to ask firms to waive claims to their assets in the case of bankruptcy. In turn, the insider who has issued bank guarantees experiences a heightened level of risk.

The insider's perception of risk will influence the firm's decisions on borrowing money. On the other hand, if "Barking Dog" prevails, then the creditor's risk is low and the creditor is unlikely to convince the firm's insiders to waive their rights to assets.

The creditor risk concept, elaborated as subconcepts of high (N5) and low degrees of risk (N6), when linked to the insider's waiver of a claim to assets, forms a nested simple FCM. The insider risk concept, also elaborated as subconcepts of high (N1) and low degrees of risk (N2), links to the waiver and the decision for the firm to issue guarantees. The nested FCMs are shown below linked to their common nodes, N3, "Insiders Waive Claim to Firm's Assets" and N4, "Firm Issues Guarantees". Note that the insider actions affect insider risk. Thus, if the firm issues guarantees or the insiders waive claims to assets, the degree of insider risk rises. Likewise, if these actions do not occur, the degree of insider risk falls.

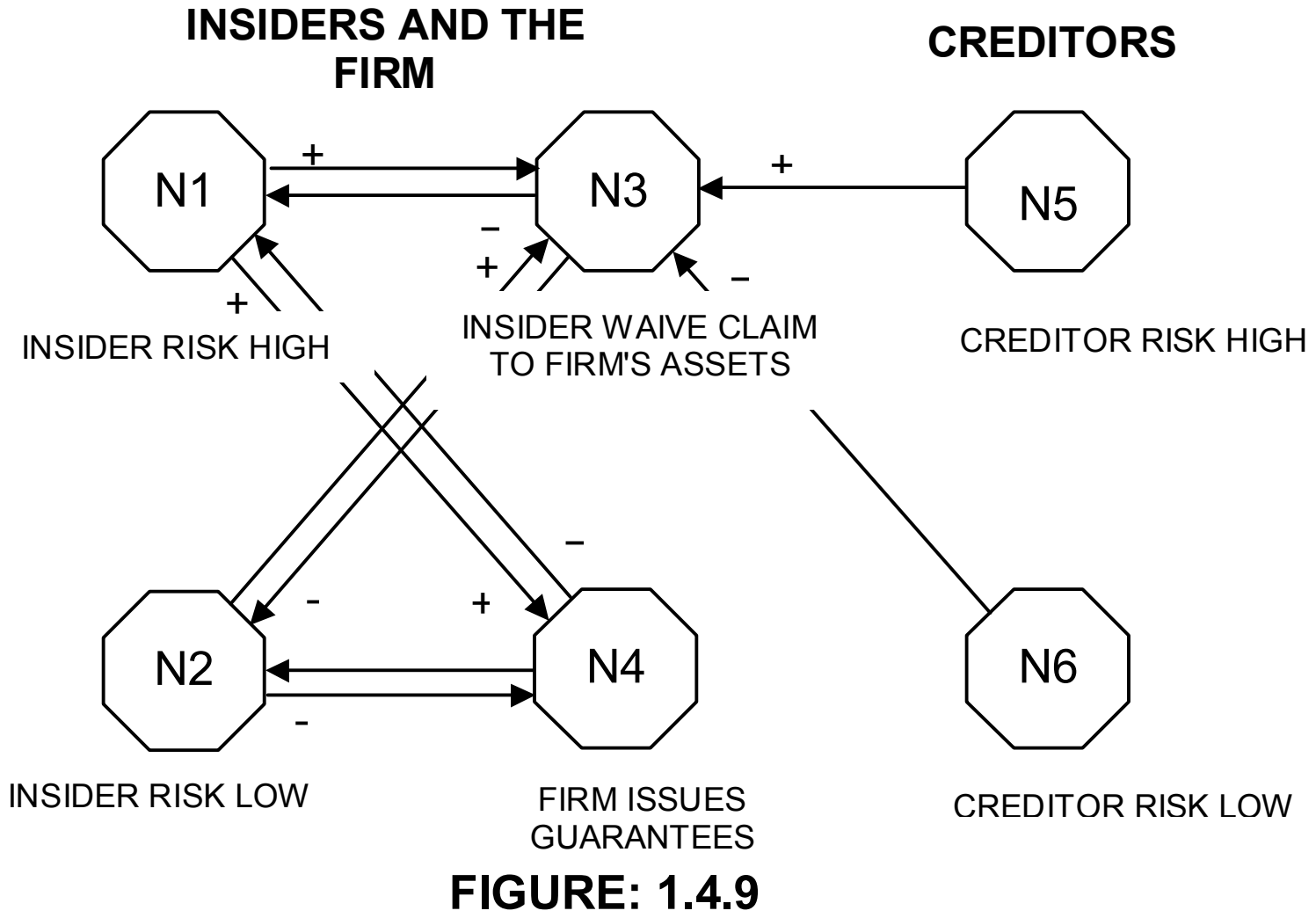

**FIGURE: 1.4.9**

As shown above, the edge weights are represented as pluses and minuses or {–1, 1}. However, this need not be the case. Any fuzzy weight can be assigned to the edges, representing the degree to which the risk rises or falls. Fuzzy reasoning models, as presented for the BFP case, can be used to determine the weights. In a similar manner, weights can be used to determine the weights. In a similar manner, weights can be assigned to the edges affecting the center nodes. Thus, one can investigate the degree to which insiders are willing to issue guarantees or the degree to which the insiders are willing to waive claims to assets. The connection matrix and the mathematical procedure for determining its properties are found in [1].

The simple risk FCM, as with the simple FCMs for the Deprizio decision, can also be represented in matrix form. However, to form a more complete model of the real world affected by the Deprizio decision, it is helpful to introduce a few more concepts. Preferential transfers to insider-backed creditors lower insider risk levels. The other two important concepts are the rules of law before and immediately after Deprizo. Depending on the state of the law, creditor risk either will be high or low and preferential transfers either will occur or will not.

As an augmented FCM, the "Real World" of creditor and firm relations appears in Figure 1.4.10.

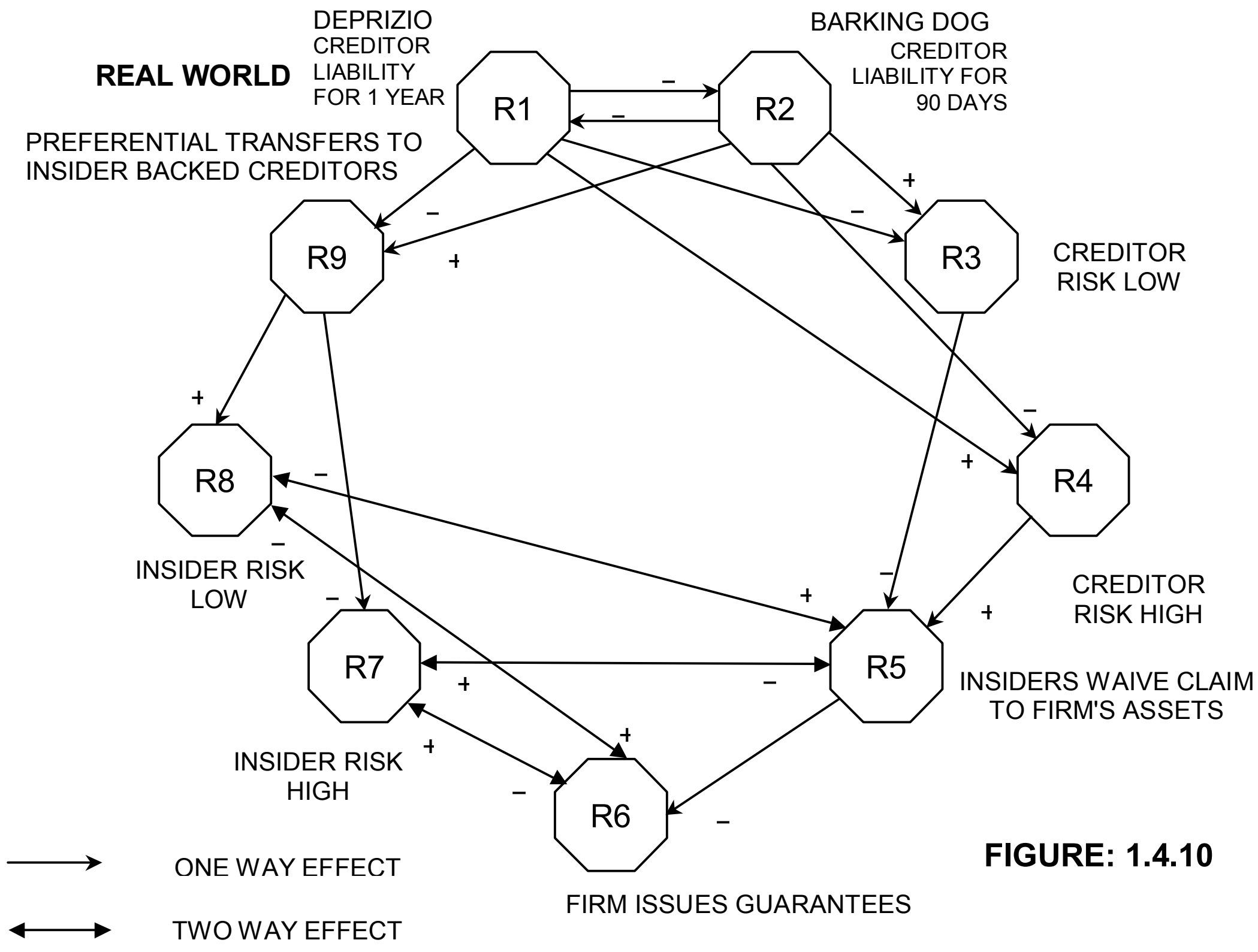

FIGURE: 1.4.10

The bridge between the judicial decision FCM and the Real world FCM in traversed by linking nodal concepts in the Trustee's and Creditor's arguments to the Real World. The links are straightforward : T8 → R1, {+1} and C6 → R2, {+1}. However, with Easterbrook's textual and economic theories intervening, only R1 will fire. The FCM for the judicial decision and the Real World is shown in Figure 1.4.11.

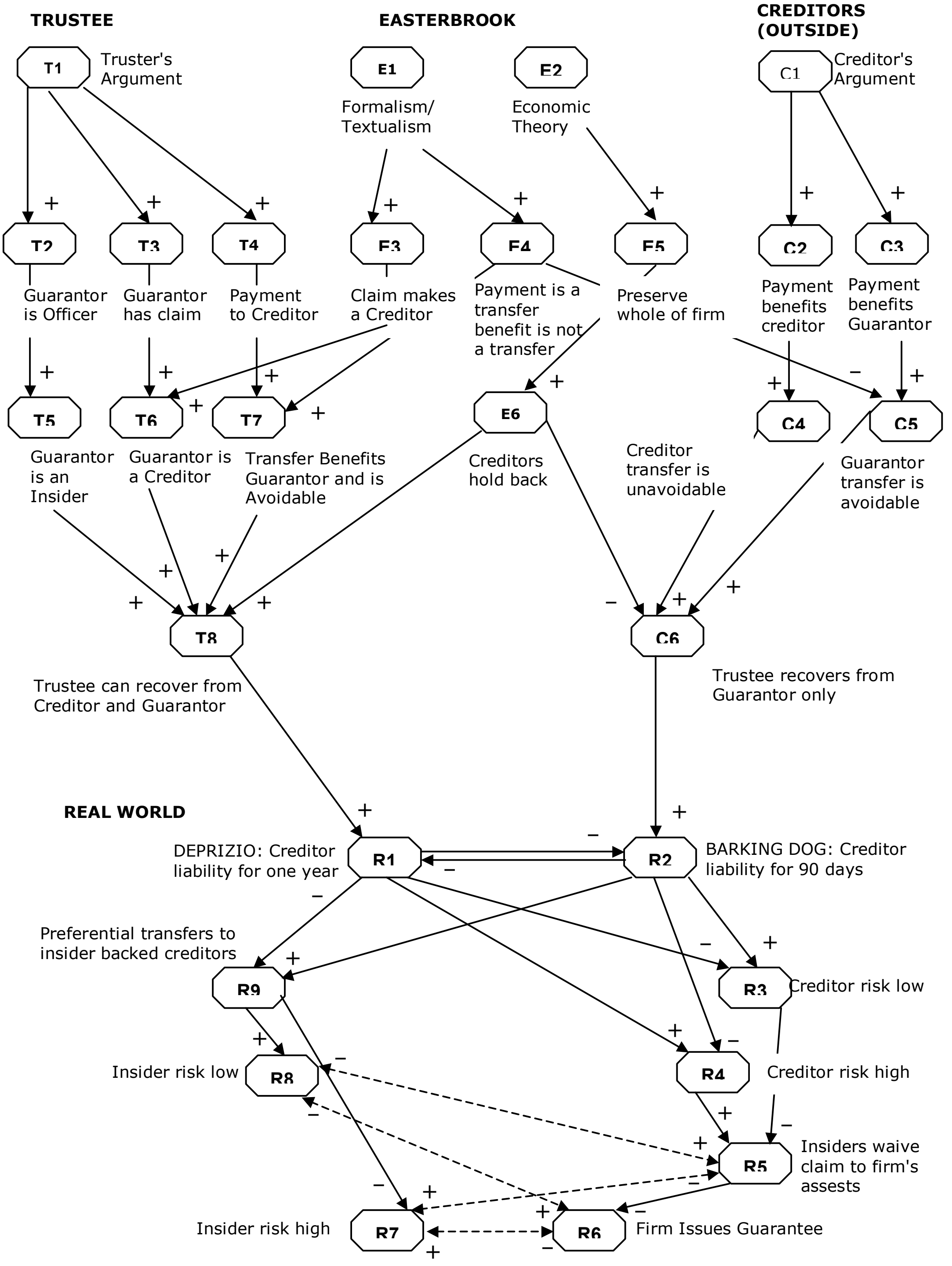

**FIGURE: 1.4.11**

The complete system, shown in Illustration 1.4.11, has been used as the basis for several simulations. The first simulation is for the actual Deprizio case and its effect on the Real World. To start, an input vector activates the Trustee's Argument (T1), the Creditor's Argument (C1), and Easterbrook's Textualism (E1) and Economic Theory (E2).

The choice of input vector leads the FCM to a limit cycle. The limit cycle has four steps. It starts with "Creditor Risk is High" (R4), a direct result of the Deprizio decision. The high degree of risk experienced by the creditors leads to "Insiders' behavior increases their degree or risk in the third step, "Insider Risk is High" (R7). In the fourth step, the insiders no longer choose to waive their claim to the firm's assets, although their risk remains high. The return to the first step, the beginning of the limit cycle, occurs next by a moderation of insider risk. In reality, one can imagine that the insiders no longer participate in the credit system. Thus, their degree of risk diminishes. The creditors' risk remains high, however, as long as Deprizio is in effect and there is need to make an acceptable rate of return on their money. At this point, it may be attractive for the insiders to once again participate and to ask again for insiders to waive their claim to the firm's assets. Overall, even as a simple non-adaptive FCM with no time dependent variables, the results are meaningful.

The next simulation examines the effect of the Creditors' Argument unimpeded by the arguments of the Trustees or Easterbrook. The initial input vector is simply "Creditors' Argument" (C1).

The creditors' argument alone produces a different limit cycle. Here, as before, it is a four step limit cycle. The first step is "Creditor Risk is Low" (R3) and "Preferential Transfers to Insider Backed Creditors" (R9). These conditions are a direct result of "Barking Dog." In the second step, the bahavior of preferential transfers acts to lower the degree of insider risk (R8). In the next step, the low insider risk encourages the firm to start issuing more insider guarantees (R6). The issuance of more insider-backed guarantees, however, acts to increase the degree of insider risk. In this last step, insider risk. In this last step, insider risk is not "high", but it is to a sufficient degree not "low".

As with the previous simulation, the results are instructive. Time dependence and variable edge weights could be used to improve the FCM. With variable edge weights, it may be possible to tract changes in the degree of insider risk and to predict at which degree their behavior will trigger the firm to issue guarantees.

The next simulation looks at a change in law from "Barking Dog" to Deprizio. In this situation, the system adjusts to the Deprizio Limit Cycle in just two iterations. The third input vector is Step III of the Deprizio Limit Cycle. The second input vector shows that the firm is still issuing guarantees even though the law has changed and preferential transfers have stopped. As shown in the previous examples, the increase in degree of creditor risk leads to insider waivers and an increase in the degree of insider risk as well.

These three simulations show some of the benefits of virtual world FCM modeling. The three examples are to a large extent focused on the "Real World" FCM. However, the main focus is indeed System. Next, it is instructive to take a closer look at the judicial decision FCM to see how a change in Easterbrook's reasoning may influence the "Real World" outcome.

For more please about this refer [1] as this paper is a very revolutionary one and the authors were impressed and attracted by it that is why all the FCM applications are given in this book. Also as several cases in the legal side have indeteterminable situation we felt it would do a lot in legal world if we apply NCM.

### APPLICATION 1.4.6: Creating Metabolic and regulatory Network Models using Fuzzy Cognitive Maps

A model of metabolic networks that uses fuzzy cognitive maps as given by Dickerson et al [32] are recollected here. Nodes of the map represent specific biochemicals such as proteins, RNA, and small molecules, or stimuli, such as light, heat, or nutrients. Edges of the map capture regulatory and metabolic relationships are established by a domain expert, the biological literature, and extracted from RNA microarray data. This work is part of the development of a software tool, FCModeler, which models and visualizes metabolic networks. A model of the metabolism of the plant hormone gibberellin in Arabidopsis is used to show the capabilities of the fuzzy model. The nodes in the FCM represent specific biochemicals such as proteins, RNA, and small molecules, or stimuli, such as light, heat, or nutrients. There are three basic types of links specified. In a conversion link (black arrow, shown as a dotted line), a node (typically a chemical(s)) is converted into another node, and used up in the process. In a regulatory link (Green and red arrows, shown as solid arrows with a plus or minus sign), the node activates or deactivates another node, and is not used up in the process. A catalytic link (blue arrows, shown as a thick line) represents an enzyme that enables a chemical conversion and does not get used up in the process. For figures please refer [32].

Thus metabolic and regulatory networks are modeled using FCM, and this kind of research has far-reaching consequences in the field of biomathematics. For more detailed discussion, refer [32].

### APPLICATION 1.4.7: Use of FCM to decide one's driving speed in a freeway

Brubaker [13-14] has modeled a way of determining one's speed when driving on a California freeway. We are extracting his column on FCMs here because of his jaunty style of narrating how exactly it feels to work with FCMs, the problems one faces and methods of arriving at greater sensitivity.

"It is a simple, six-concept model for determining one's speed when driving on a California freeway. The concepts, or the FCM's nodes, are bad weather, freeway congestion, auto accidents, patrol frequency, and own risk aversion, and the system output is own driving speed. Fuzzy values weight the causal edges. I simulated this FCM with a spreadsheet to see how it ran. Kosko defines the concept bad weather to be the only system input, although all other concepts necessarily have biases that depend driven by accidents, is also based on some patrol-scheduling algorithm that could itself be the output of another FCM.

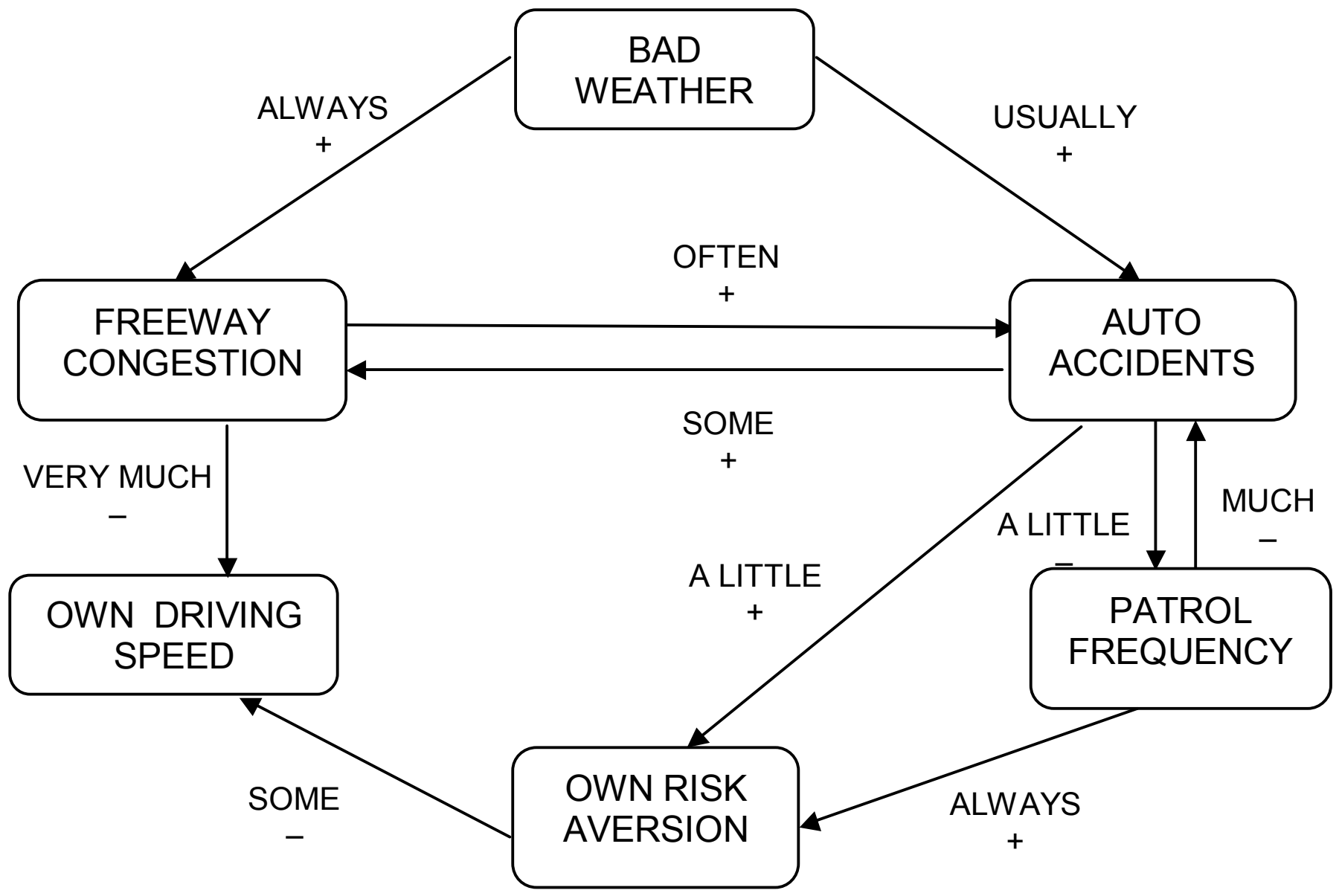

**FIGURE: 1.4.12**

"When first running the model, I discovered I could never get my own driving speed above its midpoint, which was zero, because I set the concept values between – 1 (zero driving speed) and +1 (the maximum speed 1 would ever drive). This midpoint value of own driving speed, which would map to half of my maximum driving speed, occurred when the FCM was running with no congestion and zero risk aversion.

"Considering possible solutions, if I define own driving speed to be a differential value added to some reference value (for example, a command change in speed, positive or negative, added to the speed limit), then own driving speed would not be stuck at its midpoint. However, during operation own driving speed value could only grow smaller. This fact is because, based on Kosko's model (Figure 1.4.12), there are no other FCM concepts positively affecting own driving speed.

"We need to re-examine initial assumptions. Why do we drive fast on the freeway? A principal reason is to promptly get to one's destination; let's call that schedule. For some drivers, other possible reasons would be impatience and attitude; an example of attitude would be the feeling of machismo associated with driving fast. Ignoring these last two, but adding the concept schedule and also adding what I feel are a number of legitimate edges to the original, I came up with a revised FCM (Figure 1.4.13). This implementation more accurately models my own decision-making process.

"At this point, I had planned to show plots of the FCM's operations. I unsuccessfully worked with a number of ways to do this. I have concluded that, because FCMs are dynamic, changing systems, it would be difficult or impossible to adequately demonstrate their operations in still figures. This bothers me, because I feel a sense of responsibility, now that I have introduced FCMs, to be thorough in describing how they work.

"I have an out-for me, but not for you. My “out” is that you implement the FCM and see how it works yourself; if you are interested, of course. This column and the last one give the necessary basics. If these basics are insufficient, look through the fuzzy literature for papers on the subject by Kosko and Rod Taber. Kosko also covers FCMs in both of his books, Neural Networks and Fuzzy Systems for technical readers and Fuzzy Thinking for non-technical readers.

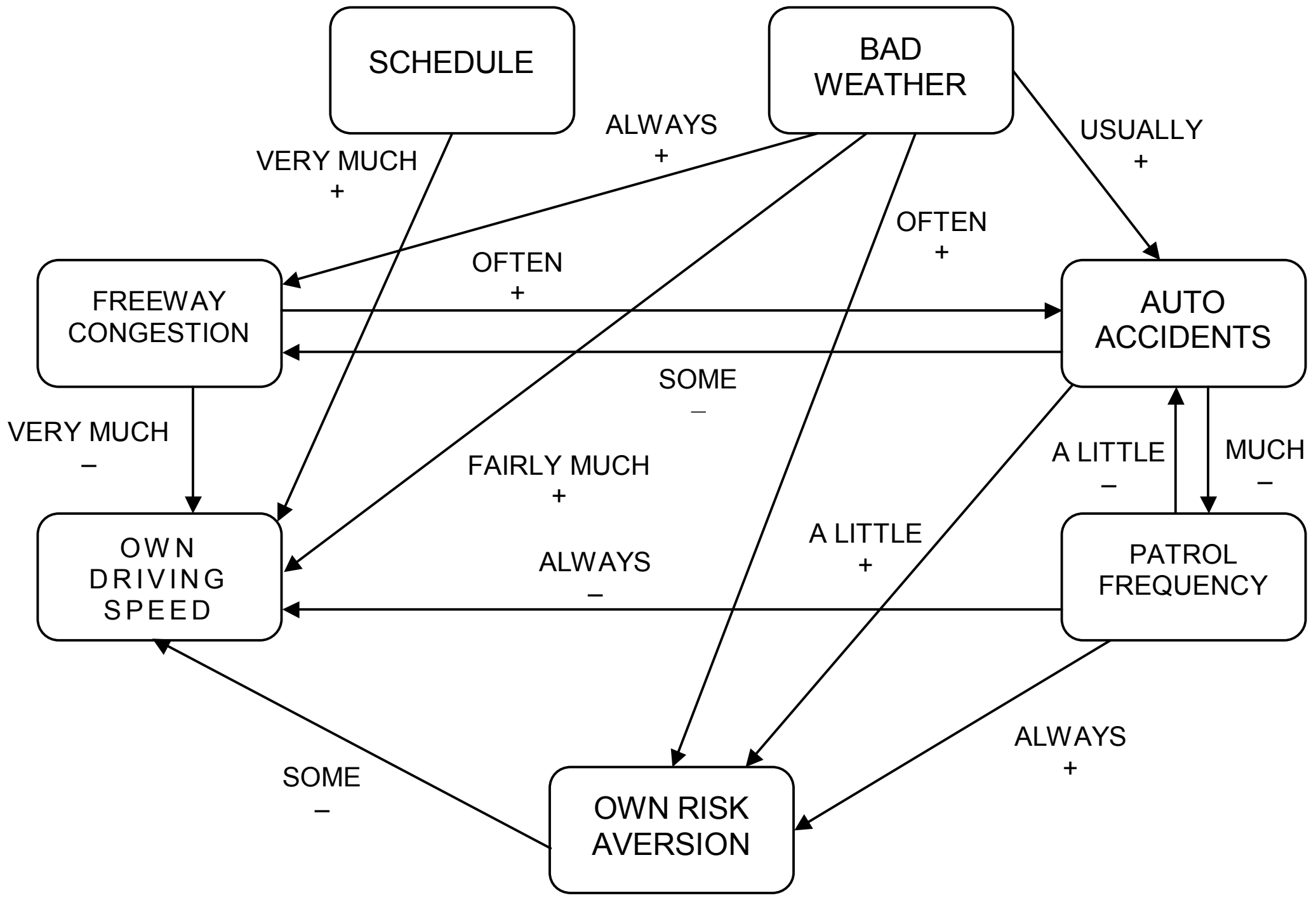

**FIGURE: 1.4.13**

"You may use a spreadsheet as a tool, as I have, or one of two commercially available FCM tools: Fuzzy Thought Amplifier from Fuzzy Systems Engineering (Poway, CA) or RTFCM (Version 2), developed by Taber and Ring Technology, (Lavale, MD). Martin McNeil wrote Fuzzy Thought Amplifier and two other Fuzzy-logic products-one for rule-base development and the other for decision making.

"Kosko suggests that the structure of an FCM lends itself well to representing an expert’s opinion. Recognizing that experts structure their thinking differently from that of other experts and that they also disagree, he has demonstrated a method of combining smaller FCMs into larger FCMs.

"Second, in Kosko’s driving-speed example, he used linguistically expressed constraint weights. In my simulation, because it is cumbersome to manipulate fuzzy variables in my spreadsheet application, I translate these variables into singleton values. This approach is a common practice. The resulting FCM is still fuzzy, because both edge and concept values are continuous. However, the example demonstrates that, if you can handle fuzzy outputs and can live with the additional computation involved in manipulating fuzzy numbers, FCM designers can express both edge and concept values as fuzzy values, represented with membership functions.

"Doing this representation would, at least in part, resolve another problem. As Kosko point out, experts tend to express causality in trivalent terms: A concept positively affects, negatively affects, or does not all affect a next concept. Fuzzy edge values would allow modeling expressions of causality in linguistic terms, such as a little positive or very negative.

"Fuzzy, as currently defined, the FCM represents causality of concepts: how concepts impact each other. Connections among concepts include weights, representing strength of causality, which are summed and normalized. Although there are likely good reasons for maintaining the purity of the structure, there is nothing to prohibit the use of additional connection operators. For example, a concept may depend on the (fuzzy) logical AND of the current values of two or more other concepts. In physical systems, both integral and derivative functions would be of value. I recently designed a system that, had I implemented it using an FCM, the arithmetic-multiply operator would have been a necessary operator.

"Armed with this information and my encouragement, go to! I will answer questions and requests for clarification set to me. However, as I have often stated the best way to learn a new technology is not to read or listen to what someone else is doing, but rather to dig in and do it ourselves. I learn far more from my mistakes than from someone else's successes."

The whole article as given by Brubaker [13] is given in the first person for not only is it inspiring in content and appealing in narration, but it also proves to us that even persons with completely non-mathematical backgrounds, can make very good use of FCMs. This article by Brubaker illustrates very candidly the application of FCMs in intelligent decision-making processes. Further, the simplicity of the tool enables anyone to make use of it for practical, day-to-day problems.

### APPLICATION 1.4.8: Use of FCM in medical diagnostics

FCMs have been very handy in medical diagnostics. Here we discuss an article by Vysoky [133]. In diagnostics of the depression, psychiatrists lack an objective knowledge about patient state. Except for very strong depressions no significant changes of the physiological parameters are measurable. Almost all information about patient mental state is obtained in verbal form on the base of conversation with the patient. Due to this fact, data are vague, uncertain and often inconsistent.

Contemporary most popular approach to the etiopathogenesis of a depression is based on cognitive-behavioral theory of depression. This theory assumes, that incorrect cognitive processing of information coming from an environment causes the depression. Patient's interpretation of events in this environment is distorted, and this misinterpretation influenced his behaviour. They are different "negative ideas", "automatic ideas", activation of the disfunction rules, etc.

A degree of the depression is usually evaluated by therapist with help of the Beck's scale. It is pseudo-quanitification based on complicated questionnaire. Patient answers expressing his attitudes to the different problems are evaluated by means of the four-degree scale and their sum, final score can be considered as the degree of depressive mood. It enables of course to estimate the degree of immediate depressive mood. But is doesn't say anything about way leading to this state (which mental states or moods were predecessors of this instantaneous mood), and it doesn't say anything about prognosis of future states or moods. Automatic and negative ideas have complicated dynamics depending on the preceding states.

The submitted approach based on the fuzzy cognitive maps enables to catch the mentioned dynamics. A first step is to classify different moods into a set of qualitatively distinguished and verbally labeled categories. Each category (specific mood), is determined with help of specific factors e.g. possibility of concentration, efficiency, feeling of meaningfulness of patients work, feeling of understanding of his environment etc. They are the same or very similar factors like in questionnaire of the Beck's scale. Now these categories may be considered like terms of linguistic variable "mood". It means that they are fuzzy sets defined on factor space.

During the time patient passes from one mood to another. Transition depends on previous state and on instantaneous influences from the environment. It can be investigated by interviewing the patient and be vividly depicted with help of the cognitive map. The nodes of the map correspond with states (moods) and edges correspond with transitions from one state to another. The cognitive map is constructed intuitively according to patient's own description of his feelings and self-evaluation. This cognitive map provides much more information about patient mental states then simple estimating of instantaneous degree of depressive mood. We can see, that some transitions are more frequent than others, and that some transitions are caused or inhibited by influences from patient's environment. On the other hand, cycles in the graph of the cognitive map can be recognized. These cycles correspond with "automatic ideas", and under the certain conditions they may behave as attractors in factor space. This description provides a deeper understanding of dynamic of the patient thinking, his misinterpretation of "inputs" leading him to the undesirable cycles. It enables to find such "inputs" which remove him safely far from attracting areas of these cycles and to prepare base for effective psychotherapy.

Any patient can distinguish a finite number of moods. According to his verbal description we can construct fuzzy cognitive map. This FCM can be considered as a graph of the state transition structure of some fuzzy finite automaton. Analyzing this FCM we could obtain its input and output alphabet, set of inner states and the state transition and output relations. Having this more formal description, we can analyse more correctly boundaries of attracting areas and we can simulate some therapeutical procedures.

Thus we see FCMs can be used to study several medical problems and especially for medical areas as important as psychotherapy where diagnostics plays a major role. For more about this, please refer [133].

### APPLICATION 1.4.9: Fuzzy Mechanisms for Causal Relations

Here we recall a method to implement Fuzzy Causal Relations (FCR) that can be used in Rule Based Fuzzy Cognitive Maps (RBFCM). This method introduces a new fuzzy operation that simulates the "accumulative" property associated with causal relations

– the Fuzzy Carry Accumulation (FCA). The FCA allows a great flexibility in the addition and removal of concepts and links among concepts while keeping compatibility with classic fuzzy operations: Carvalho and Tome [20].

There is a causal relation between two given concepts whenever a change in one of those concepts affects the other one. For example, there is a causal relation between police vigilance and robbery: a major increase in police vigilance will probably cause a decrease in robbery.

Causal associations appear to be the most widely used in Cognitive Maps; the following statements can easily explain this fact:

- Causal relations are the major way in which understanding about the world is organized;
- Causality is the primary form of post hoc explanation of events;
- Choice among alternative actions (Decision processes) involve causal evaluation.

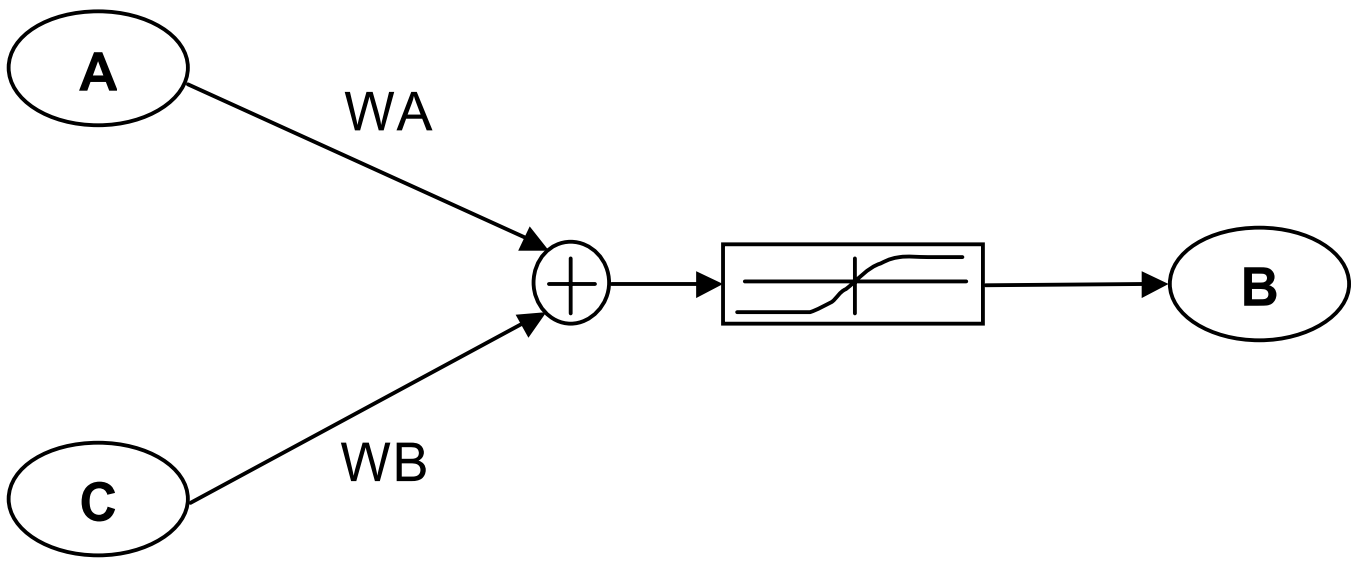

**FIGURE: 1.4.14**

Usually, causal relations in causal maps always involve change. The result of a causal effect is always a variation in one or more concepts.

The idea behind the property of fuzzy causal relations can be easily explained using singletons. Let us imagine that the fuzzy set increase of concepts Robbery is a singleton at x = 0.5 (Figure 1.4.15) and we have the following rules:

- If Police_Vigilance Decreases Then Robbery Increases;
- If Wealth_of_Residents Increases_Much Then Robbery Increases

If the consequent of the application of the two rules is Increase ($\mu = 0.7$) and Increase ($\mu = 0.5$), how should these consequents be combined in order to produce the accumulative effect (as we saw above, the result should represent a variation larger than Increase)?

In a traditional fuzzy system, the result would always be (after de-fuzzification) x = 0.5, even when the sum of the beliefs is greater than 1. Here, we introduce the concept of the Fuzzy Carry Accumulation (FCA):

- If the sum of the beliefs is lower or equal to 1, then we have a standard fuzzy operation. If we have Increase ($\mu = 0.3$) and Increase ($\mu = 0.5$), then the result is Increase with μ = 0.8, since Carvalho and Tome [20] do not fully believe that the result should be more than Increase.

However

- If the sum of the beliefs is greater than 1, then there is an overflow of the reminder (just like a carry in a sum operation) towards a value representing a larger variation.

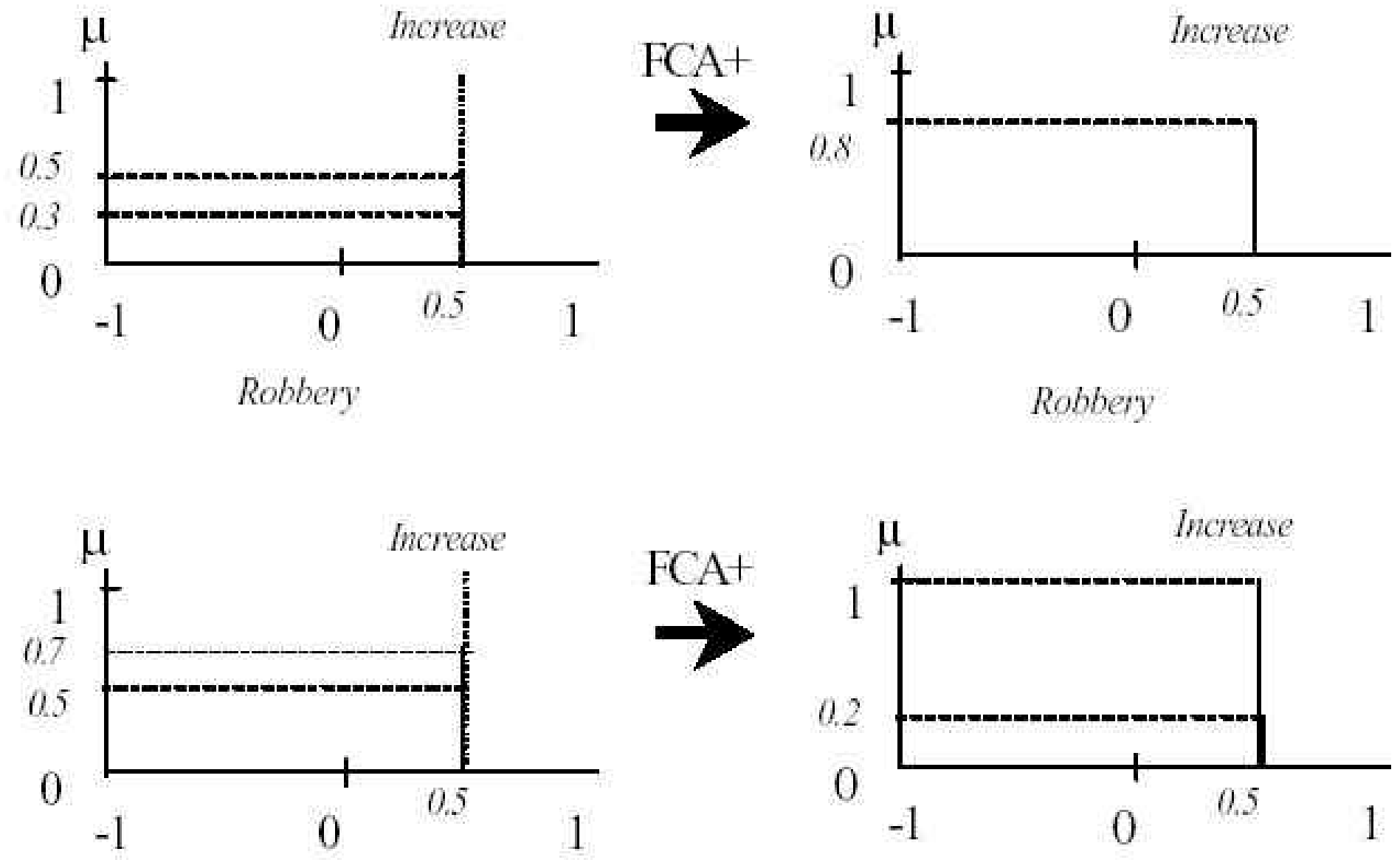

Figure 1.4.15 represents this operation. With this overflow of the excess, after de-fuzzication, the result of Robbery will be larger than Increase.

If the rules involved represent a decrease, then the carry is performed in the opposite direction.

If there are several rules that cause large reminder, whenever the reminder exceeds 1 there is another overflow towards a larger variation.

When the consequents involve different fuzzy sets-for example Increase and Increase_Much_then since we know that the result must be larger than the largest consequent, the solution is to shift the smaller consequent towards the former.

This shift operation arises several problems and constraints that will affect the implementation of fuzzy causality, since there must be a way of retaining and constraints, fuzzy causal operations can not be implemented with singletons. Singletons were used a simpler way to show the FCA principles.

Several more interesting results can be had from [20].

### 1.4.10: Use of Fuzzy Cognitive Maps to Simulate the Information Systems Strategic Planning Process

Kardaras and Karakostas [47] have given the use of FCM in Information Technology (IT). As IT plays a critical role in all organizations we have taken several factors from this paper to substantiate that FCMs are applied in study and analysis of IT. In the early 1980s articles began to focus on Strategic Planning of Information Systems (SISP) and to argue the critical importance of IT in today's organizations. Since then, a large number of models were presented in order to analyse IT from a strategic point of view and suggest new IT projects. However, researchers urge for alternative approaches to SISP, as current ones fall short in taking into consideration both the business and IT perspectives as well as they fail to tackle the complexity of the domain and suggest specific Information Systems (IS) opportunities. This article suggests Fuzzy Cognitive Maps (FCM) as an alternative modeling approach and describes how they can be developed and used to simulate the SISP process.

In this article Kardaras and Karakostas suggest the development of FCMs as an approach to alleviating the problems associated with SISP models. A list of the main problems follows:

1. Despite the large number of models available for SISP, they fail to provide a statement of the IS organizational objectives and to identify specific new projects.

2. Currently available models fail to tackle the complexity of the SISP process. Most of the models are oversimplified 2 × 2 tables with each dimension marked high and low. These unnecessary simplifications involve reductions of reality, but can also involve distortions. Knowledge can be represented in much richer ways than in matrix. Artificial intelligence techniques were proposed such as production rules, frames and semantic nets, object-oriented systems, etc. In this article the appropriateness of FCM is argued.

3. A criticism was also developed with regard to models applicability to different industrial sectors. As an example a particular model may be suitable for the chemical industry but appropriate to service industry.

4. Models as described in the literature make the implicit assumption that they exist outside time. It is assumed that the models will be as true tomorrow as they are today. Further, they are not flexible enough to consider unanticipated changed in the organization and its environment.

5. Models fail to trigger brainstorming. They fail to mobilise an organisation's biggest asset: its collective brain-power.

6. The lack of computer based models affects the effectiveness of the SISP process.

The proposed FCM considers both organizational and Information Technology (IT) related variables and links them with their interrelationships. Each inter-relationship is associated with a fuzzy linguistic weight, which indicates the Belief of Truth (BoT) that planners share for the relationship. The model contains 165 variables and 210 relationships which were extracted after thorough analysis of the SISP literature, i.e. case studies and theoretical frameworks and consideration of relevant practical experience. It is particularly important though, on one hand that planners opinion is also taken into account in an attempt to facilitate the development of a consensus with respect to concepts modeled in the FCM, and on the other to minimize planners' bias with regard to their suggestions. Both objectives are achieved with the analysis of the available SISP relevant material, which ensures the reliability and validity of the FCM and with the participation of the planners who should customize the resulting FCM to their organizational context.

The FCM Kardaras and Karakostas [47] propose allows planners to develop scenarios and assess alternative ways of applying IT in order to order to improve organizational performance. Planners and other stakeholders involved in an SISP study, usually have differing perceptions and assumptions about their organization and its environment, mainly because of constantly changing conditions, turbulence, uncertainty, and personal interests. As a result in several versions of the FCM, because of alternative suggestions as to which factors, relationships and strength of the relationships, should be modeled. The assumptions about how to conceive the organization and its environment, as presented in different versions of the FCM, constitute a source for developing scenarios. Scenarios were proposed and extensively used as an approach to trigger and structure strategic thinking. They are used to surface assumptions which consequently enable the stakeholders to get an insight and understanding of the issues involved in the planning process and the causal-effect relationships of the factors modeled. Scenarios also offer alternative views of the future and enable planners to express their beliefs with regard to future fluctuations of the factors. The proposed FCM provides the basis for developing scenarios in order to:

1. Assess the general organizational environment,
2. Evaluate the organizational information systems effectives.
3. Analyze the impact of alternative business and IT strategic options.

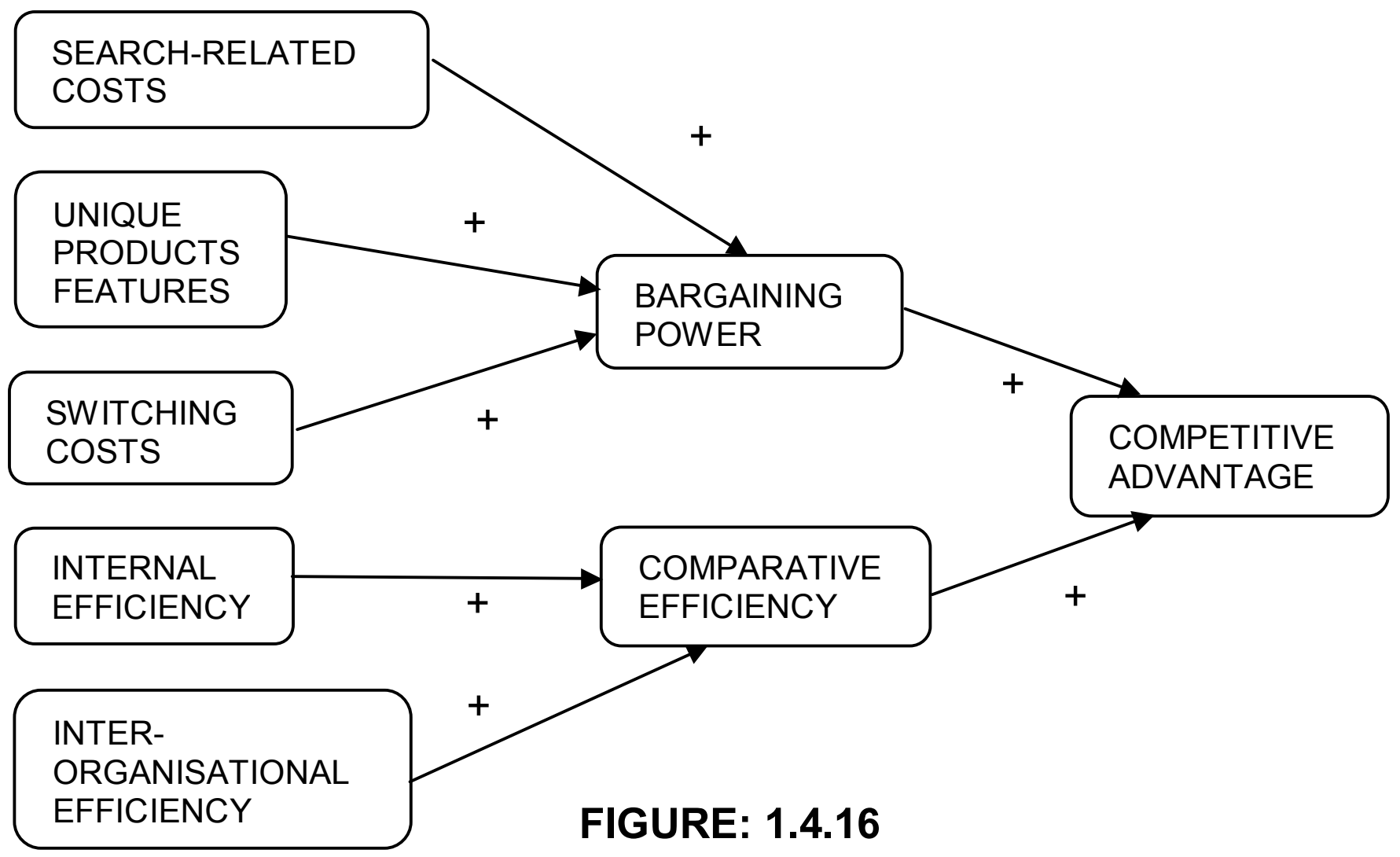

**FIGURE: 1.4.16**

The example, which is shown in Figure 1.4.16, considers the case where the corporate objective is to increase the organizational competitiveness through IT. The planners rely on the above model and analyze several strategic options, i.e. to increase comparative efficiency or to increase their business bargaining power, and then other more specific strategic options, i.e. to increase comparative efficiency or to increase their business bargaining power, to improve internal efficiency, or inter organizational efficiency, etc.

Planners may exclude or include other concepts in the model, by asking questions such as “What other options should be considered?” It is expected that individual planners will express their own personal views. There are techniques available to combine individual cognitive maps by assigning weighting factors for each map, which indicate the credibility of the individual planner. However, it is essential that a consensus about the strategic options will be developed among the members of the planning team. If opposing views go unrecognized and unresolved, agreement on the strategy may be difficult to achieve and the success of the strategic plan ill be jeopardized. It is important for both the richness of the FCM and the success of the planning that the valuable diversity in the views of the team is modeled and taken into consideration.

By continuing the analysis of the map and for example the consideration of actions, which improve the inter-organizational efficiency strategic options, the planners use the map shown in Figure 1.4.17.

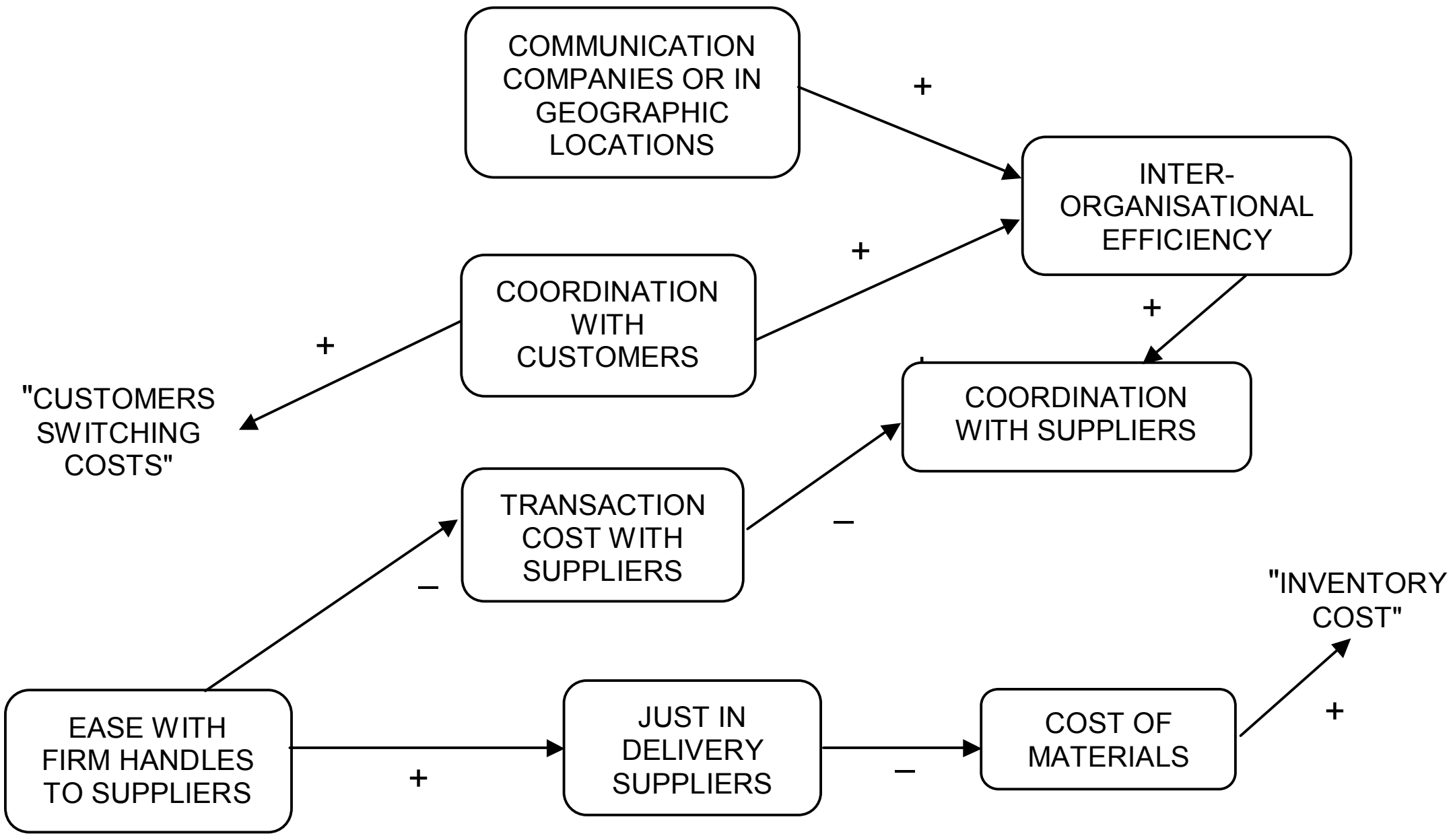

**FIGURE: 1.4.17**

Further analysis of the possible strategic options will become more detailed and focused on specific IT and organizational issues, until the train of thought of the planners is exhausted. Figure 1.4.17 shows several options such as improving the co-ordination with the customers or improving co-ordination with suppliers with a different set of more detailed suggestions for each option. Planners develop scenarios in order to assess the expected results of the alternative strategies in terms of their impact, i.e. Increase or decrease, on the variables in the model. The impact analysis relies on the change propagation through the causal relationships of the FCM. The scenarios are combined then in order to identify the business areas to which IT can contribute the most. Which alternative to consider as the final strategic direction depends on

1. The extent to which it is generally accepted by the planners, i.e. the number of planners who support a particular strategy.

2. The degree it is believed it will yield the anticipated results, as this is shown by the BoT degrees associated with each relationship and the total effect of a scenario.

The most well-accepted suggestions are considered first, then the less well-accepted suggestions are considered in turn, etc.

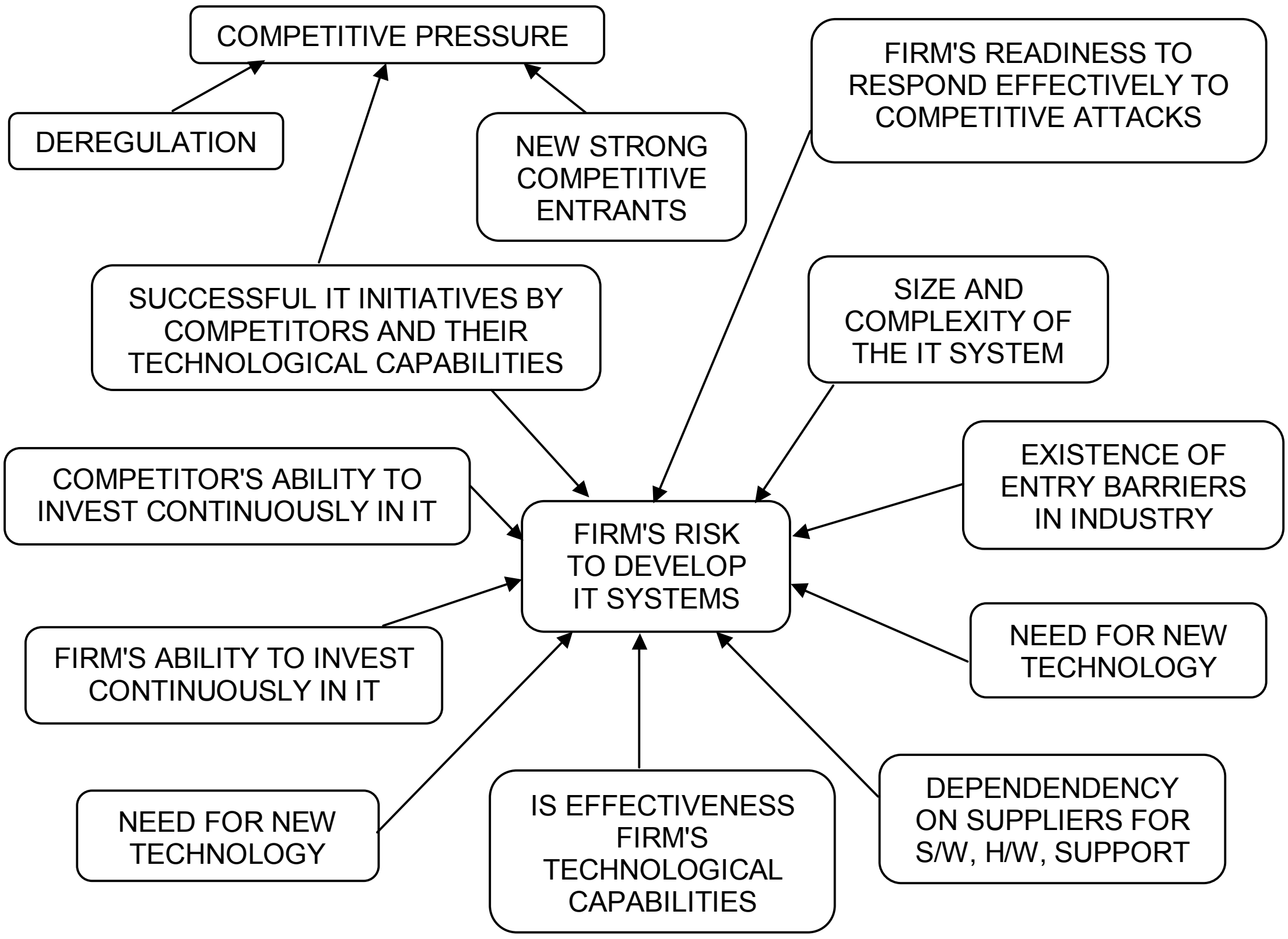

**FIGURE: 1.4.18**

Planners assign fuzzy weights to every relationship in the FCM, in order to express their belief in the strength of a certain causal relationship, in a certain organizational environment. Relationship weights are therefore customized to specific organizations and industry sectors. Fuzzy weights in our FCM show the belief, which planners share with regard to the existence of a certain relationship, and not the magnitude of change that a variable may undergo because of its causal relationship with other variables. Planners during modeling should answer the following question for each relationship:

How strongly do we believe that there exists a causal relationship between variable (X) and variable (Y)?

The quantity space of the relationships' weights Q(w), is the following set: Q(w) {Undefined, Weak, Moderate, Strong}. It is assumed that the following ordering applies for Q(w):

It is suggested that linguistic fuzzy weights are used as it is easier for the planners to express their beliefs than it is when they have to suggest real values for weights or calculate the average of different responses as proposed. Moreover, it is particularly difficult to try to measure concepts within the SISP domain, if it is to adopt the approach proposed. An attempt to measure highly aggregate and qualitative in nature concepts such as "environment rich the IT opportunities", or "business readiness to respond effectively to competitive attacks", etc. could shift the focus of the planners from a creativity oriented process to a process with a narrow focus on measurement procedures. That would affect the effectiveness of the planning session and would hinder the brainstorming of the planners. Moreover, the participation of the planners or experts in a particular domain was widely suggested. Especially in SISP the development of consensus among the planning team is vital for its success and weights as with concepts should be imposed to planners.

Certain types of causal relationships in our FCM are proposed in order to deal with spurious relationships and inferences as well as to control the change propagation. A relationship between two or more variables is spurious, when it is not an actual cause-effect relationship. Consequently, spurious relationships may lead to spurious inferences and incorrect decisions. The types of the relationships are presented in the following set:

Q(r) = {Affects, Requires, Multiple_relationship, Stop-relationship}

The notation used for each type of relationship (causal) is shown in the following examples:

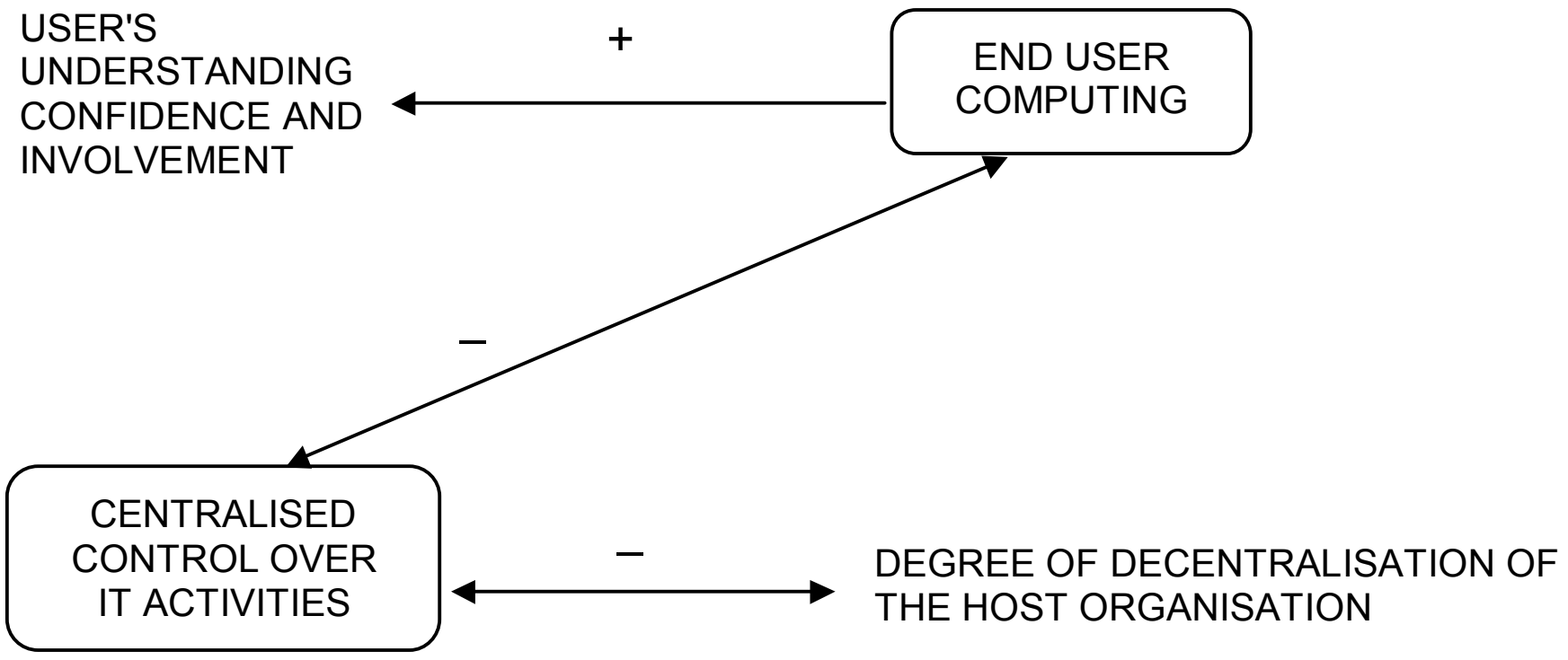

**FIGURE: 1.4.19**

These types control the direction of each relationship, thus affecting the propagation of the causality and prevent incorrect decisions to be made. Their definition and rules for modeling relationships follow:

**DEFINITION OF AFFECTS (A) [47]:** *If an {increase or decrease} in a variable "X" results in an {increase or decrease} of a variable "Y" then X affects Y.*

***Example [47]:*** Figure 1.4.18 shows a part of the FCM with affects An increase in the cause variable Volume of Transactions with Customers/Suppliers affects positively (as the sign of relationship is +) the effect variable Environment of Realize IT Opportunities of Threats. This type of relationship is used to allow the estimation of the effects of a change to be made following the causal paths.

**DEFINITION OF REQUIRES (R) [47]:**
*If an {increase or decrease} in a variable "X" does not result in an {increase or decrease}of a variable "Y" AND if an {increase or decrease} of "Y" is required for an {increase or decrease} of variable "X" then X requires Y.*

***Example [47]:*** In Figure 1.4.18 the required type is also shown and is denoted with an arrow. This type is used to answer questions on how to attain specific goals. Therefore, in order for a company to increase its IT opportunities (or minimize the threats) it is required to increase its readiness to respond effectively to competitive attacks. In this case a certain objective is set, and ways of attaining it are specified through causal paths and by exploiting the required relationship type.

**DEFINITION OF MULTIPLE RELATIONSHIP (M) [47]:**
*If an {increase or decrease} in a variable "X" affects a variable "Y" AND if an {increase or decrease} of the variable "X" requires an {increase or decrease} of the variable "Y" then X multiple Y.*

**Note [47]:** The operator "AND" in the rule, implies that both affects and requires types are generally true for the two variables. However, they cannot be true simultaneously. Therefore under certain circumstances only one type holds true. Which type is in effect, depends on which variable is considered to be a cause variable and which an effect one.

The example which follows clarifies the differences.

***Example [47]:*** Figure 1.4.19 shows a multiple relationship (denoted with arrow heads on both sides of the link), between End User Computing and Centralized Control over IT activities. A multiple relationship type can be interpreted either as affect or as requires, depending on which is the cause variable in each case. If End User Computing is the cause variable and it increases then it causes (Affects) Centralized Control over IT activities, which is then decreased (negative sign). It the End User Computing is defined as objective then this multiple type is interpreted as requires the Centralized Control over IT activities to decreases. The same holds for the Centralized Control over IT activities variable.

| **Type of relationship of cause variable** | **Compatible type of relationship of effect variable** |
|---|---|
| Affects | Affects or multiple_relationship (interpreted as Affects). |
| Requires | Requires or multiple_relationship (interpreted as Requires). |
| Multiple_relationship | Affects or requires or multiple relationship. |
| Stop_relationship | none |

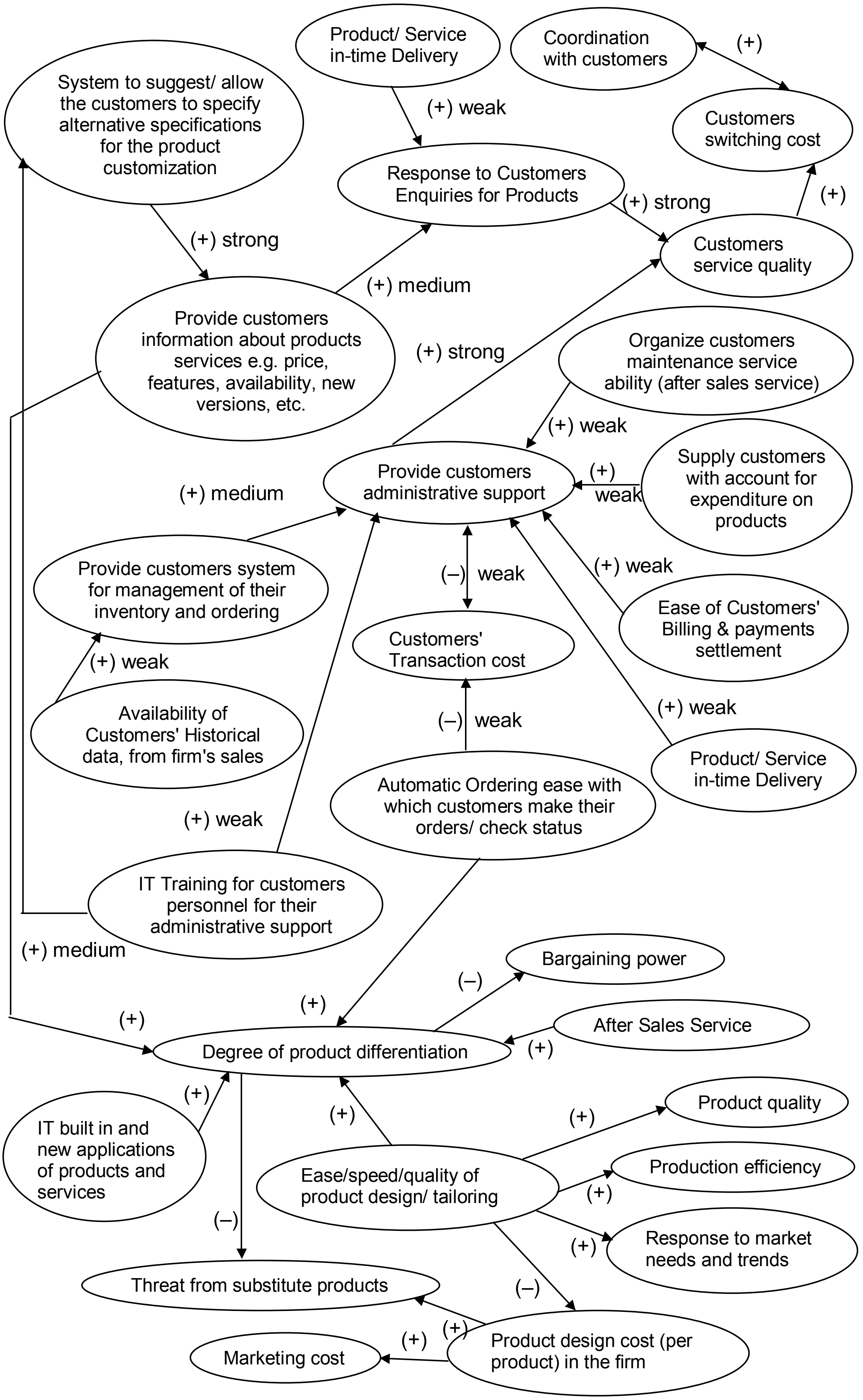
Product/ Service in-time Delivery
Coordination with customers
(+)
System to suggest/ allow the customers to specify alternative specifications for the product customization
(+) weak
Customers switching cost
Response to Customers Enquiries for Products
(+) strong
(+)
(+) strong
Customers service quality
(+) medium
Provide customers information about products services e.g. price, features, availability, new versions, etc.
(+) strong
Organize customers maintenance service ability (after sales service)
(+) weak
Provide customers administrative support
(+) weak
Supply customers with account for expenditure on products
(+) medium
Provide customers system for management of their inventory and ordering
(–) weak
(+) weak
Ease of Customers' Billing & payments settlement
Customers' Transaction cost
(+) weak
Availability of Customers' Historical data, from firm's sales
(–) weak
(+) weak
Product/ Service in-time Delivery
Automatic Ordering ease with which customers make their orders/ check status
(+) weak
IT Training for customers personnel for their administrative support
(+) medium
Bargaining power
(–)
(+)
(+)
After Sales Service
Degree of product differentiation
(+)
(+)
IT built in and new applications of products and services
(+)
(+)
Product quality
Ease/speed/quality of product design/ tailoring
Production efficiency
(+)
(–)
(+)
Response to market needs and trends
Threat from substitute products
(–)
(+)
(+)
Marketing cost
Product design cost (per product) in the firm

**FIGURE: 1.4.20**

Figure 1.4.20 gives the customer service FCM as given by Kardaras and Karakostas in [47]. We discuss in Chapter two how NCMs can be applied to this problem.

### 1.4.11 Application of FCM to Diagnosis of Specific Language Impairment

Despite the numerous studies that have been conducted since the first half of the 19th century, Specific Language Impairment (SLI) remains a language disorder that cannot be easily diagnosed and described due to its similarity to other disorders. Research has shown that almost 160 factors can be taken into account in the diagnosis. Therefore it is necessary to develop a model of differential diagnosis of SLI which will aid the specialist in the diagnosis by suggesting a possible diagnosis. SLI is a significant disorder of spoken language ability that is not accompanied by mental retardation, frank neurological damage or hearing impairment. Georgopoulous et al [37] have tackled the differential diagnosis of SLI using Fuzzy Cognitive Maps. The diagnosis of SLI is vital and imperative, because children suffering from the impairment face a wide variety of problems both on language and cognitive levels. Also, the diagnosis becomes tougher for children recognized with SLI constitute a heterogeneous group, and because SLI is congenital and large number of factors are involved.

Georgopoulos, V.C. et al have used FCMs to analyse this problem in a critical way using several experts' opinions. For more information, please refer [37].

### 1.4.12: FCM Approach to Web-mining Inference Amplification

We recall the definition of web-mining as given by Lee K.C. et al [63]. Definition of web-mining is slightly different from data mining such that web-mining is an activity of discovering knowledge or patterns from a collection of web-documents, which are usually composed of directed graph of document nodes and hypothesis. Web mining is sharply different from data-mining in that the former is based on ill-structured web documents with little machine readable semantic. In contrast data mining derives patterns from structural data in database.

An extensive study in this direction is made by [63], and it can be referred for more about this. The web content mining is related with using web search engines, main role of which is to discover the web contents according to the users requirements and constraints. Results of the web usage mining provide decision makers with crucial information about the lifetime value of customers cross marketing strategies across product and effectiveness of promotion campaigns. Further other factors like targeting ads or web page, categorize user preferences and restructure a website to create a more effective management of workgroup communication etc.

This paper by K.C. Lee et al [63] suggests a new concept of web-mining procedure called WEMIA which is aimed at enriching the interpretation of the draft results and association rule mining. They further suggest that the basic nature of FCMs used for this study needs to be improved so that more complicated web-mining results can be analyzed enough to yield meaningful inference results. This book will deal with this topic in Chapter two and suggest a better tool i.e. NCMs.

## 1.5 Definition and Illustration of Fuzzy Relational Maps (FRMs)

In this section, we introduce the notion of Fuzzy relational maps (FRMs); they are constructed analogous to FCMs described and discussed in the earlier sections. In FCMs we promote the correlations between causal associations among concurrently active units. But in FRMs we divide the very causal associations into two disjoint units, for example, the relation between a teacher and a student or relation between an employee or employer or a relation between doctor and patient and so on. Thus for us to define a FRM we need a domain space and a range space which are disjoint in the sense of concepts. We further assume no intermediate relation exists within the domain elements or node and the range spaces elements. The number of elements in the range space need not in general be equal to the number of elements in the domain space.

Thus throughout this section we assume the elements of the domain space are taken from the real vector space of dimension n and that of the range space are real vectors from the vector space of dimension m (m in general need not be equal to n). We denote by R the set of nodes $R_1,\ldots, R_m$ of the range space, where $R = \{(x_1,\ldots, x_m) \mid x_j = 0 \text{ or } 1 \}$ for j = 1, 2,…, m. If $x_i = 1$ it means that the node $R_i$ is in the on state and if $x_i = 0$ it means that the node $R_i$ is in the off state. Similarly D denotes the nodes $D_1, D_2,\ldots, D_n$ of the domain space where $D = \{(x_1,\ldots, x_n) \mid x_j = 0 \text{ or } 1\}$ for i = 1, 2,…, n. If $x_i = 1$ it means that the node $D_i$ is in the on state and if $x_i = 0$ it means that the node $D_i$ is in the off state.

Now we proceed on to define a FRM.

**DEFINITION 1.5.1:** *A FRM is a directed graph or a map from D to R with concepts like policies or events etc, as nodes and causalities as edges. It represents causal relations between spaces D and R .*

*Let $D_i$ and $R_j$ denote that the two nodes of an FRM. The directed edge from $D_i$ to $R_j$ denotes the causality of $D_i$ on $R_j$ called relations. Every edge in the FRM is weighted with a number in the set $\{0, \pm 1\}$. Let $e_{ij}$ be the weight of the edge $D_iR_j$, $e_{ij} \in \{0, \pm 1\}$. The weight of the edge $D_i R_j$ is positive if increase in $D_i$ implies increase in $R_j$ or decrease in $D_i$ implies decrease in $R_j$ ie causality of $D_i$ on $R_j$ is 1. If $e_{ij} = 0$, then $D_i$ does not have any effect on $R_j$ . We do not discuss the cases when increase in $D_i$ implies decrease in $R_j$ or decrease in $D_i$ implies increase in $R_j$.*

**DEFINITION 1.5.2:** *When the nodes of the FRM are fuzzy sets then they are called fuzzy nodes. FRMs with edge weights $\{0, \pm 1\}$ are called simple FRMs.*

**DEFINITION 1.5.3:** *Let $D_1, \ldots, D_n$ be the nodes of the domain space D of an FRM and $R_1, \ldots, R_m$ be the nodes of the range space R of an FRM. Let the matrix E be defined as $E = (e_{ij})$ where $e_{ij}$ is the weight of the directed edge $D_iR_j$ (or $R_jD_i$), E is called the relational matrix of the FRM.*

***Note:*** It is pertinent to mention here that unlike the FCMs the FRMs can be a rectangular matrix with rows corresponding to the domain space and columns

corresponding to the range space. This is one of the marked difference between FRMs and FCMs.

**DEFINITION 1.5.4:** *Let $D_1$, ..., $D_n$ and $R_1,\ldots, R_m$ denote the nodes of the FRM. Let $A = (a_1,\ldots,a_n)$, $a_i \in \{0, \pm 1\}$. A is called the instantaneous state vector of the domain space and it denotes the on-off position of the nodes at any instant. Similarly let $B = (b_1,\ldots, b_m)$ $b_i \in \{0, \pm 1\}$. B is called instantaneous state vector of the range space and it denotes the on-off position of the nodes at any instant $a_i = 0$ if $a_i$ is off and $a_i = 1$ if $a_i$ is on for i= 1, 2,…, n Similarly, $b_i = 0$ if $b_i$ is off and $b_i = 1$ if $b_i$ is on, for i= 1, 2,…, m.*

**DEFINITION 1.5.5:** *Let $D_1$, ..., $D_n$ and $R_1,\ldots, R_m$ be the nodes of an FRM. Let $D_iR_j$ (or $R_j D_i$) be the edges of an FRM, j = 1, 2,…, m and i= 1, 2,…, n. Let the edges form a directed cycle. An FRM is said to be a cycle if it posses a directed cycle. An FRM is said to be acyclic if it does not posses any directed cycle.*

**DEFINITION 1.5.6:** *An FRM with cycles is said to be an FRM with feedback.*

**DEFINITION 1.5.7:** *When there is a feedback in the FRM, i.e. when the causal relations flow through a cycle in a revolutionary manner, the FRM is called a dynamical system.*

**DEFINITION 1.5.8:** *Let $D_i R_j$ (or $R_j D_i$), $1 \leq j \leq m$, $1 \leq i \leq n$. When $R_i$ (or $D_j$) is switched on and if causality flows through edges of the cycle and if it again causes $R_i$ (or$D_j$), we say that the dynamical system goes round and round. This is true for any node $R_j$ (or $D_i$) for $1 \leq i \leq n$, (or $1 \leq j \leq m$). The equilibrium state of this dynamical system is called the hidden pattern.*

**DEFINITION 1.5.9:** *If the equilibrium state of a dynamical system is a unique state vector, then it is called a fixed point. Consider an FRM with $R_1$, $R_2,\ldots$, $R_m$ and $D_1$, $D_2,\ldots$, $D_n$ as nodes. For example, let us start the dynamical system by switching on $R_1$ (or $D_1$). Let us assume that the FRM settles down with $R_1$ and $R_m$ (or $D_1$ and $D_n$) on, i.e. the state vector remains as (1, 0, …, 0, 1) in R (or 1, 0, 0, … , 0, 1) in D), This state vector is called the fixed point.*

**DEFINITION 1.5.10:** *If the FRM settles down with a state vector repeating in the form*

$$A_1 \to A_2 \to A_3 \to \ldots \to A_i \to A_1 \text{ (or } B_1 \to B_2 \to \ldots \to B_i \to B_1)$$

*then this equilibrium is called a limit cycle.*

**METHODS OF DETERMINING THE HIDDEN PATTERN**

Let $R_1$, $R_2,\ldots$, $R_m$ and $D_1$, $D_2,\ldots$, $D_n$ be the nodes of a FRM with feedback. Let E be the relational matrix. Let us find a hidden pattern when $D_1$ is switched on i.e. when an input is given as vector $A_1 = (1, 0, \ldots, 0)$ in $D_1$, the data should pass through the relational matrix E. This is done by multiplying $A_1$ with the relational matrix E. Let $A_1E = (r_1, r_2,\ldots, r_m)$, after thresholding and updating the resultant vector we get $A_1 E \in R$. Now let $B = A_1E$ we pass on B into $E^T$ and obtain $BE^T$. We update and threshold the vector $BE^T$ so that $BE^T \in D$. This procedure is repeated till we get a limit cycle or a fixed point.

**DEFINITION 1.5.11:** *Finite number of FRMs can be combined together to produce the joint effect of all the FRMs. Let $E_1, \ldots, E_p$ be the relational matrices of the FRMs with nodes $R_1, R_2, \ldots, R_m$ and $D_1, D_2, \ldots, D_n$, then the combined FRM is represented by the relational matrix $E = E_1 + \ldots + E_p$.*

Now we give a simple illustration of a FRM, for more about FRMs please refer [124, 125, 136].

***Example 1.5.1:*** Let us consider the relationship between the teacher and the student. Suppose we take the domain space as the concepts belonging to the teacher say $D_1, \ldots, D_5$ and the range space denote the concepts belonging to the student say $R_1$, $R_2$ and $R_3$.

We describe the nodes $D_1, \ldots, D_5$ and $R_1$, $R_2$ and $R_3$ as follows:

Domain Space

$D_1$ – Teaching is good
$D_2$ – Teaching is poor
$D_3$ – Teaching is mediocre
$D_4$ – Teacher is kind
$D_5$ – Teacher is harsh [or rude]

(We can have more concepts like teacher is non-reactive, unconcerned etc.)

Range Space

$R_1$ – Good Student
$R_2$ – Bad Student
$R_3$ – Average Student

The relational directed graph of the teacher-student model is given in Figure 1.5.1.

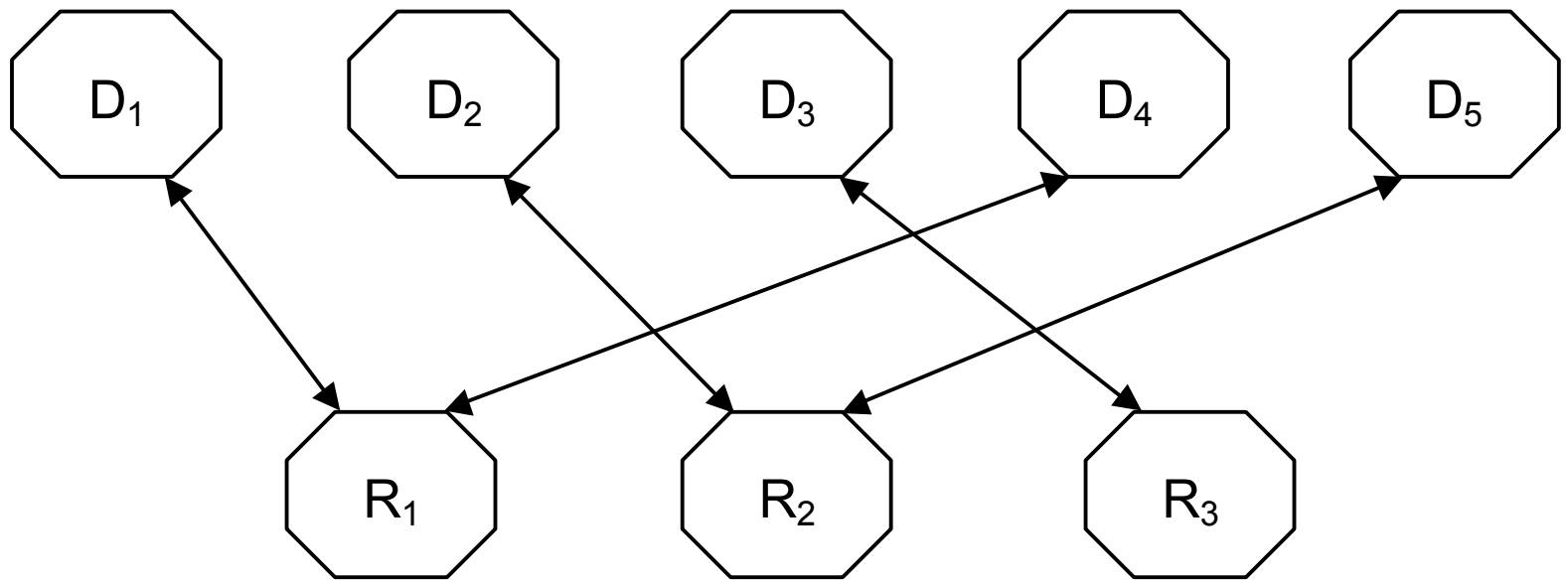

**FIGURE: 1.5.1**

The relational matrix E got from the above map is

$$E = \begin{bmatrix} 1 & 0 & 0 \\ 0 & 1 & 0 \\ 0 & 0 & 1 \\ 1 & 0 & 0 \\ 0 & 1 & 0 \end{bmatrix}$$

If A = (1 0 0 0 0) is passed on in the relational matrix E, the instantaneous vector, AE = (1 0 0) implies that the student is a good student . Now let AE = B, $BE^T$ = (1 0 0 1 0) which implies that the teaching is good and he / she is a kind teacher. Let $BE^T$ = $A_1$, $A_1E$ = (2 0 0) after thresholding we get $A_1E$ = (1 0 0) which implies that the student is good, so on and so forth.

Now we have worked on the real world problem, which studies the Employee-Employer Relationship model in the following section. The reader is expected to develop several other model relationships and illustrations.

## 1.6 Models Illustrating FRM and combined FRMs

Now we in this section give two illustrations, which are data, collected from the real world problems analyzed using both FRMs and combined FRMs.

### 1.6.1: Use of FRMs in Employee - Employer Relationship model

The employee-employer relationship is an intricate one. For, the employers expect to achieve performance in quality and production in order to earn profit, on the other hand employees need good pay with all possible allowances. Here we have taken three experts opinion in the study of Employee and Employer model. The three experts whose opinions are taken are the Industry Owner, Employees' Association Union Leader and an Employee. The data and the opinion are taken only from one industry. Using the opinion we obtain the hidden patterns.

The following concepts are taken as the nodes relative to the employee. We can have several more nodes and also several experts' opinions for it a clearly evident theory which professes that more the number of experts the better is the result.

We have taken as the concepts / nodes of domain only 8 notions which pertain to the employee.

| | | |
|---|---|---|
| $D_1$ | – | Pay with allowances and bonus to the employee |
| $D_2$ | – | Only pay to the employee |
| $D_3$ | – | Pay with allowances (or bonus) to the employee |
| $D_4$ | – | Best performance by the employee |
| $D_5$ | – | Average performance by the employee |
| $D_6$ | – | Poor performance by the employee |
| $D_7$ | – | Employee works for more number for hours |
| $D_8$ | – | Employee works for less number of hours. |

$D_1, D_2, \ldots, D_8$ are elements related to the employee space which is taken as the domain space.

We have taken only 5 nodes / concepts related to the employer in this study.

These concepts form the range space which is listed below.

| | | |
|---|---|---|
| $R_1$ | – | Maximum profit to the employer |
| $R_2$ | – | Only profit to the employer |
| $R_3$ | – | Neither profit nor loss to the employer |
| $R_4$ | – | Loss to the employer |
| $R_5$ | – | Heavy loss to the employer |

The directed graph as given by the employer is given in Figure 1.6.1. The associated relational matrix $E_1$ of the employer as given by following.

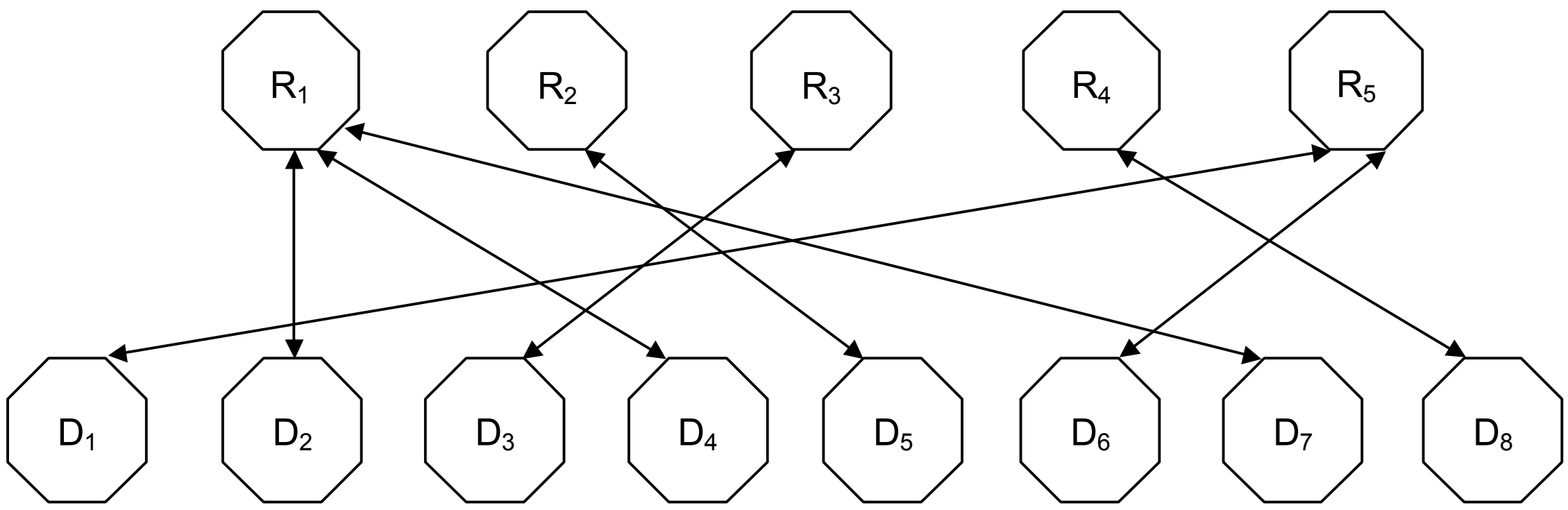

**FIGURE: 1.6.1**

$$E_1 = \begin{bmatrix} 0 & 0 & 0 & 0 & 1 \\ 1 & 0 & 0 & 0 & 0 \\ 0 & 0 & 1 & 0 & 0 \\ 1 & 0 & 0 & 0 & 0 \\ 0 & 1 & 0 & 0 & 0 \\ 0 & 0 & 0 & 0 & 1 \\ 1 & 0 & 0 & 0 & 0 \\ 0 & 0 & 0 & 1 & 0 \end{bmatrix}$$

Suppose we consider the node $D_1$ to be in the on state and rest of the nodes in the off state. That is the employee is paid with allowance and bonus i.e. $G_1 = (1\ 0\ 0\ 0\ 0\ 0\ 0\ 0)$ The effect of $G_1$ on the expert systems $E_1$ is $G_1E_1 = (0\ 0\ 0\ 0\ 1)$, it belongs to the range space R.

Let $G_1E_1 = H_1 = (0\ 0\ 0\ 0\ 1)$, $H_1E_1^T = (1\ 0\ 0\ 0\ 0\ 1\ 0\ 0) \in D$, after updating and thresholding the instantaneous vector at each stage we obtain the following chain

$$G_1 \to H_1 \to G_2 \to H_1$$

i.e., $G_1$ is a fixed point and according to the opinion of the employer who is taken as an expert we see if the employee is paid with allowances and bonus the company suffers a heavy loss due to the poor performance of the employee.

Suppose we input the vector $G_3$ = (0 0 0 1 0 0 0 0) which indicates that the node $D_4$ "best performance by the employee" is in the on state, Effect of $G_3$ on the system, $G_3 E_1$ = (1 0 0 0 0) = $H_3 \in$ R. $H_3E_1^T$ = (0 1 0 1 0 0 1 0) $\in$ D.

After updating and thresholding the instantaneous vector at each stage we obtain the following chain:

$$G_3 \rightarrow H_3 \rightarrow G_4 \rightarrow H_3$$

We see from the above the resultant is also a fixed point.

According to the first expert we see the company enjoys maximum profit by giving only pay to the employee in spite of his best performance and putting in more number of working hours.

Now to analyze the effect of employee on the employer let us input the vector $S_1$ = (1 0 0 0 0) indicating the on state of the node $R_1$ (maximum profit to the employer)

$$S_1E^T{}_1 = (0\ 1\ 0\ 1\ 0\ 0\ 1\ 0) = T_1 \in D.$$

After updating and thresholding the instantaneous vector at each stage we obtain the following chain

$$S_1 \rightarrow T_1 \rightarrow S_1.$$

From the above chain we see that when $S_1$ is passed on, it is a limit cycle.

The company enjoys maximum profit by getting best performance and more number of working hours from the employee and by giving only pay to the employee.

Suppose we input the vector $S_2$ = (0 0 0 0 1) indicating the on state of the node $R_5$.

$$S_2\, E_1^T = T_2 = (1\ 0\ 0\ 0\ 0\ 1\ 0\ 0) \in D.$$

After updating and thresholding the instantaneous vector at each stage we obtain the following chain

$$S_2 \rightarrow T_2 \rightarrow S_2.$$

Thus in his industry the employer feels that the employee performs poorly inspite of getting pay with allowances and bonus, company suffers a heavy loss.

The union leader of the same company was asked to give his opinion keeping the same nodes for the range space and the domain space i.e. as in case of the first expert. The directed graph given by the union leader of the company is given by the following diagram.

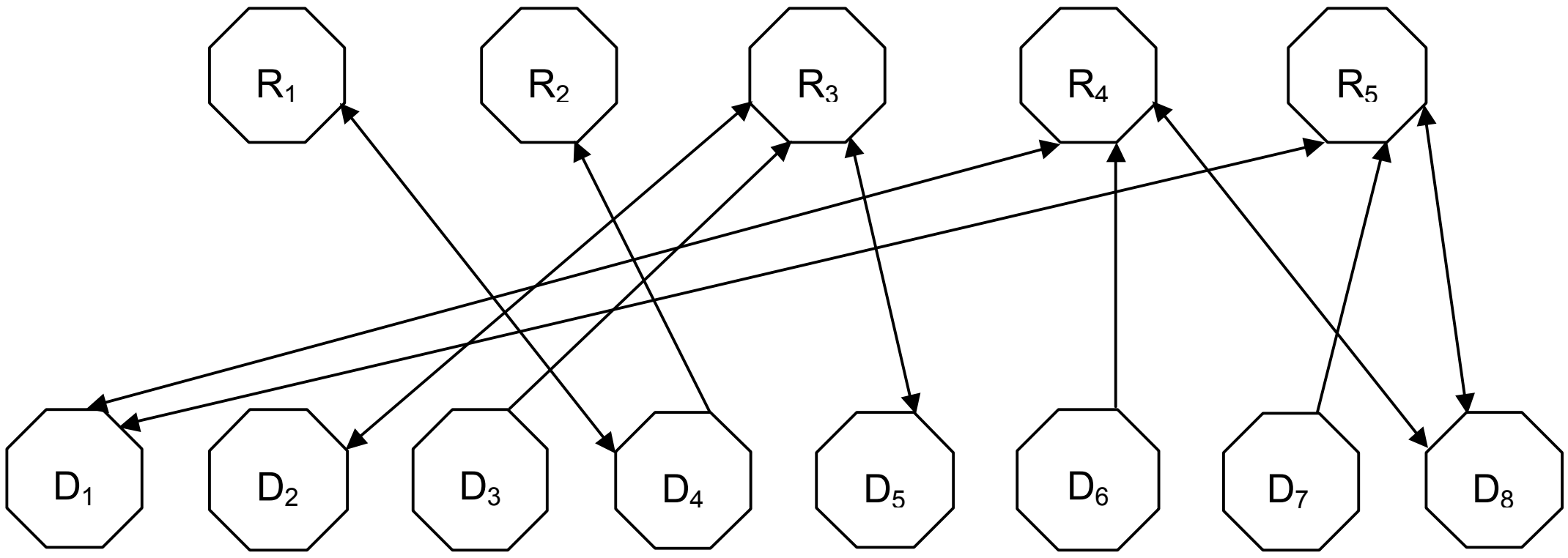

**FIGURE: 1.6.2**

The related matrix of the directed graph given by the second expert is given by $E_2$

$$E_2 = \begin{bmatrix} 0 & 0 & 0 & 1 & 1 \\ 0 & 0 & 1 & 0 & 0 \\ 0 & 0 & 1 & 0 & 0 \\ 1 & 1 & 0 & 0 & 0 \\ 0 & 0 & 1 & 0 & 0 \\ 0 & 0 & 0 & 1 & 0 \\ 0 & 0 & 0 & 0 & 1 \\ 0 & 0 & 0 & 1 & 1 \end{bmatrix}.$$

Let us input the vector
$L_1$ = (1 0 0 0 0 0 0 0) which indicates that the node $D_1$ viz. employee is paid with allowance and bonus is in the on state.

$L_1E_2 \rightarrow$ (1 0 0 1 1) = $N_1 \in R$
$N_1E_2^T \rightarrow$ (1 0 0 1 0 1 1 1) = $L_2$

$L_2$ is got only after updating and thresholding. Now $L_2E_2 \rightarrow N_2$ = (11011) which is got after updating and thresholding

$N_2 E_2^T \rightarrow$ (1 0 0 1 0 1 1 1) = $L_3$

$$L_1 \rightarrow N_1 \rightarrow L_2 \rightarrow N_2 \rightarrow L_3 = L_2 \rightarrow N_2$$

is a fixed point.

Thus the union leader's viewpoint is in a way very balanced though seemingly contradictory. When the company pays an employee with pay, allowances and bonus the company may have maximum profit or only profit. Or on the contrary if the employee does not put forth more number of working hours or puts in only the required number of working hours and if his/her performance is poor or average certainly the company will face loss and or heavy loss. Thus both can happen depending highly on the nature of the employees.

Now we proceed on to study the same industry using the third expert as an employee of the same industry.

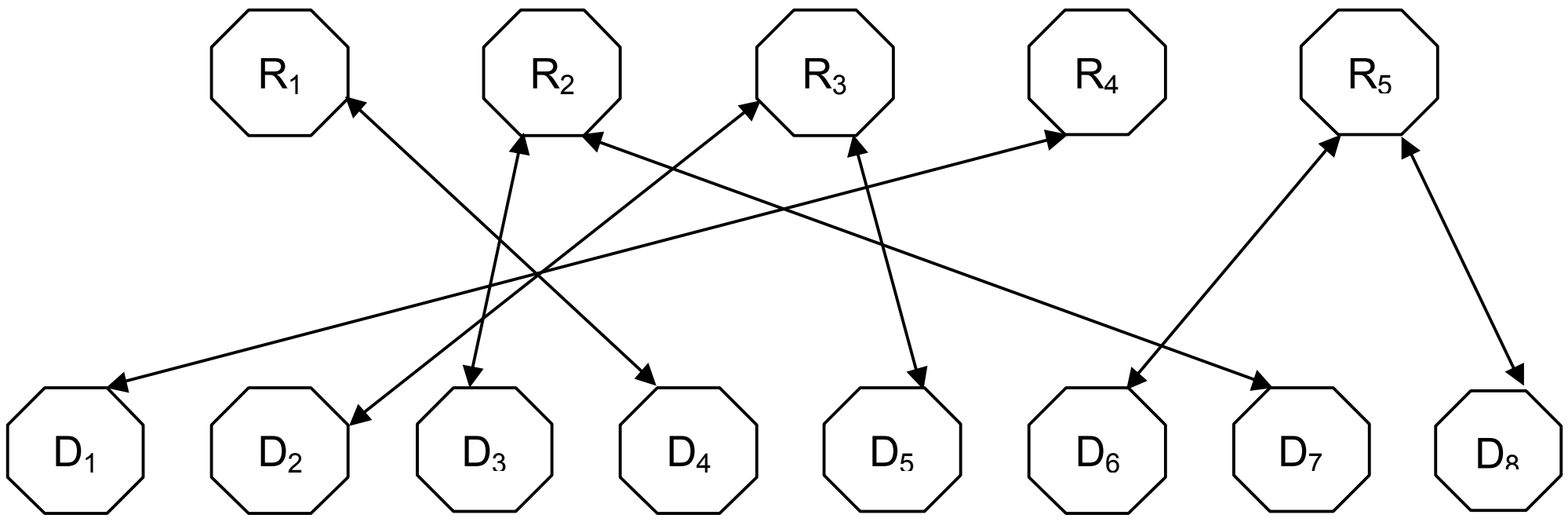

**FIGURE: 1.6.3**

The associated relational matrix $E_3$ of the third expert's opinion got from the directed graph given in Figure 1.6.3 is as follows:

$$E_3 = \begin{bmatrix} 0 & 0 & 0 & 1 & 0 \\ 0 & 0 & 1 & 0 & 0 \\ 0 & 1 & 0 & 0 & 0 \\ 1 & 0 & 0 & 0 & 0 \\ 0 & 0 & 1 & 0 & 0 \\ 0 & 0 & 0 & 0 & 1 \\ 0 & 1 & 0 & 0 & 0 \\ 0 & 0 & 0 & 0 & 1 \end{bmatrix}$$

Let us input the vector $U_1$ = (1 0 0 0 0 0 0 0) indicating the on state of the note $D_1$, effect of $U_1$ on the system $E_3$ is

$$U_1 E_3 = V_1 = (0\ 0\ 0\ 1\ 0) \in R.$$

After thresholding and updating the vector. Now $V_1E^T_3$ = (1 0 0 0 0 0 0 0) = $U_1$. So

$$U_1 \to V_1 \to U_1.$$

Hence when $U_1$ is passed on it is a limit cycle. So according to the employee when the company gives pay, bonus/ perks and allowances the company suffers a loss.

Suppose we input vector $U_2$ = (0 0 0 1 0 0 0 0) which indicates the on state of the node $D_4$, employee puts in his best performance; effect of $U_2$ on the system $E_3$ is $U_2E_3$ = (1 0 0 0 0) $\in$ R.

After updating and thresholding the instantaneous vector at each state we get the following chain:

$$U_2 \rightarrow V_2 \rightarrow U_2$$

when $U_2$ is passed on, it is a limit cycle.

It is left for the reader to input any instantaneous state vector and obtain the resultant.

Now we proceed on to find the combined FRMs. We take the opinion of the three experts discussed above and find their opinions.

We first draw the directed graph of all the three experts, which is given by the Figure 1.6.4.

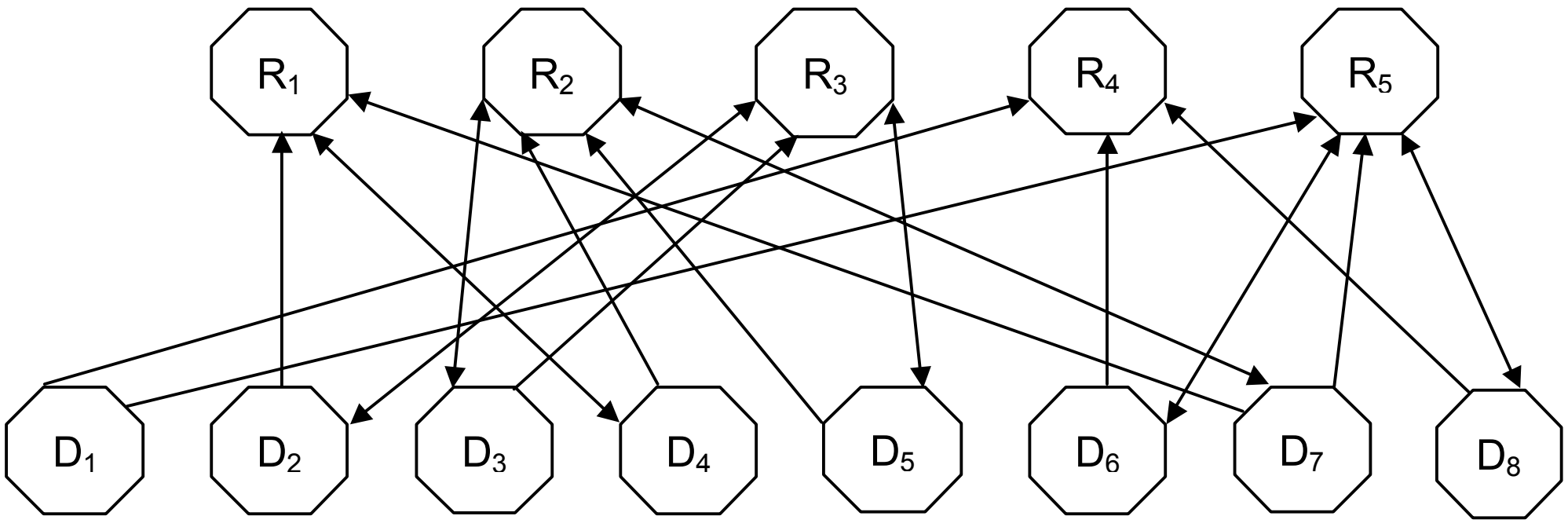

**FIGURE: 1.6.4**

The corresponding fuzzy relational matrix is given as the sum of the three fuzzy relational matrices $E_1$ , $E_2$ and $E_3$. Let

$$E = E_1 + E_2 + E_3 = \begin{bmatrix} 0 & 0 & 0 & 2 & 2 \\ 1 & 0 & 2 & 0 & 0 \\ 0 & 1 & 2 & 0 & 0 \\ 3 & 1 & 0 & 0 & 0 \\ 0 & 1 & 2 & 0 & 0 \\ 0 & 0 & 0 & 1 & 2 \\ 1 & 1 & 0 & 0 & 1 \\ 0 & 0 & 0 & 2 & 2 \end{bmatrix}.$$

Let us now input the state vector $X_1 = (1\ 0\ 0\ 0\ 0\ 0\ 0\ 0)$, indicating the on state of the node $D_1$, effect of $X_1$ on the combined system is

$$X_1 E = (0\ 0\ 0\ 2\ 2) \rightarrow (0\ 0\ 0\ 1\ 1) \in R.$$

Let $Y_1 = (0\ 0\ 0\ 1\ 1)$,

$$Y_1 E^T = (4\ 0\ 0\ 0\ 0\ 3\ 1\ 4) \rightarrow (1\ 0\ 0\ \ 0\ 0\ 1\ 1\ 1) \in D.$$

($\rightarrow$ denotes after updating and thresholding the vector).

Let $X_2 = (1\ 0\ 0\ 0\ 0\ 1\ 1\ 1)$; $X_2E = (1\ 1\ 0\ 5\ 7) \rightarrow (1\ 1\ 0\ 1\ 1) = Y_2 \in R$.

$$Y_2 E^T = (4\ 1\ 1\ 4\ 1\ 3\ 2\ 4) \rightarrow (1\ 1\ 1\ 1\ 1\ 1\ 1\ 1) \in D.$$

Let $X_3 = (1\ 1\ 1\ 1\ 1\ 1\ 1\ 1)$

$$X_3 E = (5\ 4\ 6\ 5\ 7) \rightarrow (1\ 1\ 1\ 1\ 1) = Y_3.$$

$$X_1 \rightarrow Y_1 \rightarrow X_2 \rightarrow Y_2 \rightarrow X_3 - Y_3 \rightarrow X_3 \rightarrow Y_3$$

Thus $X_1$ is a fixed point.

In fact a lot of discussion can be made on such fixed points, for the opinion of the employee and the employer or the union leader happens to be very contradictory that is why it is visibly seen that all the nodes both in the domain space and the range space becomes on, at the very advent of seeing the effect of only one node of the domain space to be in the on state i.e. $D_1$ – pay, bonus and allowances. The interpretations are also to be carefully given. On one hand when the employee is paid with pay bonus and allowances, they may work so well to see that the company runs with maximum profit which will automatically turn on all other linguistic nodes like the notions of only profit to the employer, neither profit nor loss to the employer etc.

Likewise in case of the domain nodes all nodes come to the on state. Thus we are not able to distinguish or give a nice interpretation of the state vectors. That is why if at certain stages when the expert is not able to give opinion, if we use the concept of neutrosophic logic when the relation is an indeterminate, we may not be facing such situations. To avoid this sort of contradictory fixed points / limit cycles we proceed on to construct and analyze using neutrosophic logic in Chapter 2.

For more about this illustration refer [124, 125, 136].

The fuzzy relational maps can also be used in the prediction using the past year data. So in this case it is not only the expert's opinion but the interpretation of the data and the opinion analyzed using other methods. We for the first time give such analysis in case of cement industries. We analyze how to maximize the production by giving the maximum satisfaction to the employees, we have used the data from Ramco Cement Industries for the years 1995-2001.

### 1.6.2: Maximizing Production by Giving Maximum Satisfaction to Employee using FRMs

At present India needs huge quantity of cement for the construction of dwellings, houses, apartments, dams, reservoirs, village roads, flyovers and cementing of the entire national highways. Hence maximizing cement production in each and every factory is essential which is directly dependent on the relationship between the employer and the employees. Employer and Employees congenial relationship (Industrial Harmony) is a most complicated one.

For example employer expects to achieve consistent production, quality product at optimum production to earn profit. However there may be profit, no loss and no profit, loss in business or heavy loss depending on various factors such as demand and supply, rent, electricity, raw materials transportation and consumables, safety, theft, medical aid, employees co-operation and enthusiasm. At the same time some of the expectations of the employees are pay, various allowances, bonus, welfare unanimity such as uniforms, hours of work and job satisfaction etc. Since it is a problem of both employer and employee we analyse taking into consideration the problem on both sides using fuzzy theory in general and FRMs in particular; as the very nature of the problem is one, which is dominated by uncertainties.

Some industries provide

1. Additional incentive for regularity in attendance
2. Additional days of vacation
3. Incentive linked with production level
4. Award for good suggestion for operation
5. Forming quality circles in various department for finding out new ways and means to carry out jobs economically, quickly and with safety
6. Provide hospital for employees and their family
7. Credit society and co-operative stores
8. Awards for highest sale made by salesman
9. Overtime wages
10. Loan facilities.

The employer's final achievement is to maximize the production level; so, when production linked incentive was introduced in a cement industry, the employees were so enthusiastic to run the kiln continuously and efficiently to ensure more than the targeted minimum production. Employees voluntarily offered their services in case of kiln stoppages and made the kiln to run again with a minimum shutdown. Small incentive to employees made the industry to earn considerable profit out of maximum production than the previous year. This study is significant because most of the cement industries have the common types of problems between employee and the employer.

Here in this analysis we are also trying to give a form of best relationship between the employees and employers, because the smooth relationship between the employers and employees is an important one for the cement industries to maximize the production. Probably, so far no one has approached the problem of finding a best form of relationship between employers and employees using FRMs. We use the connection fuzzy relational matrix to find the best form of relationship between the employees and employers and aim to maximize the production level.

Here, we approach the employee and employer problem of finding a best form of relationship between employee and employer. Thus this study tries to improve maximize the production level in cement industries. A good relationship between the employee and employer is very essential to run the industry successfully with least friction.

For this study the raw data is obtained from the cement industrialists and the employees; which is converted into a relational map. The relational map takes edge values as positive real numbers. If the number of concepts in the domain space is M and that of the range space is N we get an $M \times N$ causal relational matrix, which we call as the relational matrix.

Let $X_i$ and $Y_j$ denote the two nodes of the relational map. The nodes $X_i$ and $Y_j$ are taken as the attributes of the employee and employer respectively. The directed edge from $X_i$ to $Y_j$ and $Y_j$ to $X_i$ denotes the causality relations. Every edge in the relational map is weighted with a positive real number . Let $X_1, X_2, \ldots, X_M$ be the nodes of the domain space X and $Y_1, Y_2, \ldots, Y_N$ be the nodes of the range space Y of the relational map.

After obtaining a relational matrix, we find an average matrix; for simplification of calculations. The average matrix is then converted into a fuzzy matrix. Finally using the different parameter $\alpha$ (membership grades) we identify the best form of relationship between employee and employer in cement industries.

We approach the employee and the employer problems in cement industries using fuzzy relational matrix. Probably, so far no one has approached the employee and employer problems via a fuzzy relational matrix method. The raw data under investigation is classified under six broad heads viz. $Y_1, Y_2, Y_3, Y_4, Y_5, Y_6$ which are nodes associated with the employer. $X_1, X_2, \ldots, X_8$ are the eight broad heads, which deal with the attributes of the employees. It is pertinent to mention here one can choose any number of nodes associated with employees or employer. We prefer to choose $8 \times 6$ matrix for spaces X and Y respectively. Using the fuzzy matrix we find a best relationship between employer and employee in two-stages.

In the first stage, we convert the data into relational map, the relational matrix obtained from the relational map is then converted into an average matrix. In this matrix, take along the columns the data related to the employer and along the rows the data related to the employee.

In the second stage, we convert the average matrix into fuzzy matrix using different parameter $\alpha$ ($\alpha$-membership grade ), and we give the graphical illustration to the fuzzy matrix, which enables one to give a better result.

**DESCRIPTION OF THE PROBLEM**

Using the data from (Ramco) cement industries, we analyse the data via a fuzzy relational matrix and obtain a best relationship between employee and employer in cement industries. Hence using these suggestions got from the study, cement industries can maximize their production. Here we find the best relationship between the employee and the employer.

Here we consider only eight attributes of the employee and its six effects on employers in production level. Hence, in a cement industry the attributes of the employee are described by the eight nodes $X_1, X_2, X_3, X_4, X_5, X_6, X_7, X_8$ and the six effects on the employers in production level are described by the nodes $Y_1, Y_2, Y_3, Y_4, Y_5, Y_6$.

In the first stage of the problem, the data obtained from the cement industry (1995-2001) are used in the illustration of the problem and is finally applied in the fuzzy matrix model to verify the validity of the method described. The relational map is obtained using the above nodes.

From the relational map we get the M × N relational matrix. Let M represent the six effects of employers. Let N represent the eight attributes of employees feelings. Production levels are treated as rows and the various attributes of employee feelings are treated as columns resulting in the M × N relational matrix into an average matrix.

In the second step we use mean, standard deviation of each of the columns of the M × N matrix and parameter $\alpha$ (membership grade $\alpha \in [0,1]$ to convert the average matrix ($a_{ij}$) into the fuzzy matrix ($b_{ij}$); where i represents the $i^{th}$ row and j represents the $j^{th}$column. We calculate the mean $\mu_j$ and the standard deviation $\sigma_j$ for each attribute j, for j = 1, 2,…, n using the average matrix ($a_{ij}$).

For varying values of the parameter $\alpha$ where $\alpha \in [0,1]$.
We determine the values of the entry $b_{ij}$ in the average matrix using the following rule:

$$b_{ij} = \begin{cases} 0 & \text{if} \quad a_{ij} \le \mu_j - \alpha * \sigma_j \\ \dfrac{a_{ij} - (\mu_j - \alpha * \sigma_j)}{(\mu_j + \alpha * \sigma_j) - (\mu_j - \alpha * \sigma_j)} & \text{if} \quad a_{ij} \in (\mu_j - \alpha * \sigma_j, \mu_j + \alpha * \sigma_j) \\ 1 & \text{if} \quad a_{ij} \ge \mu_j + \alpha * \sigma_j \end{cases}$$

Here '*' denotes the usual multiplication.
Thus for different value of $\alpha$, we obtain different fuzzy matrices. Finally we add up the rows of each fuzzy matrix and we define the new fuzzy membership function and allocate a value between [0, 1] to each row sum. Here the highest membership grade gives the best form of relationship between employee and employer, which maximizes the production level.

**STAGE 1**

The cement industry is having eight types of attributes related to the employee, which are as follows:

$X_1$ – Salaries and wages
$X_2$ – Salaries with wages and bonus to the employee
$X_3$ – Bonus to the employee
$X_4$ – Provident fund (PF) to the employee
$X_5$ – Employee welfare medical
$X_6$ – Employee welfare LTA (Leave travel allowances)
$X_7$ – Employee welfare others
$X_8$ – Staff Training expenses.

Now we take the six effects on the employer when he proposes to pay and get the work done by the employee to maximize the production level, which are as follows:

$Y_1$ – 1995 – 1996 Production level is 8,64,685 tonnes
$Y_2$ – 1996 – 1997 Production level is 7,74,044 tonnes
$Y_3$ – 1997 – 1998 Production level is 7,22,559 tonnes
$Y_4$ – 1998 – 1999 Production level is 8,19,825 tonnes

$Y_5$ – 1999 – 2000 Production level is 7,98,998 tonnes
$Y_6$ – 2000 – 2001 Production level is 7,79,087 tonnes.

The problem now is to find to best relationship between the employee and employer, so we convert the $X_1$, $X_2$, $X_3$, $X_4$, $X_5$, $X_6$, $X_7$, $X_8$ and $Y_1$, $Y_2$, $Y_3$, $Y_4$, $Y_5$, $Y_6$ into the relational map. The relational map is shown in the Figure 1.6.5.

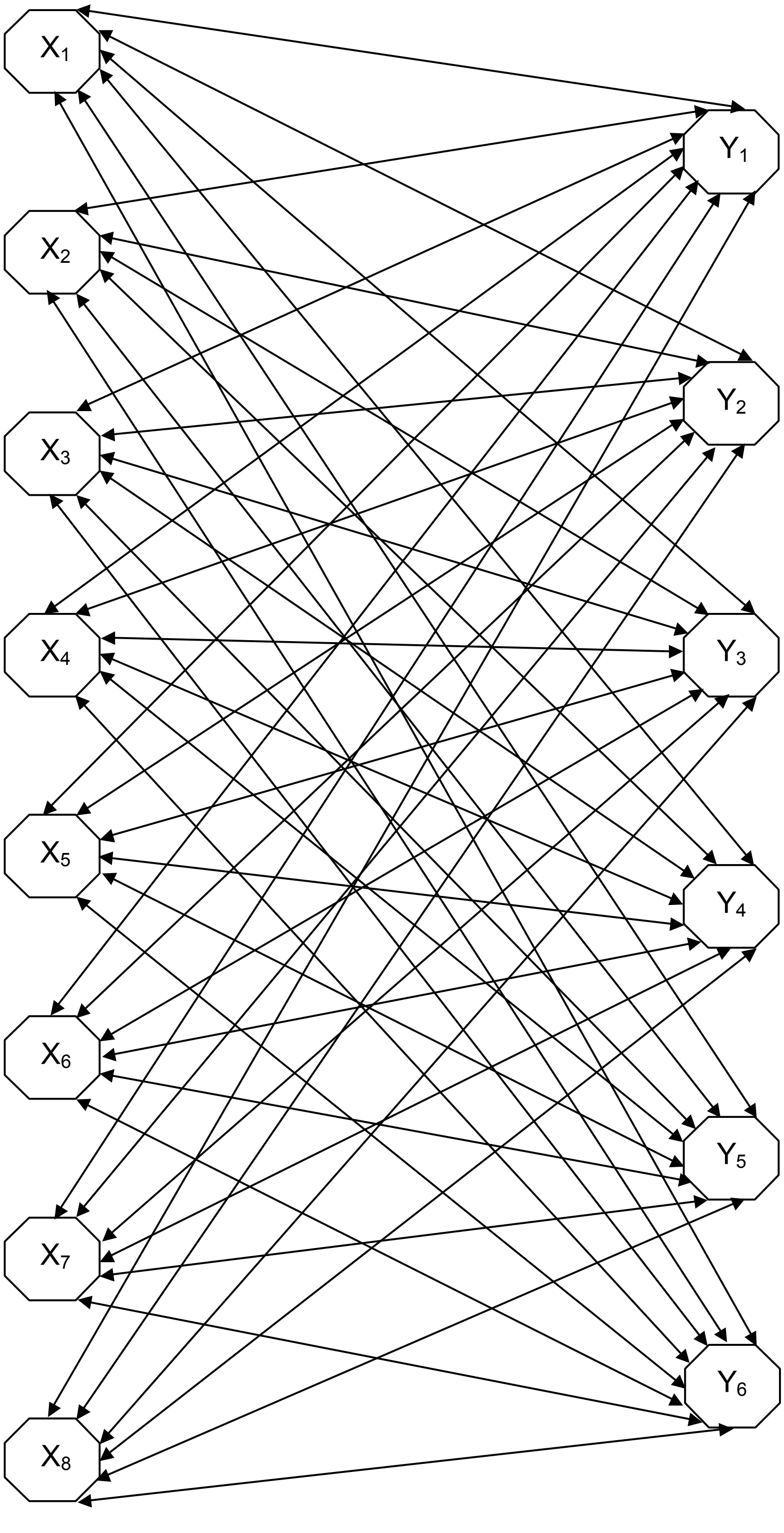

**FIGURE:1.6.5**

We obtain the initial 8 × 6 relational matrix from the relational map using weights from the data.

$$\begin{bmatrix} 70.99 & 82.61 & 11.61 & 5.856 & 2.239 & 1.23 & 7.19 & 0.867 \\ 74.33 & 85.05 & 10.71 & 5.558 & 2.058 & 0.88 & 6.02 & 0.427 \\ 74.25 & 84.29 & 10.04 & 5.969 & 2.704 & 0.769 & 6.068 & 0.144 \\ 71.08 & 81.42 & 10.34 & 7.070 & 2.386 & 1.197 & 7.411 & 0.498 \\ 70.65 & 80.85 & 10.19 & 7.102 & 2.589 & 1.112 & 8.093 & 0.238 \\ 72.03 & 80.06 & 8.028 & 6.297 & 3.276 & 1.373 & 8.563 & 0.413 \end{bmatrix}.$$

Now for simplification of calculation, we convert the relational matrix into the average matrix which is given by ($a_{ij}$).

$$(a_{ij}) = \begin{bmatrix} 35.49 & 41.30 & 5.805 & 2.928 & 1.119 & 0.615 & 3.595 & 0.433 \\ 37.16 & 42.53 & 5.355 & 2.779 & 1.029 & 0.44 & 3.01 & 0.213 \\ 37.13 & 42.15 & 5.02 & 2.985 & 1.352 & 0.385 & 3.034 & 0.072 \\ 35.54 & 40.71 & 5.17 & 3.535 & 1.193 & 0.598 & 3.076 & 0.249 \\ 35.33 & 40.43 & 5.095 & 3.551 & 1.295 & 0.556 & 4.046 & 0.119 \\ 36.02 & 40.03 & 4.014 & 3.148 & 1.638 & 0.687 & 4.281 & 0.207 \end{bmatrix}.$$

To find a best form of relationship between the employee and employer, we use mean, standard deviation and the parameter $\alpha \in [0, 1]$ and proceed to the second stage of the problem to convert the average matrix into a fuzzy matrix.

**STAGE 2**

In the second stage we use mean ($\mu$), standard deviation ($\sigma$) and the parameter $\alpha \in [0, 1]$ to find the best form of relationship between employer and employee. To convert the above average matrix into a fuzzy matrix ($b_{ij}$).

$$b_{ij} \in \left[ 0, \quad \frac{(a_{ij} - \mu_i - \alpha * \sigma_j)}{(\mu_j + \alpha * \sigma_j) - (\mu_j - \alpha * \sigma_j)}, \quad 1 \right]$$

where i represents the $i^{th}$ row and j represents the $j^{th}$ column. The value of the entry $b_{ij}$ corresponding to each intersection is determined from this interval. This interval is obtained strictly by using the average and standard deviation calculated from the raw data.

The calculations are as follows: First we calculate the $\mu_j$ corresponding to each column of the matrix. $\mu_1 = 36.1116$, $\mu_2 = 41.191$, $\mu_3 = 5.076$, $\mu_4 = 3.154$, $\mu_5 = 1.271$, $\mu_6 = 0.547$, $\mu_7 = 3.612$, $\mu_8 = 0.2155$, where $\mu_j$ are the average of each column respectively for j = 1, 2,…, 8.

Now the standard deviation $\sigma_j$ is calculated as follows.

**TABLE 1: MEAN AND STANDARD DEVIATION OF COLUMN 1**

When $\mu_1 = 36.1116$

| d | $d^2$ |
|---|---|
| 0.6216 | 0.386386 |
| 1.048 | 1.09830 |
| 1.0184 | 1.03713 |
| 0.5716 | 0.32672 |
| 0.7816 | 0.6108985 |
| 0.0916 | 0.008390 |
| $\Sigma d = 4.13280$ | $\Sigma d^2 = 3.4678245$ |

$$\sigma_1 = \sqrt{0.577970 - 0.4744454}$$

$\sigma_1 = 0.32175$.

In a similar way, the values of $\sigma_j$'s are as follows:

$\sigma_2 = 0.41368$ when $\mu_2 = 41.191$,
$\sigma_3 = 0.390285$ when $\mu_3 = 5.076$,
$\sigma_4 = 0.141763$ when $\mu_4 = 3.154$,
$\sigma_5 = 0.1162292$ when $\mu_5 = 1.271$,
$\sigma_6 = 0.05250$ when $\mu_6 = 0.547$,
$\sigma_7 = 0.253748$ when $\mu_7 = 3.612$,
$\sigma_8 = 0.07790$ when $\mu_8 = 0.2155$.

Now we calculate $\mu_j - \alpha * \sigma_j$ and $\mu_j + \alpha * \sigma_j$ when $\alpha = 0.1, 0.2, 0.3, 0.4, 0.5, 0.6, 0.7, 0.8, 0.9$ and 1.

**TABLE 2: FOR THE MEMBERSHIP GRADE $\alpha = 0.1$**

| j | $\mu_j - \alpha * \sigma_j$ | $\mu_j + \alpha * \sigma_j$ |
|---|---|---|
| 1 | 36.0794250 | 36.1437750 |
| 2 | 41.1496320 | 41.2323680 |
| 3 | 5.03697150 | 5.11502850 |
| 4 | 3.13982370 | 3.16817630 |
| 5 | 1.25937700 | 1.28262200 |
| 6 | 0.54175000 | 0.55225000 |
| 7 | 3.58662500 | 3.63737400 |
| 8 | 0.20771000 | 0.22329000 |

From the table the corresponding fuzzy relational matrix for the parameter value $\alpha = 0.1$ is as follows:

$$(b_{ij}) = \begin{array}{c} \\ Y_1 \\ Y_2 \\ Y_3 \\ Y_4 \\ Y_5 \\ Y_6 \end{array}\begin{array}{c} \begin{array}{cccccccc} X_1 & X_2 & X_3 & X_4 & X_5 & X_6 & X_7 & X_8 \end{array} \\ \begin{bmatrix} 0 & 1 & 1 & 0 & 0 & 1 & 0.165 & 1 \\ 1 & 1 & 1 & 0 & 0 & 0 & 0 & 0.339 \\ 1 & 1 & 0 & 0 & 1 & 0 & 0 & 0 \\ 0 & 0 & 1 & 1 & 0 & 1 & 1 & 1 \\ 0 & 0 & 0.743 & 1 & 1 & 1 & 1 & 0 \\ 0 & 0 & 0 & 0.288 & 1 & 1 & 1 & 0 \end{bmatrix} \end{array}.$$

The corresponding row sum of the above matrix is given by

| | | |
|---|---|---|
| $R_1$ | = | The first row sum = 4.165 |
| $R_2$ | = | The second row sum = 3.3390 |
| $R_3$ | = | The third row sum = 3 |
| $R_4$ | = | The fourth row sum = 5 |
| $R_5$ | = | The fifth row sum = 4.743 |
| $R_6$ | = | The sixth row sum = 3.288 |

where $R_1, R_2 \ldots R_6$ are nothing but the corresponding years of $Y_1, Y_2, \ldots Y_6$.

Now we define the new fuzzy membership function for graphical illustration, which converts the row sum to take values in the interval [0, 1].

We define the new fuzzy membership function as follows:

$$\mu_x(R_i) = \begin{cases} 1 & \text{if Row sum of max value} \\ \dfrac{R_i - \text{Rowsum of min value}}{\text{Rowsum of max value} - \text{Rowsum of min value}} & \text{if } \begin{array}{l}\text{Row sum of min value} \le R_i \\ \le \text{Row sum of max value}\end{array} \\ 0 & \text{if Row sum of min value} \end{cases}$$

Using the membership function the corresponding values of

$\mu_x(R_i)$ is 0.5825, 0.1695, 0, 1, 0.87, 0.144

respectively where i = 1, 2, 3, 4, 5, 6.

Graphical illustration for $\alpha = 0.1$ is shown in graph 1.

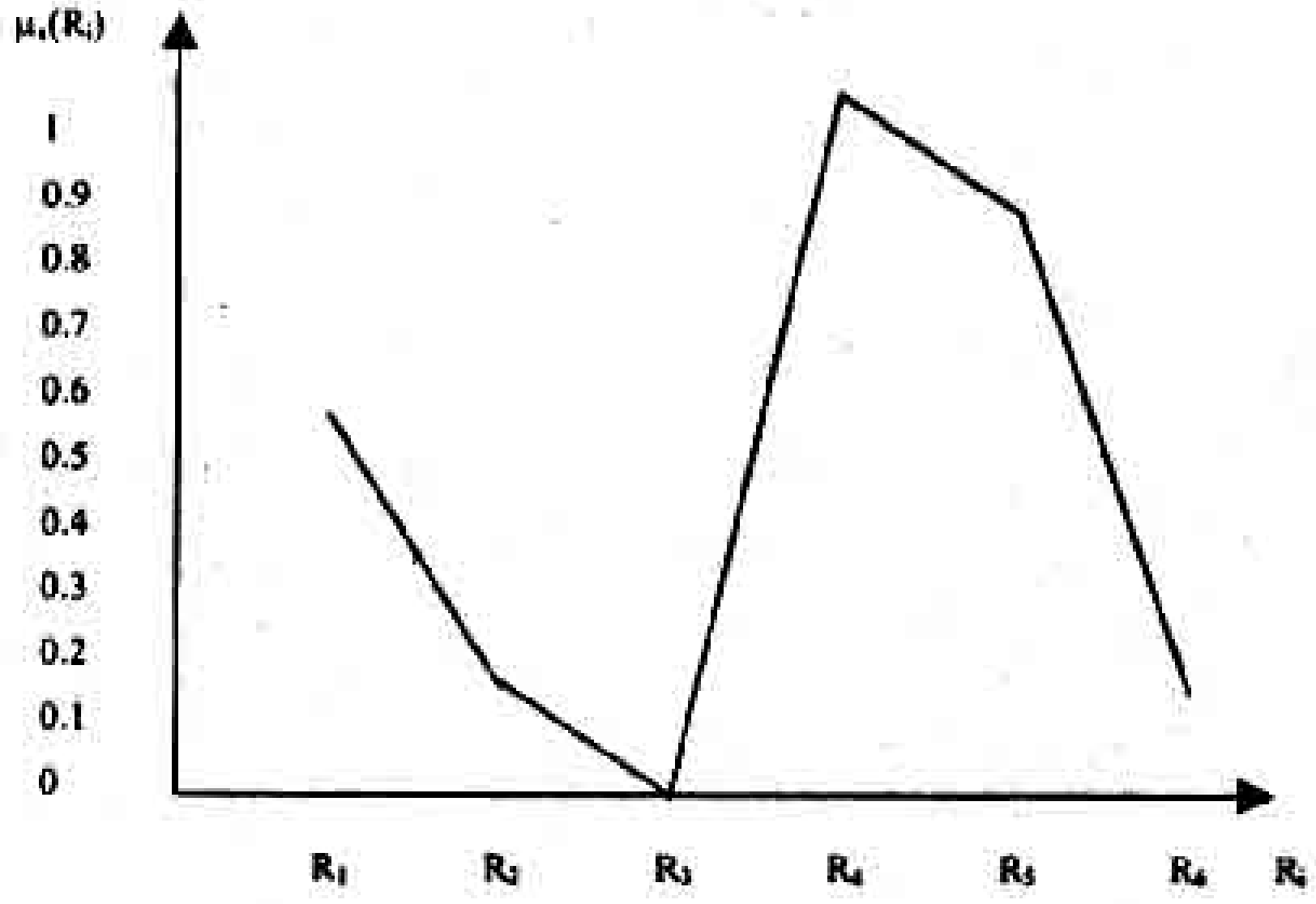

GRAPH 1: GRAPHICAL ILLUSTRATION FOR α = 0.1

Here the fourth row sum $R_4$ is getting the highest membership grade that is the membership grade is 1 and the fifth row sum $R_5$ is getting the next highest membership grade that is the membership grade is 0.87.

**TABLE 3: FOR THE MEMBERSHIP GRADE α = 0.2**

| J | $\mu_j - \alpha * \sigma_j$ | $\mu_j + \alpha * \sigma_j$ |
|---|---|---|
| 1 | 36.04720 | 36.175950 |
| 2 | 41.10826 | 41.273736 |
| 3 | 4.997900 | 5.1540500 |
| 4 | 3.125647 | 3.1823520 |
| 5 | 1.247754 | 1.2942450 |
| 6 | 0.536500 | 0.5575000 |
| 7 | 3.561250 | 3.6627496 |
| 8 | 0.199920 | 0.2310800 |

From the table the fuzzy matrix corresponding to the parameter value of α = 0.2 is as follows:

$$(b_{ij}) = \begin{array}{c} \\ Y_1 \\ Y_2 \\ Y_3 \\ Y_4 \\ Y_5 \\ Y_6 \end{array}\begin{array}{c} \begin{array}{cccccccc} X_1 & X_2 & X_3 & X_4 & X_5 & X_6 & X_7 & X_8 \end{array} \\ \begin{bmatrix} 0 & 1 & 1 & 0 & 0 & 1 & 0.333 & 1 \\ 1 & 1 & 1 & 0 & 0 & 0 & 0 & 0.419 \\ 1 & 1 & 0.1415 & 0 & 1 & 0 & 0 & 0 \\ 0 & 0 & 1 & 1 & 0 & 1 & 1 & 1 \\ 0 & 0 & 0.6218 & 1 & 1 & 0.9285 & 1 & 0 \\ 0 & 0 & 0 & 0.394 & 1 & 1 & 1 & 0.227 \end{bmatrix} \end{array}$$

The corresponding row sum of the above matrix is given by

| | | |
|---|---|---|
| $R_1$ | = | The first row sum = 3.333 |
| $R_2$ | = | The second row sum = 3.419 |
| $R_3$ | = | The third row sum = 3.1415 |
| $R_4$ | = | The fourth row sum = 5 |
| $R_5$ | = | The fifth row sum = 4.5503 |
| $R_6$ | = | The sixth row sum = 3.6210. |

Using the membership function the corresponding values of $\mu_x$ ($R_i$) is 0.1030, 0.149, 0, 1, 0.7580, 0.258 respectively where i = 1, 2, 3, 4, 5, 6. Graphical illustration for $\alpha$ = 0.2 is shown in graph 2.

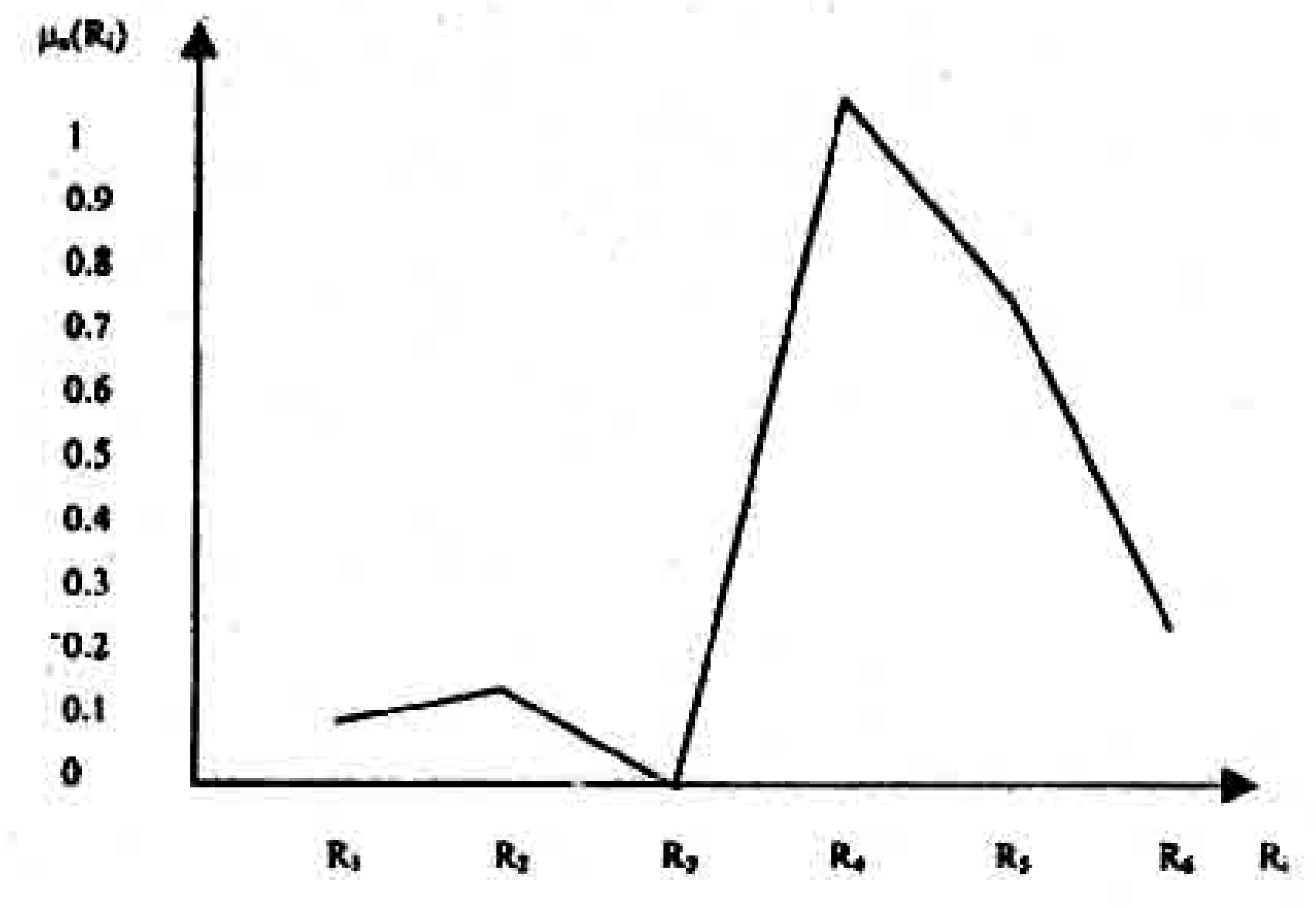

GRAPH 2: GRAPHICAL ILLUSTRATION FOR $\alpha$ = 0.2

Here the fourth row sum $R_4$ is getting the highest membership grade that is the membership grade is 0.7580.

Using the same method described for the $\alpha$ cut values 0.1 and 0.2, we see when $\alpha$ takes the values 0.3, 0.5, 0.6 the corresponding row sum of the membership grades are as follows:

$$\mu_x(R_i)_{I=1,2,\ldots 6} = \begin{cases} 0.649, 0.113,0,1,0.5791,0.312 & \text{when} \quad \alpha = 0.3 \\ 0.710,0.094,0,1,0.481,0.597 & \text{when} \quad \alpha = 0.5 \\ 0.79,0.09,0,1,0.492,0.767 & \text{when} \quad \alpha = 0.6 \end{cases}.$$

Here the fourth row sum $R_4$ is getting the highest membership grade that is the membership grade is 1 and the first row sum $R_1$ is getting the next highest membership grade. Graphical illustration for $\alpha$ = 0.3, 0.5, 0.6 is shown in graph 3, graph 5 and graph 6.

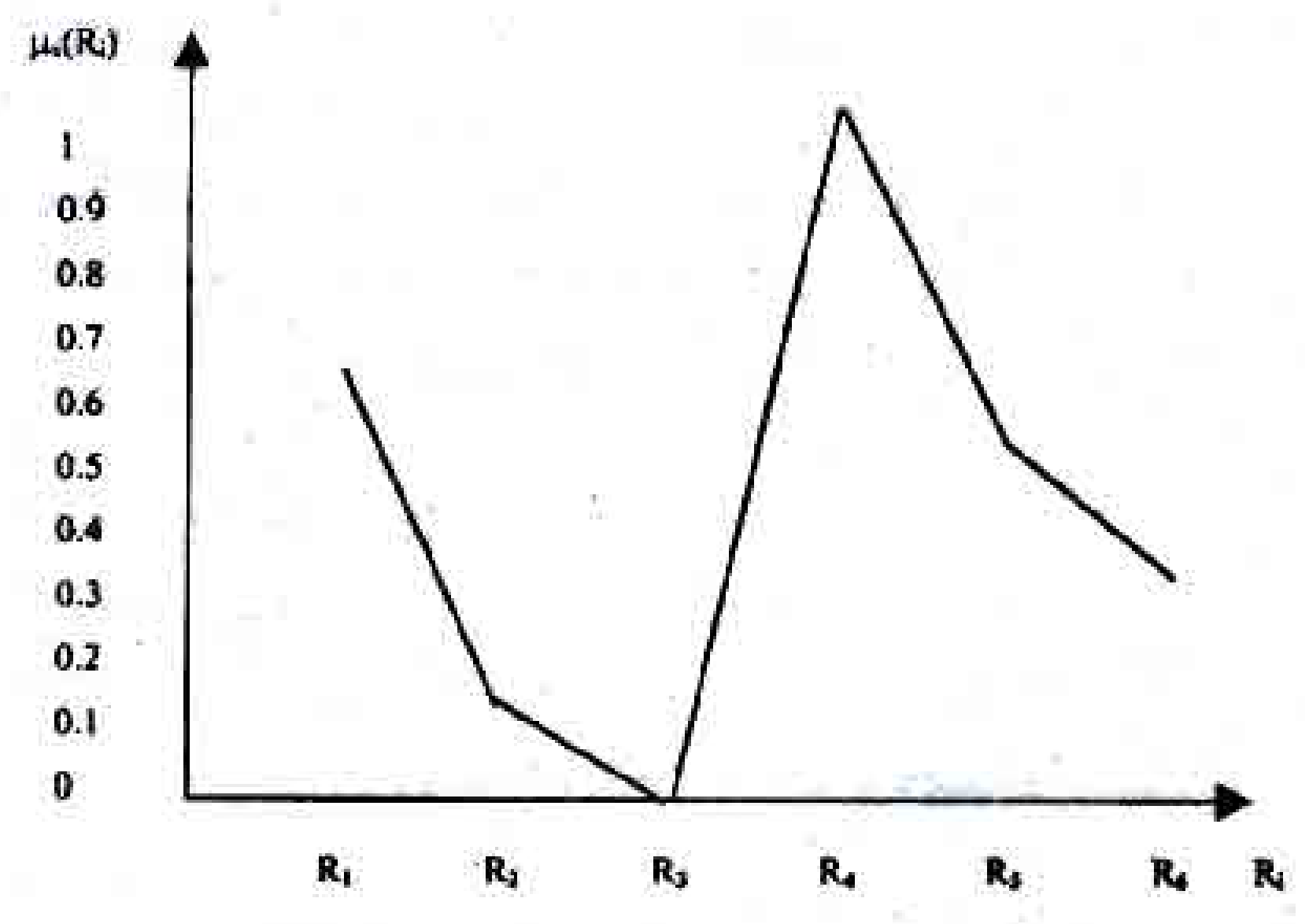

GRAPH 3: GRAPHICAL ILLUSTRATION FOR $\alpha$ = 0.3

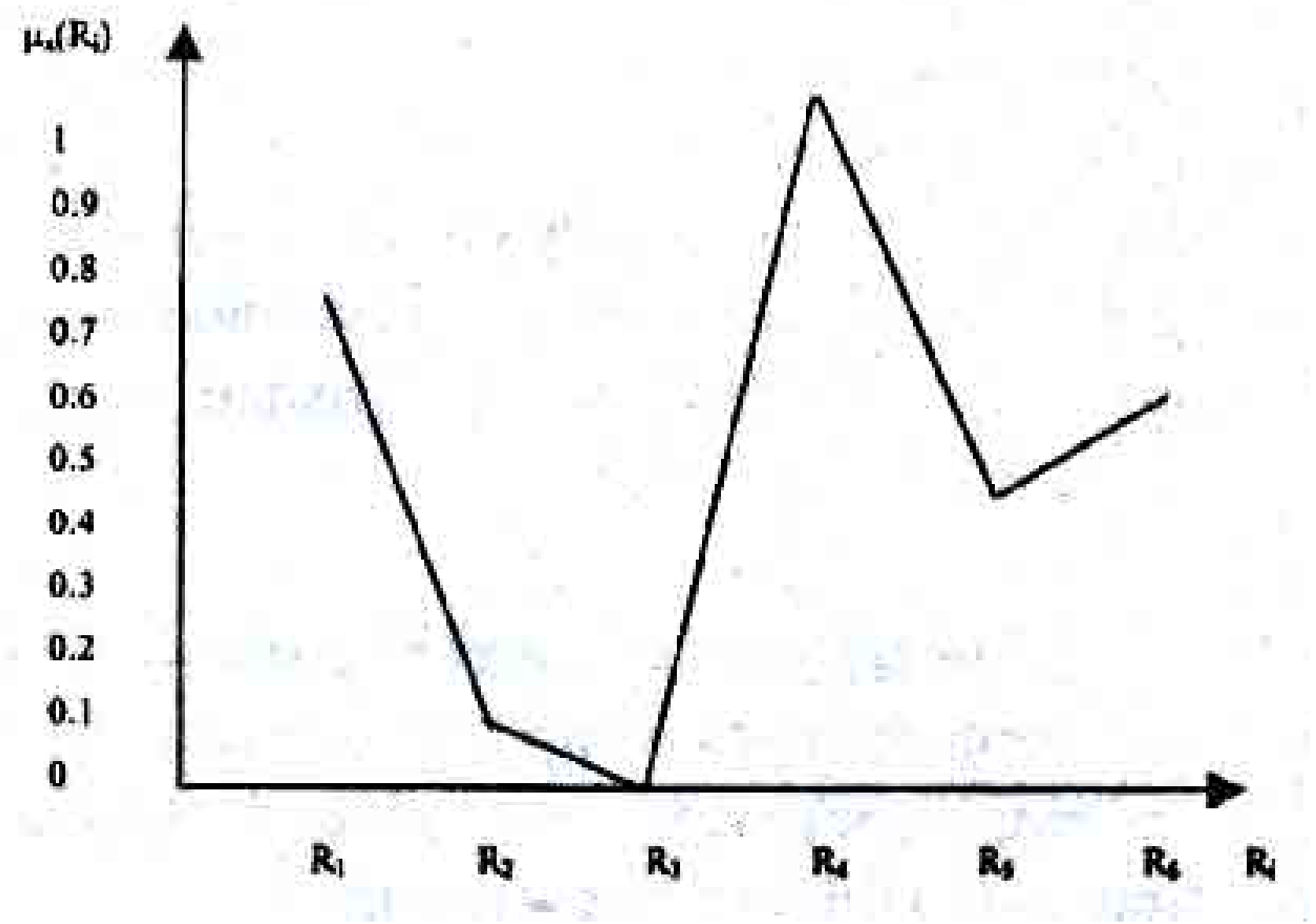

GRAPH 5: GRAPHICAL ILLUSTRATION FOR $\alpha$ = 0.5

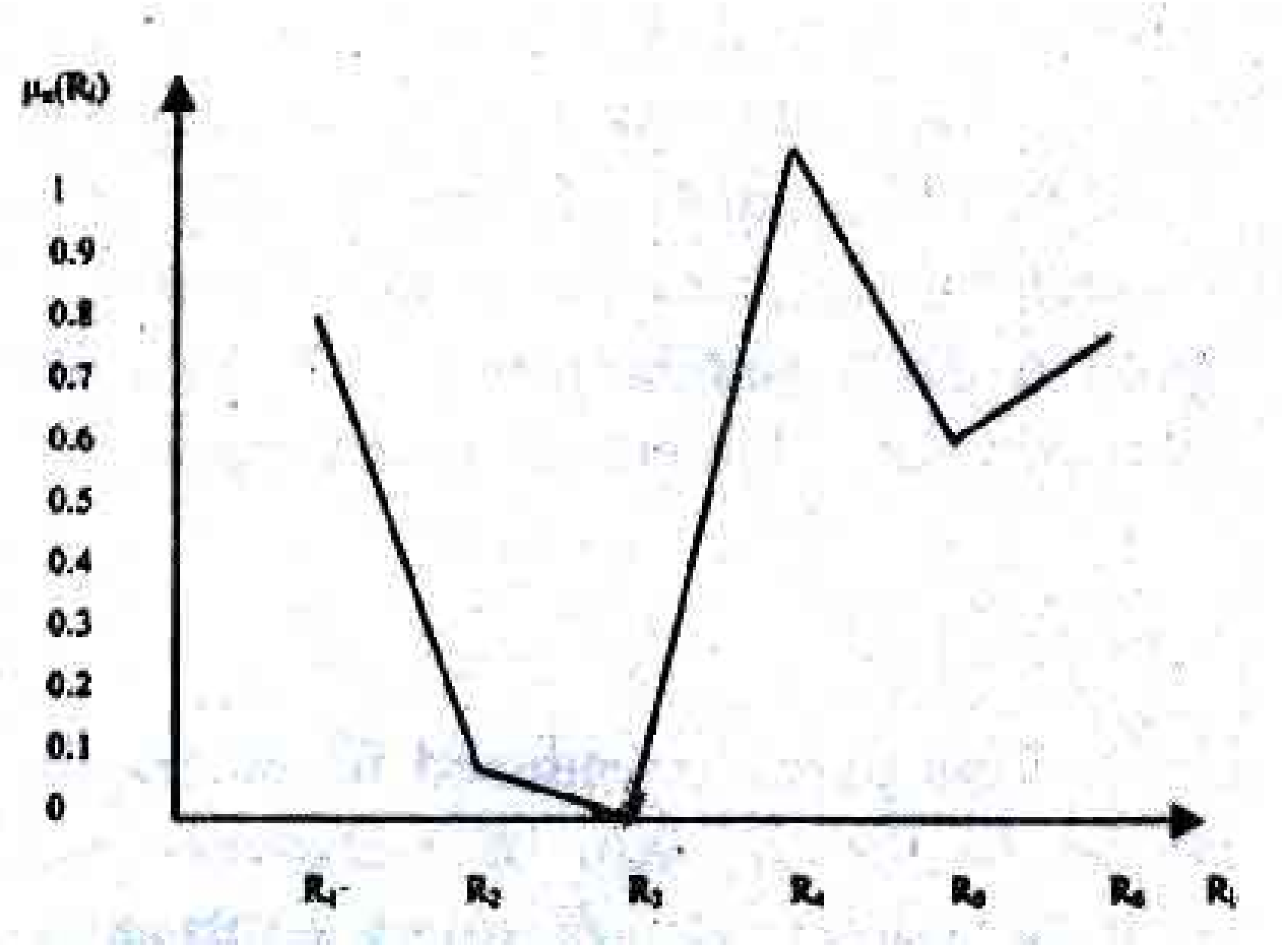

GRAPH 6: GRAPHICAL ILLUSTRATION FOR α = 0.6

When α takes the value 0.4 the corresponding row sum of the membership grades are as follows:

$\mu_x (R_i)_{i=1,2,\ldots 6} = \{0.411, 0.096, 0, 1, 0.493, 0.439$ when $\alpha = 0.4$.

Here the fourth row sum $R_4$ is getting the highest membership grade that is the membership grade is 1 and the fifth row sum $R_5$ is getting the next highest membership grade. Graphical illustration for α = 0.4 is shown in graph 4.

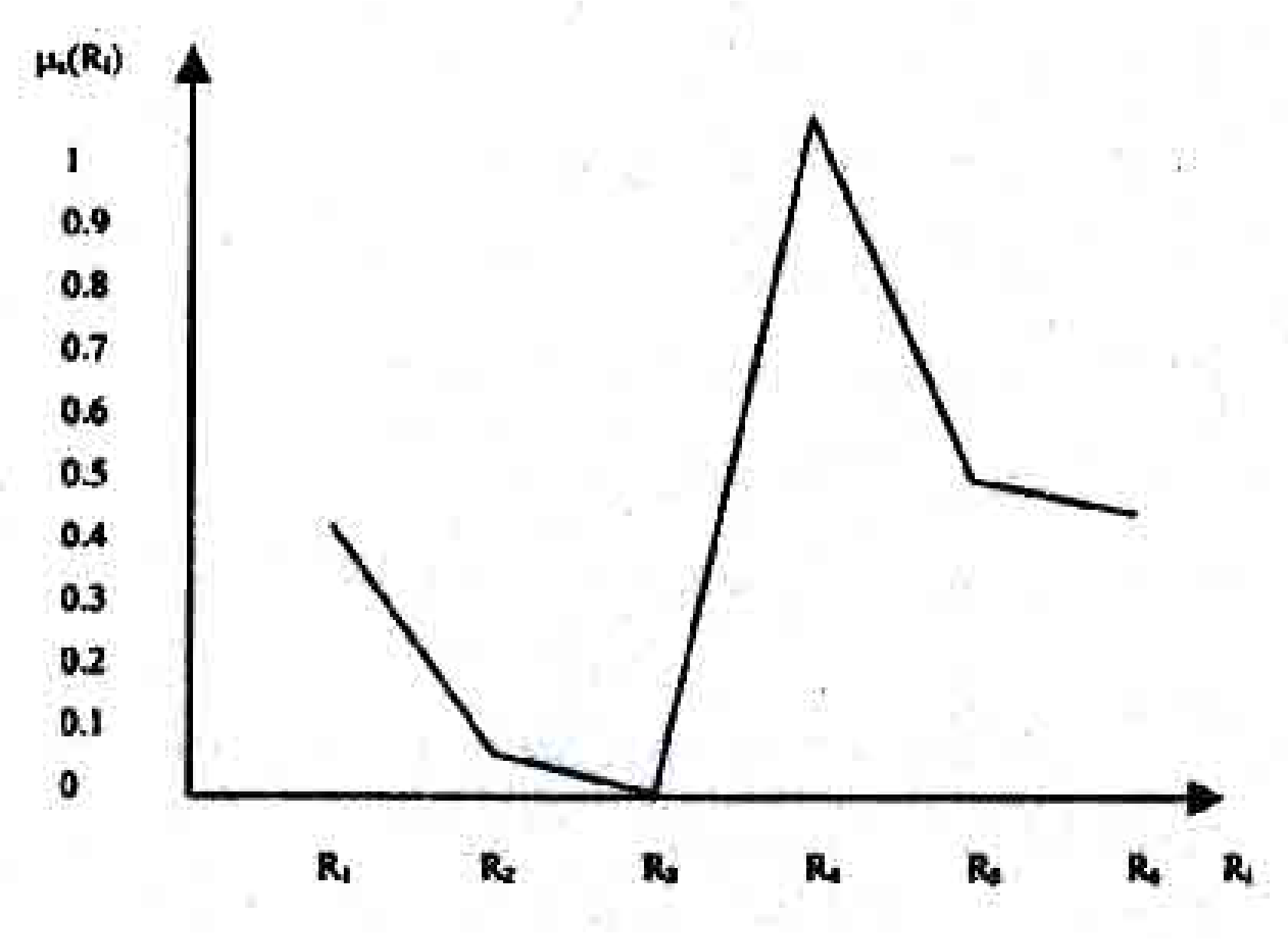

GRAPH 4: GRAPHICAL ILLUSTRATION FOR α = 0.4

When α takes the values 0.7, 0.8, 0.9, 1 the corresponding row sum of the membership grades are as follows:

$$\mu_x\,(R_i)_{i=1,2,\ldots 6} = \begin{cases} 0.521, 0.09, 0, 1, 0.478, 0.9 & \text{when } \alpha = 0.7 \\ 0.878, 0.092, 0, 1, 0.475, 0.994 & \text{when } \alpha = 0.8 \\ 0.892, 0.081, 0, 1, 0.483, 0.989 & \text{when } \alpha = 0.9 \\ 0.91, 0.07, 0, 1, 0.484, 0.99 & \text{when } \alpha = 1 \end{cases}.$$

Here the fourth row sum $R_4$ is getting the highest membership grade that is the membership grade is 1 and the sixth row sum $R_6$ is getting the next highest membership grade. Graphical illustration for $\alpha$ = 0.7, 0.8 is shown in graph 7 and graph 8.

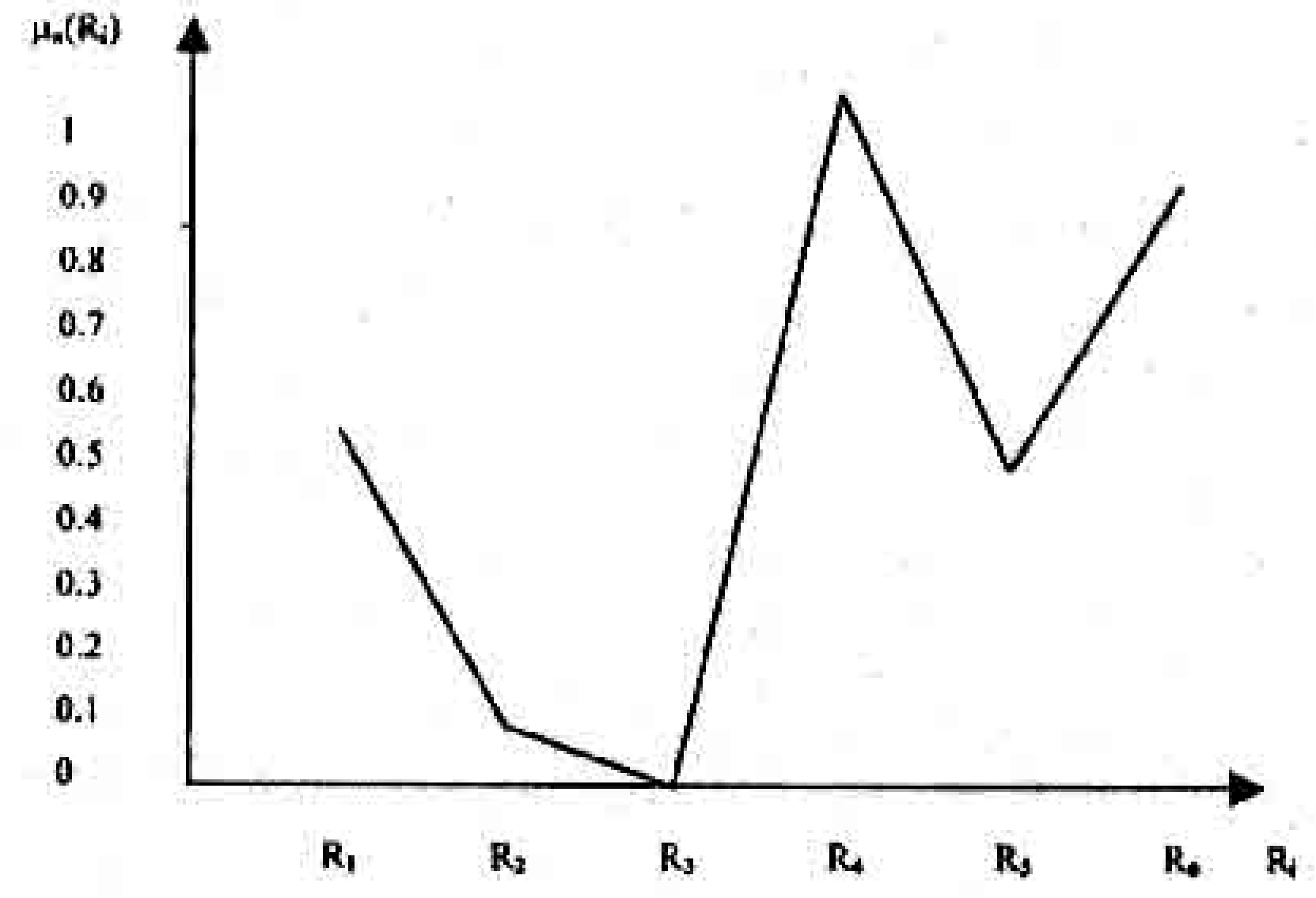

GRAPH 7: GRAPHICAL ILLUSTRATION FOR $\alpha$ = 0.7

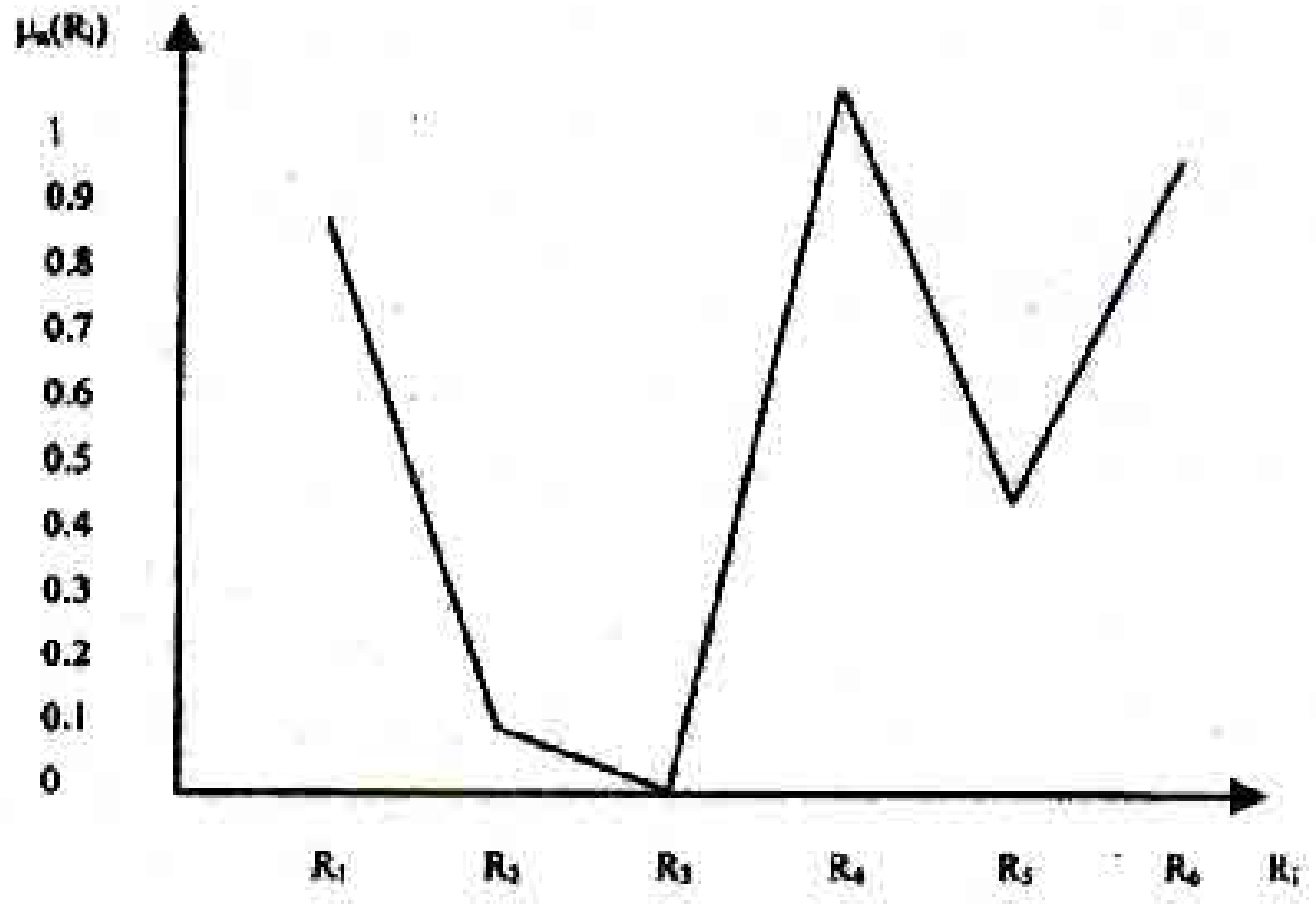

GRAPH 8: GRAPHICAL ILLUSTRATION FOR $\alpha$ = 0.8

Graphical illustration for α = 0.9 and 1 is shown in graph 9 and graph 10 which is given below:

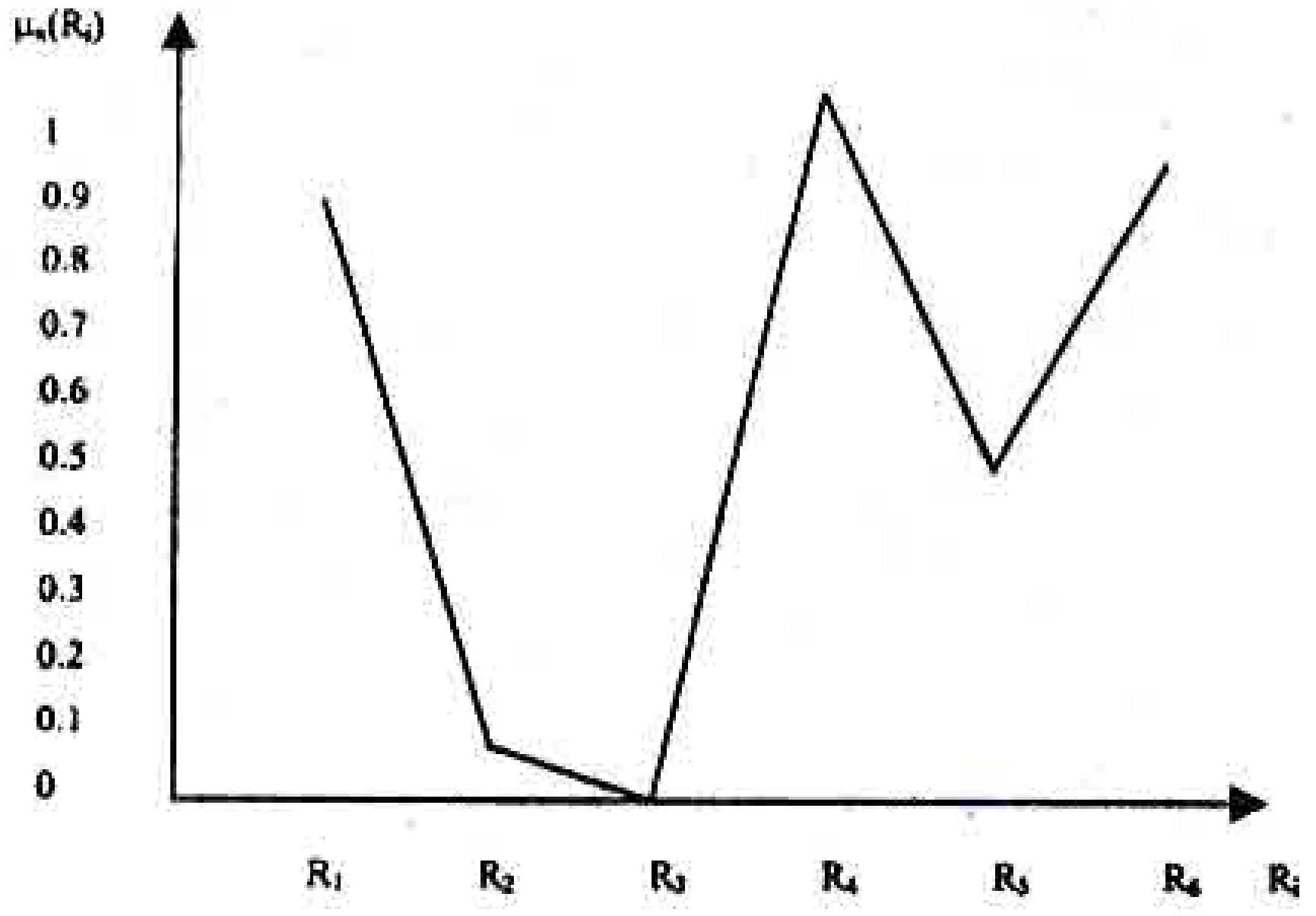

GRAPH 9: GRAPHICAL ILLUSTRATION FOR α = 0.9

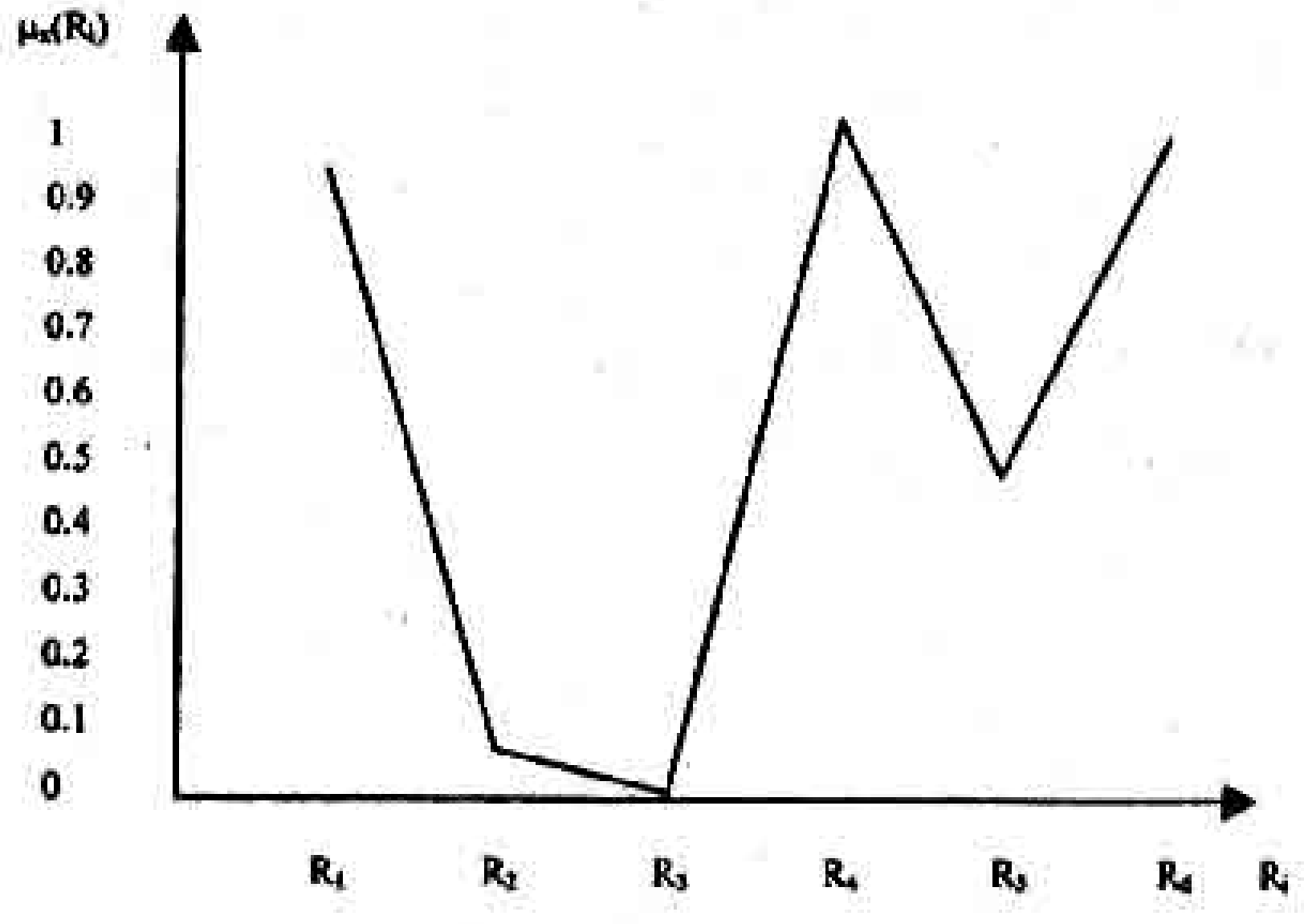

GRAPH 10: GRAPHICAL ILLUSTRATION FOR α =1

The combined fuzzy matrix for all the values of α ∈ [0, 1] is given below:

$(b_{ij})$ =

| | $X_1$ | $X_2$ | $X_3$ | $X_4$ | $X_5$ | $X_6$ | $X_7$ | $X_8$ |
|---|---|---|---|---|---|---|---|---|
| $Y_1$ | 0 | 7.596 | 10 | 0 | 0 | 10 | 3.689 | 10 |
| $Y_2$ | 10 | 10 | 9.702 | 0 | 0 | 0 | 0 | 4.533 |
| $Y_3$ | 10 | 10 | 3.116 | 0 | 9.671 | 0 | 0 | 0 |
| $Y_4$ | 0 | 0 | 7.724 | 10 | 0.399 | 9.99 | 8.531 | 8.819 |
| $Y_5$ | 0 | 0 | 5.713 | 10 | 7.480 | 7.174 | 10 | 0 |
| $Y_6$ | 1.93 | 0 | 0 | 4.345 | 9.98 | 9.98 | 9.98 | 3.413 |

The corresponding row sum of the above matrix is given by

| | | | | |
|---|---|---|---|---|
| $R_1$ | = | The first row sum | = | 41.285 |
| $R_2$ | = | The second row sum | = | 34.235 |
| $R_3$ | = | The third row sum | = | 32.787 |
| $R_4$ | = | The fourth row sum | = | 45.463 |
| $R_5$ | = | The fifth row sum | = | 40.367 |
| $R_6$ | = | The sixth row sum | = | 39.628. |

Using the membership function the corresponding values of $\mu_x$ ($R_i$) is 0.679, 0.144, 0, 1, 0.598, 0.539 respectively where i = 1, 2, 3, 4, 5, 6. Graphical illustration for the values of α ∈ [0,1] is shown in graph 11.

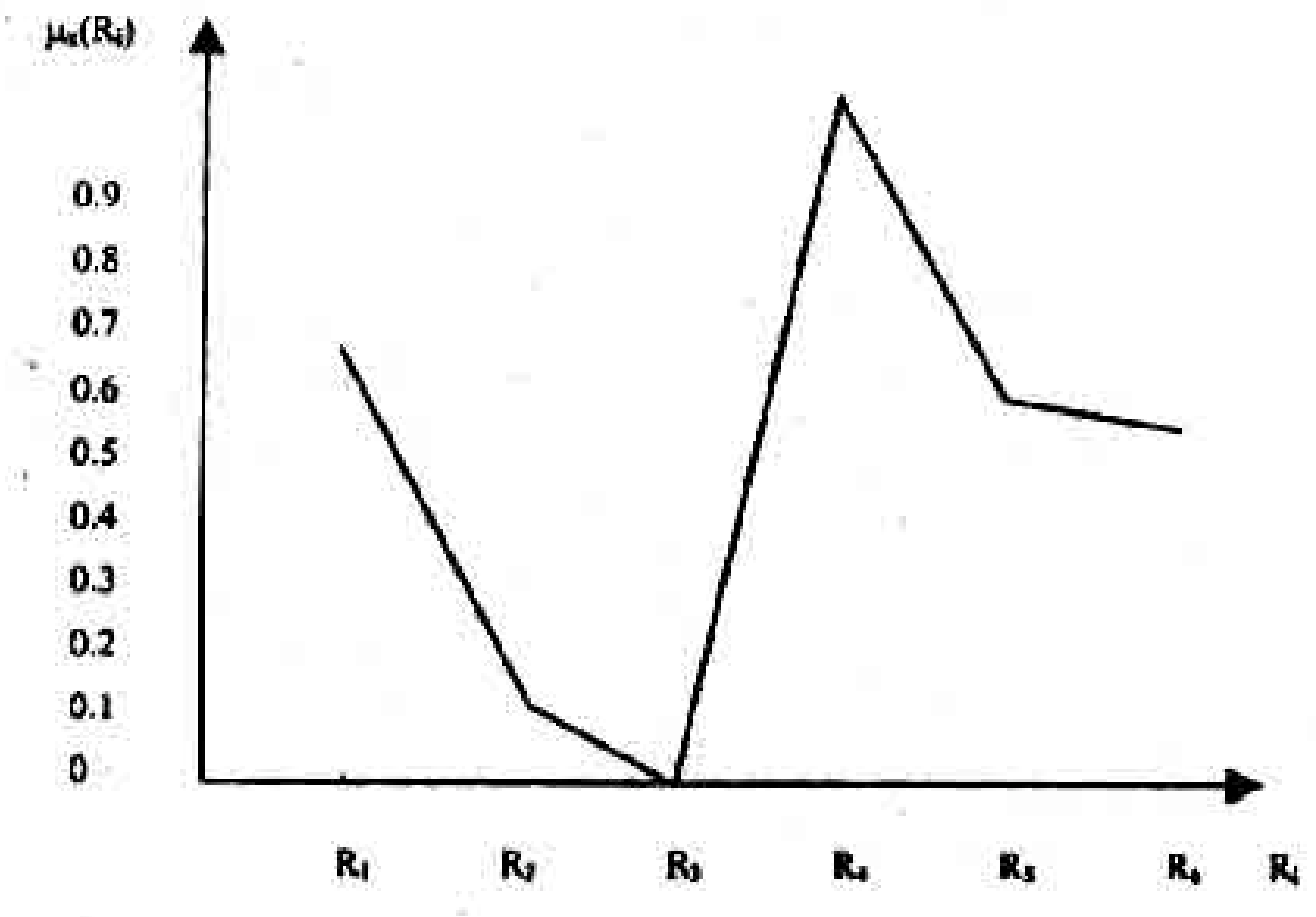

GRAPH 11: GRAPHICAL ILLUSTRATION FOR α ∈ [0, 1]

Here the fourth row sum $R_4$ is getting the highest membership grade that is the membership grade is 1, first row sum $R_1$ is getting the next highest membership grade that is the membership grade is 0.670 and the row sum $R_3$ gets the lowest membership grade which is zero.

## CONCLUSIONS

Here we introduce a fuzzy relational matrix model to the problem of employees and employers in cement industries and this finds best form of relationship between the employees and the employers. The raw data is obtained by the expert's opinion from the cement industrialists as well as feelings of the employees. Using the data from the cement industry, we carry out the analysis of the data via fuzzy relational matrix model and found the best form of relationship between employees and employers. This has been done in two ways.

In the first stage the data was converted into a relational map and the relational matrix was obtained form the relational map, and using this the average relational matrix is obtained.

In the second stage the average relational matrix is converted into fuzzy relational matrix and using different parameter $\alpha$, $\alpha \in [0, 1]$ the best form of relationship is obtained. This is also explicitly described by the graphs. Here we consider eight attributes of the employees and six effects of employers. The production levels of each year are treated as rows and the various attributes of employee are treated as columns in the fuzzy matrix.

The highest membership grade gives the best form of relationship between employers and employees, which maximizes the production level of the cement industries. Using the different parameter $\alpha$ = 0.1, 0.2, 0.3, 0.4, 0.5, 0.6, 0.7, 0.8, 0.9, 1 with the row sum of the fuzzy matrix, the fourth row sum that is the year 1998-1999 gives maximum production and employer gets maximum satisfaction.

From the combined fuzzy matrix, we observe that the row sum $R_4$ gets the highest membership grade that is the value is 1, which is in the year 1998-1999 giving the maximum production with maximum satisfaction of employees and the row sum $R_1$ gets the next highest membership grade 0.670 that is in the year 1995-1996. Now from the analysis of combined fuzzy matrix the row sum $R_3$ gets the lowest membership grade which is zero that is, in the year 1997-1998 the employees satisfaction was poor with minimum production.

For more please refer [127].

Thus we have seen how best the FRMs can be used when we do not have any data but only the opinions and when we have the data using the past experience how best we can give the predictions. Here we give some justifications to state why the use of FRMs are sometimes better than the FCMs.

The first marked difference between FCMs and FRMs is that FCMs cannot directly give the effect of one group on the other. But FRMs can give the effect of one group on the other group and vice versa .

FCMs cannot give any benefit when the nodes or causalities are mutually exclusive ones. But in the case of FRMs since we divide them into two groups and relational maps are sent from one group to other, it gives the maximum benefit.

FRMs give the direct effect of one node from space to other node or nodes of the other space very precisely. We see also in case of FRMs when more than one node is in the on state the hidden pattern ends in a limit cycle and when only one node is on the hidden pattern happens to be a fixed points. Further FRMs depends so much on the experts' opinion some times the resultants happen to be contradictory. So to overcome this problem we in the next chapter introduce the neutrosophic logic of Florentin Smarandache [90-94] which also includes the indeterminacy of relation between nodes.

Another positive point about the FRMs is when the data can be divided disjointly the size of the matrix is considerably and significantly reduced.

For if in FCMs say we have just 12 nodes then we have a 12 × 12 matrix with 144 entries. But if 12 nodes are divided into 7 and 5 we get only a 7 × 5 matrix with 35 entries. Likewise if it is a 8 × 4 matrix we deal only with 32 entries thus a three digit number is reduced to a two digit number.

Thus we now proceed on to define yet a new method of analysis of FRMs which we will be carrying out in the next section.

## 1.7 Linked Fuzzy Relational Maps

In this section we introduce yet another new technique using FRMs which is impossible in case of FCMs. This method is more adaptable in case of data when we are not in a position to inter-relate two systems but we know they are inter-related indirectly. Such study is possible when we adopt linked fuzzy relational maps. First we give the definition of pair wise linked FRMs.

**DEFINITION 1.7.1:** *Let us assume that we are analyzing some concepts which are divided into 3 disjoint units. Suppose we have 3 spaces say P, Q and R we say some m set of nodes in the space P, some n set of nodes in the space Q and r set of nodes in the space R. We can directly find FRMs or directed graphs relating P and Q, FRMs or directed graphs relating Q and R.*

*But we are not in a position to link or get a relation between P and R directly but in fact there exists a hidden link between them which cannot be easily weighted; in such cases we use linked FRMs. Thus pairwise linked FRMs are those FRMs connecting three distinct spaces P, Q and R in such a way that using the pair of FRMs we obtain a FRM relating P and R.*

If $E_1$ is the connection matrix relating P and Q then $E_1$ is a m × n matrix and $E_2$ is the connection matrix relating Q and R which is a n × r matrix. Now $E_1E_2$ is a m × r matrix which is the connection matrix relating P and R. and $E_2^T E_1^T$ matrix relating R and P, when we have such a situation we call it the pair wise linked FRMs. We illustrate this by the following example:

***Example 1.7.1:*** The child labor problem in India is one of the major problems which is given the least attention. Here we analyze the child labor problem. We have three sets of conceptual nodes in three spaces. The spaces under study are G, C and P where G - The concepts attributes associated with government policies preventing / helping child labor. C - Attributes or concepts associated with the children working as laborers and P - Attributes associated with public awareness and support of child labor. The attributes / concepts associated with government policies preventing / helping such child labor are given below:

G – The attributes associated with government.

| | | |
|---|---|---|
| $G_1$ | – | Children don't form vote bank |
| $G_2$ | – | Business men/industrialists who practice are the main source of vote bank |
| $G_3$ | – | Free and compulsory education for children |
| $G_4$ | – | No proper punishment by the government for the practice of child labor. |

C – The attributes associated with the children working as child laborer.

| | | |
|---|---|---|
| $C_1$ | – | Abolition of child labor |
| $C_2$ | – | Uneducated parents |
| $C_3$ | – | School dropouts / never attended any school |
| $C_4$ | – | Social status of child laborers |
| $C_5$ | – | Poverty / source of livelihood |
| $C_6$ | – | Orphans, Runaways, and Parents are beggars, Fathers in prison. |
| $C_7$ | – | Habits like Cinema, smoking drink etc. |

P – The attributes associated with public awareness in support of child labor.

| | | |
|---|---|---|
| $P_1$ | – | Cheap and long hours of labor from children |
| $P_2$ | – | Children as domestic servants |
| $P_3$ | – | Sympathetic public |
| $P_4$ | – | Motivation by teachers to children to pursue education |
| $P_5$ | – | Perpetuating slavery and caste bias. |

Taking the experts opinions we first give the directed graph relating the child labor and the government policies in Figure 1.7.1.

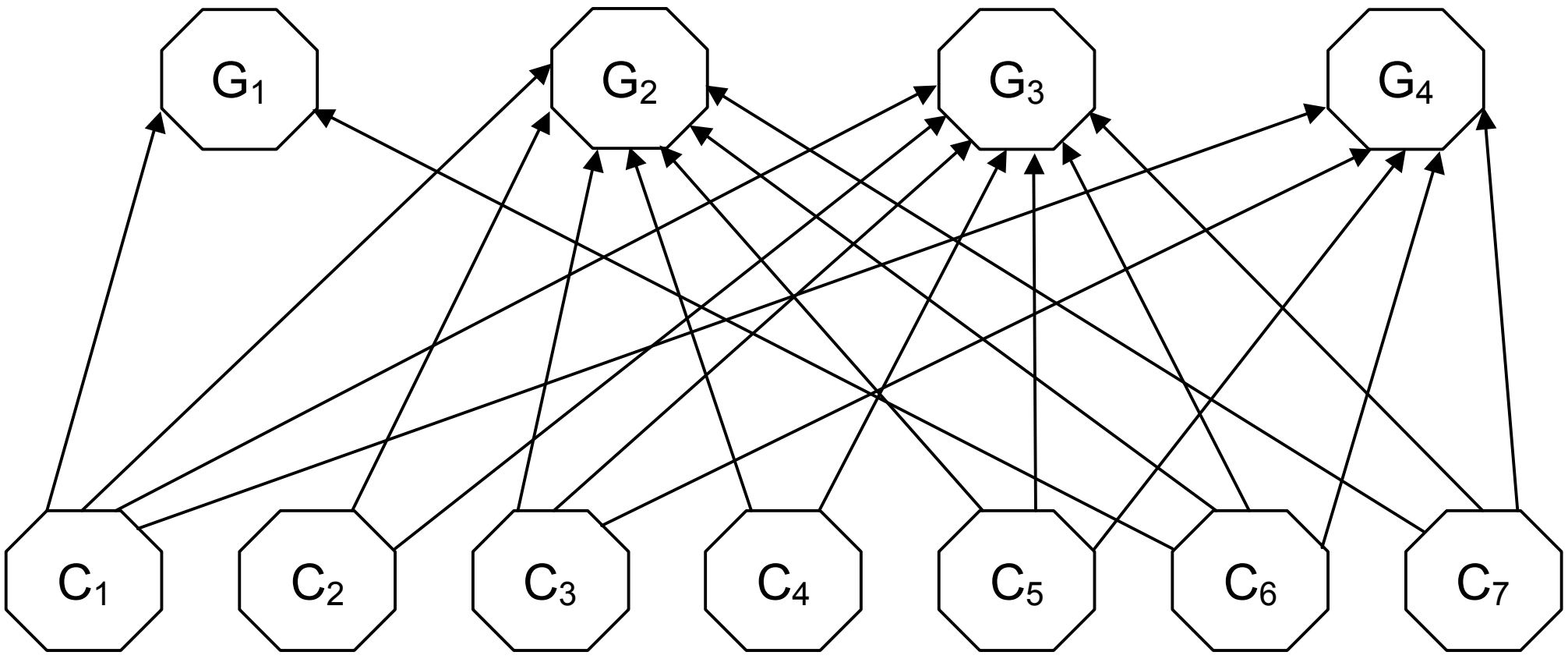

**FIGURE: 1.7.1**

The corresponding relational or connection matrix is as follows:

$$E_1 = \begin{bmatrix} 1 & -1 & 1 & -1 \\ 0 & 1 & -1 & 0 \\ 0 & 1 & -1 & 1 \\ 0 & 1 & -1 & 0 \\ 0 & 1 & 1 & 0 \\ 1 & 1 & -1 & 1 \\ 0 & 1 & -1 & 1 \end{bmatrix}.$$

We do not illustrate the effect of each state vector on the system; here we are more interested in the illustration of how the model interconnects two spaces, which have no direct relation. The same experts opinion is given by the directed graph which is given in the following Figure 1.7.2.:

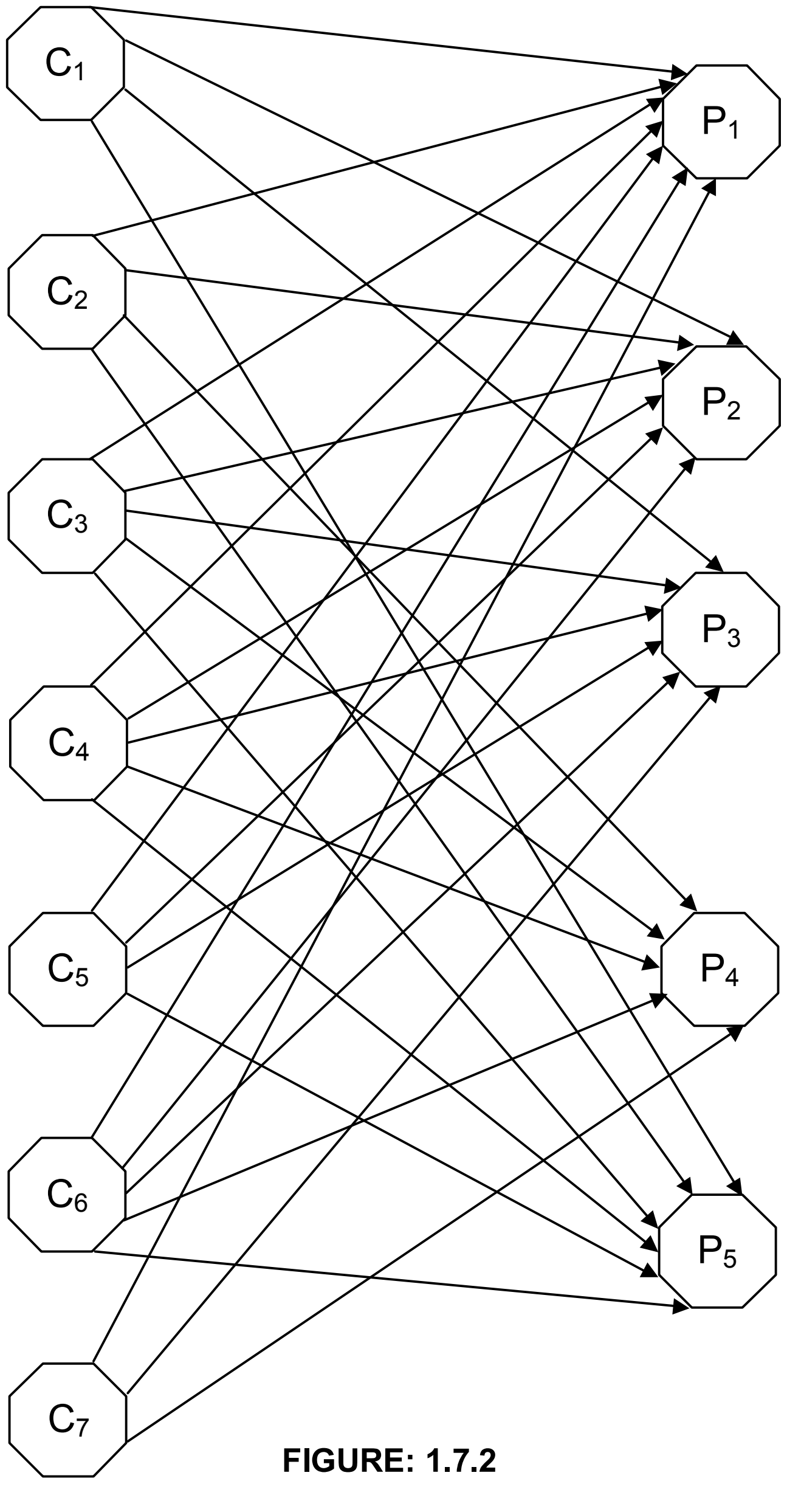

**FIGURE: 1.7.2**

The same experts opinion is sought connecting the following child labor and the role played by the public.

The connection matrix $E_2$ related to the directed graph and the weightages are seen from the matrix.

$$E_2 = \begin{bmatrix} -1 & -1 & 1 & 0 & -1 \\ -1 & 1 & 0 & -1 & 1 \\ 1 & 1 & -1 & -1 & 1 \\ -1 & -1 & 1 & 1 & -1 \\ 1 & 1 & 1 & 0 & 1 \\ 1 & 1 & -1 & -1 & 1 \\ 1 & 0 & -1 & -1 & 0 \end{bmatrix}.$$

Thus $E_1$ is a $7 \times 4$ matrix and $E_2$ is a $7 \times 5$ matrix. Without the aid of the expert we can now obtain the FRM relating the public and the government.

$E_2^T \bullet E_1$ gives the relational matrix which is a $5 \times 4$ matrix say $E_3$

$$E_3 = \begin{bmatrix} 0 & 1 & -1 & 1 \\ 0 & 1 & -1 & 1 \\ 0 & -1 & 1 & -1 \\ -1 & -1 & 1 & 1 \\ 0 & 1 & -1 & 1 \end{bmatrix}.$$

Thus by this method even if we are not in a position to get directed graph using the expert we can indirectly obtain the FRMs relating them. Now using these three FRMs and their related matrices we derive the conclusion by studying the effect of each state vector.

Now we can define 3-linked FRMs, 4-linked FRMs and so on – in general an n-linked FRM on similar lines and obtain the hidden relation which are not explicitly calculable.

Now we proceed on to define the hidden fuzzy relational map and the hidden connection matrix in case of pairwise linked FRMs.

**DEFINITION 1.7.2:** *Let P, Q and R be three spaces with nodes related to some problem. Suppose P and Q are related by a FRM and Q and R are related by an FRM then the indirectly calculated FRM between P and R (or R and P) got as the product of the connection matrices related with the FRMs. P and Q, and Q and R is called as the Hidden connection matrix and the directed graph drawn using the Hidden connection matrix is called the Hidden Directed graph of the pairwise linked FRMs.*

We now proceed on to define on similar lines the three-linked FRMs, four linked FRMs and in general n-linked FRMs.

**DEFINITION 1.7.3:** *Suppose we are analyzing data for which the nodes / concepts are divided into four disjoint classes say A, B, C and D; where A has n concepts, B has m concepts, C has p concepts and D has q nodes or concepts. Suppose the data under study is an unsupervised one and the expert is able to inter-relate, A and B, B and C and C and D through directed graphs that is the FRMs and suppose the expert is not in a clear way to inter-relate the existing relation between A and D then we can get the connection matrix between A and D using the product of the three known matrix resulting in a* $n \times q$ *matrix which will be known as the hidden connection matrix and the directed graph given by this* $n \times q$ *matrix will be known as the hidden directed graph of the 3-linked FRM.*

*On similar lines we can define 4-linked FRMS. So if we are given unsupervised data which can be partitioned into n + 1 sets say* $A_1, \ldots, A_n\ A_{n+1}$ *which will be the spaces having their respective nodes. Then if the expert is able to relate* $A_1A_2$ *,* $A_2A_3$*, ...,* $A_nA_{n+1}$ *then we can inter-relate* $A_1A_{n+1}$ *which will be known as the n-linked FRMs their hidden connection matrix and the hidden directed graph can be obtained. Even at each stage using this method one can always find the hidden relational matrix between any* $A_i$ *and* $A_j$ *and their respective hidden FRMs.*

The main advantage of FRMs is that unlike FCMs we are able to get even hidden FRMs and hidden connection matrix when the expert is not in a position to relate or draw relation between nodes but interrelation between nodes exists in a hidden way.

We give the main theorem relating FCMs and FRMs in the language of graphs.

**THEOREM 1.7.1:** *A FCM can be made into a FRM if and only if the directed graph is a bigraph.*

*Proof:* If the FCMs is made into an FRM then we know the directed graph is a bigraph for we can partition the nodes/concepts into two subsets. Conversely a FRM is a bigraph for the nodes are already made into two disjoint subsets. Hence the claim.

**Chapter Two**

# ON NEUTROSOPHIC COGNITIVE MAPS AND NEUTROSOPHIC RELATIONAL MAPS – PROPERTIES AND APPLICATIONS

This chapter introduces the concept of Neutrosophic Cognitive Maps (NCMs). NCMs are a generalization of Fuzzy Cognitive Maps (FCMs). To study or even to define Neutrosophic Cognitive Maps we have to basically define the notions of Neutrosophic Matrix. Further analogous to Fuzzy Relational Maps (FRMs) we define the notion of Neutrosophic Relational Maps, which are clearly a generalization of FRMs. Study of this type is for the first time experimented in this book. When the data under analysis has indeterminate concepts we are not in a position to give it a mathematical expression. Because of the recent introduction of neutrosophic logic by Florentine Smarandache [90-94] this problem has a solution. Thus we have introduced the additional notion of Neutrosophy in place of Fuzzy theory.

This chapter has ten sections. In the first section we just give a very brief introduction to neutrosophy. In section two, we, for the first time, define the concept of neutrosophic fields, neutrosophic matrices and two types of neutrosophic vector spaces. Section three is completely devoted to the introduction of neutrosophic graphs. The indeterminacy of a path connecting two vertices was never in vogue in mathematical literature. This section gives introduction to basic neutrosophic graphs. The definition and some real world illustrations about neutrosophic cognitive maps is carried out in section four. Section five gives some more illustrations of NCMs. Section six deals with the applications of NCMs. Section seven deals with the comparison of NCMs with FCMs. The notion of Neutrosophic Relational Maps (NRMs) is defined and described in section eight. Section nine gives some of the probable applications of NRM to real world problems. The tenth section defines the linked NRMs, and in this final section we also compare NRMs and FRMs.

## 2.1 An Introduction to Neutrosophy

In this section we introduce the notion of neutrosophic logic created by Florentine Smarandache [90-94], which is an extension / combination of the fuzzy logic in which indeterminacy is included. It has become very essential that the notion of neutrosophic logic play a vital role in several of the real world problems like law, medicine, industry, finance, IT, stocks and share etc. Use of neutrosophic notions will be illustrated/ applied in the later sections of this chapter. Fuzzy theory only measures the grade of membership or the non-existence of a membership in the revolutionary way but fuzzy theory has failed to attribute the concept when the relations between notions or nodes or concepts in problems are indeterminate. In fact one can say the inclusion of the concept of indeterminate situation with fuzzy concepts will form the neutrosophic logic.

As in this book the concept of only fuzzy cognitive maps are dealt which mainly deals with the relation / non-relation between two nodes or concepts but it fails to deal the relation between two conceptual nodes when the relation is an indeterminate one. Neutrosophic logic is the only tool known to us, which deals with the notions of indeterminacy, and here we give a brief description of it. For more about Neutrosophic logic please refer Smarandache [90-94].

**DEFINITION 2.1.1:** *In the neutrosophic logic every logical variable x is described by an ordered triple x = (T, I, F) where T is the degree of truth, F is the degree of false and I the level of indeterminacy.*

(A). To maintain consistency with the classical and fuzzy logics and with probability there is the special case where T + I + F = 1.

(B). But to refer to intuitionistic logic, which means incomplete information on a variable proposition or event one has T + I + F < 1.

(C). Analogically referring to Paraconsistent logic, which means contradictory sources of information about a same logical variable, proposition or event one has T + I + F > 1.

Thus the advantage of using Neutrosophic logic is that this logic distinguishes between relative truth that is a truth is one or a few worlds only noted by 1 and absolute truth denoted by $1^+$. Likewise neutrosophic logic distinguishes between relative falsehood, noted by 0 and absolute falsehood noted by $^-0$.

It has several applications. One such given by [90-94] is as follows:

***Example 2.1.1:*** From a pool of refugees, waiting in a political refugee camp in Turkey to get the American visa, a% have the chance to be accepted – where a varies in the set A, r% to be rejected – where r varies in the set R, and p% to be in pending (not yet decided) – where p varies in P.

Say, for example, that the chance of someone Popescu in the pool to emigrate to USA is (between) 40-60% (considering different criteria of emigration one gets different percentages, we have to take care of all of them), the chance of being rejected is 20-25% or 30-35%, and the chance of being in pending is 10% or 20% or 30%. Then the neutrosophic probability that Popescu emigrates to the Unites States is

NP (Popescu) = ((40-60) (20-25) ∪ (30-35), {10,20,30}), closer to the life.

This is a better approach than the classical probability, where 40 P(Popescu) 60, because from the pending chance – which will be converted to acceptance or rejection – Popescu might get extra percentage in his will to emigrating and also the superior limit of the subsets sum

$$60 + 35 + 30 > 100$$

and in other cases one may have the inferior sum < 0, while in the classical fuzzy set theory the superior sum should be 100 and the inferior sum µ 0. In a similar way, we could say about the element Popescu that Popescu ((40-60), (20-25) ∪ (30-35), {10, 20, 30}) belongs to the set of accepted refugees.

***Example 2.1.2:*** The probability that candidate C will win an election is say 25-30% true (percent of people voting for him), 35% false (percent of people voting against him), and 40% or 41% indeterminate (percent of people not coming to the ballot box, or giving a blank vote – not selecting any one or giving a negative vote cutting all candidate on the list). Dialectic and dualism don't work in this case anymore.

***Example 2.1.3:*** Another example, the probability that tomorrow it will rain is say 50-54% true according to meteorologists who have investigated the past years weather, 30 or 34-35% false according to today's very sunny and droughty summer, and 10 or 20% undecided (indeterminate).

***Example 2.1.4:*** The probability that Yankees will win tomorrow versus Cowboys is 60% true (according to their confrontation's history giving Yankees' satisfaction), 30-32% false (supposing Cowboys are actually up to the mark, while Yankees are declining), and 10 or 11 or 12% indeterminate (left to the hazard: sickness of players, referee's mistakes, atmospheric conditions during the game). These parameters act on players' psychology.

As in this book we use mainly the notion of neutrosophic logic with regard to the indeterminacy of any relation in cognitive maps we are restraining ourselves from dealing with several interesting concepts about neutrosophic logic. As FCMs deals with unsupervised data and the existence or non-existence of cognitive relation, we do not in this book elaborately describe the notion of neutrosophic concepts.

However we just state, suppose in a legal issue the jury or the judge cannot always prove the evidence in a case, in several places we may not be able to derive any conclusions from the existing facts because of which we cannot make a conclusion that no relation exists or otherwise. But existing relation is an indeterminate. So in the case when the concept of indeterminacy exists the judgment ought to be very carefully analyzed be it a civil case or a criminal case. FCMs are deployed only where the existence or non-existence is dealt with but however in our Neutrosophic Cognitive Maps we will deal with the notion of indeterminacy of the evidence also. Thus legal side has lot of Neutrosophic (NCM) applications. Also we will show how NCMs can be used to study factors as varied as stock markets, medical diagnosis, etc.

## 2.2 Some basic Neutrosophic structures

In this section we define some new neutrosophic algebraic structures like neutrosophic fields, neutrosophic spaces and neutrosophic matrices and illustrate them with examples. For these notions are used in the definition of neutrosophic cognitive maps which is dealt in the later sections of this chapter.

Throughout this book by 'I' we denote the indeterminacy of any notion/ concept/ relation. That is when we are not in a position to associate a relation between any two concepts then we denote it as indeterminacy.

Further in this book we assume all fields to be real fields of characteristic 0 all vector spaces are taken as real spaces over reals and we denote the indeterminacy by 'I' as i will make a confusion as i denotes the imaginary value viz $i^2 = -1$ that is $\sqrt{-1} = i$.

**DEFINITION 2.2.1:** *Let K be the field of reals. We call the field generated by K $\cup$ I to be the neutrosophic field for it involves the indeterminacy factor in it. We define $I^2 = I$, $I + I = 2I$ i.e., $I + ... + I = nI$, and if $k \in K$ then $k.I = kI$, $0I = 0$. We denote the neutrosophic field by K(I) which is generated by K $\cup$ I that is $K(I) = \langle K \cup I \rangle$.*

***Example 2.2.1:*** Let R be the field of reals. The neutrosophic field is generated by $\langle R \cup I \rangle$ i.e. R(I) clearly $R \subset \langle R \cup I \rangle$.

***Example 2.2.2:*** Let Q be the field of rationals. The neutrosophic field is generated by Q and I i.e. $Q \cup I$ denoted by Q(I).

**DEFINITION 2.2.2:** *Let K(I) be a neutrosophic field we say K(I) is a prime neutrosophic field if K(I) has no proper subfield which is a neutrosophic field.*

***Example 2.2.3:*** Q(I) is a prime neutrosophic field where as R(I) is not a prime neutrosophic field for $Q(I) \subset R(I)$.

It is very important to note that all neutrosophic fields are of characteristic zero. Likewise we can define neutrosophic subfield.

**DEFINITION 2.2.3:** *Let K(I) be a neutrosophic field, $P \subset K(I)$ is a neutrosophic subfield of P if P itself is a neutrosophic field. K(I) will also be called as the extension neutrosophic field of the neutrosophic field P.*

Now we proceed on to define neutrosophic vector spaces, which are only defined over neutrosophic fields. We can define two types of neutrosophic vector spaces one when it is a neutrosophic vector space over ordinary field other being neutrosophic vector space over neutrosophic fields. To this end we have to define neutrosophic group under addition.

**DEFINITION 2.2.4:** *We know Z is the abelian group under addition. Z(I) denote the additive abelian group generated by the set Z and I, Z(I) is called the neutrosophic abelian group under '+'.*

Thus to define basically a neutrosophic group under addition we need a group under addition. So we proceed on to define neutrosophic abelian group under addition. *Suppose G is an additive abelian group under '+'. $G(I) = \langle G \cup I \rangle$, additive group generated by G and I, G(I) is called the neutrosophic abelian group under '+'.*

***Example 2.2.4:*** Let Q be the group under '+'; $Q(I) = \langle Q \cup I \rangle$ is the neutrosophic abelian group under addition; '+'.

***Example 2.2.5:*** R be the additive group of reals, $R(I) = \langle R \cup I \rangle$ is the neutrosophic group under addition.

***Example 2.2.6:*** $M_{n \times m}(I) = \{(a_{ij}) \mid a_{ij} \in Z(I)\}$ be the collection of all $n \times m$ matrices under '+' $M_{n \times m}(I)$ is a neutrosophic group under '+'.

Now we proceed on to define neutrosophic subgroup.

**DEFINITION 2.2.5:** *Let G(I) be the neutrosophic group under addition. P ⊂ G(I) be a proper subset of G(I). P is said to be neutrosophic subgroup of G(I) if P itself is a neutrosophic group i.e. P = ⟨$P_1$ ∪ I⟩ where $P_1$ is an additive subgroup of G.*

***Example 2.2.7:*** Let Z(I) = ⟨Z ∪ I⟩ be a neutrosophic group under '+'. ⟨2Z ∪ I⟩ = 2Z(I) is the neutrosophic subgroup of Z (I).

In fact Z(I) has several neutrosophic subgroups.

Now we proceed on to define the notion of neutrosophic quotient group.

**DEFINITION 2.2.6:** *Let G (I) = ⟨G ∪ I⟩ be a neutrosophic group under '+', suppose P (I) be a neutrosophic subgroup of G (I) then the neutrosophic quotient group*

$$\frac{G(I)}{P(I)} = \{a + P(I) \mid a \in G\,(I)\}.$$

***Example 2.2.8:*** Let Z (I) be a neutrosophic group under addition, Z the group of integers under addition P = 2Z(I) is a neutrosophic subgroup of Z(I) the neutrosophic subgroup of Z(I), the neutrosophic quotient group

$$\frac{Z(I)}{P} = \{a + 2Z(I) \mid a \in Z(I)\} = \{(2n+1) + (2n+1)\,I \mid n \in Z\}.$$

Clearly $\frac{Z(I)}{P}$ is a group. For P = 2Z (I) serves as the additive identity. Take a, b ∈ $\frac{Z(I)}{P}$. If a, b ∈ Z(I) \ P then two possibilities occur.

a + b is odd times I or a + b is odd or a + b is even times I or even if a + b is even or even times I then a + b ∈ P. if a + b is odd or odd times I a + b ∈ $\frac{Z(I)}{P = 2Z(I)}$.
It is easily verified that P acts as the identity and every element in

$$a + 2Z\,(I) \in \frac{Z(I)}{2Z(I)}$$

has inverse. Hence the claim.

Now we proceed on to define the notion of neutrosophic vector spaces over fields and then we define neutrosophic vector spaces over neutrosophic fields.

**DEFINITION 2.2.7:** *Let G(I) by an additive abelian neutrosophic group. K any field. If G(I) is a vector space over K then we call G(I) a neutrosophic vector space over K.*

Now we give the notion of strong neutrosophic vector space.

**DEFINITION 2.2.8:** *Let G(I) be a neutrosophic abelian group. K(I) be a neutrosophic field. If G(I) is a vector space over K(I) then we call G(I) the strong neutrosophic vector space.*

**THEOREM 2.2.1:** *All strong neutrosophic vector space over K(I) are a neutrosophic vector space over K; as $K \subset K(I)$.*

*Proof:* Follows directly by the very definitions.

Thus when we speak of neutrosophic spaces we mean either a neutrosophic vector space over K or a strong neutrosophic vector space over the neutrosophic field K(I). By basis we mean a linearly independent set which spans the neutrosophic space.

Now we illustrate with an example.

***Example 2.2.9:*** Let R(I) × R(I) = V be an additive abelian neutrosophic group over the neutrosophic field R(I). Clearly V is a strong neutrosophic vector space over R(I). The basis of V are {(0,1), (1,0)}.

***Example 2.2.10:*** Let V = R(I) × R(I) be a neutrosophic abelian group under addition. V is a neutrosophic vector space over R. The neutrosophic basis of V are {(1,0), (0,1), (I,0), (0,I)}, which is a basis of the vector space V over R.

A study of these basis and its relations happens to be an interesting form of research.

**DEFINITION 2.2.9:** *Let G(I) be a neutrosophic vector space over the field K. The number of elements in the neutrosophic basis is called the neutrosophic dimension of G(I).*

**DEFINITION 2.2.10:** *Let G(I) be a strong neutrosophic vector space over the neutrosophic field K(I). The number of elements in the strong neutrosophic basis is called the strong neutrosophic dimension of G(I).*

We denote the neutrosophic dimension of G(I) over K by $N_k$ (dim) of G (I) and that the strong neutrosophic dimension of G (I) by $SN_{K(I)}$ (dim) of G(I).

Now we define the notion of neutrosophic matrices.

**DEFINITION 2.2.11:** *Let $M_{nxm} = \{(a_{ij}) \mid a_{ij} \in K(I)\}$, where K (I), is a neutrosophic field. We call $M_{nxm}$ to be the neutrosophic matrix.*

***Example 2.2.11:*** Let Q(I) = ⟨Q ∪ I⟩ be the neutrosophic field.

$$M_{4x3} = \begin{pmatrix} 0 & 1 & I \\ -2 & 4I & 0 \\ 1 & -I & 2 \\ 3I & 1 & 0 \end{pmatrix}$$

is the neutrosophic matrix, with entries from rationals and the indeterminacy I. We define product of two neutrosophic matrices whenever the production is defined as follows:

Let

$$A = \begin{pmatrix} -1 & 2 & -I \\ 3 & I & 0 \end{pmatrix}_{2\times 3} \quad \text{and} \quad B = \begin{pmatrix} -I & 1 & 2 & 4 \\ 1 & I & 0 & 2 \\ 5 & -2 & 3I & -I \end{pmatrix}_{3x4}$$

$$AB = \begin{bmatrix} -6I+2 & -1+4I & -2-3I & I \\ -2I & 3+I & 6 & 12+2I \end{bmatrix}_{2x4}$$

(we use the fact $I^2 = I$).

To define Neutrosophic Cognitive Maps we direly need the notion of Neutrosophic Matrices. We use square neutrosophic matrices for Neutrosophic Cognitive Maps and use rectangular neutrosophic matrices for Neutrosophic Relational Maps.

## 2.3 Some Basic Notions about Neutrosophic Graphs

In this section we for the first time introduce the notion of neutrosophic graphs, illustrate them and give some basic properties. We need the notion of neutrosophic graphs basically to obtain neutrosophic cognitive maps which will be nothing but directed neutrosophic graphs. Similarly neutrosophic relational maps will also be directed neutrosophic graphs.

It is no coincidence that graph theory has been independently discovered many times since it may quite properly be regarded as an area of applied mathematics. The subject finds its place in the work of Euler. Subsequent rediscoveries of graph theory were by Kirchhoff and Cayley. Euler (1707-1782) became the father of graph theory as well as topology when in 1936 he settled a famous unsolved problem in his day called the Konigsberg Bridge Problem.

Psychologist Lewin proposed in 1936 that the life space of an individual be represented by a planar map. His view point led the psychologists at the Research center for Group Dynamics to another psychological interpretation of a graph in which people are represented by points and interpersonal relations by lines. Such relations include love, hate, communication and power. In fact it was precisely this approach which led the author to a personal discovery of graph theory, aided and abetted by psychologists L. Festinger and D. Cartwright. Here it is pertinent to mention that the directed graphs of an FCMs or FRMs are nothing but the psychological inter-relations or feelings of different nodes; but it is unfortunate that in all these studies the concept of indeterminacy was never given any place, so in this chapter for the first time we will be having graphs in which the notion of indeterminacy i.e. when two vertex should be connected or not is never dealt with. If graphs are to display human feelings then this point is very vital for in several situations certain relations between concepts may certainly remain an indeterminate.

So this section will purely cater to the properties of such graphs the edges of certain vertices may not find its connection i.e., they are indeterminates, which we will be defining as neutrosophic graphs.

The world of theoretical physics discovered graph theory for its own purposes. In the study of statistical mechanics by Uhlenbeck the points stands for molecules and two adjacent points indicate nearest neighbor interaction of some physical kind, for example magnetic attraction or repulsion. But it is forgotten in all these situations we may have molecule structures which need not attract or repel but remain without action or not able to predict the action for such analysis we can certainly adopt the concept of neutrosophic graphs.

In a similar interpretation by Lee and Yang the points stand for small cubes in Euclidean space where each cube may or may not be occupied by a molecule. Then two points are adjacent whenever both spaces are occupied.

Feynmann proposed the diagram in which the points represent physical particles and the lines represent paths of the particles after collisions. Just at each stage of applying graph theory we may now feel the neutrosophic graph theory may be more suitable for application.

Now we proceed on to define the neutrosophic graph.

**DEFINITION 2.3.1:** *A neutrosophic graph is a graph in which atleast one edge is an indeterminacy denoted by dotted lines.*

**NOTATION***:* The indeterminacy of an edge between two vertices will always be denoted by dotted lines.

***Example 2.3.1:*** The following are neutrosophic graphs:

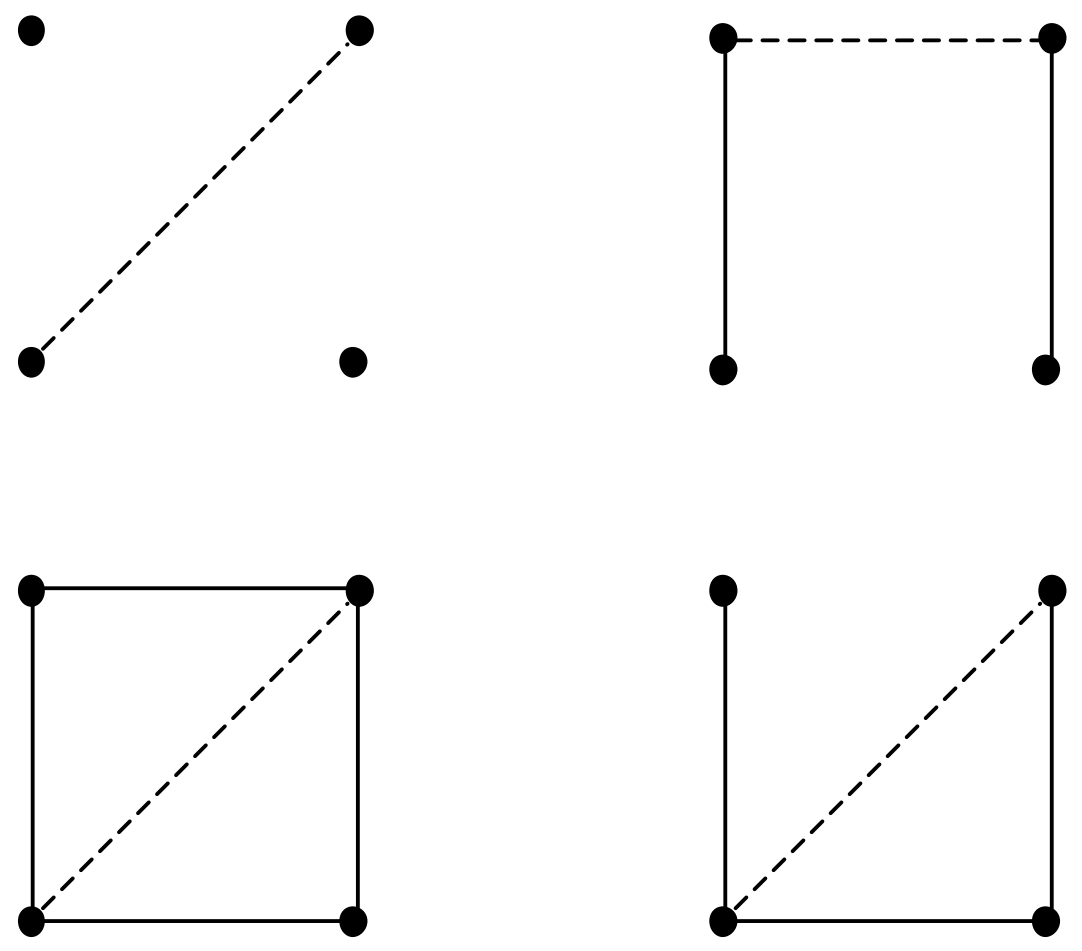

FIGURE:

All graphs in general are not neutrosophic graphs.

***Example 2.3.2:*** The following graphs are not neutrosophic graphs given in Figure 2.3.2:

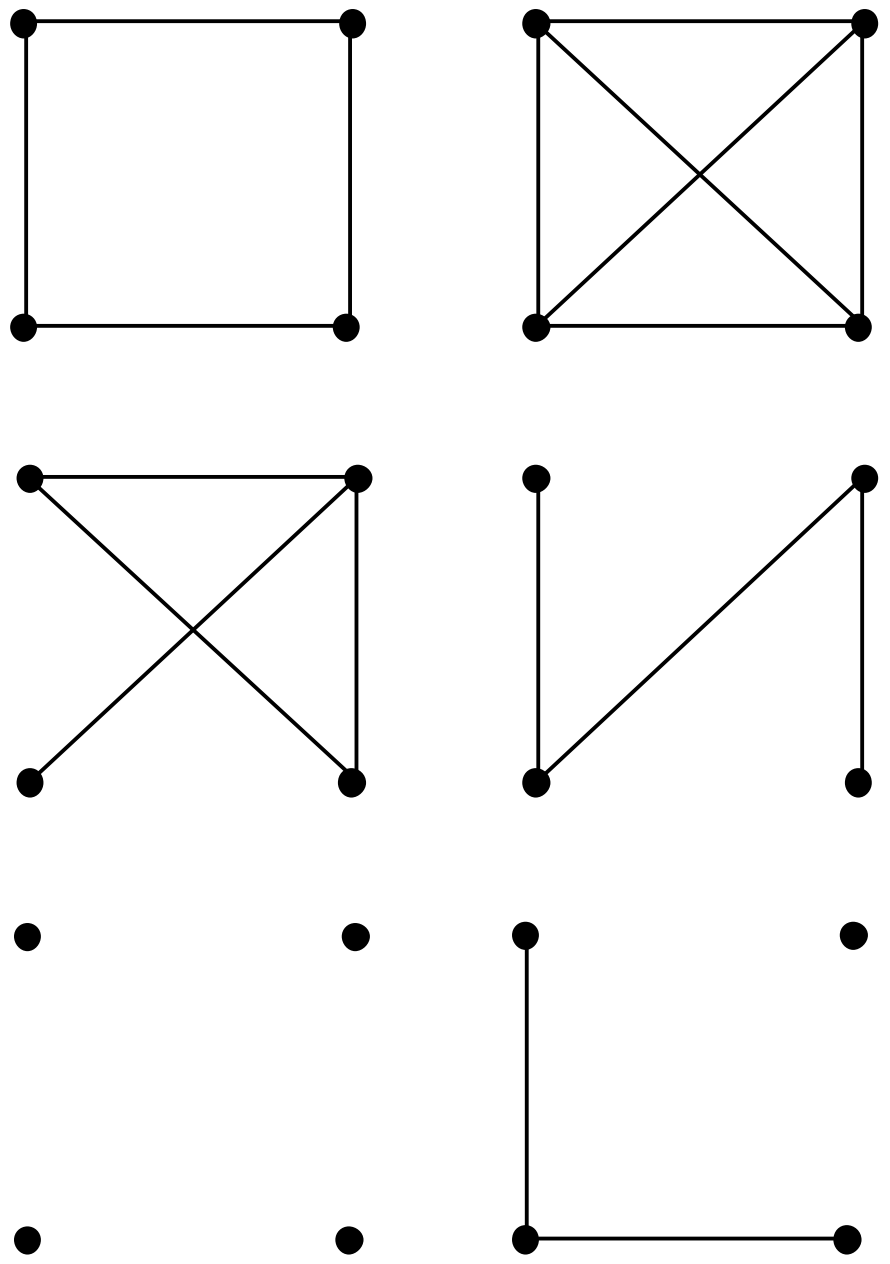

FIGURE:

**DEFINITION 2.3.2:** *A neutrosophic directed graph is a directed graph which has atleast one edge to be an indeterminacy.*

**DEFINITION 2.3.3:** *A neutrosophic oriented graph is a neutrosophic directed graph having no symmetric pair of dircted indeterminacy lines.*

**DEFINITION 2.3.4:** *A neutrosophic subgraph H of a neutrosophic graph G is a subgraph H which is itself a neutrosophic graph.*

**THEOREM 2.3.1:** *Let G be a neutrosophic graph. All subgraphs of G are not neutrosophic subgraphs of G.*

*Proof:* By an example. Consider the neutrosophic graph given in Figure 2.3.3.

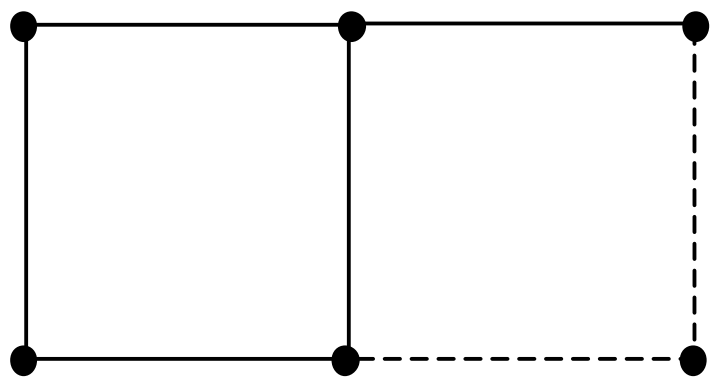

**FIGURE: 2.3.3**

This has a subgraph given by Figure 2.3.4.

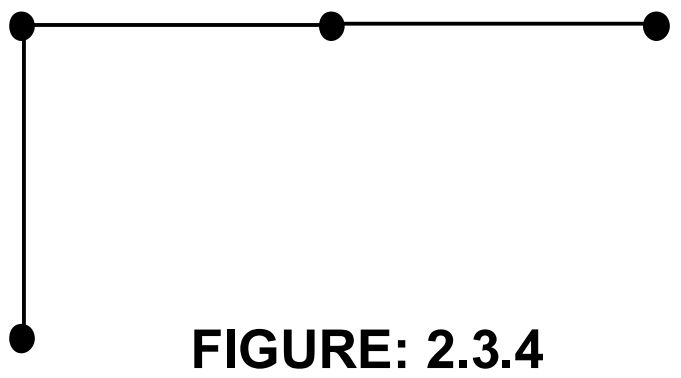

**FIGURE: 2.3.4**

which is not a neutrosophic subgraph of G.

**THEOREM 2.3.2:** *Let G be a neutrosophic graph. In general the removal of a point from G need be a neutrosophic subgraph.*

*Proof:* Consider the graph G given in Figure 2.3.5.

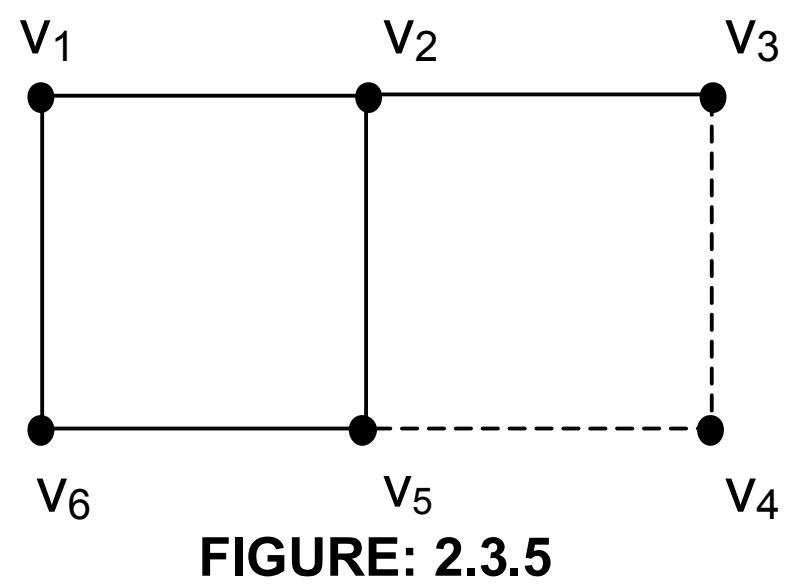

**FIGURE: 2.3.5**

G \ $v_4$ is only a subgraph of G but is not a neutrosophic subgraph of G.

Thus it is interesting to note that this a main feature by which a graph differs from a neutrosophic graph.

**DEFINITION 2.3.5:** *Two graphs G and H are neutrosophically isomorphic if*

i. *They are isomorphic*
ii. *If there exists a one to one correspondence between their point sets which preserve indeterminacy adjacency.*

**DEFINITION 2.3.6:** *A neutrosophic walk of a neutrosophic graph G is a walk of the graph G in which at least one of the lines is an indeterminacy line. The neutrosophic walk is neutrosophic closed if $v_0 = v_n$ and is neutrosophic open otherwise.*

*It is a neutrosophic trial if all the lines are distinct and atleast one of the lines is a indeterminacy line and a path if all points are distinct (i.e. this necessarily means all lines are distinct and atleast one line is a line of indeterminacy). If the neutrosophic walk is neutrosophic closed then it is a neutrosophic cycle provided its n points are distinct and $n \geq 3$. A neutrosophic graph is neutrosophic connected if it is connected and atleast a pair of points are joined by a path. A neutrosophic maximal connected neutrosophic subgraph of G is called a neutrosophic connected component or simple neutrosophic component of G. Thus a neutrosophic graph has atleast two neutrosophic components then it is neutrosophic disconnected. Even if one is a component and another is a neutrosophic component still we do not say the graph is neutrosophic disconnected.*

Neutrosophic disconnected graphs are given in Figure 2.3.6.

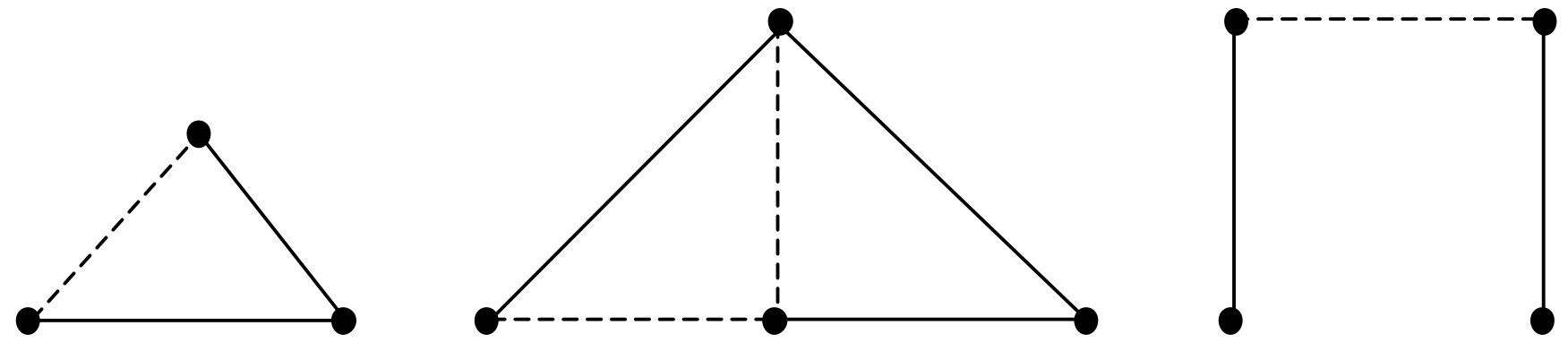

**FIGURE: 2.3.6**

Graph which is not neutrosophic disconnected is given by Figure 2.3.7.

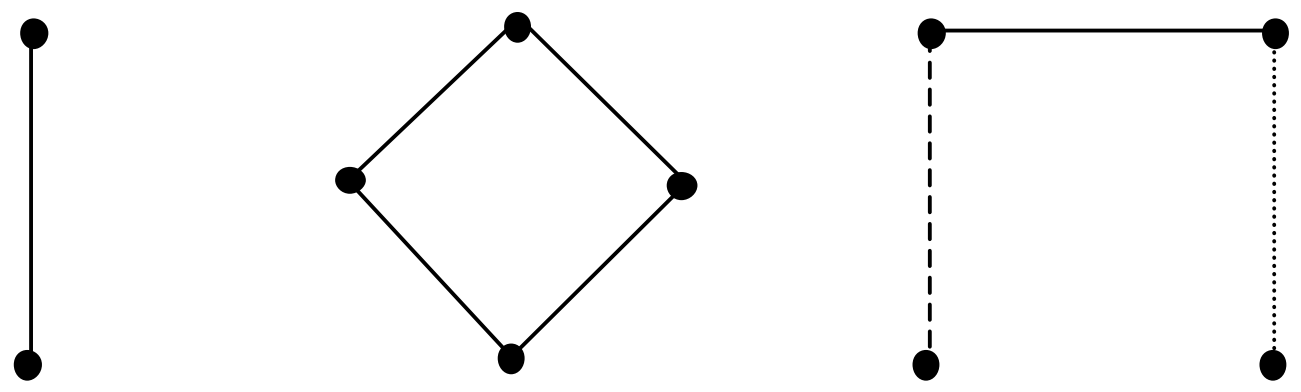

**FIGURE: 2.3.7**

Several results in this direction can be defined and analysed.

**DEFINITION 2.3.7:** *A neutrosophic bigraph, G is a bigraph, G whose point set V can be partitioned into two subsets* $V_1$ *and* $V_2$ *such that at least a line of G which joins* $V_1$ *with* $V_2$ *is a line of indeterminacy.*

*This neutrosophic bigraphs will certainly play a role in the study of FRMs and FCMs and in fact give a method of conversion of data from FRMs to FCMs.*

As both the models FRMs and FCMs work on the adjacency or the connection matrix we just define the neutrosophic adjacency matrix related to a neutrosophic graph G given by Figure 2.3.8.

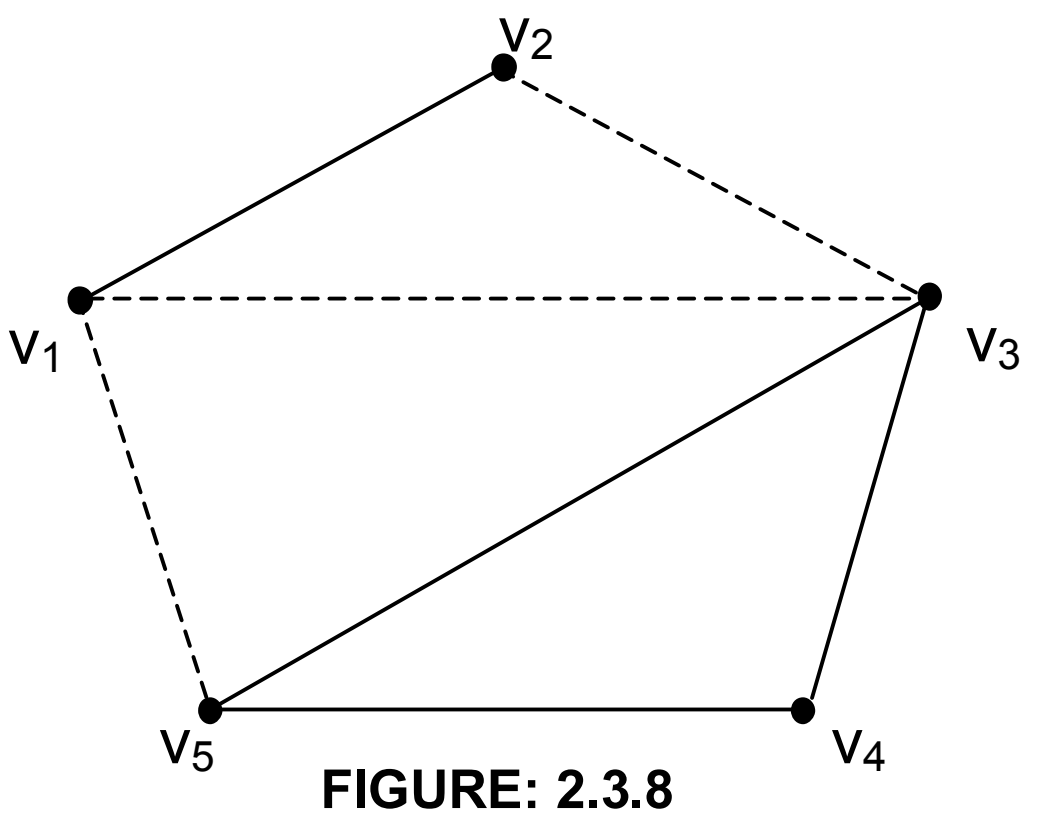

**FIGURE: 2.3.8**

*The neutrosophic adjacency matrix is N(A)*

$$N(A) = \begin{bmatrix} 0 & 1 & I & 0 & I \\ I & 0 & I & 0 & 0 \\ I & I & 0 & 1 & 1 \\ 0 & 0 & 1 & 0 & 1 \\ I & 0 & 1 & 1 & 0 \end{bmatrix}.$$

*Its entries will not only be 0 and 1 but also the indeterminacy I.*

**DEFINITION 2.3.8:** *Let G be a neutrosophic graph. The adjacency matrix of G with entries from the set (I, 0, 1) is called the neutrosophic adjacency matrix of the graph.*

Now as our main aim is the study of Neutrosophic Cognitive Maps we do not divulge into a very deep study of Neutrosophic Graphs or its properties but have given only the basic and the appropriate notions which are essential for studying of this book.

## 2.4 On Neutrosophic Cognitive Maps with Examples

The notion of Fuzzy Cognitive Maps (FCMs) which are fuzzy signed directed graphs with feedback are discussed and described in Chapter 1. The directed edge $e_{ij}$ from causal concept $C_i$ to concept $C_j$ measures how much $C_i$ causes $C_j$. The time varying concept function $C_i(t)$ measures the non negative occurrence of some fuzzy event, perhaps the strength of a political sentiment, historical trend or opinion about some topics like child labor or school dropouts etc. FCMs model the world as a collection of classes and causal relations between them.

The edge $e_{ij}$ takes values in the fuzzy causal interval $[-1, 1]$ ($e_{ij} = 0$ indicates no causality, $e_{ij} > 0$ indicates causal increase; that $C_j$ increases as $C_i$ increases and $C_j$ decreases as $C_i$ decreases, $e_{ij} < 0$ indicates causal decrease or negative causality $C_j$ decreases as $C_i$ increases or $C_j$, increases as $C_i$ decreases. Simple FCMs have edge value in {-1, 0, 1}. Thus if causality occurs it occurs to maximal positive or negative degree.

It is important to note that $e_{ij}$ measures only absence or presence of influence of the node $C_i$ on $C_j$ but till now any researcher has not contemplated the indeterminacy of any relation between two nodes $C_i$ and $C_j$. When we deal with unsupervised data, there are situations when no relation can be determined between some two nodes. So in this section we try to introduce the indeterminacy in FCMs, and we choose to call this generalized structure as Neutrosophic Cognitive Maps (NCMs). In our view this will certainly give a more appropriate result and also caution the user about the risk of indeterminacy.

Now we proceed on to define the concepts about NCMs.

**DEFINITION 2.4.1:** *A Neutrosophic Cognitive Map (NCM) is a neutrosophic directed graph with concepts like policies, events etc. as nodes and causalities or indeterminates as edges. It represents the causal relationship between concepts.*

Let $C_1, C_2, \ldots, C_n$ denote n nodes, further we assume each node is a neutrosophic vector from neutrosophic vector space V. So a node $C_i$ will be represented by $(x_1, \ldots, x_n)$ where $x_k$'s are zero or one or I (I is the indeterminate introduced in Sections 2.2 and 2.3 of the chapter 2) and $x_k = 1$ means that the node $C_k$ is in the on state and $x_k = 0$ means the node is in the off state and $x_k = I$ means the nodes state is an indeterminate at that time or in that situation.

Let $C_i$ and $C_j$ denote the two nodes of the NCM. The directed edge from $C_i$ to $C_j$ denotes the causality of $C_i$ on $C_j$ called connections. Every edge in the NCM is weighted with a number in the set {-1, 0, 1, I}. Let $e_{ij}$ be the weight of the directed edge $C_iC_j$, $e_{ij} \in \{-1, 0, 1, I\}$. $e_{ij} = 0$ if $C_i$ does not have any effect on $C_j$, $e_{ij} = 1$ if increase (or decrease) in $C_i$ causes increase (or decreases) in $C_j$, $e_{ij} = -1$ if increase (or decrease) in $C_i$ causes decrease (or increase) in $C_j$ . $e_{ij}$ = I if the relation or effect of $C_i$ on $C_j$ is an indeterminate.

**DEFINITION 2.4.2:** *NCMs with edge weight from {-1, 0, 1, I} are called simple NCMs.*

**DEFINITION 2.4.3:** *Let $C_1$, $C_2$, ..., $C_n$ be nodes of a NCM. Let the neutrosophic matrix N(E) be defined as N(E) = ($e_{ij}$) where $e_{ij}$ is the weight of the directed edge $C_i$ $C_j$, where $e_{ij} \in \{0, 1, -1, I\}$. N(E) is called the neutrosophic adjacency matrix of the NCM.*

**DEFINITION 2.4.4:** *Let $C_1$, $C_2$, ..., $C_n$ be the nodes of the NCM. Let A = ($a_1$, $a_2$,..., $a_n$) where $a_i \in \{0, 1, I\}$. A is called the instantaneous state neutrosophic vector and it denotes the on – off – indeterminate state position of the node at an instant*

| | | |
|---|---|---|
| $a_i$ | = | *0 if $a_i$ is off (no effect)* |
| $a_i$ | = | *1 if $a_i$ is on (has effect)* |
| $a_i$ | = | *I if $a_i$ is indeterminate(effect cannot be determined)* |

*for i = 1, 2,..., n.*

**DEFINITION 2.4.5:** *Let $C_1$, $C_2$, ..., $C_n$ be the nodes of the FCM. Let $\overrightarrow{C_1C_2}$, $\overrightarrow{C_2C_3}$, $\overrightarrow{C_3C_4}$, … , $\overrightarrow{C_iC_j}$ be the edges of the NCM. Then the edges form a directed cycle. An NCM is said to be cyclic if it possesses a directed cyclic. An NCM is said to be acyclic if it does not possess any directed cycle.*

**DEFINITION 2.4.6:** *An NCM with cycles is said to have a feedback. When there is a feedback in the NCM i.e. when the causal relations flow through a cycle in a revolutionary manner the NCM is called a dynamical system.*

**DEFINITION 2.4.7:** *Let $\overrightarrow{C_1C_2}, \overrightarrow{C_2C_3}, \cdots, \overrightarrow{C_{n-1}C_n}$ be cycle, when $C_i$ is switched on and if the causality flow through the edges of a cycle and if it again causes $C_i$, we say that the dynamical system goes round and round. This is true for any node $C_i$, for i = 1, 2,..., n. the equilibrium state for this dynamical system is called the hidden pattern.*

**DEFINITION 2.4.8:** *If the equilibrium state of a dynamical system is a unique state vector, then it is called a fixed point. Consider the NCM with $C_1$, $C_2$,..., $C_n$ as nodes. For example let us start the dynamical system by switching on $C_1$. Let us assume that the NCM settles down with $C_1$ and $C_n$ on, i.e. the state vector remain as (1, 0,..., 1) this neutrosophic state vector (1,0,..., 0, 1) is called the fixed point.*

**DEFINITION 2.4.9:** *If the NCM settles with a neutrosophic state vector repeating in the form*

$$A_1 \rightarrow A_2 \rightarrow \ldots \rightarrow A_i \rightarrow A_1,$$

*then this equilibrium is called a limit cycle of the NCM.*

**METHODS OF DETERMINING THE HIDDEN PATTERN:**

Let $C_1$, $C_2$,…, $C_n$ be the nodes of an NCM, with feedback. Let E be the associated adjacency matrix. Let us find the hidden pattern when $C_1$ is switched on when an input is given as the vector $A_1$ = (1, 0, 0,…, 0), the data should pass through the neutrosophic matrix N(E), this is done by multiplying $A_1$ by the matrix N(E). Let $A_1N(E) = (a_1, a_2,…, a_n)$ with the threshold operation that is by replacing $a_i$ by 1 if $a_i >$ k and $a_i$ by 0 if $a_i < k$ (k – a suitable positive integer) and $a_i$ by I if $a_i$ is not a integer. We update the resulting concept, the concept $C_1$ is included in the updated vector by making the first coordinate as 1 in the resulting vector. Suppose $A_1N(E) \to A_2$ then consider $A_2N(E)$ and repeat the same procedure. This procedure is repeated till we get a limit cycle or a fixed point.

**DEFINITION 2.4.10:** *Finite number of NCMs can be combined together to produce the joint effect of all NCMs. If $N(E_1)$, $N(E_2)$,…, $N(E_p)$ be the neutrosophic adjacency matrices of a NCM with nodes $C_1$, $C_2$,…, $C_n$ then the combined NCM is got by adding all the neutrosophic adjacency matrices $N(E_1)$,…, $N(E_p)$. We denote the combined NCMs adjacency neutrosophic matrix by $N(E) = N(E_1) + N(E_2) + … + N(E_p)$.*

**NOTATION:** Let $(a_1, a_2, … , a_n)$ and $(a'_1, a'_2, … , a'_n)$ be two neutrosophic vectors. We say $(a_1, a_2, … , a_n)$ is equivalent to $(a'_1, a'_2, … , a'_n)$ denoted by $((a_1, a_2, … , a_n) \sim (a'_1, a'_2, …, a'_n)$ if $(a'_1, a'_2, … , a'_n)$ is got after thresholding and updating the vector $(a_1, a_2, … , a_n)$ after passing through the neutrosophic adjacency matrix N(E).

The following are very important:

***Note 1:*** The nodes $C_1$, $C_2$, …, $C_n$ are nodes are not indeterminate nodes because they indicate the concepts which are well known. But the edges connecting $C_i$ and $C_j$ may be indeterminate i.e. an expert may not be in a position to say that $C_i$ has some causality on $C_j$ either will he be in a position to state that $C_i$ has no relation with $C_j$ in such cases the relation between $C_i$ and $C_j$ which is indeterminate is denoted by I.

***Note 2:*** The nodes when sent will have only ones and zeros i.e. on and off states, but after the state vector passes through the neutrosophic adjacency matrix the resultant vector will have entries from {0, 1, I} i.e. they can be neutrosophic vectors.

The presence of I in any of the coordinate implies the expert cannot say the presence of that node i.e. on state of it after passing through N(E) nor can we say the absence of the node i.e. off state of it the effect on the node after passing through the dynamical system is indeterminate so only it is represented by I. Thus only in case of NCMs we can say the effect of any node on other nodes can also be indeterminates. Such possibilities and analysis is totally absent in the case of FCMs.

***Note 3:*** In the neutrosophic matrix N(E), the presence of I in the $a_{ij}$ the place shows, that the causality between the two nodes i.e. the effect of $C_i$ on $C_j$ is indeterminate. Such chances of being indeterminate is very possible in case of unsupervised data and that too in the study of FCMs which are derived from the directed graphs.

Thus only NCMs helps in such analysis.

Now we shall represent a few examples to show how in this set up NCMs is preferred to FCMs. At the outset before we proceed to give examples we make it clear that all unsupervised data need not have NCMs to be applied to it. Only data which have the relation between two nodes to be an indeterminate need to be modeled with NCMs if the data has no indeterminacy factor between any pair of nodes one need not go for NCMs; FCMs will do the best job.

*Example 2.4.1:* **The child labor problem prevalent in India is modeled in this example using NCMs.**

**Let us consider the child labor problem with the following conceptual nodes**

| | | |
|---|---|---|
| $C_1$ | - | Child Labor |
| $C_2$ | - | Political Leaders |
| $C_3$ | - | Good Teachers |
| $C_4$ | - | Poverty |
| $C_5$ | - | Industrialists |
| $C_6$ | - | Public practicing/encouraging Child Labor |
| $C_7$ | - | Good Non-Governmental Organizations (NGOs) |

$C_1$ - Child labor, it includes all types of labor of children below 14 years which include domestic workers, rag pickers, working in restaurants / hotels, bars etc. (It can be part time or fulltime).

$C_2$ - We include political leaders with the following motivation: Children are not vote banks so political leaders are not directly concerned with child labor but they indirectly help in the flourishing of it as industrialists who utilize child laborers or cheap labor are the decision makers for the winning or losing of the political leaders. Also industrialists financially control political interests. So we are forced to include political leaders as a node in this problem.

$C_3$ - Teachers are taken as a node because mainly school dropouts or children who have never attended the school are child laborers. So if the motivation by the teacher is very good, there would be less school dropouts and therefore there would be a decrease in child laborers.

$C_4$ - Poverty which is the most responsible reason for child labor.

$C_5$ - Industrialists – when we say industrialists we include one and all starting from a match factory or beedi factory, bars, hotels etc.

$C_6$ - Public who promote child labor as domestic servants, sweepers etc.

$C_7$ - We qualify the NGOs as good for some NGOs may not take up the issue fearing the rich and the powerful. Here "good NGOs" means NGOs who try to stop or prevent child labor.

Now we give the directed graph as well as the neutrosophic graph of two experts in the following Figures 2.4.1 and 2.4.2:

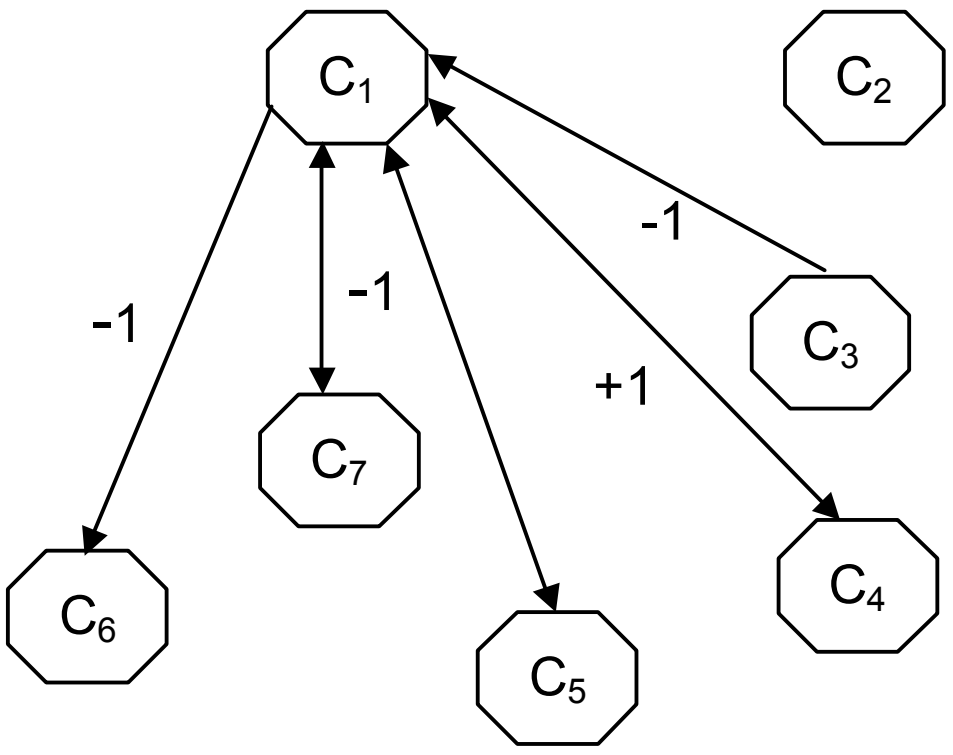

FIGURE:

Figure 2.4.1 gives the directed graph with $C_1$, $C_2$, …, $C_7$ as nodes and Figure 2.4.2 gives the neutrosophic directed graph with the same nodes.

The connection matrix E related to the graph in Figure 2.4.1. is given below:

$$E = \begin{bmatrix} 0 & 0 & 0 & 1 & 1 & 1 & -1 \\ 0 & 0 & 0 & 0 & 0 & 0 & 0 \\ -1 & 0 & 0 & 0 & 0 & 0 & 0 \\ 1 & 0 & 0 & 0 & 0 & 0 & 0 \\ 1 & 0 & 0 & 0 & 0 & 0 & 0 \\ 0 & 0 & 0 & 0 & 0 & 0 & 0 \\ -1 & 0 & 0 & 0 & 0 & 0 & 0 \end{bmatrix}.$$

According to this expert no connection however exists between political leaders and industrialists.

Now we reformulate a different format of the questionnaire where we permit the expert to give answers like the relation between certain nodes is indeterminable or not known. Now based on the expert's opinion also about the notion of indeterminacy we obtain the following neutrosophic directed graph:

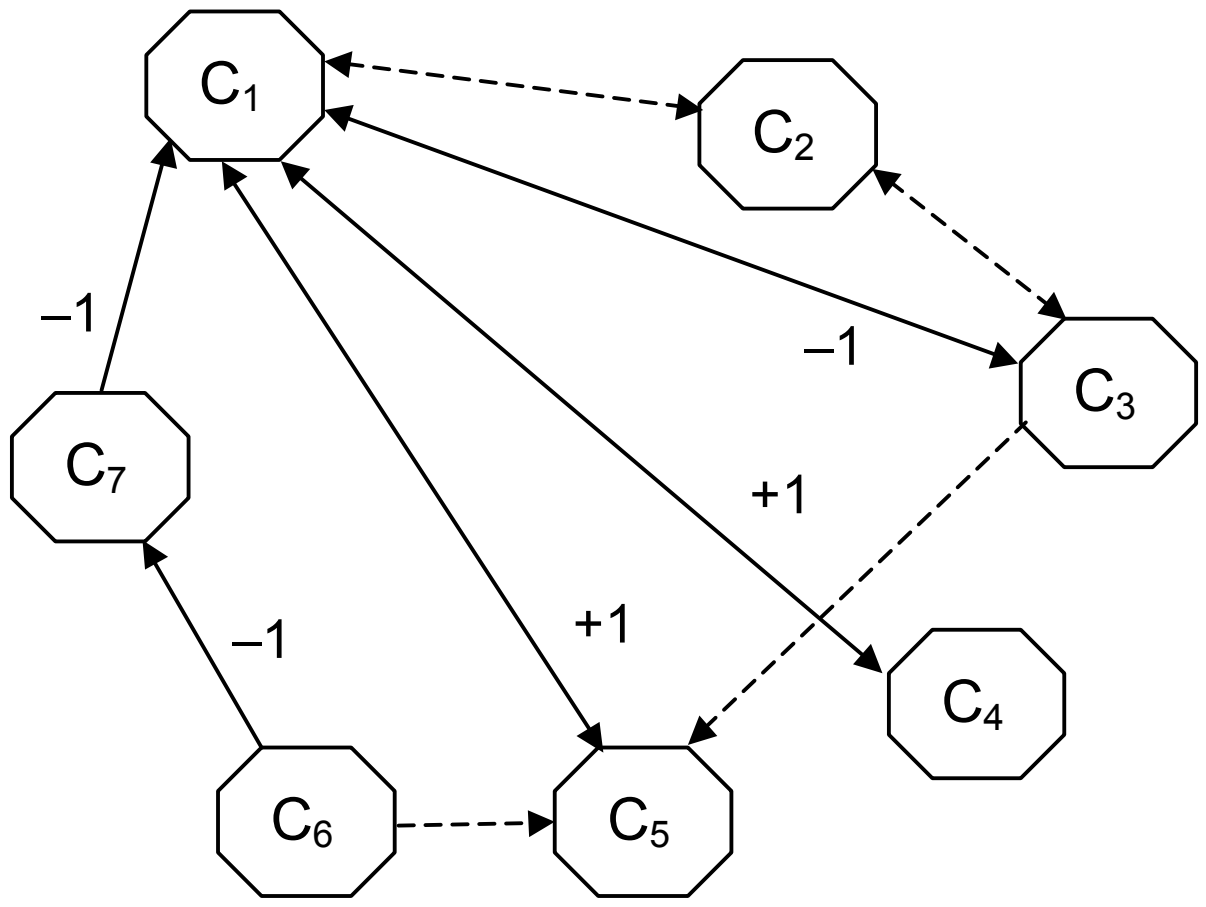

**FIGURE: 2.4.2**

The corresponding neutrosophic adjacency matrix N(E) related to the neutrosophic directed graph (Figure 2.4.2.) is given below:

$$N(E) = \begin{bmatrix} 0 & I & -1 & 1 & 1 & 0 & 0 \\ I & 0 & I & 0 & 0 & 0 & 0 \\ -1 & I & 0 & 0 & I & 0 & 0 \\ 1 & 0 & 0 & 0 & 0 & 0 & 0 \\ 1 & 0 & 0 & 0 & 0 & 0 & 0 \\ 0 & 0 & 0 & 0 & I & 0 & -1 \\ -1 & 0 & 0 & 0 & 0 & 0 & 0 \end{bmatrix}$$

Suppose we take the state vector $A_1$ = (1 0 0 0 0 0 0). We will see the effect of $A_1$ on E and on N(E).

| | | | | | | |
|---|---|---|---|---|---|---|
| $A_1E$ | = | (0 0 0 1 1 1 –1) | → | (1 0 0 1 1 1 0) = | $A_2$ |
| $A_2$ E | = | (2 0 0 1 1 1 0) | → | (1 0 0 1 1 1 0) = | $A_3 = A_2$. |

Thus child labor flourishes with parents' poverty and industrialists' action. Public practicing child labor also flourish but good NGOs are absent in such a scenario. The state vector gives the fixed point.

Now we find the effect of $A_1$ = (1 0 0 0 0 0 0) on N(E).

| | | | | | |
|---|---|---|---|---|---|
| $A_1$ N(E) | = | (0 I –1 1 1 0 0) | → | (1 I 0 1 1 0 0) = | $A_2$ |
| $A_2$ N(E) | = | (I + 2, I, –1+ I, 1 1 0 0) | → | (1 I 0 1 1 0 0) = | $A_2$ |

Thus $A_2$ = (1 I 0 1 1 0 0), according to this expert the increase or the on state of child labor certainly increases with the poverty of parents and other factors are indeterminate to him. This mainly gives the indeterminates relating to political leaders and teachers in the neutrosophic cognitive model and the parents poverty and Industrialist become to on state.

However, the results by FCM give as if there is no effect by teachers and politicians for the increase in child labor. Actually the increase in school dropout increases the child labor hence certainly the role of teachers play a part. At least if it is termed as an indeterminate one would think or reflect about their (teachers) effect on child labor.

Also the node the role played by political leaders has a major part; for if the political leaders were stern about stopping the child labor, certainly it cannot flourish in the society. They are ignored for two reasons: First, if children were vote banks certainly their position would be better. The second reason is, industrialists who practice child labor, are a main source of help to politicians, and their victory/defeat depends on their (financial) support so the causes for politicians ignoring child labor is two-fold.

Now we seek the opinion of another expert who is first asked to give a FCM model and then a provocative questionnaire discussing about the indeterminacy of relation between nodes is suggested and he finally gives a neutrosophic version of his ideas.

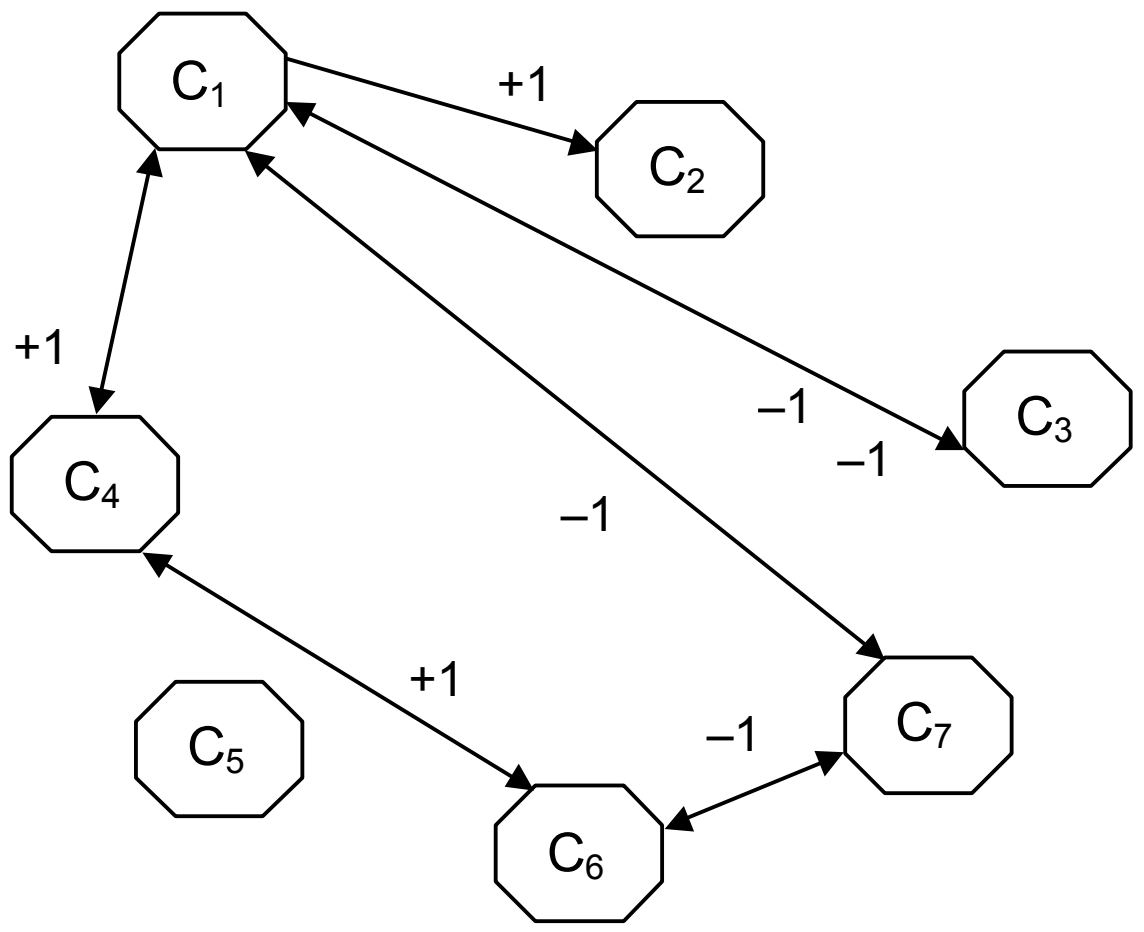

**FIGURE: 2.4.3**

Figure 2.4.3 is the directed graph of the expert. The related connection matrix $E_l$ is as follows:

$$E_l = \begin{bmatrix} 0 & 1 & -1 & 1 & 0 & 0 & -1 \\ 0 & 0 & 0 & 0 & 0 & 0 & 0 \\ -1 & 0 & 0 & 0 & 0 & 0 & 0 \\ 1 & 0 & 0 & 0 & 0 & 1 & 0 \\ 0 & 0 & 0 & 0 & 0 & 0 & 0 \\ 0 & 0 & 0 & 1 & 0 & 0 & -1 \\ -1 & 0 & 0 & 0 & 0 & -1 & 0 \end{bmatrix}.$$

Take $A_l$ = (1 0 0 0 0 0 0) the effect of $A_l$ on the system $E_l$ is

$$A_lE_l = (0\ 1\ -1\ 1\ 0\ 0\ -1) \to (1\ 1\ 0\ 1\ 0\ 0\ 0) = A_2$$
$$A_2E_2 = (1\ 1\ -1\ 1\ 0\ 1\ -1) \to (1\ 1\ 0\ 1\ 0\ 1\ 0) = A_3$$
$$A_3E_2 = (1\ 1\ -1\ 2\ 0\ 1\ -2) \to (1\ 1\ 0\ 1\ 0\ 1\ 0) = A_4 = A_3.$$

Thus according to this expert child labor has direct effect on political leaders, no effect on good teachers, effect on poverty and industrialists and no-effect on the public who encourage child labor; and good NGOs.

The same person was now put with the neutrosophic questions i.e. terms like "can you find any relation between the nodes or are you not in apposition to decide any relation between two nodes and so on"; so that a idea of indeterminacy is introduced to them.

Now the neutrosophic directed graph is drawn using this experts opinion given in Figure 2.4.4.

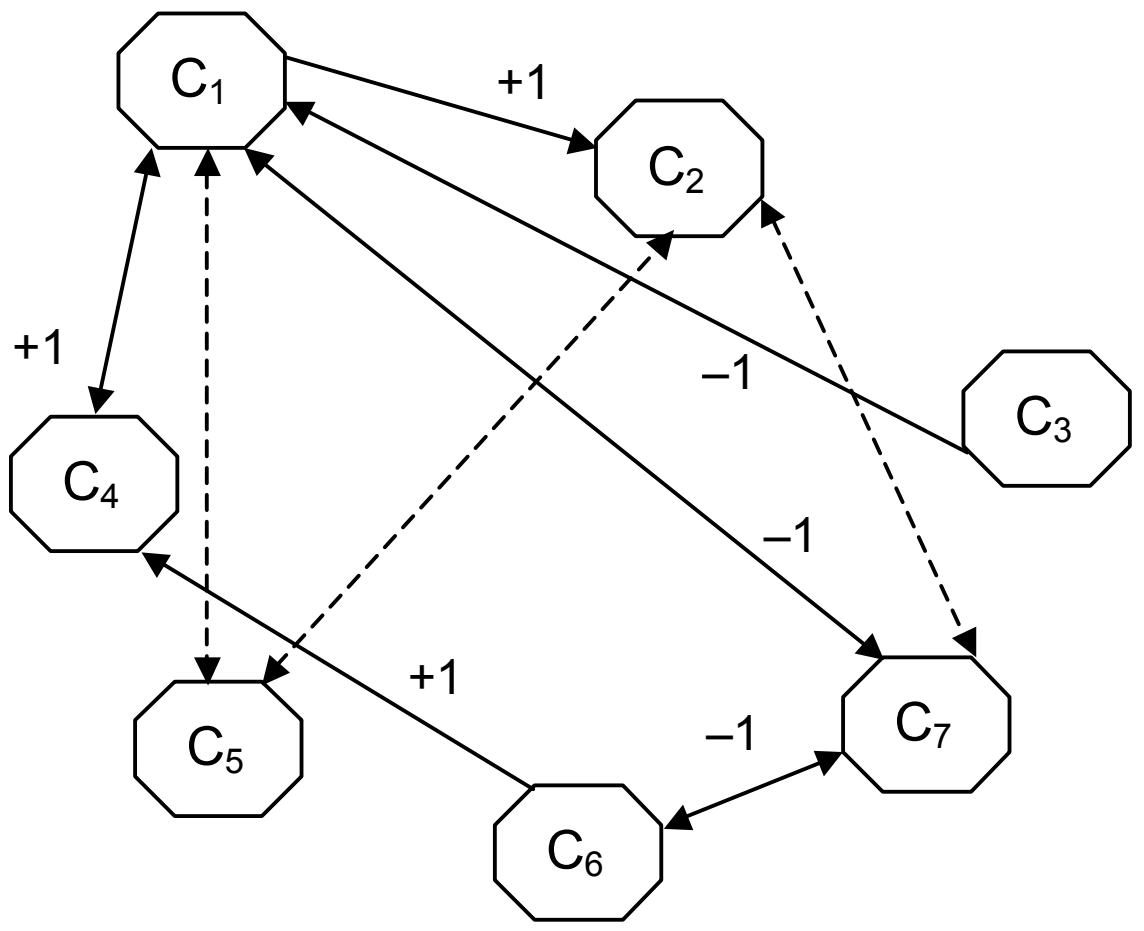

**FIGURE: 2.4.4**

The corresponding neutrosophic connection matrix $N(E_1)$ is as follows:

$$N(E_1) = \begin{bmatrix} 0 & 1 & -1 & 1 & I & 0 & -1 \\ 0 & 0 & 0 & 0 & I & 0 & I \\ -1 & 0 & 0 & 0 & 0 & 0 & 0 \\ 1 & 0 & 0 & 0 & 1 & 0 & 0 \\ I & I & 0 & 0 & 0 & 0 & 0 \\ 0 & 0 & 0 & 1 & 0 & 0 & -1 \\ -1 & 0 & 0 & 0 & 0 & -1 & I \end{bmatrix}.$$

Suppose $A_1 = (1, 0, 0, 0, 0, 0, 0)$ is the state vector whose effect on the neutrosophic system $N(E_1)$ is to be considered.

$A_1N(E_1) = (0\ 1\ -1\ 1\ I\ 0\ -1) \rightarrow (1\ 1\ 0\ 1\ I\ 0\ 0) = A_2$

$A_2N(E_1) = (1+I, 1+I, -1, 1, 2I+1, 0\ -1+I) \rightarrow (1\ 1\ 0\ 1\ 1\ 0\ 0) = A_3$

$A_3N(E_1) = (1+I, 1+I, -1, 2\ I+1\ 0\ -1 + I) \rightarrow (1\ 1\ 0\ 1\ 1\ 0\ 0) = A_4.$

We see $A_2 = A_3$.

But according to the NCM when the conceptual node child labor is on it implies that the cause of it is political leaders, poverty and industrialists participation by employing children as laborers.

The reader is expected to compare this NCM with FCMs for the same problem which is dealt earlier in this book as we have now indicated how a NCM works.

***Example 2.4.2:*** Application of NCM to study the Hacking of e-mail by students. One of the major problems of today's world of information technology that is faced by one and all is; How safe are the messages that are sent by e-mail? Is there enough privacy? For if a letter is sent by post one can by certain say that it cannot be read by any other person, other than the receiver. Even tapping or listening (over hearing) of phone calls from an alternate location / extension is only a very uncommon problem.

However compared to these modes of communication even though e-mail guarantees a lot of privacy it is a highly common practice to hack e-mail. Hacking is legally a cyber crime but is also one of the crimes that does not leave any trace. Hacking of another persons e-mail account can be carried out for a variety of purposes to study the factors, which are root-causes of such crimes we use NCM to analyse them. The following nodes are taken as the conceptual nodes.

| | | |
|---|---|---|
| $C_1$ | - | Curiosity |
| $C_2$ | - | Professional rivalry |
| $C_3$ | - | Jealousy/ enmity |
| $C_4$ | - | Sexual satisfaction |
| $C_5$ | - | Fun/pastime |
| $C_6$ | - | To satisfy ego |
| $C_7$ | - | Women students |
| $C_8$ | - | Breach of trust. |

However more number of conceptual nodes can be added as felt by the expert or the investigator. The neutrosophic directed graph as given by an expert is given in Figure 2.4.5.

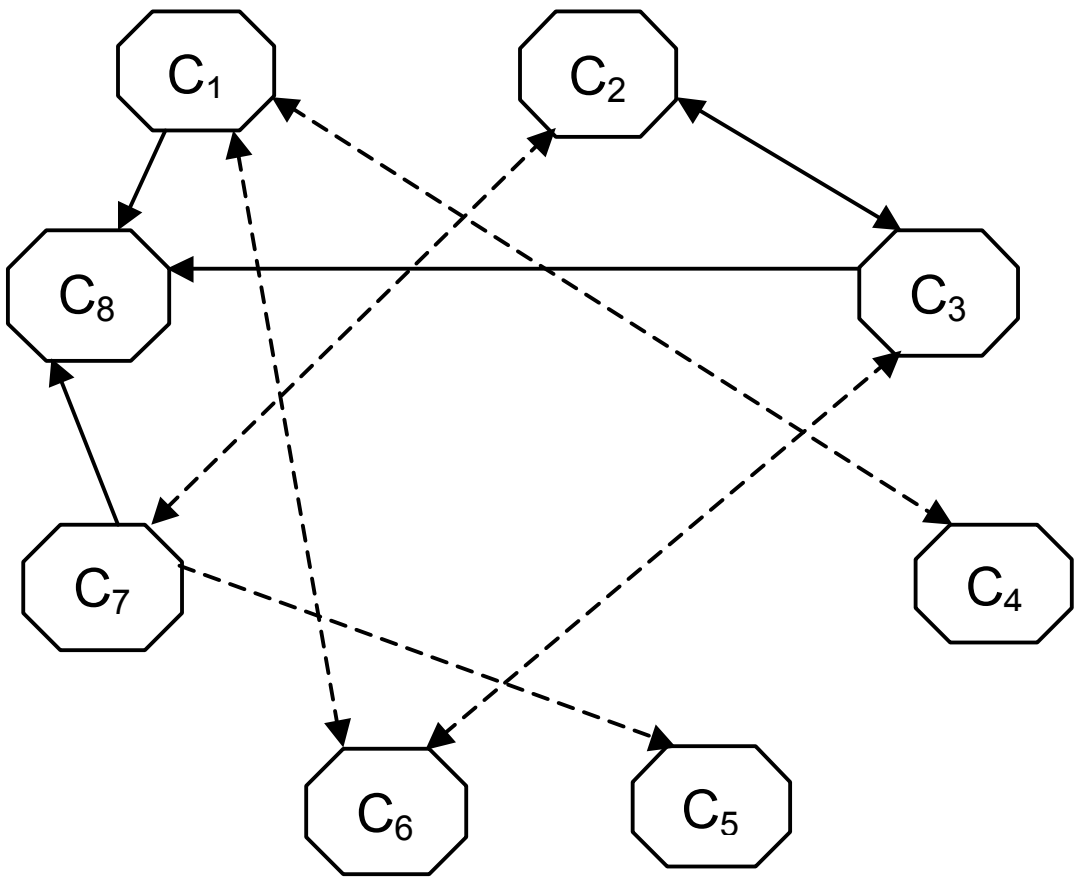

**FIGURE: 2.4.5**

The corresponding neutrosophic connection matrix N(E) is as follows:

$$N(E) = \begin{bmatrix} 0 & 0 & 0 & I & 0 & I & 0 & 1 \\ 0 & 0 & 1 & 0 & 0 & 0 & I & 0 \\ 0 & 1 & 0 & 0 & 0 & I & 0 & 1 \\ I & 0 & 0 & 0 & 0 & 0 & 0 & 0 \\ 0 & 0 & 0 & 0 & 0 & 0 & 0 & 0 \\ I & 0 & I & 0 & 0 & 0 & 0 & 0 \\ 0 & I & 0 & 0 & I & 0 & 0 & 1 \\ 0 & 0 & 0 & 0 & 0 & 0 & 0 & 0 \end{bmatrix}.$$

Suppose we take the instantaneous state vector $A_1$ = (0 0 0 0 0 0 1 0), women students node is in the on state then the effect of $A_1$ on the neutrosophic system N(E) is given by

| | | | | | | |
|---|---|---|---|---|---|---|
| $A_1N(E)$ | = | (0 I 0 0 I 0 0 1) | → | (0 I 0 0 I 0 1 1) | = | $A_2$ |
| $A_2N(E)$ | = | (0 I I 0 I 0 I 1) | → | (0 I I 0 I 0 1 1) | = | $A_3$ |
| $A_3N(E)$ | = | (0 I I 0 I I I 1+I) | → | (0 I, I, 0, I, I, 1, 1) | = | $A_4$ |
| $A_4N(E)$ | = | (I, I I 0 I I I I +1) | → | (I I I 0 I I 1 1) | = | $A_5$ |
| $A_5N(E)$ | = | (I I I 0 I I I 1) | → | (I I I 0 I I 1 1) | = $A_6$ = | $A_5$. |

So in case the node "women students" is in the on state node we see curiosity is an indeterminate, professional rivalry is an indeterminate, jealousy/ enmity is an indeterminate sexual satisfaction is in the off state, fun/ pastime is an indeterminate, to satisfy ego is an indeterminate and breach of trust is in the on state whereas if the 'I's are removed and N(E) is used as a usual FCM matrix then the effect $B_1$ = (0 0 0 0 0 0 1 0) in the on state when passed through the system we get B = (0 0 0 0 1 0 1 1) implies to satisfy ego becomes the on state and Breach of trust is in the on state. Thus we see other sates are in the off state.

The reader is expected to work with other coordinates and compare with FCMs which is got by replacing all I's in the neutrosophic connection matrix N(E) by 0.

Several other examples can be shown using the method of NCM.

## 2.5 Some more Illustrations of NCMs

The study of illustrations of FCMs are given in Sections 1.3. of this book. Here in this section we apply NCMs only to some of the illustrations mentioned in that section.

The main motivation for doing so is that any reader can compare the FCM and NCM on the same model and draw conclusions in that regard.

### 2.5.1: NCM for Decision Support in an Intelligent Intrusion Detective System

A complete introduction of the problem of intelligent intrusion detective system [89] with description is given in Section 1.3.1. of this book.

Now in case of FCMs, the edge values come from expert knowledge and experience. These are functions of the expert's common sense and engineering judgments. Moreover the parameters are tunable in FCMs flexible structure. Now when we replace the FCMs by NCMs, we make it possible for the expert to express his feelings of indeterminacy of some relation between certain nodes. If FCM is used, these edges cannot find any place except a 'zero' but if NCM is used certainly they can be weighed by 'I' an element of indeterminacy.

The forged signature will be one of concepts of indeterminacy as no one can easily detect the forged signature. Likewise, several other concepts might possibly be indeterminate.

The task of modeling this problem as an NCM is left as a research exercise for the reader. For more on Intelligent Intrusion Detection Systems the reader can refer [89].

### 2.5.2: Analysis of Strategic Planning Simulation based on NCMs Knowledge and Differential Game

FCM has been used in the study of differential game [65] we use same map of FCM but after discussing with an expert convert it into an NCM by adjoining the edges which are indeterminate, and this is mainly carried out for easy comparison.

Now according to this expert, competitiveness and market demand is an indeterminate. Also sales price and economic condition is an indeterminate. Also according to him the productivity and market share is an indeterminate whether a relation exists directly cannot be said but he is not able to state that there is no relation between these concepts so he says let it be an indeterminate.

Also according to him quality control and market share is an indeterminate. Thus on the whole the market share is an FCM with a lot of indeterminacy so is best fit with an NCM model.

The reader is expected to study the two tables given in Section 1.3.2 (the Table 1 and Table 2 related to the figure 1.3.4) with NCM using Figure 2.5.1 and obtain analysis and conclusions using NCMs and compare it with FCM. Thus obtain the initial version of NCM matrix and refined version of NCM matrix, also give the corresponding comment. Study the factor of indeterminacy and prove the result is nearer to truth for finding solutions to the market share problem. Compare FCM and NCM in the case of market share problem.

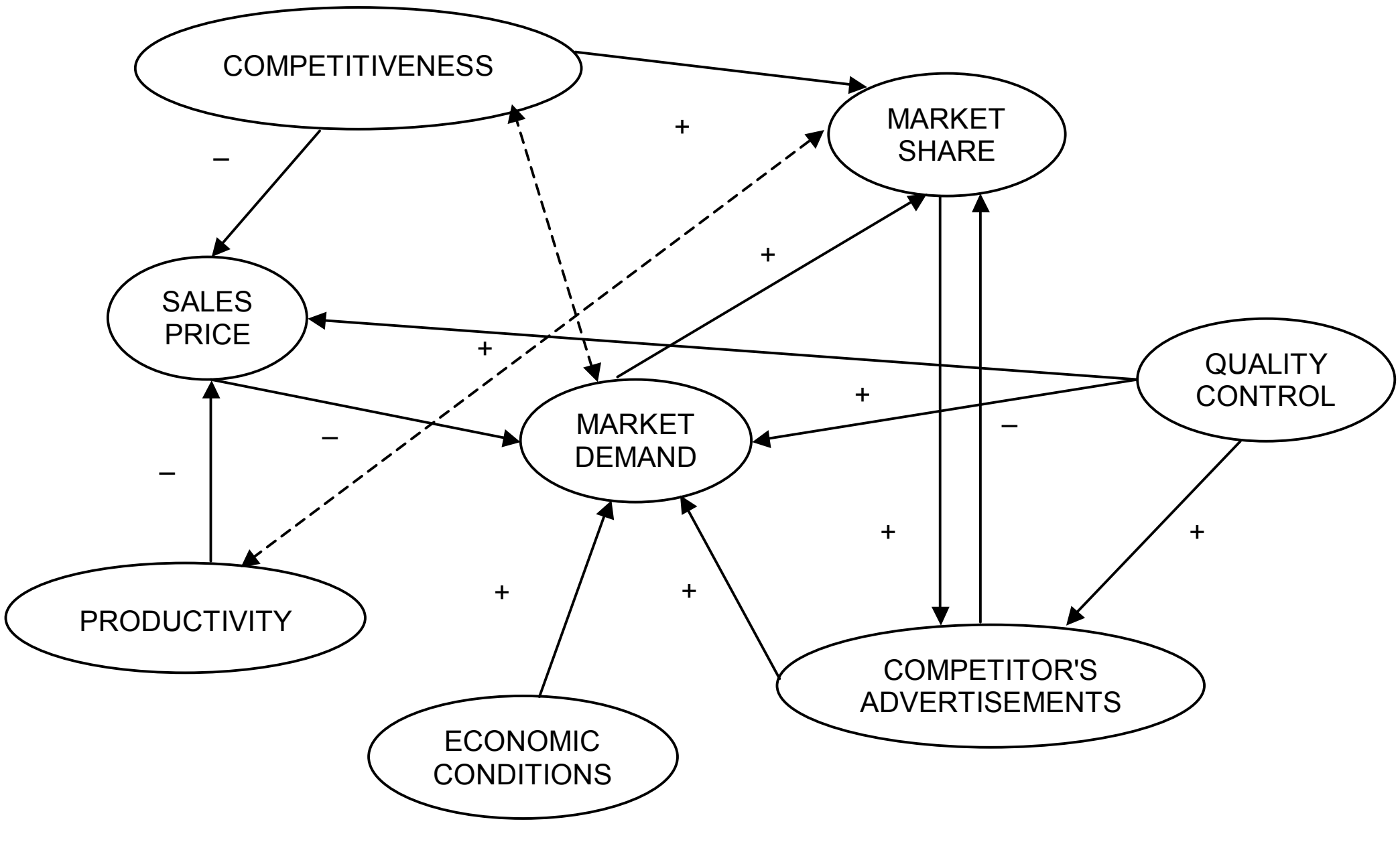

**FIGURE: 2.5.1**

It is pertinent to mention that K.C. Lee et al [65] have suggested that for better results more refined FCM with edge values can be used. In our opinion NCM may be a better solution to Lee's problems suggested in [65].

## 2.5.3: Adaptive Neutrosophic Cognitive maps for hyper knowledge representation in strategy formation process

The application or use of adaptive fuzzy cognitive maps for hyper knowledge representation in strategy formation process was already recalled in Chapter 1; Section 1.3.3. For more about it please refer Carlsson and Fuller [17].

Adaptive Fuzzy Cognitive Maps can learn the weight from historical data. Once the FCM is trained, it lets us play what-if games (eg. what if demand goes up and prices remain stable? i.e. we improve our market position). Likewise adaptive neutrosophic cognitive maps are fuzzy cognitive maps with an addition concept between two nodes when the relation between them is an indeterminate. Once the NCM and the related FCM which is got when I's are replaced by zeros is trained, let us play what-if games to predict the future in a realistic way.

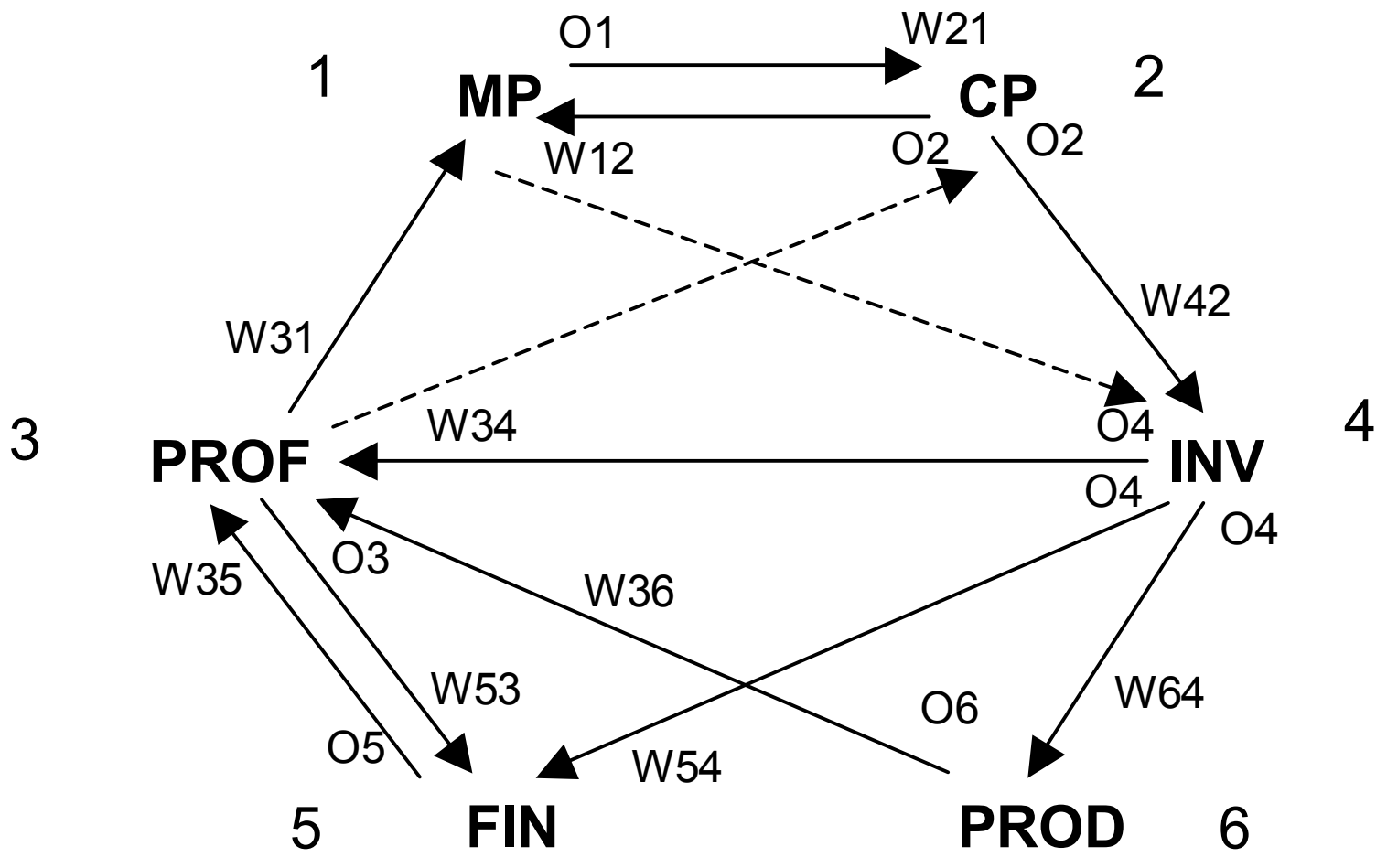

**FIGURE: 2.5.2**

Adaptive Neutrosophic Cognitive Map for the strategy formation process. The neutrosophic matrix N(w) is obtained using Figure 2.5.2. The reader is advised to study analyse and compare the results of NCM with FCM, for more about Adaptive FCM please refer [17]. From the Figure 2.5.2 we see according to the expert the relations between MP and INV is indeterminate also that between CP and PROF is an indeterminate. Using this and results in Figure 1.3.7, discuss the differences.

## 2.5.4: A New Balance Degree for Neutrosophic Cognitive Maps (NCMs)

We first define when is an NCM imbalanced. In the opposite case, we say the NCM is balanced analogous to balance degree of FCM given by Tsadiras et al [107].

**DEFINITION 2.5.4.1:** *An NCM is imbalanced if we can find two paths between the same two nodes that create causal relations of different sign. In the opposite case the NCM is balanced. The term 'balanced' neutrosophic digraph is used in the following sense that is in a imbalanced NCM we cannot determine the sign or the presence of indeterminacy of the total effect of a concept to another.*

*Now on similar lines based on the idea that as the length of the path increases, the indirect causal relations become weakened the total effect should have the sign of the shortest path between two nodes.*

The degree to which a neutrosophic digraph of the NCM is balanced or imbalanced is given by the balance degree of the neutrosophic digraph. There are as in case of graph various types of balance degree in case of neutrosophic graphs also defined purely in an analogous way. An interested reader can obtain nice results on Balanced Degree of Neutrosophic Digraphs using results from [107].

### 2.5.5: Rule based Neutrosophic cognitive maps

NCM is a FCM in which the relation between some of the causal nodes are indeterminates. Rule-based FCMs have been studied by [21] described in Section 1.3.5 of this book. We can on similar lines analyse and formulate the Rule Based Neutrosophic Cognitive Maps (RB-NCM).

Rule Based Neutrosophic Cognitive Maps (RB-NCM) provide a representation of the dynamics of complex real world qualitative systems with feedback and allow the simulation of the occurrences of events and their influence in the system. They are fuzzy directed neutrosophic graphs with feedback which are composed of fuzzy nodes and both neutrosophic links and fuzzy links. In RB-FCM unlike concepts of NCM, concepts of neutrosophic variables defined using neutrosophic logic and the neutrosophic rule bases are used to express relations; i.e. in short RB-NCM are RB-FCMs together with the indeterminate relations 'I'.

The reader is expected to apply and illustrate it in real models using RB-NCMs and compare it with the methods of RB-FCMs.

### 2.5.6: Neutrosophic Causal Relationship and Neutrosophic Partially Causal Relationship

Kim and Lee have studied and defined in the notion of Fuzzy Causal Relationship and Fuzzy Partial Causal Relationship in [50]. The definitions were recalled in Section 1.3.6. Here we just define Neutrosophic Causal Relation. Using the usual notation for $Q_i$ and $M_i$ from [50], the definition follows:

**DEFINITION 2.5.6.1:**

1. *$C_i$ causes $C_j$ if and only if*
   a. $(Q_i \cap M_i) \subset (Q_j \cap M_j)$ *and*
   b. $(\tilde{}Q_i \cap M_i) \subset (\tilde{}Q_j \cap M_j)$

2. *$C_i$ causally decreases $C_j$ if and only if*
   a. $(Q_i \cap M_i) \subset (\tilde{}Q_j \cap M_j)$
   b. $(\tilde{}Q_i \cap M_i) \subset (Q_j \cap M_j)$

3. *$C_i$ indetermines $C_j$ (i.e. one has doubt whether any relation exists i.e., relation between $C_i$ and $C_j$ is indeterminable) if and only if one of the relations 1a or 1b or 2a or 2b holds good and other behaves in a chaotic way then we say the relation between $C_i$ and $C_j$ is an indeterminate. ('or' in the mutually exclusive sense the occurrence of the other).*

*Now we call this relationship a Neutrosophic Causal Relation denoted by NCR. Thus if NCR is true then one of the condition 1a, 1b, 2a or 2b is true and the other relations remain indeterminate which can occur in the real world problem of stock market.*

On similar lines, one can define Neutrosophic Partially Causal Relations analogous to Fuzzy Partially Causal Relations and study them.

### 2.5.7: Automatic construction of NCMs

Schneider et al have introduced the notion of automatic construction of FCMs. It is briefly described in Section 1.3.7. Here we describe automatic construction of NCM.

The method of automatically constructing Neutrosophic Cognitive Maps based on the user provided data:

This method consists of finding the degrees of similarity between any two variables (represented by numerical neutrosophic vectors), finding whether the relation between is direct or inverse or neutrosophic which is determined by finding out the causality among variables (neutrosophic vectors).

Using the same type of formulas except replacing the usual vectors by neutrosophic vectors all the equations described by Schneider and et al recalled in Section 1.3.7 can be derived with proper or appropriate changes as the case may be. It can be achieved as a matter of routine. The readers are advised to find automatic construction of NCMs for the economic model given in [86] and compare it with automatic construction of FCMs given by [86].

Finally we proceed to describe the notion of Neutrosophic Cognitive State Maps.

### 2.5.8 Illustration of Neutrosophic Cognitive State maps of users web behavior

Searching for information in general is complex, with lot of indeterminacies and it is an uncertain process for it depends on the search engine; number of key words, sensitivity of search, seriousness of search etc. Hence we can see several of the factors will remain as indeterminates for the $C_1, C_2, \ldots, C_7$ as given by [71] in 1.3.8 we can remodel using NCM.

The NCM modeling of the users web behavior is given by the following neutrosophic digraph (Figure 2.5.3) and the corresponding N(E) built using an expert opinion is given by the following neutrosophic matrix:

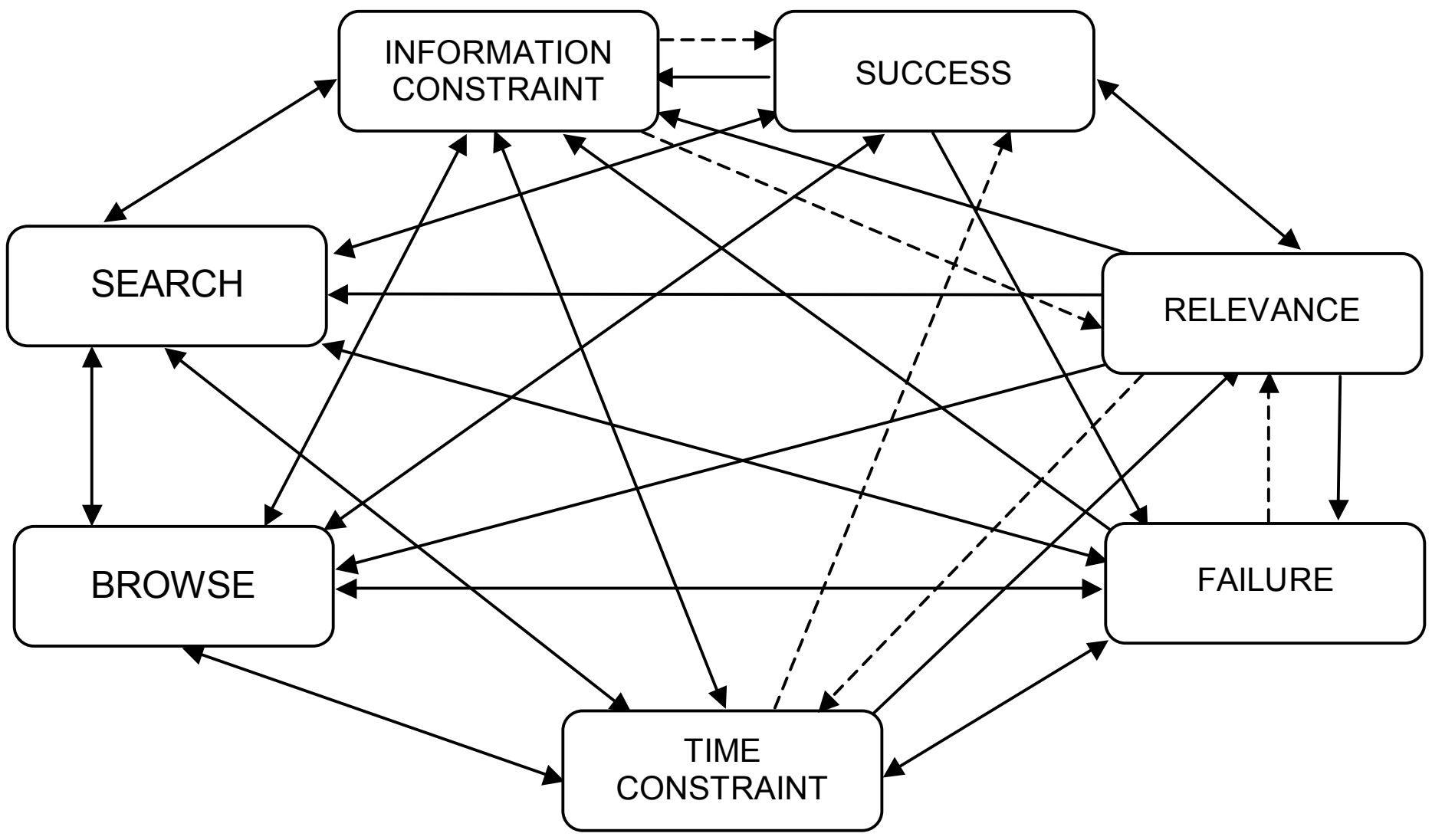

**FIGURE: 2.5.3**

$$N(E) = \begin{bmatrix} 0 & -1 & -1 & 1 & -1 & -1 & 1 \\ 1 & 0 & -1 & -1 & -1 & -1 & 1 \\ -1 & -1 & 0 & -1 & I & 1 & 1 \\ -1 & -1 & -1 & 0 & 0 & I & I \\ 1 & 1 & 1 & 1 & 0 & 1 & -1 \\ 1 & 1 & I & 1 & 1 & 0 & -1 \\ 1 & 1 & 1 & 1 & 0 & I & 0 \end{bmatrix}$$

Several results and conclusions can be derived for each of the state vectors. The reader is given the work of comparing Fuzzy Cognitive State Map given in Section 1.3.8 with the Neutrosophic Cognitive State Map given in Figure 2.5.3.

Now we proceed on to give some more applications of NCM.

### 2.5.9 Neutrosophic mapping and uncertainty Neuron neutrosophic cognitive maps

Tsadiras et al have defined and studied in [108] the notion of cognitive maps and Certainty Neuron FCMs. In Section 1.3.9 we briefly recall certain important concepts from [108]. Here, we just define the concept of Neutrosophic mapping and the Certainty Neuron Neutrosophic Cognitive Maps. We expect the readers to develop real world models using them.

**DEFINITION 2.5.9.1:** *Neutrosophic mapping is a cognitive map in which we associate value 'I' an indeterminacy between at least two nodes i.e., to be more formal a neutrosophic map (NM) can be represented by means of signed, directed neutrosophic graph, where concepts are represented by nodes and causal relation and indeterminable casual relations are also included.*

*A Neutrosophic Cognitive Map using the public health model is given in Figure 1.3.12(a). As far as India is concerned when we cannot determine any direct relation between migration into city and sanitation facility (possibly owing to ethnicity, type of migrants, and period of settlement), the concept of indeterminacy comes in handy.*

*Similarly migration into city and number of diseases per 1000 residents is an indeterminate for some experts are of the opinion that if the migrants are not used to city life they can develop different types of new diseases which will be found in only the migrants. Thus relation between these two concepts can be treated as an indeterminate. Likewise different experts when given an option, can, depending on their place of study, give the indeterminacy between nodes.*

Thus we can use NM in places where CM are used whenever we have a factor of the indeterminacy involved. CM is an FCM when it is weighed with values between [–1, 1]. Likewise NM will be a NCM if the edges take values from the set {0, 1, –1, I}. Certainty in any prediction cannot exist in reality, for, we cannot be certain of results that too in an unsupervised data.

For neurons that are used as transfer functions between two variable functions, which are dependent on the new activation level, can in certain cases be such that the very existence of any sort of relation is indeterminate or unpredictable or uncertain. In such cases we built the notion of Uncertainty Neuron Neutrosophic Cognitive Maps. These instances do occur in the study or analysis of unsupervised data.

Keeping these arguments in mind the reader is expected to build NM and Uncertainty Neuron Neutrosophic Maps.

### 2.5.10: NCM considering Time Relationships

We have recalled the definition of FTCM in Section 1.3.10 of this book from [78].

Now we define Neutrosophic Time Cognitive Map (NTCM).

**DEFINITION 2.5.10.1:** *A Neutrosophic Time Cognitive Map (NTCM) is a neutrosophic directed graph G = (X, E) consisting of a finite set X of N nodes* $X = \{i\}^{n}_{i=1}$ *and a set E of edges (which include neutrosophic edges also)* $e_{ij}$ , $i, j \in X$ *i.e., the edge connecting some i and j can also be an indeterminate edge of the nodes i and j. Each edge* $e_{ij}$ *of the NTCM has two kinds of relative causalities the strength* $s_{ij} \in [-1, 1, I]$ *and the time lag* $t_{ij} \in [a, b]$ *where a and b are constants* $0 < a \leq b$*. In short NTCM is identified by X,* $S = ||s_{ij}||_{N\times M}$ $s_{ij}$ *can be –1, 1 or I and* $T = ||t_{ij}||_{N\times N}$.

The study of NTCM can be made in an analogous way as in case of FTCM , except that in the case of NTCM we use Neutrosophic Directed Graphs and the matrices are Neutrosophic Matrices.

The reader is expected to construct real world models using NTCM. For more about FTCM, please refer [78].

### 2.5.11: The study of balanced differential learning algorithm using NCM

The table of linguistic labels used for the balanced differential algorithm in FCM as given by [131] is discussed in 1.3.11. In study of NCM we have the following mapping between labels and values:

| SYMBOLIC VALUES | NUMERIC VALUES |
|---|---|
| Affects a lot | 1.0 |
| Affects | 0.5 |
| Does not affect | 0 |
| Affects negatively | -0.5 |
| Affects negatively a lot | –1 |
| Cannot measure how it affects | I - indeterminacy |

Using this neutrosophic mapping between labels and values, adaptive NCM can be constructed in an analogous way as adaptive FCM is constructed. The reader is given the task of building the adaptive NCM.

## 2.6 Application of NCMs

Here we will give the applications of NCMs in the place of FCMs wherever possible. The applications of FCMs, which have been dealt with in Section 1.4 of this book, are also discussed here with reference to NCMs. The reader is advised to work with any new application of his/her own using some real world problems using NCMs.

### 2.6.1: Modeling Supervisory Systems using NCM

The concept of NCM can be used in modeling of supervisory systems in the design of hybrid models for complex systems wherever the concept of indeterminacy play a role in that model. This model "The challenge of modeling supervisory systems using fuzzy cognitive maps" as carried out by C. Stylios et al [100] does not lend itself for the application of NCM as we see their does not exist an indeterminacy when the conceptual nodes are taken as the state of the valve closed open or partially closed.

However when we assume the valve is closed it may happen that there is some leakage in the valve or it is not closed properly in such case certainly we suggest the FCM structure can be modified and the NCM model can be implemented so that the results can be more accurate and when indeterminacy occurs we can have a doubt about the specific gravity of the liquid that is produced during the mixing when the measured specific gravity lies in a specified range.

Another indeterminacy that may occur is that in the events 1 to 8 listed in [100] it may so happen when the variation /amount of liquid in the tank is taken the six liquid filling the pipe wastage of liquid while operating the valves etc should be given enough representation or all these can be kept as an indeterminacy and work should be carried out so that the results obtained under this model happens to be better than the existing ones. The reader is expected to develop this model using NCM.

### 2.6.2: Design of Hybrid Models for Complex Systems: a NCM Approach

FCM is applied in the Hybrid Models of Complex Systems [39]. We now apply NCM to model the Complex Systems. By incorporating neutrosophic principles into a neural network we feel that the edge weight $e_{ij}$ can also be indeterminate. So even the simple NCMs (i.e. edge weights $e_{ij} \in \{-1, 1, 0\ I\}$ act as asymmetrical networks of threshold or continuous neurons and converge to limit cycles.

It is suggested to involve NCM which best utilize existing experience in the operation of the system and are capable in modeling the behaviour of complex systems, as NCMs seems to be a useful method in complex system modeling and control which will help the designer of a system in decision analysis. The reader is expected to model complex systems using NCM working in an analogous way of [39] which is direct and easily any reader can do it as a matter of routine; and it is suggested to the reader as a piece of research.

### 2.6.3: Use of NCM in Robotics

While it has been argued that FCMs are preferred for usage in robotics and applications of intimate technologies, owing to their ability to handle contradictory inputs, NCMs would be the more viable tool, for not only are they capable of handling contradictory inputs, but they can also handle indeterminacy.

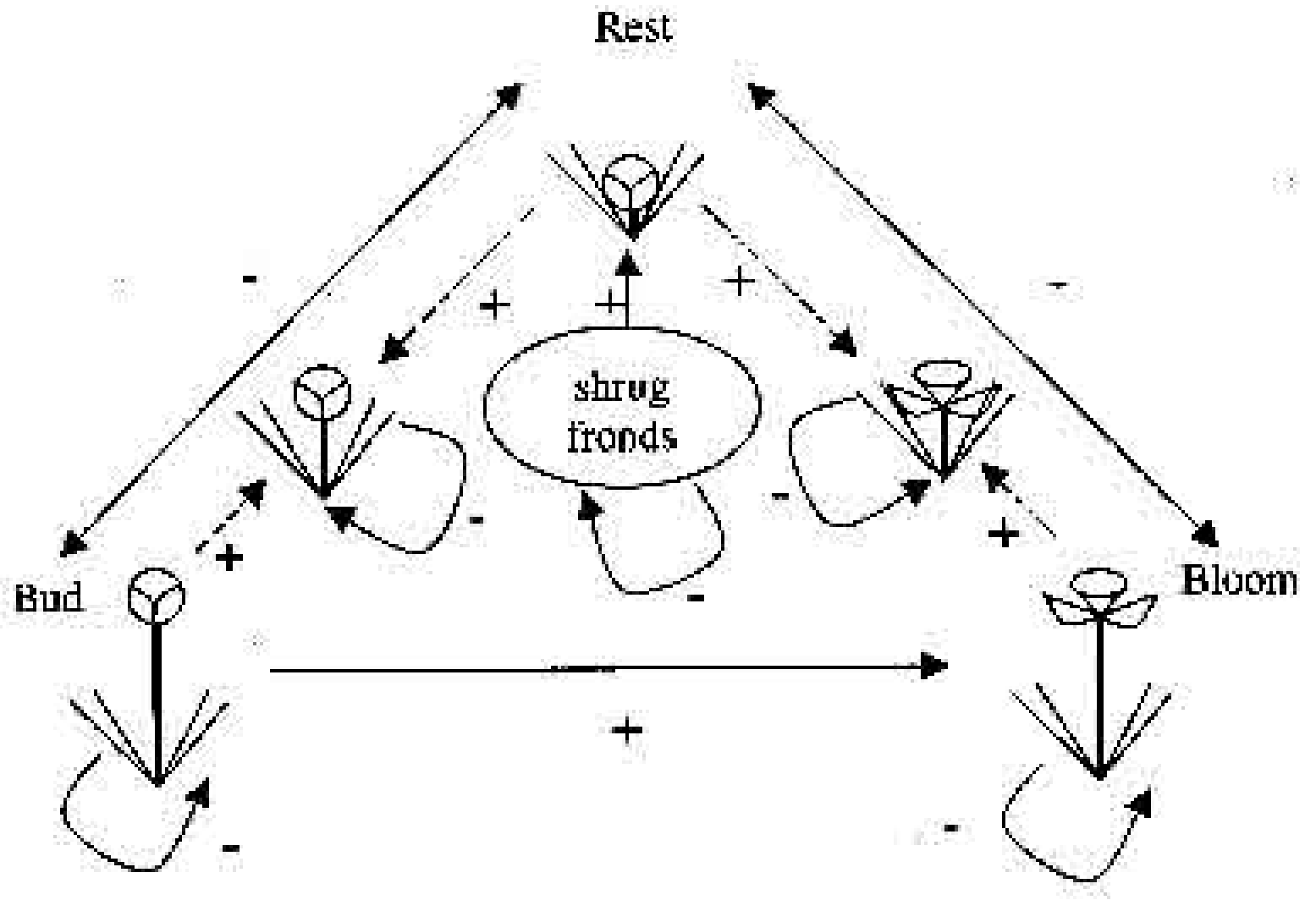

**FIGURE: 2.6.1**

Further FCMs have been used to model the Office Plant #1 to analyse the types of e-mails [9]. It is once again left as an exercise for the reader to use NCM in the place of FCM in this study. For categorically one cannot always divide the e-mails as official / non official, friendly / business like and so on for some can be termed as indeterminate, semi-friendly and semiofficial or so on and so forth.

So NCM can be adopted in mobile robots like Office Plant #1 and the study can be carried out as a maiden effort. A description of use of FCM is given in Section 1.4.3.

For example we see at each stage the relation would be indeterminate if the email received has an over-lapping attributes in which case the section of the node may be indeterminate. Thus in the behaviour of the office plant, the dotted arrows ought to be adopted in situations where there is indeterminacy.

### 2.6.4: Adaptation of NCMs to Model and Analyse Business Performance Assessment

The use of FCMs to model and analyse business performance assessment [48] was dealt by us in Section 1.4 and subsection 1.4.4.

For the study of NCM, introduce on the FCM the NCM structure so that one is in a position to analyse and compare it with the FCM in Figure 1.4.5. Using Figure 2.6.3 we can get the NCM.

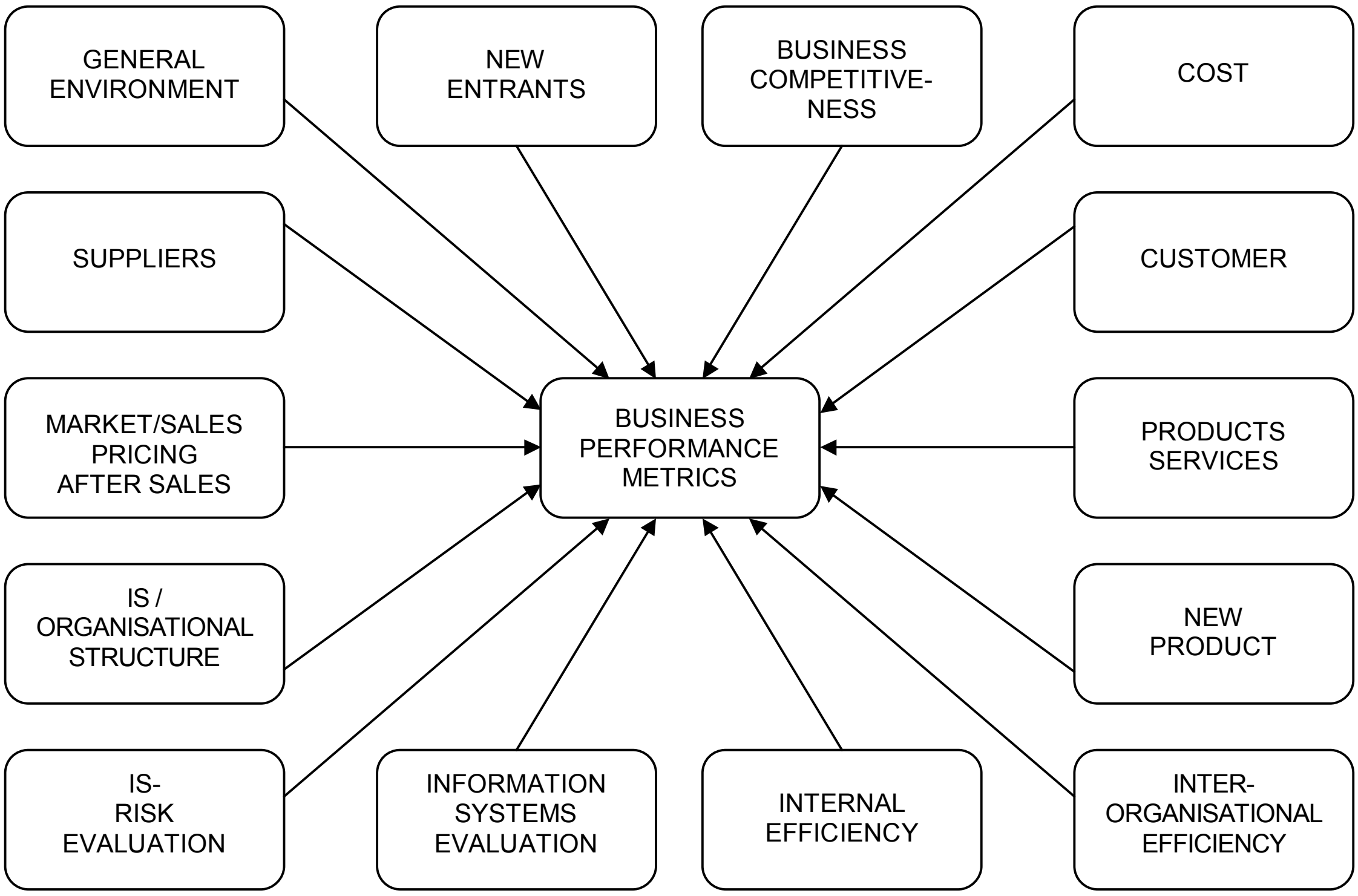

**FIGURE: 2.6.2**

The reader to requested to implement NCMs and draw conclusions based on the introduction of NCM to this model. Figure 2.6.2 describes the Business performance metrics given by [48].

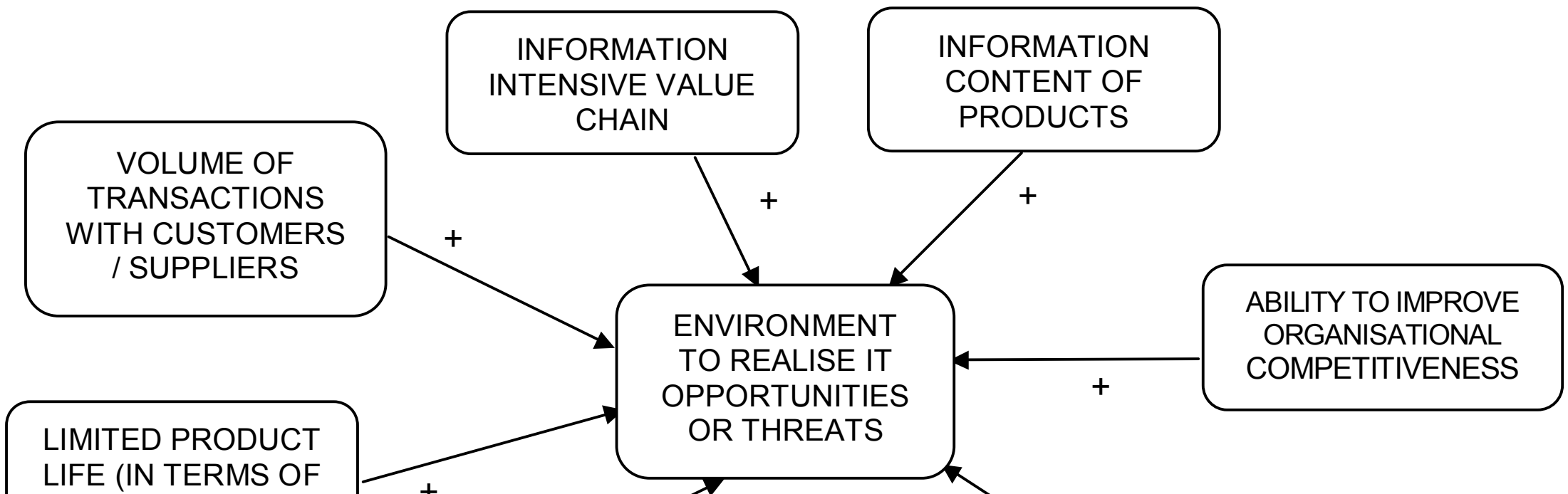

### 2.6.5: Use of NCMs in legal rules

The implementation of legal rules using FCMs are very well studied by [1] and we have discussed it in Section 1.4.5 of this book.

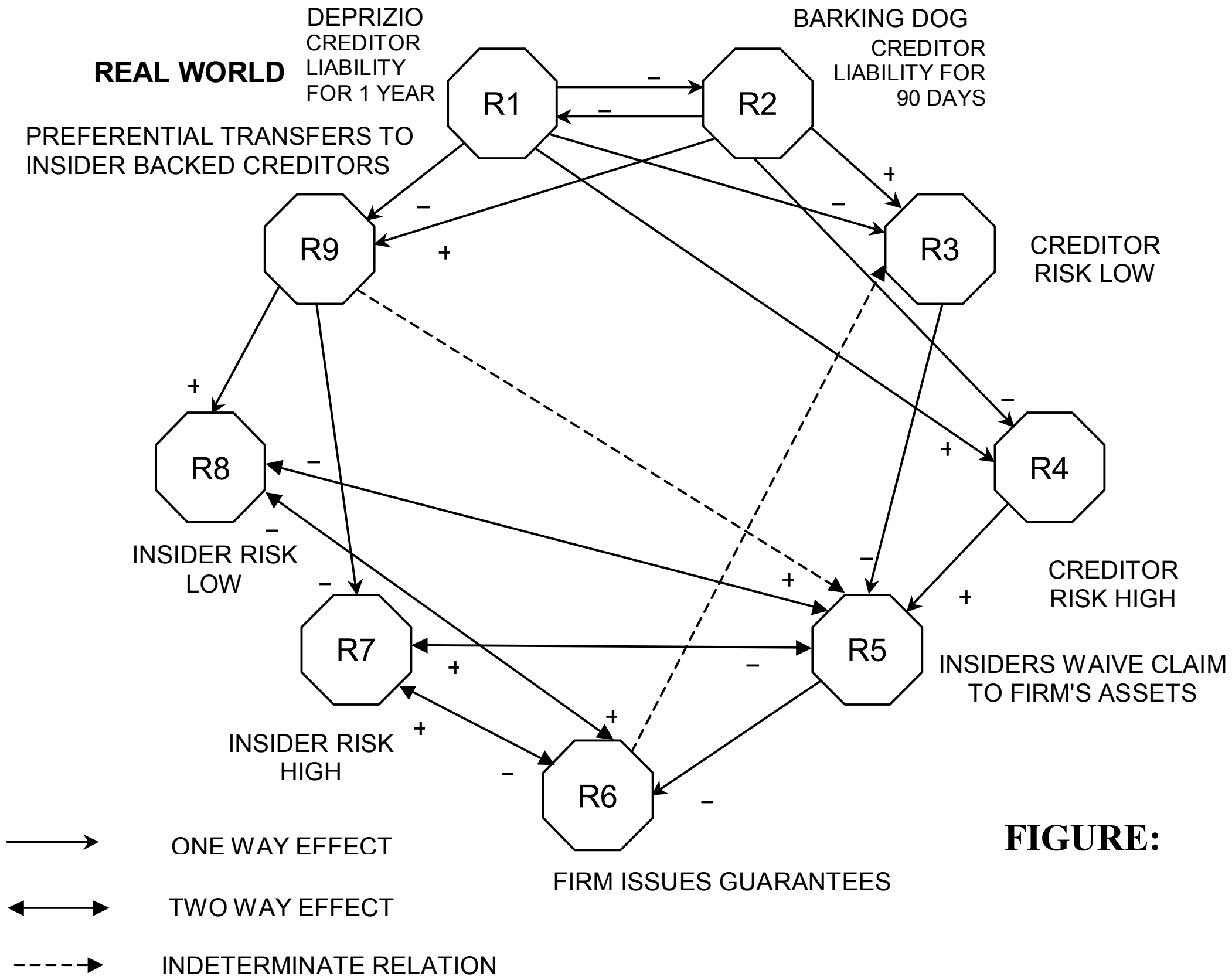

**FIGURE:**

Study NCM using Figure 2.6.4 and analyse and compare it with FCM given by [1]. Now we first show NCMs are better tools than the FCMs as in FCM we do not have the concept of indeterminacy. Only in case of NCMs we can say the relation between two concepts / attributes / nodes can be indeterminable also. For especially in criminal cases the concerned may not be able draw conclusions based on the data provided to him. Very many relations can be indeterminable so NCMs should be a better fit than the FCMs. Further instead of saying no relation between two nodes exist but still if the feelings exist with some doubt we cannot represent it in terms of FCMs but can easily implement NCMs so that while spelling out the judgment, the court can be very careful and give due weightages to the indeterminable relations.

The reader is given the task of implementing NCMs to the problem discussed in [1] using FCMs. They can also compare and contrast the FCM and NCM for this problem by studying where FCM has already been used and compare it with NCM and derive conclusions.

Now we illustrate how NCMs will play a major role in legal sides in several cases in India. We mention here a few types of them.

1. Encounter deaths with police.
2. Murder of political personalities which is very common in India
3. Custodial deaths in prison (Beaten to death, hanging, suicide etc.)
4. Undue delay in compensation cases in which government / private body is involved.
5. Means to punish intellectual harassment and torture and view it as more than the physical harassment and torture.

In all these five types of cases lot of indeterminacy is involved, so when NCM is applied certainly it will lead to better results.

### 2.6.6: Use of NCM in Creating Metabolic and Regulatory Network Models

The FCM model has been used by Dickerson et al [32] in creating metabolic and regulatory network models. In Section 1.4.6 we have provided a brief description of the problem that Dickerson et al [32] have handled. Can FCM be replaced by NCM in creating metabolic and regulatory network models?

The very fact that FCMs are used implies that there is lack of quantitative information on how different variables interact. The second is that the direction of causality is at least partly known and can be articulated by a domain expert. The third is that they link concepts from different domains together using arrows of causality.

These features are shared by the problem of modeling metabolic networks. Further FCModeler tool is intended to capture the intuitions of biologists, help test hypotheses and provide a modeling framework for assessing the large amounts of data captured by RNA, micro-arrays and other high-throughout experiments.

Since mainly the FCModeler tool is intended to capture intuitions of biologists, so in several or few cases they may not be able to say anything, for it may be also an indeterminate for the biologists so we are justified in applying NCMs and certainly NCMs can replace FCMs for effective results. The reader can surely attempt to build an NCModeler using NCMs.

### 2.6.7: Use of NCMs to find the driving speed in any one in freeway

Brubaker [13-14] used FCMs to create a model to find ones speed when driving in a California freeway. Thus FCM plays a major role in the study and analysis of transportation problem of all kind, for transportation problems are basically problems of decision-making. In our opinion, NCM can also be used to arrive at better results.

The concepts or nodes of the FCM are bad weather, freeway congestion, auto accidents, patrol frequency, own risk aversion, impatience and attitude. Now if these are taken as nodes certainly we can have pairs of nodes for which the relation is indeterminate, for the concept of impatience and attitude with other concepts like bad weather or free way congestion and speed of others is an indeterminate in itself. For fearing the bad weather one may be impatient and drive fast due to the fear that the weather may become worst or some other may fear bad weather (and consequent accidents) and be obsessed with fear and drive slow; the minute the nature of impatience or fear dominates a person certainly one cannot predict the speed, hence a lot of uncertainty and indeterminacy is involved. So the adaptation of NCM may yield a better understanding and modeling of the problem than FCMs.

Thus we request the reader to model this problem using NCM and compare it when only FCM is applied; complete discussions using FCM is given in Section 1.4 as taken from [13]. For a slight suggestion of how to go about with this NCM, we have for the reader's sake provided a possible NCM graph in Figure 2.6.5.

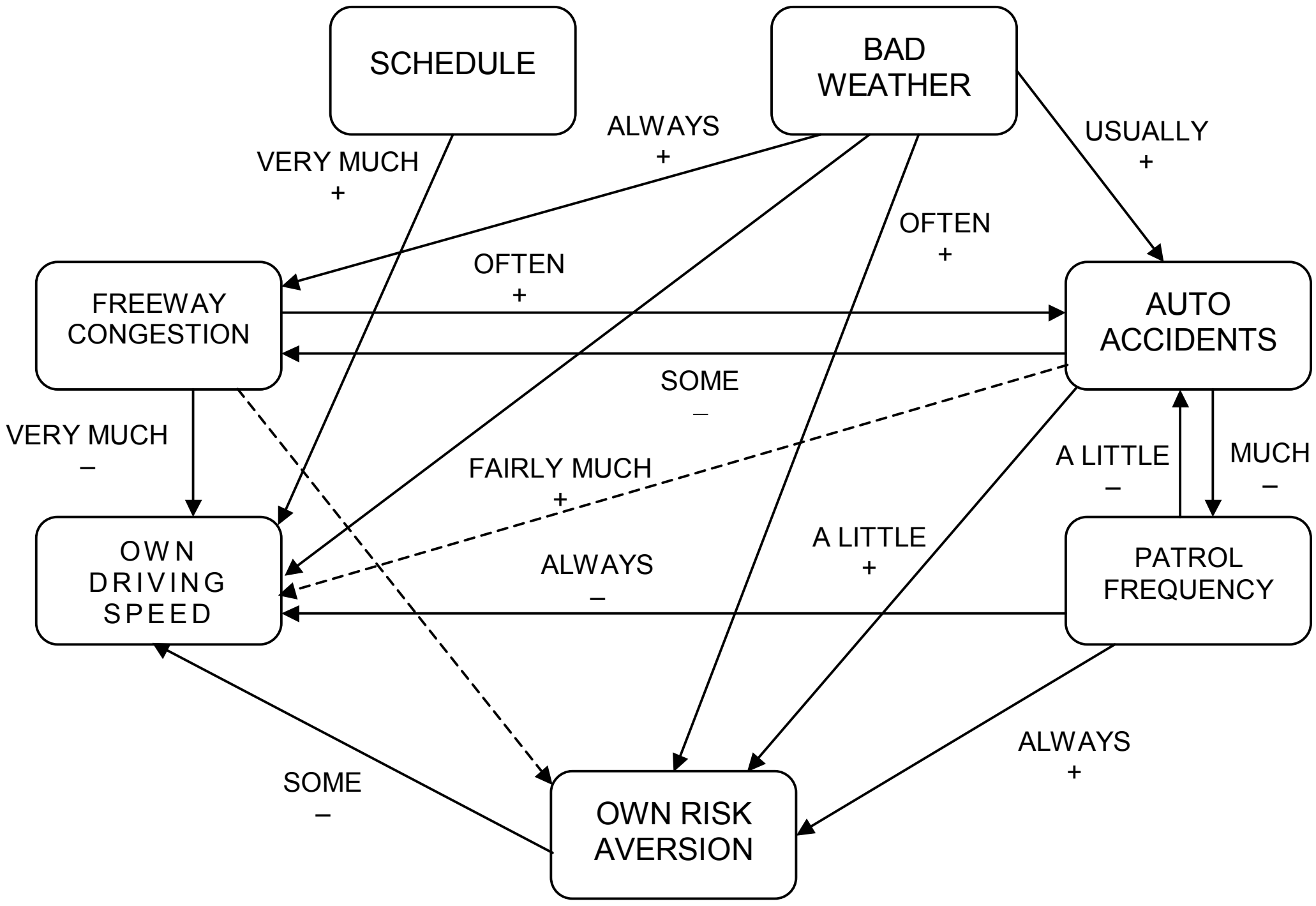

**FIGURE: 2.6.5**

As our main aim is to motivate researchers to use NCMs in place of FCMs whenever applicable and apply them to real world problems, we just give justification for the use of NCM and leave the work of constructing an NCM model to the reader. Even in this problem it is not only speed of one who drives but also several others factors like the speed of others, congestion etc. may or may not play an indeterminable role.

Accidents are very common in the countries like India, where many other factors like bad roads; reckless driving by others; drunken driving etc. wreak havoc on the number of accident deaths. We can say in conclusion that problems related to traffic and transportation can be very efficiently handled with tools like the FCM and the NCM.

Next we shall see how best NCM can be applied in medical diagnosis. FCMs have found applications in many medical diagnostics including symptom disease model (in homeopathy), studying the depression of terminally ill patients and studies like death wish of terminally ill patients, etc. Since we can have these models where indeterminacy can exist between two nodes we can as well apply NCM in the place of FCMs.

### 2.6.8: Application of NCM in Medical Diagnostics

In diagnostics of the depression psychiatrist has a lack of an objective knowledge about patient state. Except for very strong depressions no significant changes of the physiological parameters are measurable. Almost all information about patient mental state is obtained in verbal form on the basis of conversation with the patient. Due to this fact data are vague, uncertain, inconsistent and also indeterminable. The most popular modern approach to the etiopathogenesis of a depression is based on cognitive-behavioral theory of depression. This theory assumes that "the depression is caused by incorrect cognitive processing of information coming from an environment". Description of this problem and use of FCM is given by Vysoky [133] and we have dealt with it in an earlier section. We feel that use of NCM in the place of FCM will yield better results. This task of working with NCM to study medical diagnostics is left for the interested reader as a lot of indeterminacy is involved: for instance it is well known that certain types of food may trigger depression, meeting some people can cause depression, and several times mood-swings are connected with changes in the weather. Above all, gender, age, social status and the type of (other) diseases may also be factors, which help us to understand depression.

### 2.6.9: Neutrosophic Mechanisms for Causal Relations

There is causal relation between two given concepts whenever a change in one of these concepts affects the other one for example, there is a causal relation between Police vigilance and robbery, causal relations in causal maps always involve change or sometimes the change in causal relations in causal maps may be indeterminate/ vague/ uncertain so that we will be only in a position to measure the degrees of indeterminance in causal relations in causal maps.

When an indeterminate causality is present in an FCM we term it as an NCM and the degrees of indeterminate take values from [0, 1, I, –1].
The reader may try the neutrosophic mechanism for causal relations in a real world problem like intelligence and mark score of students in their school finals, social conditions and intelligence level, crime and socio-economic conditions etc.

### 2.6.10: Application of NCMs in Diagnosis and Study of Specific Language Impairment

Georgopoulous et al [37] have used FCMs in the diagnosis of specific language impairment and we have discussed it in brief in Section 1.4.11. For entire literature in this direction please refer [37]. Now in this study if we use the concept of NCMs i.e. indeterminacy of nodes and the corresponding neutrosophic matrix certainly the result would be more sensitive for the following reasons: The language impairment may at times be indeterminate when the child communicates with its parents, but refuses to communicate with strangers, or in cases where owing to fear the child may not communicate with teachers and in such cases of behavioural peculiarities, the degrees of indeterminacy become greater. Another factor is the language factor, in which the child is taught to communicate. For in India, owing to a multiplicity of languages, children often are brought up to speak a language at home which is different from the *lingua franca* of the world outside. In such cases, the relation between nodes may remain as indeterminate. Other such indeterminacy can also occur which may involve the circumstance and time.

Thus it is left for the reader to adopt NCM and model the specific language impairment problem and compare it when FCM is applied to the same problem by replacing I by 0.

### 2.6.11: The use of NCM approach in web mining inference

In Section 1.4.12 we have just sketched the FCM approach to webmining inference as given by K.C. Lee et al [63]. In conclusion of their paper, they have mentioned that "the basic nature of FCMs used for this study needs to be improved so the more complicated web-mining results can be analysed enough to yield meaningful inferences." At this juncture we would like to mention that if FCMs are replaced by NCMs then certainly we shall arrive at more effective, useful and sensitive conclusions and information as the inference would also depend on the indeterminates so the experts can be warned of the problems and deal with it in an appropriate manner. Any motivated reader can develop a new concept of web-mining procedures analogous to WEMIA in which FCMS are replaced by NCMs.

## 2.7 Neutrosophic Cognitive Maps versus Fuzzy Cognitive Maps

The great scientist Albert Einstein said, "*So far as the laws of mathematics refer to reality, they are not certain. And so far as they are certain, they do not refer to reality*".

Thus when we are analyzing the unsupervised data we cannot say anything for certain, at several times we are forced to face the indeterminacy of facts; so the only powerful tool which helps us to understand and apply the concept of indeterminacy in the analysis of data is the notion of neutrosophy. We have used in the place of fuzzy theory the concept of neutrosophy.

Here (in neutrosophy) we use the fact that between any two concepts/ nodes the existing relation may be an indeterminate (as) in reality, however as the notion of fuzzy cognitive maps do not help us to study indeterminacy in the analysis of the unsupervised data.

1. FCMs measure the existence of causal relation between two concepts and if no relation exists it is denoted by 0.

2. NCMs measure not only the existence of causal relation between concepts or the absence of causal relations between two concepts but also gives representation to the indeterminacy of relations between any two concepts.

3. We cannot apply NCMs for all unsupervised data. NCMs will have meaning only when relation between at least two concepts $C_i$ and $C_j$ are indeterminate.

4. The class of FCMs is strictly contained in the class of NCMs. (All NCMs can be made into FCMs by replacing I in the connection matrix by 0).

5. The directed graphs in case of NCMs are termed as neutrosophic graphs. i.e. in the graph we have at least two edges which are related by the dotted lines, means the edge between those two vertices is an indeterminate.

6. All connection matrices of the NCM are neutrosophic matrices i.e. they have in addition to the entries 0, 1, –1, the symbol I.

7. The resultant vectors i.e. the hidden pattern resulting in a fixed point or a limit cycle of a NCM can also be a neutrosophic vector; i.e. signifying the state of certain conceptual nodes of the system to be an indeterminate i.e. it is not off i.e. ‘0’ not on i.e.'1' and indeterminate relation is signified by I.

8. Because NCMs measure the indeterminate, the expert of the model can give due careful representation while implementing the results of the model.

9. In case of simple FCMs, we have the number of instantaneous state vectors to be the same as the number of resultant vectors but in case of NCMs we see the number of instantaneous state vectors is from the set {0, 1} where as the resultant vectors are from the bigger set {0, 1, I}. This is also one of the major differences between NCMs and FCMs.

<u>On Causal inferences in FCMs and NCMs:</u> Here we introduce from Yuan Miao and Zhi-Qiang Liu[136] some notions about FCMs using graph theory.

**DEFINITION [136]:** *A vertex function $f_{T_i}$ of $v_i$ is defined as*

$$x_i = f_{T_i}(\mu_i) = \begin{cases} 1, & \mu_i > T_i \\ 0, & \mu_i \leq T_i \end{cases}$$

*where $\mu_i$ is the total input of $v_i$, i.e. $\mu_i = \sum_k w_{ik}.x_k = W_i . x$ . $W_i$ is the $i^{th}$ row of W, $x = (x_1, \ldots ,x_n)^{Ti}$.*

**DEFINITION [136]:** *A path $P (v_1, \ldots, v_r)$ is a subdigraph of FCM, $\upsilon$ that contains k vertices $v_i$ $(1 \leq i \leq r;\ v_i \neq v_j$ if $i \neq j)$ and arcs $a_{i,i+1}$ $(1 \leq i < k)$ of $\upsilon$, where $a_{i,i+1}$ is an arc from vertex $v_i$ to $v_{i+1}$.*

**DEFINITION [136]:** *The length of the path is defined as the number of arcs it contains.*

$$L(P (v_1, \ldots, v_r)) = r.$$

*Similar to paths in graph theory, a path of the FCM does not allow any repeated vertex: $v_i \neq v_j$ if $i \neq j$, which is indicated in the earlier definition. However, for the sake of simplicity, we allow $v_1 = v_r$, i.e. in case of a circle, we consider $P (v_1, \ldots, v_{r-1}, v_1)$ as a path from $v_1$ to $v_1$ via $v_2, \ldots, v_{r-1}$. A path $P(v_1, \ldots, v_r)$ can also be denoted as $P_{v_1, v_r}$ or simply $P_{1,r}$ if it causes no confusion. In case of $v_r = v_1$, the path is simplified as $P_{v_1, v_1}$ or $P_{1,\ 1}$. Note that in FCM, we do not consider the case of multiple edges between two vertices, which is illustrated in Figure 2.7.1.*

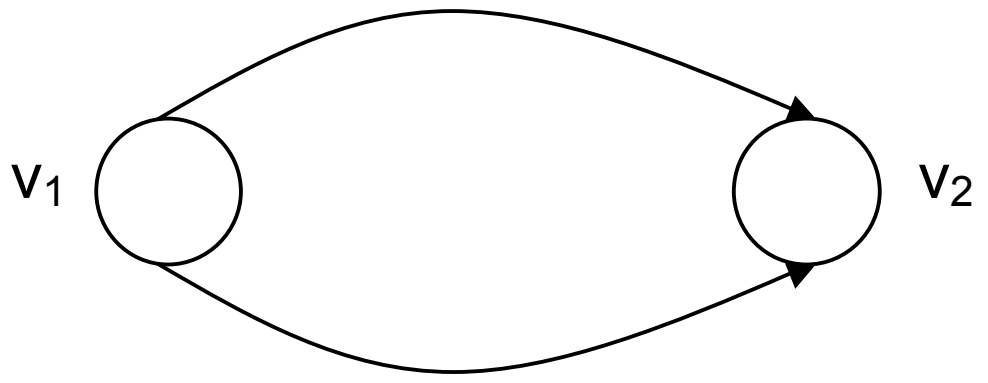

**FIGURE: 2.7.1**

*The path represents a basic structure of the FCM. As the FCM is different from a usual digraph; for instance, every vertex is capable of deciding how to respond to its causal inputs. It is useful and necessary to define different aspects of the path.*

**DEFINITION [136]:** *A path $P(v_1, \ldots, v_r)$ is called as a normal path if $v_i$ $(1 < i < r)$ has only one input.*

From this definition, it is now possible to define the total effect of a normal path.

**DEFINITION [136]:** *The function of a normal path $P(v_1, \ldots, v_r)$ is defined as*

$$f_{P_{1,k}} = f_{T_k, v_r}\left(\cdots\left(w_{3,2}.f_{T_2, v_2}\left(w_{2,1}.x_1\right)\cdots\right)\right)$$

*it $T_i = 0$, a normal path function can only be positive or negative and denoted as follows:*

$$f^1(\mu) = \begin{cases} 1, & \mu > 0 \\ 0, & \mu \leq 0 \end{cases}$$

$$f^{-1}(\mu) = \begin{cases} 0, & \mu > 0 \\ 1, & \mu \leq 0 \end{cases}$$

*if there is no path from $v_1$ to $v_r$ via $v_2, \ldots, v_{r-1}$, it is denoted as $P^{<>}(v_1, \ldots, v_r)$. Because a path transfers the causal effect from $v_1$ to $v_r$, we have the following function:*

$$f_{P^{<>}(v_1,\ldots,v_r)}(.) = 0$$

**DEFINITION [136]:** *A circle $C(v_1, \ldots, v_r)$ consists of a path $P(v_1, \ldots, v_r)$ and an arc pointing from $v_r$ to $v_1$.*

$$C(v_1, \ldots, v_r) = P(v_1, \ldots, v_r) \cup a_{1,r}.$$

An FCM without circles is relatively simple. It will be proved later that such FCMs will remain static after certain inference steps. On the other hand, an FCM with circles may have very complicated "hidden" patterns. Circles enable the feedback mechanism in FCM inference process, which will have very important consequences and will be useful in many real-world applications.

**DEFINITION [136]:** *A path $P(v_1, \ldots, v_r)$ is an input path of $\upsilon$ if*

i. *$v_1$ has no input or has only an external input sequence;*
ii. *$v_i$, $1 < i < r$ do not belong to any circle of $\upsilon$; $v_1$ is called an input vertex.*

The FCM is a causal model of the real system. The input vertex is where the external stimulus can affect the system and the input path is the channel that carries the causal effect. For example, in Fig 2.7.2, $P(v_2, v_1, v_{501})$ is an input path.

**DEFINITION [136]:** *An input path $P(v_1, \ldots, v_r)$ is a standardized input path if: 1) $v_i$, $(1 < i < k)$ has only one output; 2) $P(v_1, \ldots, v_r)$ is a normal path; and 3) $v_r$ is a vertex of a circle. Figure 2.7.2 illustrates the standardization of input paths. $P(v_{12}, v_2, v_4)$ and $P(v_{11}, v_2, v_3)$ in Figure 2.7.2(a) are all input paths but are not standardized input paths. Figure 2.7.2(b) shows the two input paths have been standardized input paths.*

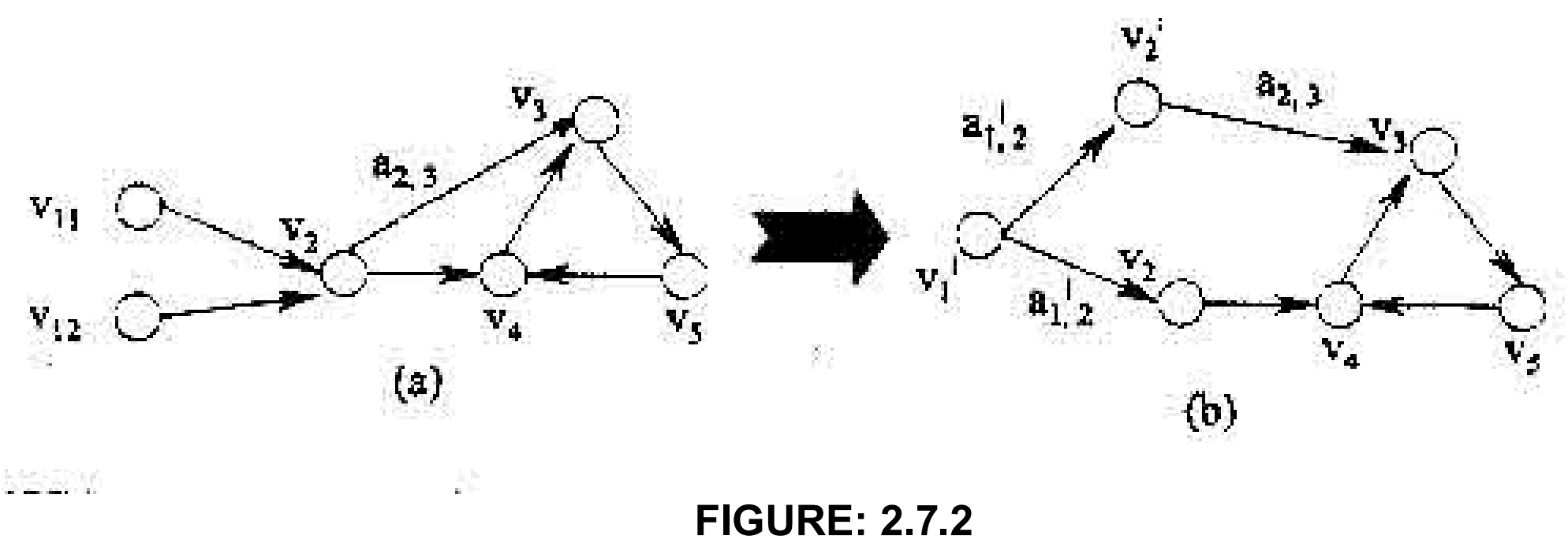

**FIGURE: 2.7.2**

Figure 2.7.2(b) shows that the two input paths have been standardized by inserting two dummy vertices $v'_1$ and $v'_2$ and $a'_{12}$ and $a'_{23}$. The new input paths $P(v'_1, v_2, v_4)$ and $P(v'_1, v'_2, v_3)$ are the standardized input paths. Obviously, any input path can be standardized by introducing some dummy vertices. After input path standardization, the FCM now has subgraphs containing circles, i.e., feedback loops. Input paths do not have feedback and can affect the other parts of the FCM. Input path standardization greatly simplifies the structure.

The FCM may also contain a part (or a subdigraph) that can be affected by other parts of the FCM, but has no effect on other parts. As far as the inference pattern of the FCM is concerned, such an affected part is passive and can be omitted without loss of generality.

**DEFINITION [136]:** *Suppose $AB(v_1,..., v_r)$ is a subdigraph of $\upsilon$ containing vertices $\{v_1,..., v_r\}$ and all arcs related to $\{v_1,..., v_r\}$. An edge $a_{i,j}$ is related to $\{v_1,..., v_r\}$ if either $v_i \in \{v_1,..., v_r\}$ or $v_j \in \{v_1,..., v_r\}$. $AB(v_1,..., v_r)$ is an affected branch of $\upsilon$ if it contains no circles and has no arcs pointing from any vertex in $\{v_1,..., v_r\}$ to any vertex in $V(\upsilon) - \{v_1,..., v_r\}$.*

An affected branch can be affected only by other parts of the FCM. As a result, the inference pattern of an affected branch is completely determined by the rest of the FCM and it has no effect on the rest of the FCM. Therefore, we can consider affected branches separately.

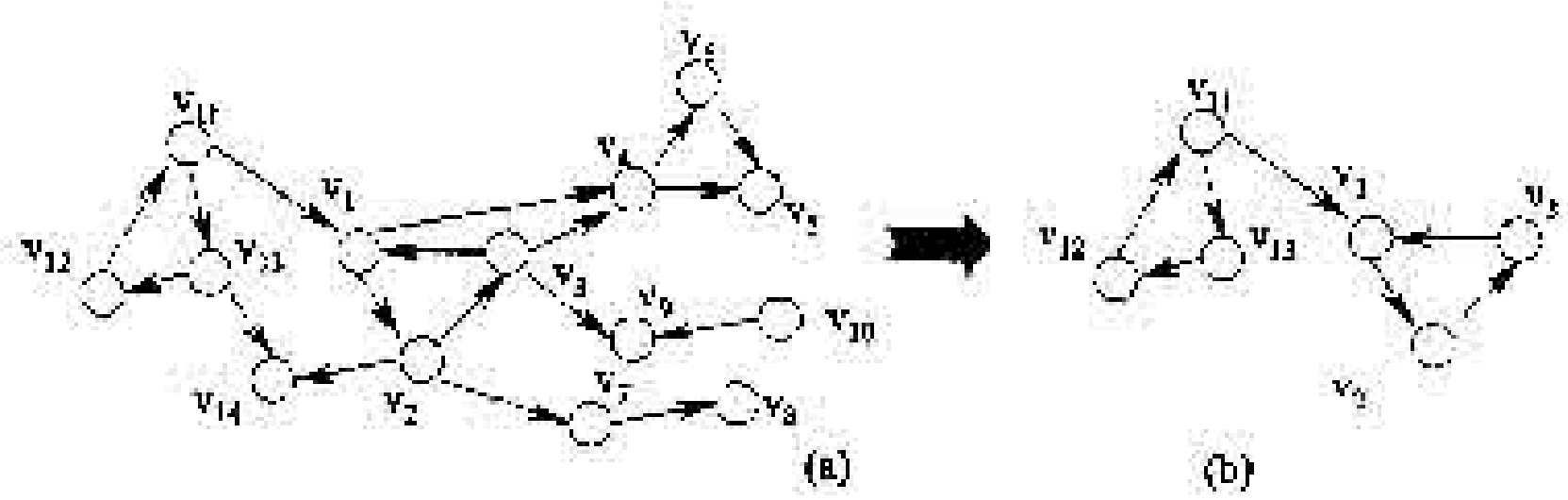

**FIGURE: 2.7.3**

Further, an affected branch is a trivial sub-FCM as it contains no circles. In the discussion that follows, they Yuan Miao and Zhi-Qiang Liu [136] first trim the affected branches from the FCM to make the in Figure 2.7.3 AB ($v_7$, $v_8$) is a trivial affected branch. AB ($v_4$, $v_5$, $v_6$) is an affected branch by two vertices that belong to the same circle. AB($v_9$, $v_{10}$) contains an input vertex. AB($v_{14}$) is an affected branch that is affected by two circles. Figure 2.7.3 is the FCM after all affected branches have been deleted.

FCM is powerful in representing structured knowledge and combining knowledge from different human experts. However, a general FCM may contain a large number of vertices with very complicated connections. It is difficult to be handled directly, if possible at all.

However, an FCM can be divided into several basic modules, which will be explicitly defined below. Every causal module is a smaller FCM. Vertices (or concepts) of a causal module infer each other and are closely connected. Basic FCM modules are the minimum FCM entities that cannot be divided further.

An FCM $\upsilon$ is "divided" as FCM $\upsilon_1$ and FCM $\upsilon_2$ if:

i. $V(\upsilon) = V(\upsilon_1) \cup V(\upsilon_2)$

ii. $A(\upsilon) = A(\upsilon_1) \cup A(\upsilon_2) \cup B(\upsilon_1, \upsilon_2)$

where $B(\upsilon_1, \upsilon_2)$

$$= \{a(\upsilon) \mid V^I(a(\upsilon)) \in V(\upsilon_1), V^o(a(\upsilon)) \in V(\upsilon_2) \text{ or } V^I(a(\upsilon)) \in V(\upsilon_2), V^o(a(\upsilon)) \in V(\upsilon_1)$$

$$\begin{aligned} V(\upsilon_1) \cap V(\upsilon_2) &= \varnothing \\ A(\upsilon_1) \cap A(,\upsilon_2) &= \varnothing \\ A(\upsilon_1) \cap B(\upsilon_1,\upsilon_2) &= \varnothing \\ A(\upsilon_2) \cap B(\upsilon_1,\upsilon_2) &= \varnothing. \end{aligned}$$

This operation is denoted as

$$\upsilon = \upsilon_1 \leftarrow \frac{B(\upsilon_1,\upsilon_2)}{\oplus} \rightarrow \upsilon_2 .$$

Particularly, we consider the causal relationships in subgraphs
$B(\upsilon_1, \upsilon_2)$

$$= \{a(\upsilon) \mid V^I(a(\upsilon)) \in V(\upsilon_1), V^0(a(\upsilon)) \in V(\upsilon_2)\}$$

this case is described as "$\upsilon_2$ is caused by $\upsilon_1$ " or "$\upsilon_1$ causes $\upsilon_2$." Such a division is called a regular division denoted as

$$\upsilon = \upsilon_1 \frac{B(\upsilon_1,\upsilon_2)}{\oplus} \rightarrow \upsilon_2 .$$

**DEFINITION [136]:** *An FCM containing no circles is called a trivial FCM. A trivial FCM has no feedback mechanism and its inference pattern is usually trivial.*

**DEFINITION [136]:** *A basic FCM is an FCM containing at least one circle, but cannot be regularly divided into two FCMs, both of which contain at least one circle.*

From the definitions, we can see that an FCM containing circles can always be regularly divided into basic FCMs. In general, the basic FCMs of $\upsilon$ can be ordered concatenatedly as follows:

$$\upsilon = \upsilon_1 \frac{B(\upsilon_1,\upsilon_2)}{\oplus} \rightarrow \upsilon_2 \frac{B(\upsilon_2,\upsilon_3)}{\oplus}$$

$$\rightarrow \upsilon_3 \ldots \frac{B(\upsilon_{m-1},\upsilon_m)}{\oplus} \rightarrow \upsilon_m .$$

We establish a formal result in the following theorem:

**THEOREM [136]:** *Suppose FCM $\upsilon$ can be regularly divided into m basic FCM's, then*

$$\upsilon = \left( \bigcup_{i=2}^{m} \upsilon_i \right) \cup \left( \bigcup_{i=2}^{m} \bigcup_{j=1}^{i} B(\upsilon_j, \upsilon_i) \right)$$

*where*

$$V(\upsilon_i) \cap V(\upsilon_j) = \varnothing \quad i \neq j$$
$$A(\upsilon_i) \cap A(,\upsilon_j) = \varnothing \quad i \neq j$$

$B(\upsilon_i \upsilon_j)$

$$= \{a(\upsilon) \,/ V^1(a(\upsilon)) \in V(\upsilon_i),\, V^0(a(\upsilon)) \in V(\upsilon_j)\}.$$

**THEOREM [136]:** *Suppose that FCM $\upsilon$ is input-path standardized and trimmed of affected branches. If a vertex $\upsilon_o$ of $\upsilon$ has at least two input arcs and does not belong to any circle, then $\upsilon$ is not a basic FCM.*

If an FCM is a trivial FCM, it will become static after L inference interactions unless it has an external input sequence, where L is the length of the longest path of the FCM.

The proof is obvious. Consequently, the following is true:

*Vertices except the end vertex of an input path will become static after L inference iterations unless it has an external input sequence, where L is the length of the path.*

**DEFINITION [136]:** *A vertex is called as a key vertex if:*

*i. It is common vertex of an input path and circle or*
*ii. It is a common vertex of two circles with at least two arcs pointing to it, which belongs to the two circles or*
*iii. It is any vertex on a circle if the circle contains no other key vertices.*

*If an FCM contains no key vertices, it contains no circle, i.e., it is a trivial FCM. So it will remain static after a certain inference steps according to the above theorem.*

Let P $\{v_1,\ldots, v_r\}$. be a normal path. If $w_{i+1,i} < 0$, then

$$x_{i+1}(k+1) = f_{i+1}(w_{i+1,i} \,.\, x_i(k))$$

$$= \begin{cases} f_{i+1}(0) = 0 & x_i(k) = 0 \\ f_{i+1}(w_{i+1,i}.1) = 0 & x_i(k) = 1 \end{cases} \quad k \geq 0$$

$$x_{i+2}(k+2) = f_{i+2}(w_{i+2,I+1} \,.\, x_{i+1}(k+1)) = f_{i+2}(0) = 0.$$

i.e., $x_{i+1} \equiv 0,\ldots, x_r \equiv 0$. This means that after certain steps of inference, $x_{i+1},\ldots, x_r$ will not affect the other concepts of the FCM. We can delete $x_{i+1},\ldots x_r$ from the FCM for further analysis and then there is no $P(v_1,\ldots, v_r)$ . A normal path containing a negative arc is a path that cannot pass causal flow after (at most) r inference steps and is denoted as $f_{-r}$ (.) = 0.

In case that an FCM is an isolated circle, C ($\nu_1, \ldots, \nu_r$). If there exists $w_{i+1,I} < 0$, the state of the FCM will remain zero after (at most) r inference steps:

If j > i ,

$$f_{P_{i,j}}(.) = f_{-(j-i)(.)} = 0$$

If j < i ,

$$f_{P_{i,j}}(.) = f_1(w_{1,r}.x_r(.))$$

$$= f_{P_{i,j}}(w_{1,r}.f_{i,r}(.)) = f_{P_{i,j}}(w_{1,r}.f_{-(r-i)}(.))$$

$$= f_{-r-I+1+j-1}(.) = f_{-r-I+j}(.) = 0$$

If all $w_{I+1,i} > 0$, select $\upsilon_1$ as the key vertex. Thus, we have $x_1(k + r) = f_{p_{1,1}}(x_1(k + r\text{-}r)$, k = 1, 2, …, where $f_{p_{1,1}}$ can only be $f^1$. Then the circle is a positive feedback.

$$x_1(k + r) = x_1(k + r\text{-}r) = x_1(k).$$

So the state of the FCM will finally be trapped in a limit cycle. The period of the limit cycle is r.

There are several more results related to FCM as discussed by [136]. But we request the reader to refer them. Once again our main motivation is to express relationships more relating to NCMs. We just mention the new definitions and results on FCMs in passing and hint whenever possible how to adopt them in NCMs. We don't provide each and every analogue for the transition between NCMs and FCMs.

The causal inferences in FCM as given by [136] is described below. We now give the analogous procedure for NCM.

We assume the existence of FCMs. Now we recall the definition of neutrosophic digraph in which vertices representing concepts and arcs between the vertices indicate causal relationships between the concepts. $N(W^T)$ denotes the adjacency neutrosophic matrix. In case of FCMs with n attributes or nodes we have utmost $2^n$ states. Thus after $2^n$ inferences steps, the state is either static state or is trapped in a limit circle. But in case of NCM even if n = 2, we get the resultant set will take its value from the set {0, 1, I}.

Thus in case of FCMs we have 4 choices of resultant state vectors on the other hand in case of NCM we have nine resultant vectors given by {(0, 0), (1, 0), (0, 1), (1, 1), (0, I), (I, 0), (1, I), (I, 1), (I, I)}. Thus when we work with n = 3 we have still more choices {(0, 0, 0), (1, 0, 0), (0, 1, 0), (0, 0, 1), (1, 1, 0), (1, 0, 1), (0, 1, 1), (1, 1, 1),

(I, 0, 0), (0, I, 0), (0, 0, I), (0, I, I), (I, 0, I), (I, I, 0), (I, 1, 1), (1, 1, I), (1, I, 1), (I, I, 1), (I, 1, I), (1, I, I), (I, I, I), (1, 0, I), (1, I, 0), (0, 1, I), (0, I, 1), (I, 0, 1), (I, 1, 0)}.

Thus we see many more resultant vectors than the instantaneous vectors. Almost all analogous definitions and results as in case of FCMs are found in the paper [136]. A NCM is divided as $(NCM)_1$ and $(NCM)_2$ if vertices V of NCM is divided into vertices of $V_i$ related to $(NCM)_i$, i = 1, 2 i.e. $V = V_1 \cup V_2$.

Suppose $v_1, v_2, \ldots, v_n$ are attributes/ nodes a circle $c(v_1, v_2, \ldots, v_r)$ consists of a path $P(v_1, v_2, \ldots, v_r)$ and an arc pointing from $v_r$ to $v_1$, $c(v_1, v_2, \ldots, v_r) = P(v_1, \ldots, v_r) \cup (\overrightarrow{v_r v_1})$.

Now in case of NCMs the path $P(v_1, \ldots, v_r)$ can also be an indeterminate paths; a circle compromising of indeterminate paths is called the neutrosophic circle. A NCM without neutrosophic circle is relatively simple. It is left for the reader to prove that such NCMs will remain static after certain inference steps.

Just as in case of FCMs we may also have in NCMs a part (or a sub-digraph) that can be affected by other parts of the NCM, but has no effect on other parts. As far as inference pattern of the NCM is concerned such an affected part is passive and can be omitted without loss of generality.

An affected branch can be affected only by other parts of the NCM. As a result the inference pattern of an affected branch is completely determined by the rest of the NCM and it has no effect on the rest of NCM. Thus in NCMs also consider the affected branches separately. Further an affected branch is a trivial sub-NCM as it contains no circles.

The reader is expected to construct real world models using NCMs and study the affected parts of the NCM. As in case of FCMs we will call a NCM to be basic if the NCM contains at least one circle but cannot be regularly divided into two NCMs both of which contain at least one circle. Keeping in mind an NCM is trivial if it has no circles we can see that NCM containing circles can always be regularly divided into basic NCMs. The reader is given the work of studying/ analyzing the concept of causal inferences in NCM.

## 2.8 Neutrosophic Relational Maps — Definition with Examples

When the nodes or concepts under study happens to be such that they can be divided into two disjoint classes and a study or analysis can be made using Fuzzy Relational Maps (FRMs) was introduced and described in the earlier chapter. Here we define a new concept called Neutrosophic Relational Maps (NRMs), analyse and study them. We also give examples of them.

**DEFINITION 2.8.1:** *Let D be the domain space and R be the range space with $D_1, \ldots, D_n$ the conceptual nodes of the domain space D and $R_1, \ldots, R_m$ be the conceptual nodes of the range space R such that they form a disjoint class i.e. $D \cap R = \phi$. Suppose there* *is a FRM relating D and R and if at least a edge relating a $D_i$ $R_j$ is an indeterminate then we call the FRM as the Neutrosophic relational maps. i.e. NRMs.*

***Note:*** In everyday occurrences we see that if we are studying a model built using an unsupervised data we need not always have some edge relating the nodes of a domain space and a range space or there does not exist any relation between two nodes, it can very well happen that for any two nodes one may not be always in a position to say that the existence or nonexistence of a relation, but we may say that the relation between two nodes is an indeterminate or cannot be decided.

Thus to the best of our knowledge indeterminacy models can be built using neutrosophy. One model already discussed is the Neutrosophic Cognitive Model. The other being the Neutrosophic Relational Maps model, which are a further generalization of Fuzzy Relational Maps.

It is not essential when a study/ prediction/ investigation is made we are always in a position to find a complete answer. This is not always possible (sometimes or many a times) it is almost all models built using unsupervised data, we may have the factor of indeterminacy to play a role. Such study is possible only by using the Neutrosophic logic.

***Example 2.8.1:*** Female infanticide (the practice of killing female children at birth or shortly thereafter) is prevalent in India from the early vedic times, as women were (and still are) considered as a property. As long as a woman is treated as a property/ object the practice of female infanticide will continue in India.

In India, social factors play a major role in female infanticide. Even when the government recognized the girl child as a critical issue for the country's development, India continues to have an adverse ratio of women to men. Other reasons being torture of the in-laws may also result in cruel death of a girl child. This is mainly due to the fact that men are considered superior to women. Also they take into account that fact that men are breadwinners for the family. Even if women work like men, parents think that her efforts is going to end once she is married and enters a new family.

Studies have consistently shown that girl babies in India are denied the same and equal food and medical care that the boy babies receive. Girl babies die more often than boy babies even though medical research has long ago established that girls are generally biologically stronger as new-borns than boys. The birth of a male child is a time for celebration, but the birth of female child is often viewed as a crisis. Thus the female infanticide cannot be attributed to single reason it is highly dependent on the feeling of individuals ranging from social stigma, monetary waste, social status etc.

Suppose we take the conceptual nodes for the unsupervised data relating to the study of female infanticide. We take the status of the people as the domain space D

$D_1$ – very rich
$D_2$ – rich
$D_3$ – upper middle class
$D_4$ – middle class
$D_5$ – lower middle class

$D_6$ – Poor
$D_7$ – Very poor.

The nodes of the range space R are taken as

$R_1$ – Number of female children - a problem
$R_2$ – Social stigma of having female children
$R_3$ – Torture by in-laws for having only female children
$R_4$ – Economic loss / burden due to female children
$R_5$ – Insecurity due to having only female children (They will marry and enter different homes thereby leaving their parents, so no one would be able to take care of them in later days.)

Keeping these as nodes of the range space and the domain space experts opinion were drawn which is given the following Figure 2.8.1:

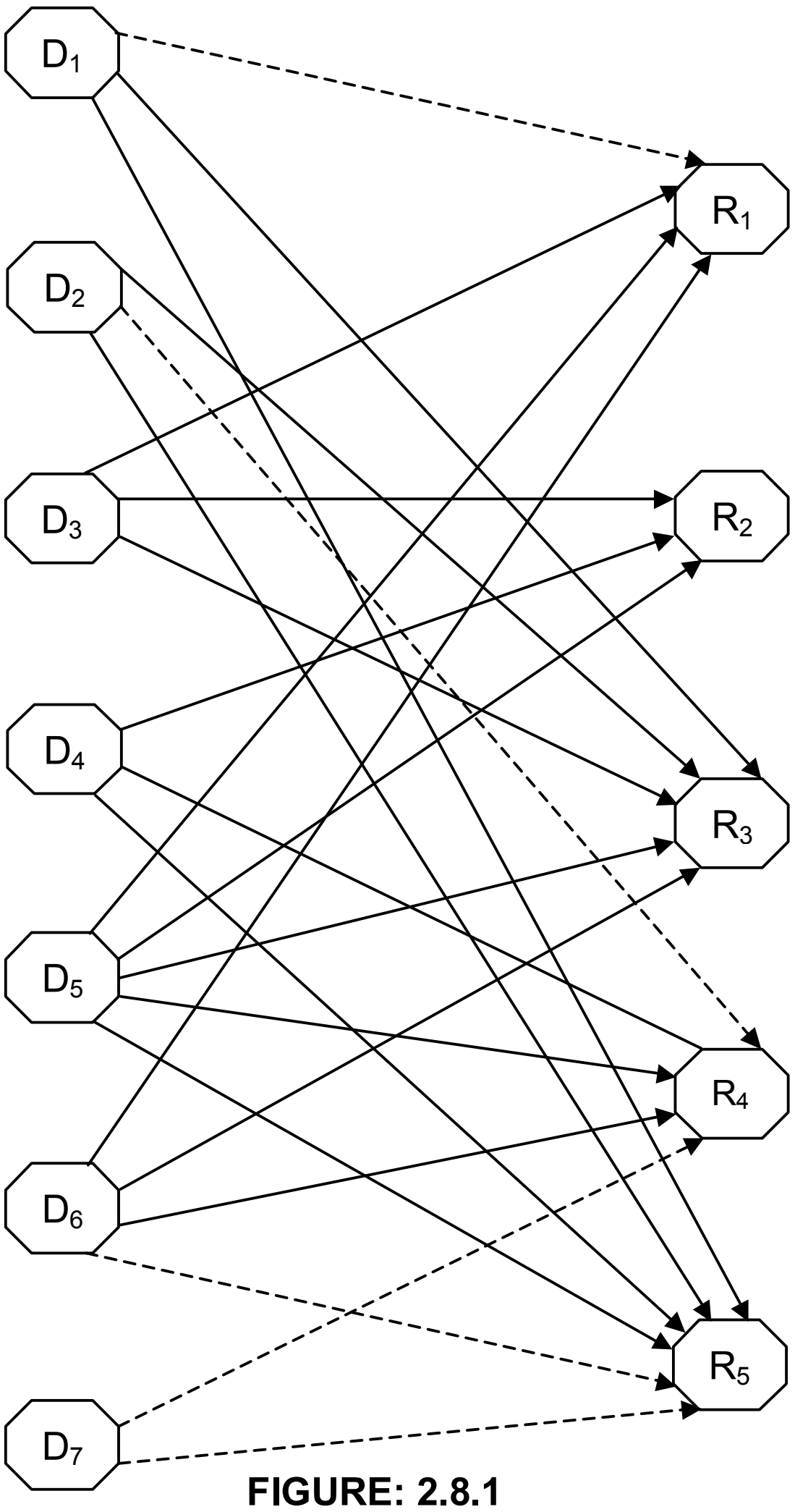

**FIGURE: 2.8.1**

Figure 2.8.1 is the neutrosophic directed graph of the NRM.

The corresponding neutrosophic relational matrix $N(R)^T$ is given below:

$$N(R)^T = \begin{bmatrix} I & 0 & 1 & 0 & 1 & 1 & 0 \\ 0 & 0 & 1 & 1 & 1 & 0 & 0 \\ 1 & 1 & 1 & 1 & 1 & 1 & 0 \\ 0 & I & 0 & 0 & 1 & 1 & I \\ 1 & 1 & 0 & 1 & 1 & I & I \end{bmatrix}$$

and

$$N(R) = \begin{bmatrix} I & 0 & 1 & 0 & 1 \\ 0 & 0 & 1 & I & 1 \\ 1 & 1 & 1 & 0 & 0 \\ 0 & 1 & 1 & 0 & 1 \\ 1 & 1 & 1 & 1 & 1 \\ 1 & 0 & 1 & 1 & I \\ 0 & 0 & 0 & I & I \end{bmatrix}.$$

Suppose $A_1$ = (0 1 0 0 0) is the instantaneous state vector under consideration i.e., social stigma of having female children.

The effect of $A_1$ on the system N(R) is

| | | | | | | |
|---|---|---|---|---|---|---|
| $A_1N(R)^T$ | = | (0 0 1 1 1 0 0) | → | (0 0 1 1 1 0 0) | = | $B_1$ |
| $B_1[N(R)]$ | = | (2, 3, 3, 1, 2) | → | (1 1 1 1 1) | = | $A_2$ |
| $A_2[N(R)]^T$ | = | (1 + I, 1 + I, 1, 1, 1, 1 + I, I) | → | (1 1 1 1 1 1 I) | = | $B_2$ |
| $B_2(N(R))$ | = | (I + 3, 3, 5, 2I + 2, 2I + 4) | → | (1 1 1 1 1) | = $A_3$ = | $A_2$. |

Thus this state vector $A_1$ = (0, 1, 0, 0, 0) gives a fixed point (1 1 1 1 1 1) indicating if one thinks that having female children is a social stigma immaterial of their status they also feel that having number of female children is a problem, it is a economic loss / burden, they also under go torture or bad treatment by in-laws and ultimately it is a insecurity for having only female children, the latter two cases hold where applicable.

On the other hand we derive the following conclusions on the domain space when the range space state vector $A_1$ = (0 1 0 0 0) is sent

| | | | | |
|---|---|---|---|---|
| $A_1[N(R)]^T$ | → | $B_1$ | | |
| $B_1[N(R)]$ | → | $A_2$ | | |
| $A_2\,[N(R)]^T$ | → | $B_2$ | | |
| $B_2[N(R)]$ | → | $A_3$ | = | $A_2$ so |
| $A_2\,[N(R)]^T$ | → | $B_2$ = (1, 1, 1, 1, 1, 1, I) | | |

leading to a fixed point. When the state vector $A_1$ = (0, 1, 0, 0, 0) is sent to study i.e. the social stigma node is on uniformly all people from all economic classes are awakened expect the very poor for the resultant vector happens to be a Neutrosophic vector hence one is not in a position to say what is the feeling of the very poor people and the "many female children are a social stigma" as that coordinate remains as an indeterminate one. This is typical of real-life scenarios, for the working classes hardly distinguish much when it comes to the gender of the child.

Several or any other instantaneous vector can be used and its effect on the Neutrosophical Dynamical System can be studied and analysed. This is left as an exercise for the reader. Having seen an example and application or construction of the NRM model we will proceed on to describe the concepts of it in a more mathematical way.

**DESCRIPTION OF A NRM:**

Neutrosophic Cognitive Maps (NCMs) promote the causal relationships between concurrently active units or decides the absence of any relation between two units or the indeterminance of any relation between any two units. But in Neutrosophic Relational Maps (NRMs) we divide the very causal nodes into two disjoint units. Thus for the modeling of a NRM we need a domain space and a range space which are disjoint in the sense of concepts. We further assume no intermediate relations exist within the domain and the range spaces. The number of elements or nodes in the range space need not be equal to the number of elements or nodes in the domain space.

Throughout this section we assume the elements of a domain space are taken from the neutrosophic vector space of dimension n and that of the range space are neutrosophic vector space of dimension m. (m in general need not be equal to n). We denote by R the set of nodes $R_1,\ldots, R_m$ of the range space, where $R = \{(x_1,\ldots, x_m) \mid x_j = 0 \text{ or } 1 \text{ for } j = 1, 2, \ldots, m\}$.

If $x_i = 1$ it means that node $R_i$ is in the on state and if $x_i = 0$ it means that the node $R_i$ is in the off state and if $x_i = I$ in the resultant vector it means the effect of the node $x_i$ is indeterminate or whether it will be off or on cannot be predicted by the neutrosophic dynamical system.

It is very important to note that when we send the state vectors they are always taken as the real state vectors for we know the node or the concept is in the on state or in the off state but when the state vector passes through the Neutrosophic dynamical system some other node may become indeterminate i.e. due to the presence of a node we may not be able to predict the presence or the absence of the other node i.e., it is indeterminate, denoted by the symbol I, thus the resultant vector can be a neutrosophic vector.

**DEFINITION 2.8.2:** *A Neutrosophic Relational Map (NRM) is a Neutrosophic directed graph or a map from D to R with concepts like policies or events etc. as nodes and causalities as edges. (Here by causalities we mean or include the indeterminate causalities also). It represents Neutrosophic Relations and Causal Relations between spaces D and R .*

*Let $D_i$ and $R_j$ denote the nodes of an NRM. The directed edge from $D_i$ to $R_j$ denotes the causality of $D_i$ on $R_j$ called relations. Every edge in the NRM is weighted with a number in the set $\{0, +1, -1, I\}$. Let $e_{ij}$ be the weight of the edge $D_i R_j$, $e_{ij} \in \{0, 1, -1, I\}$. The weight of the edge $D_i R_j$ is positive if increase in $D_i$ implies increase in $R_j$ or* *decrease in $D_i$ implies decrease in $R_j$ i.e. causality of $D_i$ on $R_j$ is 1. If $e_{ij} = -1$ then increase (or decrease) in $D_i$ implies decrease (or increase) in $R_j$. If $e_{ij} = 0$ then $D_i$ does not have any effect on $R_j$. If $e_{ij} = I$ it implies we are not in a position to determine the effect of $D_i$ on $R_j$ i.e. the effect of $D_i$ on $R_j$ is an indeterminate so we denote it by I.*

**DEFINITION 2.8.3:** *When the nodes of the NRM take edge values from {0, 1, –1, I} we say the NRMs are simple NRMs.*

**DEFINITION 2.8.4:** *Let $D_1$, ..., $D_n$ be the nodes of the domain space D of an NRM and let $R_1$, $R_2$,..., $R_m$ be the nodes of the range space R of the same NRM. Let the matrix N(E) be defined as N(E) = ($e_{ij}$) where $e_{ij}$ is the weight of the directed edge $D_i R_j$ (or $R_j D_i$) and $e_{ij} \in$ {0, 1, –1, I}. N(E) is called the Neutrosophic Relational Matrix of the NRM.*

The following remark is important and interesting to find its mention in this book.

**Remark**: Unlike NCMs, NRMs can also be rectangular matrices with rows corresponding to the domain space and columns corresponding to the range space. This is one of the marked difference between NRMs and NCMs. Further the number of entries for a particular model which can be treated as disjoint sets when dealt as a NRM has very much less entries than when the same model is treated as a NCM.

Thus in many cases when the unsupervised data under study or consideration can be spilt as disjoint sets of nodes or concepts; certainly NRMs are a better tool than the NCMs.

**DEFINITION 2.8.5:** *Let $D_1$, ..., $D_n$ and $R_1$,..., $R_m$ denote the nodes of a NRM. Let A = ($a_1$,..., $a_n$), $a_i \in$ {0, 1, –1} is called the Neutrosophic instantaneous state vector of the domain space and it denotes the on-off position of the nodes at any instant. Similarly let B = ($b_1$,..., $b_n$) $b_i \in$ {0, 1, –1}, B is called instantaneous state vector of the range space and it denotes the on-off position of the nodes at any instant, $a_i = 0$ if $a_i$ is off and $a_i = 1$ if $a_i$ is on for i = 1, 2, ..., n. Similarly, $b_i = 0$ if $b_i$ is off and $b_i = 1$ if $b_i$ is on for i = 1, 2,..., m.*

**DEFINITION 2.8.6:** *Let $D_1$,..., $D_n$ and $R_1$, $R_2$,..., $R_m$ be the nodes of a NRM. Let $D_i R_j$ (or $R_j D_i$) be the edges of an NRM, j = 1, 2,..., m and i = 1, 2,..., n. The edges form a directed cycle. An NRM is said to be a cycle if it possess a directed cycle. An NRM is said to be acyclic if it does not possess any directed cycle.*

**DEFINITION 2.8.7:** *A NRM with cycles is said to be a NRM with feedback.*

**DEFINITION 2.8.8:** *When there is a feedback in the NRM i.e. when the causal relations flow through a cycle in a revolutionary manner the NRM is called a Neutrosophic dynamical system.*

**DEFINITION 2.8.9:** *Let $D_i R_j$ (or $R_j D_i$) $1 \leq j \leq m$, $1 \leq i \leq n$, when $R_j$ (or $D_i$) is switched on and if causality flows through edges of a cycle and if it again causes $R_j$ (or $D_i$) we say that the Neutrosophical dynamical system goes round and round. This is true for any node $R_j$ ( or $D_i$ ) for $1 \leq j \leq m$ (or $1 \leq i \leq n$). The equilibrium state of this Neutrosophical dynamical system is called the Neutrosophic hidden pattern.*

**DEFINITION 2.8.10:** *If the equilibrium state of a Neutrosophical dynamical system is a unique Neutrosophic state vector, then it is called the fixed point. Consider an NRM with $R_1, R_2, \ldots, R_m$ and $D_1, D_2, \ldots, D_n$ as nodes. For example let us start the dynamical system by switching on $R_1$ (or $D_1$). Let us assume that the NRM settles down with $R_1$ and $R_m$ (or $D_1$ and $D_n$) on, or indeterminate on, i.e. the Neutrosophic state vector remains as (1, 0, 0,..., 1) or (1, 0, 0,...I) (or (1, 0, 0,...1) or (1, 0, 0,...I) in D), this state vector is called the fixed point.*

**DEFINITION 2.8.11:** *If the NRM settles down with a state vector repeating in the form $A_1 \to A_2 \to A_3 \to \ldots \to A_i \to A_1$ (or $B_1 \to B_2 \to \ldots \to B_i \to B_1$) then this equilibrium is called a limit cycle.*

**Methods of determining the Hidden pattern in a NRM**

Let $R_1, R_2, \ldots, R_m$ and $D_1, D_2, \ldots, D_n$ be the nodes of a NRM with feedback. Let N(E) be the Neutrosophic Relational Matrix. Let us find the hidden pattern when $D_1$ is switched on i.e. when an input is given as a vector; $A_1 = (1, 0, \ldots, 0)$ in D; the data should pass through the relational matrix N(E). This is done by multiplying $A_1$ with the Neutrosophic relational matrix N(E). Let $A_1N(E) = (r_1, r_2, \ldots, r_m)$ after thresholding and updating the resultant vector we get $A_1E \in R$, Now let $B = A_1E$ we pass on B into the system $(N(E))^T$ and obtain $B(N(E))^T$. We update and threshold the vector $B(N(E))^T$ so that $B(N(E))^T \in D$.

This procedure is repeated till we get a limit cycle or a fixed point.

**DEFINITION 2.8.12:** *Finite number of NRMs can be combined together to produce the joint effect of all NRMs. Let $N(E_1), N(E_2), \ldots, N(E_r)$ be the Neutrosophic relational matrices of the NRMs with nodes $R_1, \ldots, R_m$ and $D_1, \ldots, D_n$, then the combined NRM is represented by the neutrosophic relational matrix $N(E) = N(E_1) + N(E_2) + \ldots + N(E_r)$.*

Now we give a simple illustration of a NRM.

***Example 2.8.2:*** Now consider the example given in the first chapter, Section 5. We take $D_1, D_2, \ldots, D_5$ and the $R_1$, $R_2$ and $R_3$ as in Example 1.5.1:

| | | |
|---|---|---|
| $D_1$ | – | Teacher is good |
| $D_2$ | – | Teaching is poor |
| $D_3$ | – | Teaching is mediocre |
| $D_4$ | – | Teacher is kind |
| $D_5$ | – | Teacher is harsh (or Rude) |

$D_1, \ldots, D_5$ are taken as the 5 nodes of the domain space, we consider the following 3 nodes to be the nodes of the range space.

| | | |
|---|---|---|
| $R_1$ | – | Good student |
| $R_2$ | – | Bad student |
| $R_3$ | – | Average student. |

The Neutrosophic relational graph of the teacher – student model is as follows:

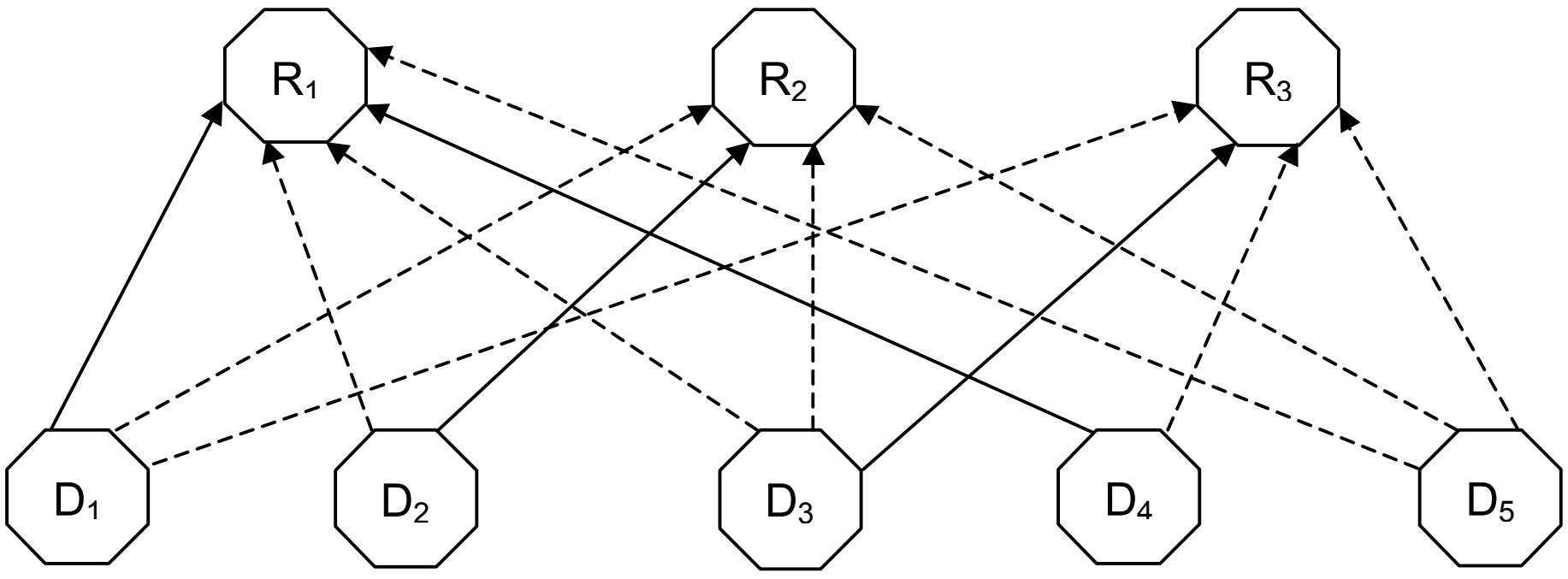

**FIGURE: 2.8.2**

$$N(E) = \begin{bmatrix} 1 & I & I \\ I & 1 & 0 \\ I & I & 1 \\ 1 & 0 & I \\ I & I & I \end{bmatrix}.$$

If $A_1 = (1, 0, 0, 0, 0)$ is taken as the instantaneous state vector and is passed on in the relational matrix N(E), $A_1N(E) = (1, I, I) = A_2$. Now

| | | | | | | |
|---|---|---|---|---|---|---|
| $A_2(N(E))^T$ | = | (1 + I, 1 + I, I 1 +I I) | → | (1, 1, I, 1, I) | = | $B_1$ |
| $B_1N$ (E) | = | (2+I, I+1, I) | → | (1 I I) | = | $A_3$ |
| $A_3N(E)$ | = | (1+I, I, I, 1 + I, I) | → | (1 I I 1 I) = $B_2 = B_1$. | | |
| $B_1N(E)$ | = | (1 I I) | | | | |

Thus we see from the NRM given that if the teacher is good it implies it produces good students but nothing can be said about bad and average students. The bad and average students remain as indeterminates. On the other hand in the domain space when the teacher is good the teaching quality of her remains indeterminate therefore both the nodes teaching is poor and teaching is mediocre remains as indeterminates but the node teacher is kind becomes in the on state and the teacher is harsh is an indeterminate, (for harshness may be present depending on the circumstances).

All other combinations of state vectors can be substituted and the reader can derive results. We give another example of the NRM.

***Example 2.8.3:*** Now we build a NRM model using the socio and psychological feeling of the AIDS / HIV patients and the feelings of the public towards the AIDS / HIV patients. As it is an example we have taken only few nodes related to HIV / AIDS patients as the domain space and that of the feelings of the public as the range space.

The nodes / concepts related with HIV / AIDS patients which is taken as the domain space D:

$D_1$ – Feeling of loneliness / aloofness

| | | |
|---|---|---|
| $D_2$ | – | Feeling of Guilt |
| $D_3$ | – | Desperation / fear in public |
| $D_4$ | – | Sufferings both mental / physical |
| $D_5$ | – | Public disgrace (feeling). |

The concepts / nodes related with public taken as the nodes of the range space

| | | |
|---|---|---|
| $R_1$ | – | Fear of getting the disease |
| $R_2$ | – | No mind to forgive the AIDS / HIV patient's "sin" |
| $R_3$ | – | Social Stigma to have HIV / AIDS patient as friend |
| $R_4$ | – | No sympathy. |

Now the NRM is obtained using an experts opinion whose Neutrosophic directed graph is given in Figure 2.8.3.

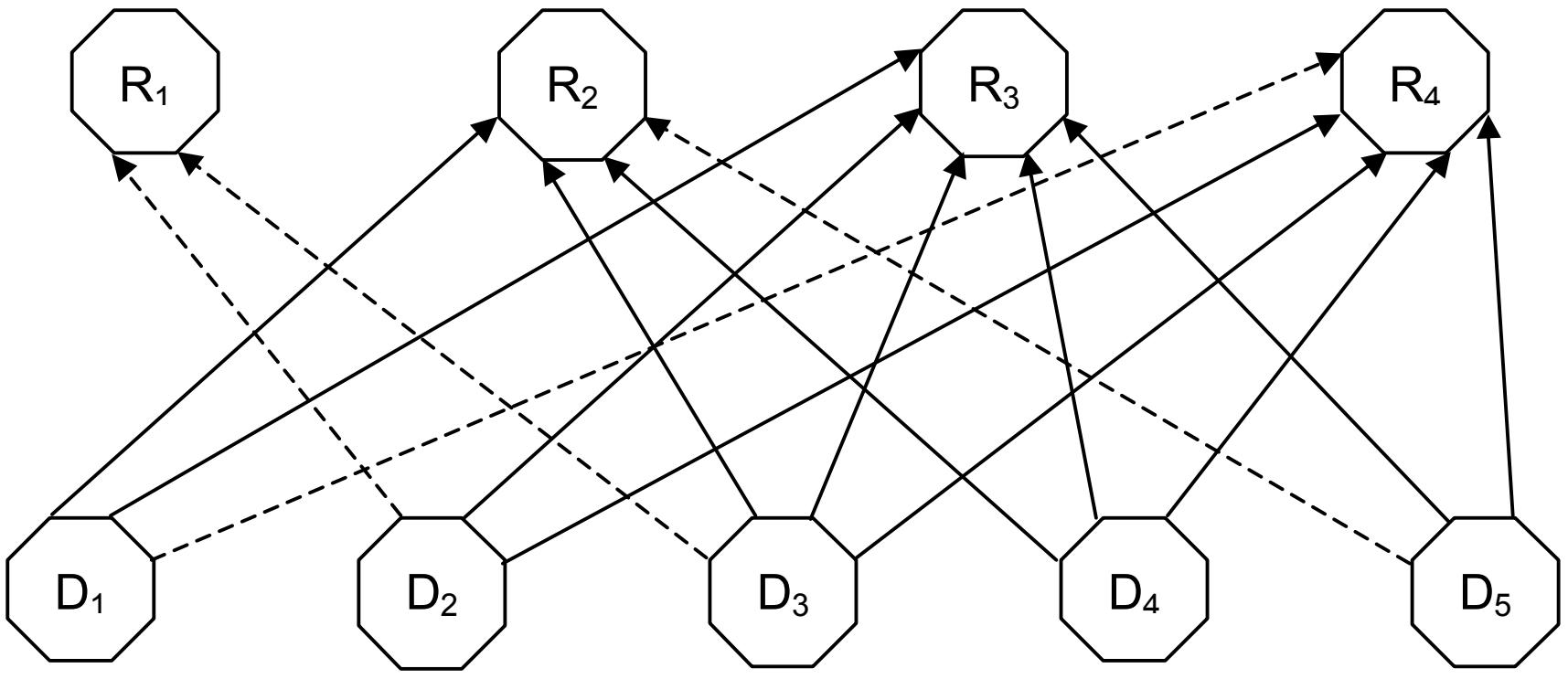

**FIGURE: 2.8.3**

The corresponding neutrosophic relational maps given by N(E).

$$N(E) = \begin{bmatrix} 0 & 1 & 1 & I \\ I & 0 & I & 1 \\ I & 1 & 1 & 1 \\ 0 & 1 & 1 & 1 \\ 0 & I & 1 & 1 \end{bmatrix}$$

Suppose we consider the instantaneous state vector $A_1$ = (0 1 0 0 0) the feeling of guilt, to be in the on state. Now the effect of $A_1$ on the Neutrosophic dynamical system N(E) is as follows:

$$A_1N(E) = (I\ 0\ I\ 1) = B_1$$
$$B_1(N(E))^T = (I, 1+I, 1+I, I+1, I+1)$$

After updating and thresholding we get

$$B_1(N(E))^T \rightarrow (I, 1, 1, 1, 1) = A_2$$
$$A_2N(E) = (I, I+1, I+1, I+1) \rightarrow (I, 1, 1, 1) = B_2$$
$$B_2(N(E))^T = (I+1, I+1, I+1, I+1) \rightarrow (1, 1, 1, 1, 1) = A_3$$
$$A_3N(E) = (I, I+1, 1+I, 1+I) = (I, 1, 1, 1) = B_2.$$

This is a fixed point.

$B_2 (N(E))^T = (I+1, I+1, I+1, I+1, I+1) \rightarrow (1, 1, 1, 1, 1) = A_4 = A_3$

which is a fixed point.

Thus we see the feeling of guilt in the AIDS / HIV patients makes all the nodes of the domain space to be in the on state i.e., feeling of guilt gives them feeling of aloofness, desperate / fear of public, suffering of both mental / physical and a feeling of public disgrace. On the other hand the feeling of guilt has the following effects on the range space. The fear of public getting the disease is an indeterminate one, where as it shows the public have no mind to forgive the AID / HIV patients sin, they still view AIDS/ HIV as a social stigma and they (public) do not have any sympathy for them. Thus we see the resultant vector in the range space is a Neutrosophic vector where as the resultant vector in the domain space is a real vector, further both the hidden patterns happen to give only a fixed point.

Having seen some examples we now proceed on to give the applications / illustrations of the NRMs, combined NRMs define the new notion of linked NRMs in the following section:

## 2.9 Application / Illustration of NRMs, combined NRMs and the introduction of linked NRMs

The notion of FRMs is itself very new and we have just introduced the concept of NRMs, however we will give some illustration / applications of NRMs in the real world problems. Further we will also study and explicitly explain the notion of combined NRMs and its application. Also we will introduce a new notion called the linked NRMs and illustrate them with examples.

### 2.9.1: NRMs in Employee and Employers relationships model

The congenial relation of the employee and employer congenial relation is an important one which reflects the industrial harmony. For example employer expects to achieve consistent production, quality product at the optimum production to earn profit. However there may be profit, no loss and no profit, loss in business or heavy loss depending on various factors such as demand and supply, rent, electricity, raw materials, transportation and consumables, safety, theft, medical aid, employees co-operation and enthusiasm. At the same time some of the expectations of the employees are pay, various allowances, bonus, welfare unanimity such as uniforms, hours of work and job satisfaction etc. Since in this relationship between employee and the employer several concepts like uncertainty, indeterminance play a major role we are justified in using NRMs to analyse them.

Opinions of three experts are taken. We have instructed the experts that the concepts can be inter related, no relation and an indeterminate relation can exist or be given between two conceptual nodes. The data and the opinion is taken using only one industry. Using the opinion we obtain the hidden pattern. The following concepts are taken as the nodes relating to the employee but we can have several more nodes and also several experts opinion. Here in our analysis only 8 notions which pertain to the employee are taken.

| | | |
|---|---|---|
| $D_1$ | – | Pay with allowances and bonus to the employee |
| $D_2$ | – | Only pay to the employee |
| $D_3$ | – | Pay with allowances (or bonus) to the employee |
| $D_4$ | – | Best performance by the employee |
| $D_5$ | – | Average performance by the employee |
| $D_6$ | – | Poor performance by the employee |
| $D_7$ | – | Employee works for more number of hours |
| $D_8$ | – | Employee works for less number of hours. |

$D_1$, $D_2$,…, $D_8$ are taken as elements of the domain space which pertain to the employee nodes.

We have taken only 5 nodes / concepts related to the employer in this study, these are taken as the nodes of the range space.

The nodes / concepts taken as the range space is given below

| | | |
|---|---|---|
| $R_1$ | – | Maximum profit to the employer |
| $R_2$ | – | Only profit to the employer |
| $R_3$ | – | Neither profit nor loss to the employer |
| $R_4$ | – | Loss to the employer |
| $R_5$ | – | Heavy loss to the employer. |

The neutrosophic directed graph as given by the employer is given in the following Figure 2.9.1.:

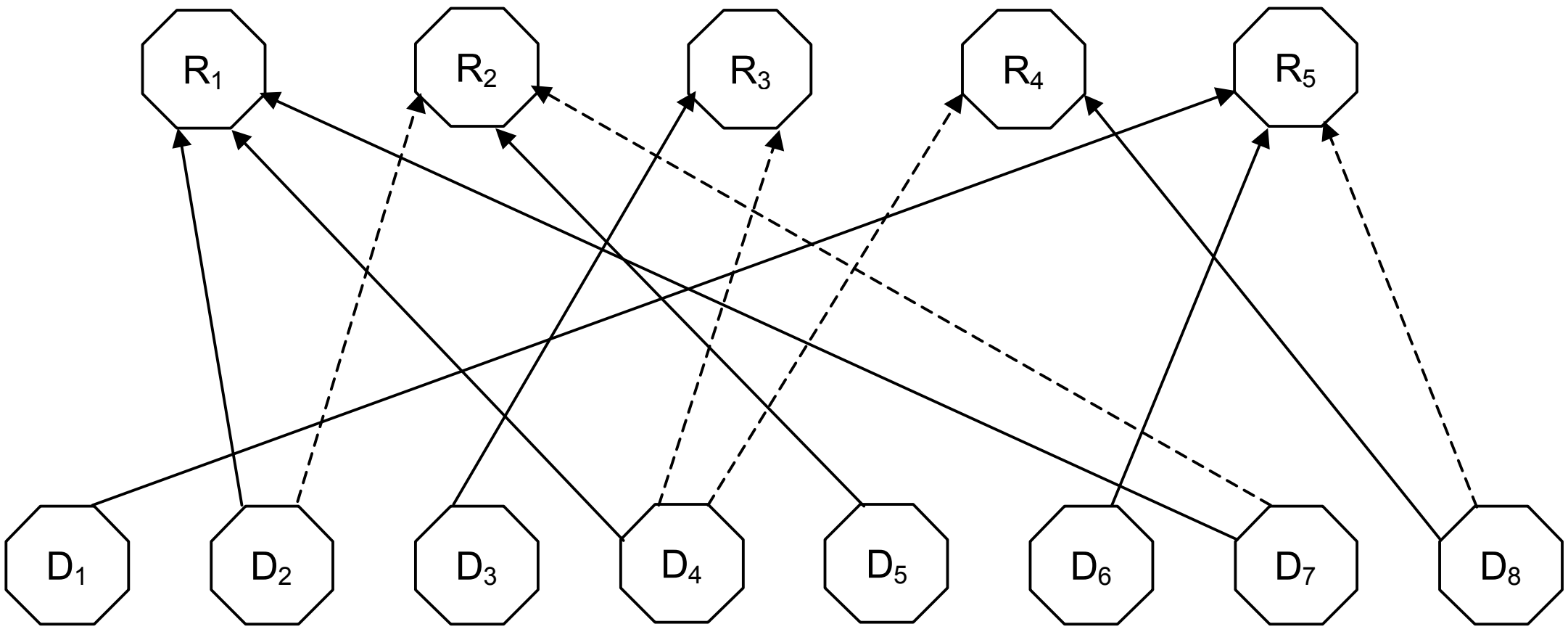

**FIGURE: 2.9.1**

The corresponding or associated neutrosophic relational matrix $N(E_1)$ for this neutrosophic directed graph is given by the following matrix:

$$N(E_1) = \begin{bmatrix} 0 & 0 & 0 & 0 & 1 \\ 1 & I & 0 & 0 & 0 \\ 0 & 0 & 1 & 0 & 0 \\ 1 & I & I & 0 & 0 \\ 0 & 1 & 0 & 0 & 0 \\ 0 & 0 & 0 & 0 & 1 \\ 1 & I & 0 & 0 & 0 \\ 0 & 0 & 0 & 1 & I \end{bmatrix}.$$

Suppose we consider the node $D_1$ to be in the on state and rest of the nodes in the off state that is the employee is paid with allowances and bonus i.e. $G_1 = (1\ 0\ 0\ 0\ 0\ 0\ 0\ 0)$, the effect of $G_1$ on the expert system $N(E_1)$ is

$$\begin{array}{lll} G_1N(E_1) & = & (0\ 0\ 0\ 0\ 1) = H_1 \\ H_1[N(E_1)]^T & = & (1\ 0\ 0\ 0\ 0\ 1\ 0\ I) = G_2 \\ G_2N(E_1) & = & (0\ 0\ 0\ 1\ 1 + I) \to (0\ 0\ 0\ I\ 1) = H_2 \\ H_2(N(E_1))]^T & = & (1\ 0\ 0\ 0\ 0\ 1\ 0\ I + I) \to (1\ 0\ 0\ 0\ 0\ 1\ 0\ I) = G_3 = G_2. \end{array}$$

Thus the hidden pattern of the NRM is a fixed point. If the employee is paid with pay and allowance we see the hidden pattern gives the feeling that "poor performance by the employee" and the concept whether the "employee works for less number of hours" is an indeterminate, suggesting that there exists some vague relationship between these concepts.

The effect of paying the employee with pay and bonus gives the following the loss to the employer is an indeterminate but it has bad effect of heavy loss to the employer.

On the other hand, the result got using FRM is that the company suffers a heavy loss due to the poor performance of the employee. The reader is expected to compare these results. Several other state vectors can be tested in this model.

The union leader of the same company was asked to give his opinion keeping the same nodes for the range space and domain space i.e. as in case of the first expert. The neutrosophic directed graph given by the union leader of the company is given in the Figure 2.9.2.

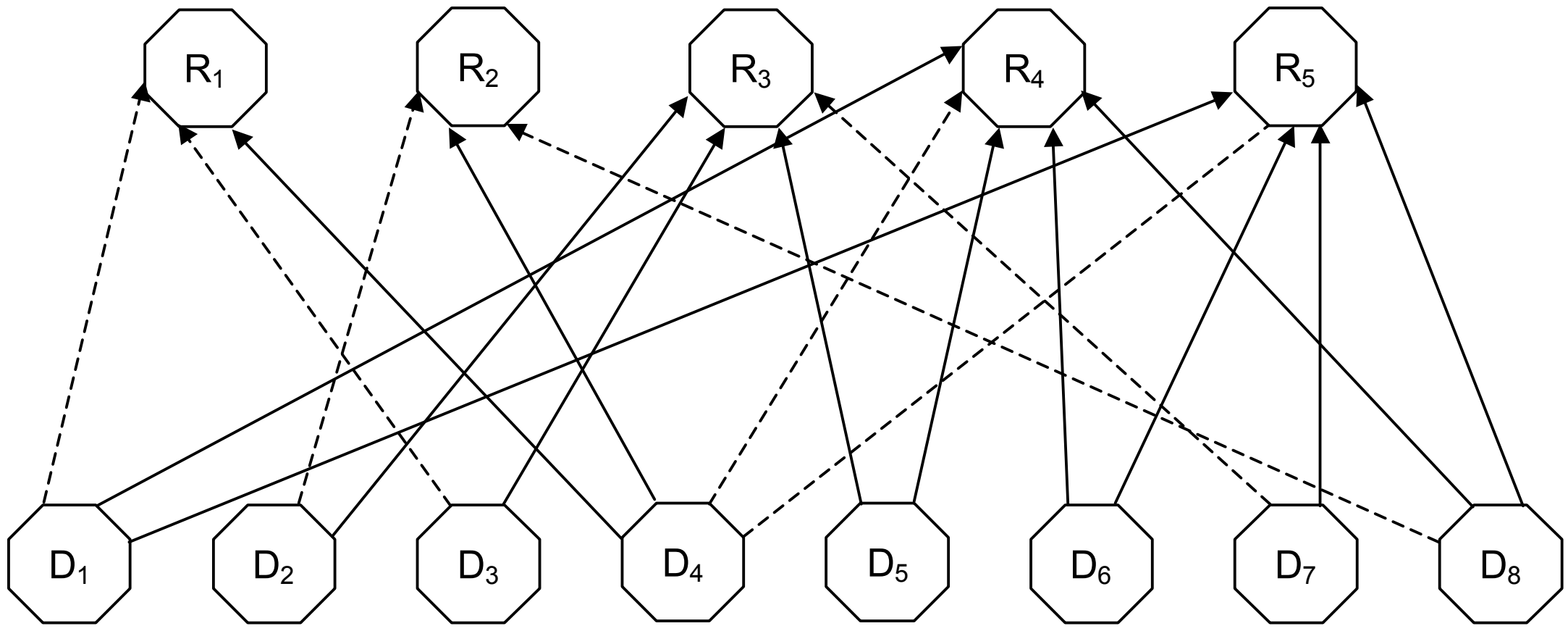

**FIGURE: 2.9.2**

The related neutrosophic connections matrix $N(E_2)$ of the neutrosophic directed graph is given below:

$$N(E_2) = \begin{bmatrix} I & 0 & 0 & 1 & 1 \\ 0 & I & 1 & 0 & 0 \\ I & 0 & 1 & 0 & 0 \\ 1 & 1 & 0 & I & I \\ 0 & 0 & 1 & I & 0 \\ 0 & 0 & 1 & 0 & I \\ 0 & 0 & I & 0 & 1 \\ 0 & I & 0 & 1 & 1 \end{bmatrix}.$$

Let us input the instantaneous vector $L_1$ = (1 0 0 0 0 0 0 0) which indicates the node $D_1$ viz. employee is paid with allowance and bonus is in the on state.

| | | | | | | |
|---|---|---|---|---|---|---|
| $L_1N(E_2)$ | = | (I 0 0 1 1) | = | $N_1$ | | |
| $N_1(N(E_2))]^T$ | = | (1+ I 0 I I I I 1 1) | → | (1, 0 I I I I, 1, 1) | = | $L_2$ |
| $L_2N(E_2)$ | = | (I I I 2 + I 3 + I) | → | (I I I, 1, 1) | = | $N_2$ |
| $N_2(N(E_2))^T$ | = | (I+1, I 2I, 4I 2I, 2I I+1 I + 2) | → | (1, I, I, I, I, I, 1, 1) | = | $L_3$ |
| $L_3N(E_2)$ | = | (2I, 3I, 4I, 1+2I 3+I) | → | (I I I 1, 1) = $N_3$ | = | $N_2$ |
| $N_2[N(E_2)]^T$ | = | (1, I, I, I, I, I, 1,1) | = | $L_4$ | = | $L_3$. |

Thus we see if the employer pays the employee with pay, bonus and allowances then the employee works for more number of hours or less number of hours, which will be interpreted after the notions of range, space and all the remaining concepts/nodes remain indeterminates. In case of range space the maximum profit, only profit or neither profit nor loss remains indeterminate whereas loss to employer and heavy loss to the employer remains in the on state, whereas the maximum profit, only profit or neither profit nor loss to the employer remains as indeterminate.

Now we proceed on to study the same industry using the third expert who is an employee of the same industry. The neutrosophic directed graph given by the third expert is given in Figure 2.9.3.

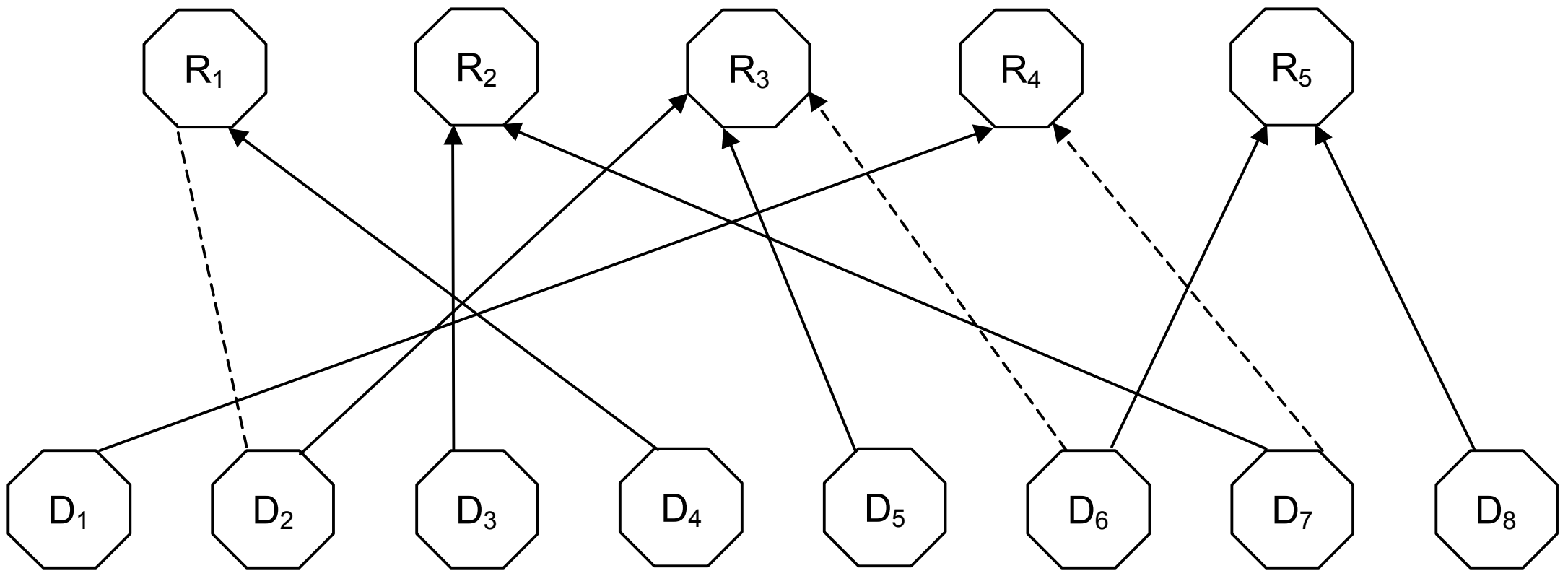

**FIGURE: 2.9.3**

The associated relational neutrosophic matrix representation of the third expert opinion is given by the following:

$$N(E_3) = \begin{bmatrix} 0 & 0 & 0 & 1 & 0 \\ I & 0 & 1 & 0 & 0 \\ 0 & 1 & 0 & 0 & 0 \\ 1 & 0 & 0 & 0 & 0 \\ 0 & 0 & 1 & 0 & 0 \\ 0 & 0 & I & 0 & 1 \\ 0 & 1 & 0 & I & 0 \\ 0 & 0 & 0 & 0 & 1 \end{bmatrix}.$$

Suppose we input the vector $V_1$ = (1 0 0 0 0 0 0 0 ) indicating the on state of the node $D_1$, effect of $V_1$ on $N(E_3)$ is

| | | | | | | | | |
|---|---|---|---|---|---|---|---|---|
| $V_1N(E_3)$ | = | (0 0 0 1 0) | = | $W_1$ | | | | |
| $W_1[N(E_3)]^T$ | = | (1 0 0 0 0 0 I 0) | = | $V_2$ | | | | |
| $V_2N(E_3)$ | = | (0 I 0 1+I 0) | → | (0 I 0 1 0) | | | = | $W_2$ |
| $W_2[N(E_3)]^T$ | = | (1 0 I 0 0 0 2I 0) | → | (1 0 I 0 0 0 I 0) | | | = | $V_3$ |
| $V_3N(E_3)$ | = | (0 2I, 0 1+I 0) | → | (0 I 0 1 0) | = | $W_3$ | = | $W_2$. |

Thus $V_3$ = (1 0 I 0 0 0 I 0) is the fixed point. The $D_1$ node in the on state implies pay with allowances is an indeterminate and employee works for more number of hours is an indeterminate. The fixed point of the range space when the $D_1$ node is in the on state we see the profit to the employer is an indeterminate and loss to the employer becomes the on state. The reader is expected to compare the NRM with FRM discussed in the Section 1.6.1 of this chapter, in order to clearly understand the greater sensitivity of NRMs.

Now using the three experts' opinion we get the combined NRM given by the relational matrix N(E).

Let $N(E) = N(E_1) + N(E_2) + N(E_3)$

$$= \begin{bmatrix} I & 0 & 0 & 2 & 2 \\ 1 & 2I & 2 & 0 & 0 \\ I & 1 & 2 & 0 & 0 \\ 3 & 1 & I & I & I \\ 0 & 1 & 2 & I & 0 \\ 0 & 0 & 1 & 0 & 1 \\ 1 & 1 & I & I & 1 \\ 0 & I & 0 & 2 & 2 \end{bmatrix}.$$

If $T_1$ = (1, 0 0 0 0 0 0 0) has $D_1$ to be in the on state then $T_1$ N(E) = (I 0 0 2 2) → (I 0 0 1 1) = $S_1$

| | | | | | |
|---|---|---|---|---|---|
| $S_1[N(E)]^T$ | = (I+4, I, I, 5 I, I, 1, I+1, 4) | → | (1, I, I, I, I, 1, 1, 1) | = | $T_2$ |
| $T_2N(E)$= | (6I+1, 6I+1, 8I+1, 4+I, 5+I) | → | (1 1 1 1 1) | = | $S_2$ |

$S_2 [N(E)]^T$ = (4+I, I+1, I+3 4+ 3I 3+I, 2, 3 +I, I + 4) → (1, I I I I I I) = $T_3$
$T_3[N(E)]^T$ → (I, I, 1, I 1) = $S_3$
$S_3N(E)$ → (1 1 1 I 1 1 1 1) = $T_4$
$T_4[N(E)]^T$ → (1 1 1 1 1).

Thus the resultant is a limit cycle.

The reader can interpret and compare the combined NRM with the combined FRMs and draw conclusions on the model. The reader can also try to work with another node in the on state and study the NRM and later compare the same in case of FRMs.

Now we will give one more illustrations of NRM.

### 2.9.2: Use of NRM to Study the Women Empowerment Relative To HIV/AIDS

The process of successful mobilization begins when a community begins to identify its concerns and takes up the vital issues. Such a widely based cost effective and beneficial community mobilization will certainly strengthen the entire social fabric. More particularly it will provide a grass-root solution to the epidemic of AIDS /HIV and it will be an immense power in reducing the spread of the disease as well as caring for the affected. The role of community mobilization with respect to the sufferings of AIDS patients is indisputable since in the long run it tends to be a continuous source of motivation and involvement. Although the methods of large-scale community mobilization vary vastly from religious commitment to social service, a better outcome can be expected only when the various groups work together.

When community mobilization is achieved by women's empowerment, the psychosocial impacts of the AIDS epidemic gather momentum. As a first step the motivation of AIDS prevention and AIDS patients care can be linked to get a better socially sustainable situation. The second factor is using women's empowerment as a tool of community mobilization increases the scope of operation to match the scope of the AIDS epidemic. Apart from the increase in the number of social workers and NGOs attending to the AIDS epidemic, the empowered women can also bring about work through organizations that already exist in these communities. It can further change the power balance in typical gender relations thereby reducing the risks of women from getting infected with HIV.

Since the very concept of women's empowerment and community mobilization in the context of AIDS epidemic is an unsupervised data having no specific statistical values to represent these concepts; we felt it is appropriate to use Neutrosophic Relational Maps to study this problem. It is the first time that such a problem is analyzed using a mathematical tool that too the tool of Neutrosophic Relational Maps. The attributes related to women's empowerment and community mobilization are taken as:

$W_1$ - Gender Balance
$W_2$ - Cost effectiveness
$W_3$ - Large-scale operation
$W_4$ - Social service.

$W_1$, $W_2$, $W_3$ and $W_4$ are taken as the nodes / attributes of the domain space. It is pertinent to mention here that certainly one can take more attributes / nodes but for the sake of simplicity of illustration we have taken only four attributes. Further it is important to mention here that we expect the reader to work with more attributes and study the same problem.

$W_1$ – Gender Balance – By this term we mainly mean that women should be given equal social status as in par with men.

$W_2$ – Cost effective. If the government or a private body should adopt or visit a village for measures to prevent the spreading of AIDS/HIV or to make a study about the treatment of AIDS/HIV etc it has to mobilize a lot of money for this purpose. If the community mobilization and women's empowerment are utilized in that village itself it will involve a very less cost and also promote a better social relationship in that village. By this we mean the phrase cost-effective.

$W_3$ – Successful Large-scale operation. Within a mobilized community identification of social concerns and taking up of the vital issues regarding AIDS especially spreading greater awareness is an easy task which is at the same time large scale. This grass-root level operation covers simultaneously large geographic area in a shorter time span which no outside agency can do with the same efficiency.

$W_4$ – Social service to AIDS patients (Giving them courage to fight the disease, wiping the social stigma, etc).

The attributes of the range space which is related to the HIV/AIDS epidemic are

$A_1$ - Care for the AIDS infected persons
$A_2$ - Prevention of spread of HIV/AIDS epidemic
$A_3$ - Creation of Awareness about AIDS
$A_4$ - Medical treatment of AIDS patients
$A_5$ - Social stigma.

The Neutrosophic Directed Graph with these attributes as given by an expert is given in Figure 2.9.4.

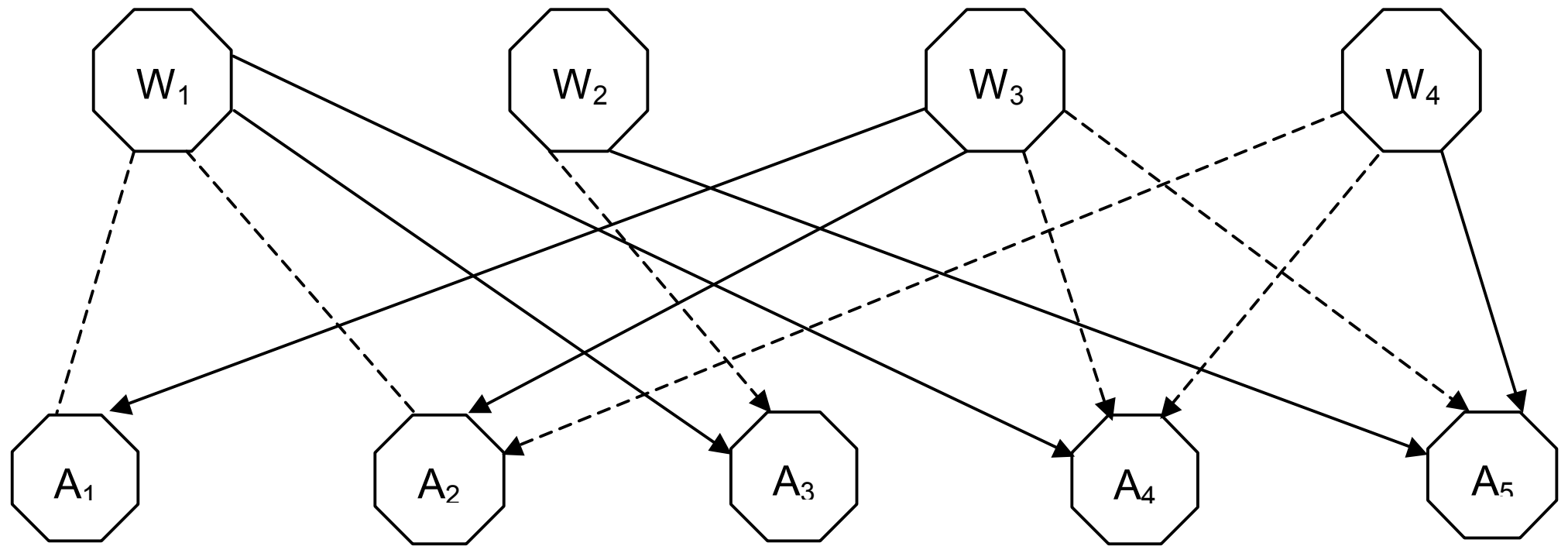

**FIGURE: 2.9.4**

The corresponding neutrosophic relational matrix N(E) is

$$N(E) = \begin{bmatrix} I & I & 1 & 1 & 0 \\ 0 & 0 & I & 0 & 1 \\ 1 & 1 & 0 & I & I \\ 0 & I & 0 & I & 1 \end{bmatrix}.$$

Suppose we take the instantaneous state vector $X_1 = (1\ 0\ 0\ 0)$ i.e. the gender balance node alone to be in the on state. The effect of the vector $X_1$ on the neutrosophic dynamical system is as follows:

| | | | | | | | | |
|---|---|---|---|---|---|---|---|---|
| $X_1N(E)$ | = | (I I 1 1 0) | = | $Y_1$ | | |
| $Y_1[N(E)]^T$ | = | (1+I I 3I I) | → | (1 I I I) | = | $X_2$ |
| $X_2N(E)$ | = | (2I, 3I 1+I 1+I 3I) | → | (I I 1 1 I) | = | ½ |
| $Y_2(N(E))^T$ | = | (2I + 2, 2I, 4I 3I) | → | (1 I I I) | = | $X_3 = X_2$. |

This $X_1$'s hidden pattern is a fixed point.

According to this expert we see that when the attribute / node, gender balance is in the on state its effect on other nodes of the domain space remain indeterminate i.e. cost effectiveness, successful large scale operation and social service relations with gender balance is an indeterminate.

The other fixed point of the range space is $Y_2$ = (I I 1 1 I) i.e. gender balance makes on the attribute / node, creation of awareness about AIDS and medical treatment of AIDS patients, but it remains indeterminate about the care for the AIDS infected patients, prevention of spread of HIV/AIDS patients and the concept of social stigma.

The reader is expected to work with other applications of on state of various other nodes from the range or domain space and draw the related results. Actually we have taken the history of over 500 AIDS/HIV patients and for this study we are deeply indebted to the encouragement and financial support given by the TNSACS (Tamil Nadu State AIDS Control Society).

The reader is expected to derive the conclusions using the model. He can also experiment with an opinion given by a different expert and see how it affects the system. This is left to reader as an exercise to make him familiar with the working techniques of NRM. Now we proceed on to define how NRMs can be used to study the emotional situation of the terminally ill patients and support of their relatives and friends. This application could serve, as a model for any nursing home and how they should take care of the terminally ill patients and how close relatives must be counseled so that the terminally ill patients can lead a peaceful life before they die.

### 2.9.3: NRM to study the Depression of Terminally ill Patients and their Outside Support

In recent medical literature several authors have expressed the views that many terminally ill patients are depressed and suicidal. Depressed patients may be irrational. In most cases doctors, nurses or close relatives miss the presence of depression and so do not realize that these patients need special management. By "outside support" we mean the support from doctors, nurses, close friends and associates and immediate relatives.

Since the very study involves feelings that too very strong feelings like depression we felt that it is essential to model this problem using NRMs. Since the concepts associated with outside support and the attributes involved with patients are disjoint it is better to use NRMs than NCMs.

The attributes related to the terminally ill patients are taken as

| | | |
|---|---|---|
| $P_1$ | - | Fear / anxiety of death |
| $P_2$ | - | Physical suffering |
| $P_3$ | - | Loss of dignity and independence |
| $P_4$ | - | Large doses of medicines [Whose side-effects itself forms a part of disease and hastens death.(For in case of cancer treatment the patient loses hair, blackened complexion etc).] |
| $P_5$ | - | Wants to die - acutely depressed. |
| $P_6$ | - | Prefers to suffer in isolation. |

These six concepts are taken as the nodes of the domain space however one can have several attributes to be associated with them.

The concepts /nodes associated with the outside support is taken as the nodes of the range space which are as follows:

| | | |
|---|---|---|
| $O_1$ | - | Unconcern of doctor / nurses |
| $O_2$ | - | Role of relatives / friends |
| $O_3$ | - | Religion as a solace |
| $O_4$ | - | Entertainment / relaxation |
| $O_5$ | - | No treatment for depression |

Taking the experts opinion the conclusions will be drawn.

Here we give the opinion of one of the experts, which is represented as the neutrosophic graph and is given in Figure 2.9.5.

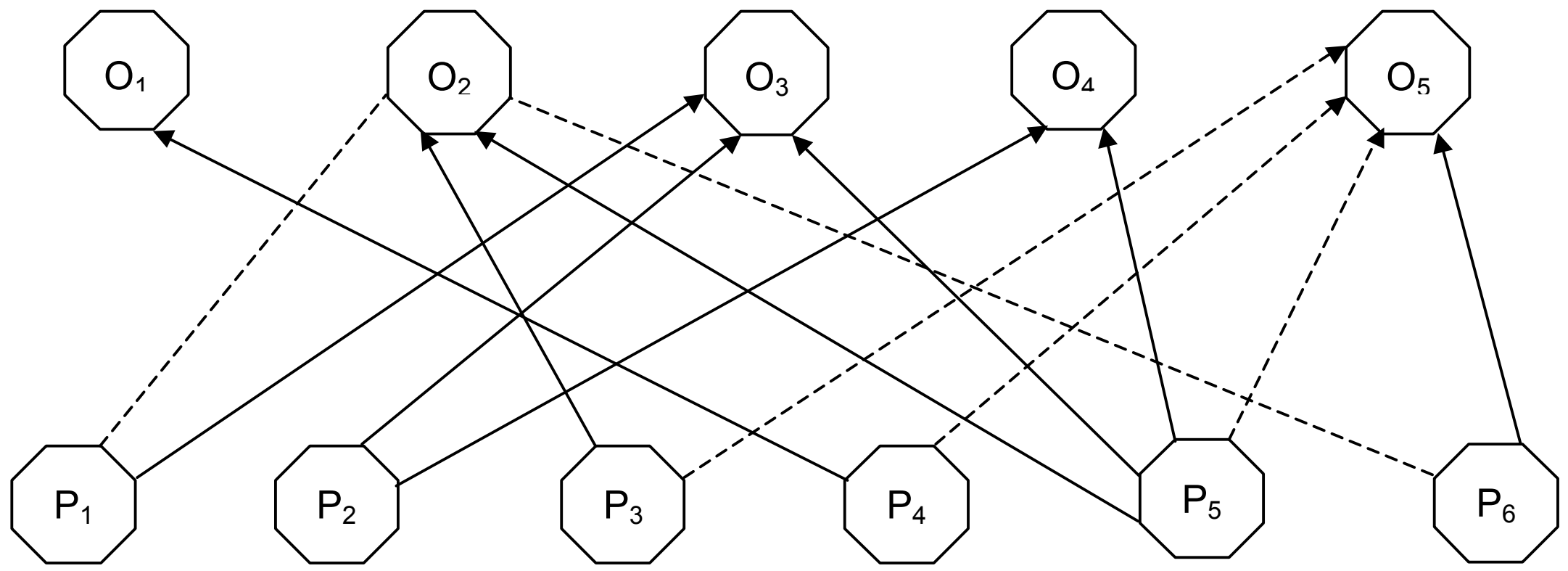

**FIGURE: 2.9.5**

This neutrosophic directed graph is used and the corresponding neutrosophic relational matrix N(E) is obtained.

$$N(E) = \begin{bmatrix} 0 & I & 1 & 0 & 0 \\ 0 & 0 & 1 & 1 & 0 \\ 0 & 1 & 0 & 0 & I \\ 1 & 0 & 0 & 0 & I \\ 0 & 1 & 1 & 1 & I \\ 0 & I & 0 & 0 & 1 \end{bmatrix}.$$

Suppose we take the vector $V_1$ = (0 0 0 0 1 0) the $P_5$ node i.e., wants to die- acutely depressed to be in the on state we study the effect of $V_1$ on the neutrosophic dynamical system N(E).

| | | | | | | |
|---|---|---|---|---|---|---|
| $V_1N(E)$ | = | (0 1 1 1 I) = $W_1$ | | | | |
| $W_1(N(E))^T$ | = | (I+1, 1, 1+I, I 3+I, I+1) | $\to$ | (1 1 1 I 1 1) | = | $V_2$ |
| $V_2N(E)$ | = | (I, 1+I, 3 3 2 I+1) | $\to$ | (I 1 1 1 1) | = | $W_2$ |
| $W_2(N(E))^T$ | = | I+1, 2, 1+I, 2I 2+I I+1) | $\to$ | (1, 1, 1, I I I) | = | $V_3$ |
| $V_3N(E)$ | = | (I, 1+3I, 2+I, 1+I 4I) | $\to$ | (I, 1 1 1 I) | = | $W_3$ |
| $W_3(N(E))^T$ | = | (1+I, 2, 1+I, 2I, 3+I, 2I) | $\to$ | (1 1 1 I 1 I) | = | $V_4$ |
| $V_4N(E)$ | = | (I, 1+I, 3, 2, I) | $\to$ | (I 1 1 1 I) | = | $W_4$ |

Thus when the terminally ill patient wants to die as a result of his/her acute depression we see the patient suffer from fear/ anxiety of death, physical suffering, loss of dignity/ independence become on states and the node which states he is given large doses of medicine remains an indeterminate, but the fact whether he wishes to suffer in isolation remains an indeterminate.

On the other hand the node that the patient wants to die is in on state it reflects the following facts from the position of nodes in the range space, unconcern of doctors remains an indeterminate, the role of relatives/friends is an indeterminate, religion serves as solace is an indeterminate, entertainment or relaxation is an indeterminate but no treatment is given for depression is in the on state.

The reader is expected to study with various nodes from the range and domain space to be in the on state and study the problem using this model and derive conclusions. From our study however we derived several conclusions using many experts' opinion. Some of them are mentioned here. Over 90% of the terminally ill patients longed for solace from religion, 8% were indeterminate about it 2% did not want any role of religion and the majority of them felt that it could reduce depression. The other major conclusion was that the presence of depression was varying in the same patient or to be more precise factually it is still unknown whether it was due to intake of drugs or certain type of foods, which created the depression. Also depression was not continuous but characterized by its erratic nature. The fact remains that almost none of the doctors treat the patients for depression. It is an unfortunate observation. [Note: Our study was carried out among the terminally ill patients in rural Tamil Nadu, India. This may differ for patients in other countries]. Few experts felt that most of the doctors were like machines with no sympathy or additional care and they didn't showing any special concern, for them it was routine. Thus a few of the experts even insisted that the doctors should be trained on how to treat a terminally ill patient, after all, medicine alone is not a full cure for any disease.

Now the reader is expected to study and analyze this problem using combined NRMs.

We have seen three applications of the NRM to the real-world problems. In fact they can be applied to a wide variety of sociological and psychological problems including the study of terminally ill cancer patient who wish for death, in study of female infanticide, suicide by the agriculturists in India, unemployment problems, health hazards due to industries and in particular the problems faced by workers of some chemical industries and so on.

Now we proceed on to define the notion of linked NRM as in the case of FRM and illustrate a few examples and application of them.

The notion of linked neutrosophic relational maps is impossible in the case of neutrosophic cognitive maps. This method is more adaptable in case of data that are not in a position to be directly related but an indirect relation exists between them and only to those data whose related conceptual nodes can be partitioned into disjoint sets. Such study is possible only by using linked NRMs. First we give the definition of pairwise linked NRMs.

**DEFINITION 2.9.1:** *Let us assume that we are analyzing some nodes / concepts which are divided into three disjoint units. Suppose we have 3 spaces say P, Q and R. We say some m nodes in the space P, some n nodes in the space Q and some r nodes in the space R.*

*We can directly find NRMs or directed neutrosophic graphs relating P and Q and Q and R. But we are not in a position to link or get a relation between P and R directly but in fact there exists a hidden link between them which cannot be easily weighted, in such cases we use linked NRM.*

*Thus pairwise linked NRMs are those NRMs connecting three distinct spaces P, Q and R in such a way that using the pair of NRM, we obtain a NRM relating P and R. Thus if $E_1$ is the connection matrix relating P and Q then $E_1$ is a $m \times n$ matrix and $E_2$ is the connection matrix relating Q and R which is a $n \times r$ matrix.*

*Now consider P and R we are not in a position to link P and R directly by any directed graph but the product matrix $E_1 E_2$ gives a neutrosophic connection matrix between P and R. Also $E^T_2 E^T_1$ gives the neutrosophic connection matrix between R and P. When we have such a situation we call the NRMs as the pairwise linked NRMs.*

We will illustrate this definition explicitly by an example.

***Example 2.9.1:*** We just recall the Example 1.7.1 where the study of child labor is carried out using linked FRM. Now instead of FRM we instruct the experts that they need not always state the presence or absence of relation between any two nodes but they can also spell out the indeterminacy of any relation between two nodes, with these additional instruction to the experts, the opinions are taken.

The spaces under study are G – the concepts / attributes associated with the government policies preventing / helping child labor. C – attributes or concepts associated with children working as child laborers and P – attributes associated with public awareness and support of child labor.

G – Concepts associated with government policies:

| | | |
|---|---|---|
| $G_1$ | - | Children do not form vote bank |
| $G_2$ | - | Business men/industrialists who practice child labor are the main source of vote bank and finance |
| $G_3$ | - | Free and compulsory education for children |
| $G_4$ | - | No proper punishment given by Government to those who practice child labor. |

Now we list out some of the attributes / concepts associated with the children working as laborers – C:

| | | |
|---|---|---|
| $C_1$ | - | Abolition of child labor |
| $C_2$ | - | Uneducated parents |
| $C_3$ | - | School dropouts / children who never attended school |
| $C_4$ | - | Social status of child laborers |
| $C_5$ | - | Poverty / sources of living |
| $C_6$ | - | Orphans runways parents are beggars or in prison. |
| $C_7$ | - | Habits like smoking cinema, drinking etc. |

Now we list out the attributes / concepts associated with public awareness or public supporting or exploiting the existence of child labor – P:

| | | |
|---|---|---|
| $P_1$ | - | Cheap and long hours of labor with less pay |
| $P_2$ | - | Children as domestic servants |
| $P_3$ | - | Sympathetic public |
| $P_4$ | - | Motivation by teachers to children to pursue education |
| $P_5$ | - | Perpetuating slavery and caste bias. |

Taking the experts opinion we first give the directed neutrosophic graph relating to child labor and the government policies in Figure 2.9.6.

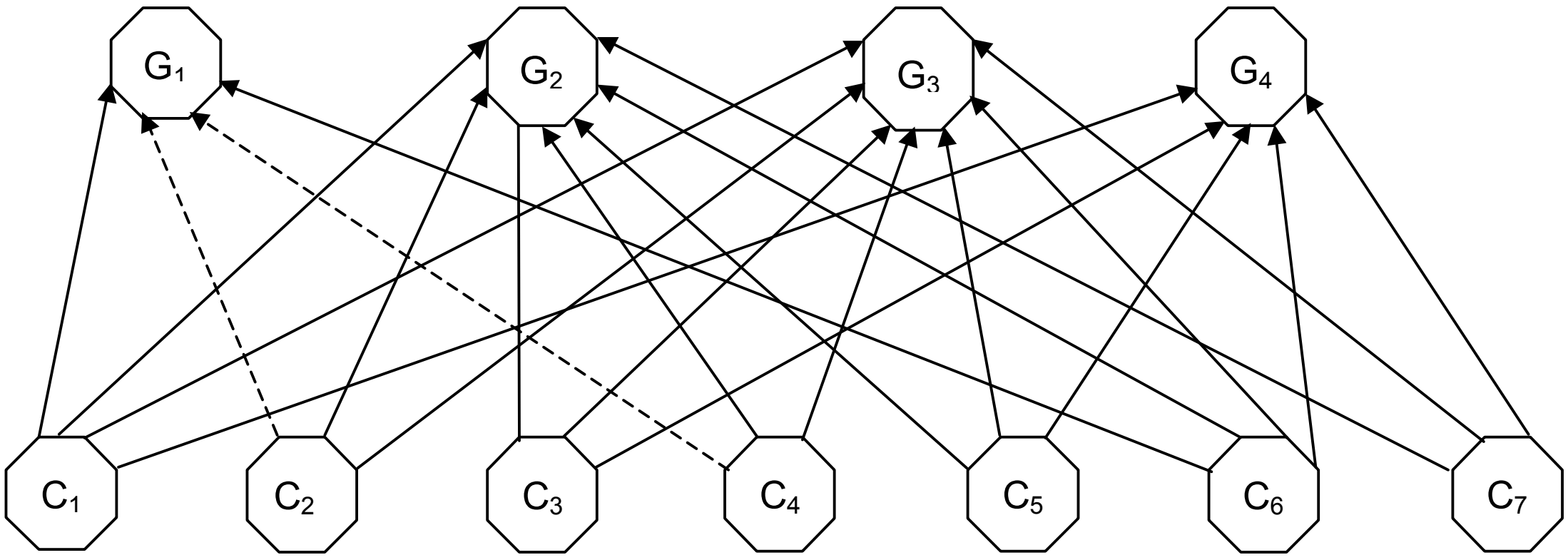

FIGURE:

The related neutrosophic connection matrix $N(E_1)$ is as follows:

$$N(E_1) = \begin{bmatrix} 1 & -1 & 1 & -1 \\ I & 1 & -1 & 0 \\ 0 & 1 & -1 & 1 \\ I & 1 & 1 & 0 \\ 0 & 1 & -1 & 1 \\ 1 & 1 & -1 & 1 \\ 0 & 1 & -1 & 1 \end{bmatrix}.$$

Now we are not interested in seeing the effect of instantaneous state vector on the neutrosophic dynamical system N($E_1$) but we are more interested in the illustration of how the model interconnects two spaces which have no direct relation, the same expert opinion is sought connecting neutrosophically the flourishing of child labor and the role played by the public. The neutrosophic directed graph relating the child labor and the public supporting child labor is given in Figure 2.9.7.

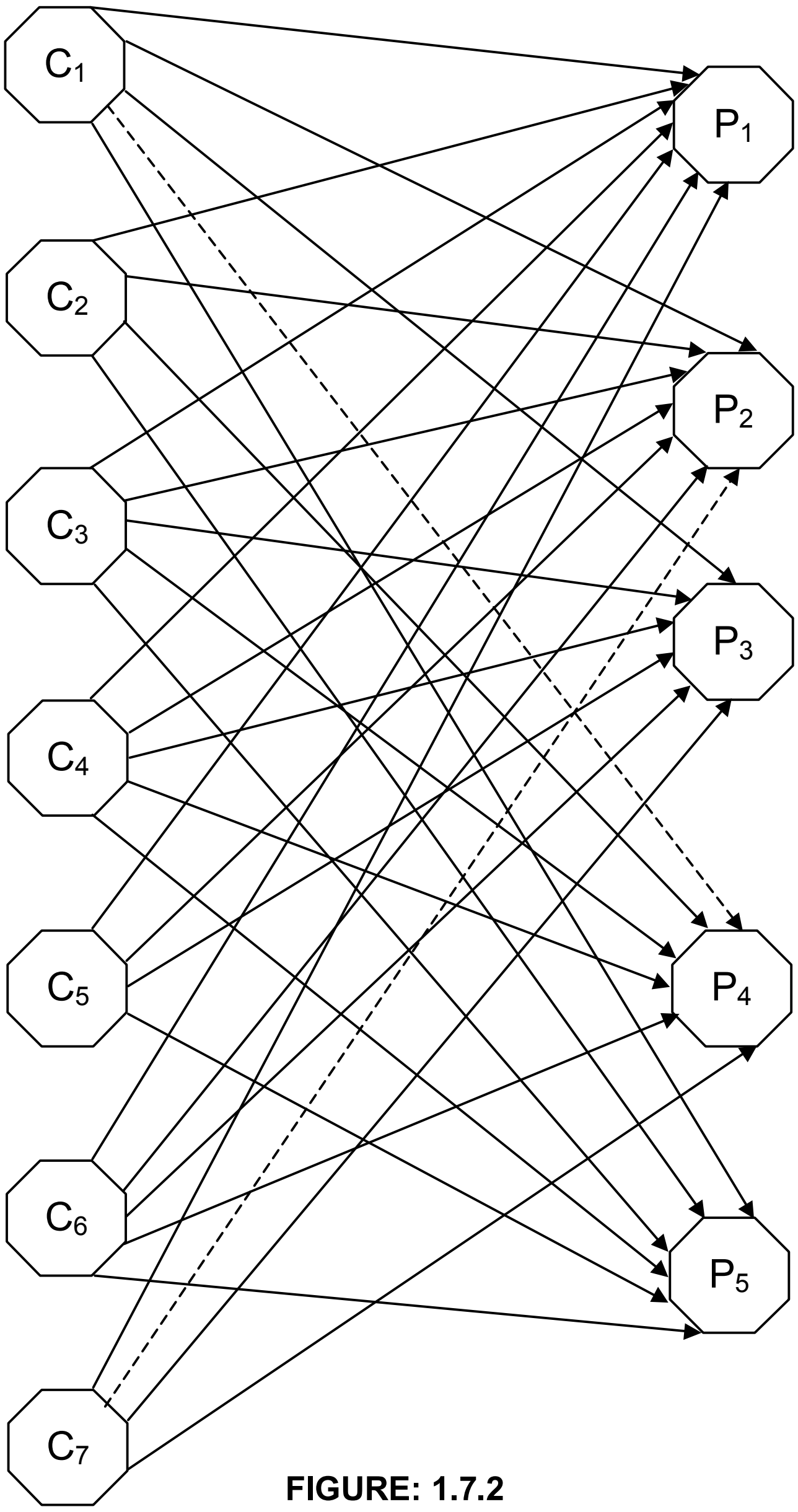

**FIGURE: 1.7.2**

The corresponding neutrosophic connection matrix related to the neutrosophic directed graph is given by $N(E_2)$.

$$N(E_2) = \begin{bmatrix} -1 & -1 & 1 & I & -1 \\ -1 & 1 & 0 & -1 & 1 \\ 1 & 1 & -1 & -1 & 1 \\ -1 & -1 & 1 & 1 & -1 \\ 1 & 1 & 1 & 0 & 1 \\ 1 & 1 & -1 & -1 & 1 \\ 1 & I & -1 & -1 & 0 \end{bmatrix}.$$

Thus $N(E_1)$ is a $7 \times 4$ matrix and $N(E_2)$ is a $7 \times 5$ matrix. Without the aid of the expert we can get the expert opinion relating the public and the government concerning child labor.

The neutrosophic relational matrix $[N(E_1)]^T N(E_2)$ is given as

$$N(E_2)]^T N(E_1) = \begin{bmatrix} 1 & 1 & 1 & -1 & 0 \\ 1 & 1 & -1 & -1 & 1 \\ -1 & -1 & 1 & 1 & -1 \\ 1 & 1 & -1 & -1 & 1 \end{bmatrix}.$$

Thus this maps variedly differs from the linked FRMs (discussed in Section 1.7).

The reader is expected to model real world problems and apply linked NRMs.

**DEFINITION 2.9.2:** *Let P, Q and R be three spaces with nodes related to the same problem. Suppose P and Q are related by an NRM and Q and R are related by an NRM then the indirectly calculated NRM relating P and R got as the product of the neutrosophic connection matrices related with the NRMs. P and Q and Q and R is called the hidden neutrosophic connection matrix and the neutrosophic directed graph drawn using the Hidden neutrosophic matrix is called the Hidden neutrosophic directed graph of the pairwise linked NRMs.*

Now we proceed on to define on similar lines the three linked NRMs four linked NRMs, and in general n-linked NRMs.

**DEFINITION 2.9.3:** *Suppose we are analyzing a data for which the nodes / concepts are divided into four disjoint classes say A, B, C, D where A has n concepts / nodes B has m nodes, C has p concepts / nodes and D has q nodes / concepts.*

*Suppose the data under study is an unsupervised one and the expert is able to relate A and B, B and C and C and D through neutrosophic directed graphs that is the related NRMs and suppose the expert is not in a clear way to interrelate the existing relation between A and D then by using the neutrosophic connection matrixes we can get the neutrosophic connection matrix between A and D using the product of the three known matrix resulting in a* $n \times q$ *neutrosophic matrix which will be known as the*

*hidden neutrosophic connection matrix and the related neutrosophic graph is called the hidden neutrosophic graph.*

*This structure is called as the 3-linked NRM.*

On similar lines we can by using 5 disjoint concepts or nodes of a problem define the 4-linked NRM. Thus if we have the nodes / concepts be divided into n + 1 disjoint classes of nodes then we can define the n-linked NRMs as in case of FRMs.

However we leave the task of constructing examples and applications of linked NRMs to the reader.

## 2.10 Neutrosophic Relational Maps versus Fuzzy Relational Maps

We know in all mathematical analysis of the unsupervised data, not only uncertainty is predominant but also the concept of indeterminacy is in abundance. But in mathematical logic no one used the concept of indeterminacy till the year 1995. Only in the year 1995 Florentine Smarandache introduced and studied the notion of indeterminacy, creating a further generalization of fuzzy logic, which was termed by him as Neutrosophic logic.

So whenever we introduce the notion of indeterminacy in our mathematical analysis we name the structure as Neutrosophic structure. Very recently (1999-2000) the notion of Fuzzy Relational Maps (FRM) was introduced. FRMs were a special particularization of FCMs when the data under study can be divided into disjoint sets. In this section we bring the comparison between FRMs and FCMs and the comparison of NRMs versus FRMs. We show that in the analysis of data, NRMs give a better and a realistic prediction than FRMs.

Here we give a comparison between FRMs and FCMs.

### 2.10.1: Fuzzy Relational Maps versus Fuzzy Cognitive Maps

FRMs are best suited when the data under consideration has its attributes or nodes to be divided into two (or more) disjoint classes.

FRMs cannot be applied when the nodes or attributes under consideration in a data cannot be divided into disjoint classes. Whenever the notion of FRMs are applied the main advantage is it minimizes labor and we need to work only with a smaller rectangular matrix which saves time and energy.

FCMs give a single hidden pattern arising from a fixed point or a limit cycle. On the other hand in the case of FRMs we have two hidden patterns for a given input vector one given by the range space and the other related to the domain space. Thus one is able to study the effect of the instantaneous vector considered in the domain space not only on the domain space but also on the range space.

The directed graph obtained from an FCM may or may not be bigraph. If the directed graph of the FCM is made into a bigraph then it implies and is implied that the FCM can be made into a FRM. All directed graphs of the FRMs are bigraphs.

If the nodes of the FCM is such that the related directed graph can never be made into a bigraph then it autionatically implies that the FCM has nodes which cannot be made into two disjoint classes; so it cannot be made into FRMs.

Thus when FCMs can be converted into FRMs, they always enjoy better and a sensitive resultant apart from being economic and time saving.

### 2.10.2: Neutrosophic Relational Maps versus Neutrosophic Cognitive Maps

The NCMs are the Neutrosophic Cognitive Maps which are FCMs in which atleast one of the directed edge which is an indeterminacy i.e. the directed graph of the NCM is a directed Neutrosophic graph. NRMs are the Neutrosophic Relational Maps i.e. the Neutrosophic directed graph is a Neutrosophic bigraph. The study of NRMs like FRMs is also economic and less time consuming.

The NRMs can be applied only in case the nodes / attributes of the unsupervised data can be divided into two disjoint sets; otherwise the data cannot have NRMs to be applied on it. The same notion in the language of Neutrosophic graphs is as follows:

If the associated Neutrosophic graph of an NCM can be made or is a Neutrosophic bigraph then certainly we can apply NRMs to the given data. If the data i.e. the Neutrosophic graph is never a Neutrosophic bigraph then certainly for that data NRMs can never be applied.

All NRMs can be easily changed into NCMs but NCMs in general cannot be always converted into NRMs. Just like FRMs, NRMs also gives two hidden patterns for a given instantaneous state vector a fixed point or a limit cycle one pertaining to the domain space and one fixed point or a limit cycle pertaining to the range space.

Also NRMs helps the influence of an instantaneous state vector on the space where it is taken and also in the other space.

Thus we can always say when the data is such that it is possible to apply NRMs it certainly yields a better conclusion than the NCMs.

Finally we give the comparison of NRMs and FRMs.

### 2.10.3: Neutrosophic Relation Maps versus Fuzzy Relational Maps

Now we recollect how NRMs are better than FRMs. In the first place in reality one cannot always say that relations between a node in the domain space is related with a node in the range space or not related with node in the range space. It may so happen that the existing relation between two nodes may not always be determinable by an expert. In FRMs there is no scope for such statement or such analysis, we can have a relation or no relation but this will not always be true in case of real world problems that too in case of unsupervised data, the relation can be an indeterminate, in such cases only NRMs are better disposed than FRMs. Thus NRMs play a better role and give a sensitive result than the FRMs.

Fuzzy world is about fuzzy data and fuzzy membership but it has no capacity to deal with indeterminate concepts, only Neutrosophy helps us to treat the notion of indeterminacy as a concept and work with it.

Thus whenever in the resultant data we get the indeterminacy i.e. the symbol I the person who analyze the data can deal with more caution their by getting sensitive results than treating the nonexistence or associating 0 to that co-ordinate.

Thus from our study we have made it very clear that NRMs and NCMs are better tools yielding sensitive and truer results than FRMs and FCMs.

The reader is given the task of constructing real world models and analyzing the data using FCMs (FRMs) and NCMs (NRMs) and compare the results and their merits.

Chapter Three

# SUGGESTED PROBLEMS

In this chapter we have just suggested 59 problems for the reader. Modeling a real world problem by NCM or NRM or linked NRM will itself form substantial research. Some of the problems tackled using FCM or FRM are given to the reader to adopt the method of NCM and NRM. Here also we once again record our research is not on the general fuzzy theory but only pertains to the recollection of the use and application of Fuzzy Cognitive Maps and introduction and analysis of Neutrosophic Cognitive Maps which is a revolutionary method in the field of applied mathematics because it specially studies uncertainty and indeterminacy. As said by Florentine Smarandache [90-94] not only uncertainty is involved but the notion of indeterminacy plays a vital role in such problems. To the best of our knowledge only neutrosophy alone can tackle such a situation Hence NCM, NRM and linked NRMs are defined using the neutrosophic notions.

1. Develop an algorithm for the FRM with m concepts in the domain space and n concepts in the range space to find the limit cycle, fixed point of the dynamical system when m and n are large (Hint:- such program has been developed in case of FCM refer [123] ).

2. Illustrate by a real world problem (i.e. real data used from any source) that FRMs serves better when concepts / nodes of the data can be divided into disjoint classes than using the FCMs.

3. Study the political situation in your hometown using

    i. FRMs

    ii. FCMs

   Justify from your study which technique is better in this case.

4. Construct a model using 3 – linked FRMs for a real world problem.

5. Build an algorithm to find a pairwise linked FRMs.

Hence or otherwise build an algorithm to find the n-linked FRMs and the hidden connection matrix.

6. Using the method of r-linked FRMs in a real world problem, show the hidden matrix is really a very difficult relation to be got as a direct opinion.

7. Consider the child labor problem. Now take 4 spaces (disjoint of course) as G – government policies, P – Public opinion, C – Child labourers, E – Educationalists (Teachers/ Heads of the School/ Educational Institutions, Correspondents / Directors of schools/ Education inspectors etc.). Now all the four spaces form disjoint notions / concepts. Form a 3 – linked FRM and study it, using data or opinions from experts.

8. It is said that the strength of the result depends on the larger number of experts. Prove or disprove it in case of FRMs using a real model.

9. Construct an example of a real model in which FRMs cannot be used and only FCMs can have a role to play.

10. If G(I) is a neutrosophic vector space over K and if G(I) is a strong neutrosophic vector space over $K(I) = \langle K \cup I \rangle$, then prove dimension of neutrosophic vector space > strong dimension of neutrosophic vector space.

11. Can we always say if G(I) is a finite dimensional neutrosophic vector space over K and G(I) is a finite strong dimensional neutrosophic vector space over $K(I) = \langle I \cup K \rangle$ and if $SN_{K(I)}(\dim)$ of G(I) is n then for what values of n we have $N_K$ (dim) of G(I) is $n^2$ ?

12. Define neutrosophic linear operator, neutrosophic linear transformations of neutrosophic vector spaces and strong neutrosophic vector spaces.

13. State and prove Ulam's conjecture in case of neutrosophic graphs.

14. Define and introduce new and interesting properties about neutrosophic graphs.

15. Can we characterize those FCMs which can never be made into FRMs by imposing conditions on the directed graphs of those nodes?

16. Give some analogues or classical results on neutrosophic graphs.

17. Illustrate by a model, which makes use of the application of NCMs, study the data both with FCM and NCM and compare the results.

18. Give an illustration of a model/example, which works best with FCMs.

19. Prove using a model/example that all NCMs are not FCMs but that FCMs are NCMs. That is, prove FCMs form a proper subclass of NCMs.

20. Prove in case of NRMs, that FRMs are a subclass of NRMs.

21. Construct a real model to show that NRMs are more closer to truth than FRMs when indeterminacy is involved in the inter-relations of concepts.

22. Can NRM be used to model a chemical plant? Justify your claim.

23. When the data is unsupervised when can one use NCMs (or NRMs)? Justify your claim.

24. Illustrate by an example that NCMs are best fit than the FCMs in case the data available for study is an unsupervised one with uncertainties and indeterminates.

25. In case of study of the mental status of a patient and the outside support to patients; can you prove NRMs play a better role than NCMs?

26. In modeling the prediction of stocks and shares which model is appropriate FCM or NCM? Justify your claim.

27. Give a model in the real world problem to illustrate working of the linked FRMs.

28. Show by means of a real world illustration that NRMs are better suited than FRMs when the nodes of the data can be divided into disjoint sets.

29. Prove by an example that all NCMs cannot be made into NRMs.

30. Give an illustrative example of a 6-linked FRM.

31. Illustrate by an example a 5-linked NRM.

32. Write an algorithm in Java or C++ to
    i. work with an n-linked NRM
    ii. work with an m-linked FRM.

33. Can a program be constructed to compare a NRM with FRM once the data is provided? Justify your claim.

34. Write a program to implement NCMs. (Hint: Use the program given for FCMs and modify it by including the indeterminacy 'I')

35. What can we apply for the death wish of the terminally ill [say, cancer] patients

    i. NRMs?
    ii. 2-linked NRMs?

    Which will be a better model? Justify your claim with real world data.

36. Compare using a real world data by applying FCMs and then NCMs. Which is better for this data? Justify your claim concretely.

37. Prove using data from any plant the modeling of supervisory systems using NCMs is better than FCMs.

38. Using the simple training set of historical data given by the following table W taken from [17],

| | MP | CP | PROF | INV | FIN | PROD |
|---|---|---|---|---|---|---|
| 1 | 3 | 3 | 3 | 3 | 3 | 3 |
| 2 | 4 | 3.5 | 3.5 | 3 | 4 | 3 |
| 3 | 4 | 4 | 3.5 | 4 | 5 | 3.5 |
| 4 | 3 | 4 | 3.5 | 4 | 4 | 3.5 |
| 5 | 3 | 3.5 | 4 | 4 | 3 | 4 |
| 6 | 2 | 3 | 4 | 4 | 2 | 4 |
| 7 | 3 | 2.5 | 4 | 5 | 1 | 4 |
| 8 | 3 | 3 | 4 | 5 | 2 | 4 |
| 9 | 4 | 3 | 4 | 5 | 3 | 3.5 |
| 10 | 3 | 3.5 | 4 | 5 | 4 | 3.5 |

    Formulate or calculate with experts opinion the weighted neutrosophic matrix N(W).

39. Obtain some nice properties about balanced degree in neutrosophic digraphs. (Recall: A neutrosophic digraph D consists of a finite set V of points and a collection of ordered pairs of distinct points. Any such pair (u, v) is called an arc or directed line or an indeterminacy (i.e. we are not certain whether they can be connected by an arc or by a directed line) and will usually be denoted by uv). In the definition of outdegree of a point υ of a neutrosophic digraph is the number of points adjacent from it and it includes also the dotted lines, which measure the indeterminacy, and in case of indegree it is the number adjacent to it (certainly all indeterminate dotted lines are taken)

40. Define different types of Balance degrees in case of neutrosophic graphs and illustrate them with examples.

41. Develop a RB-NCM syntax in language to describe real world qualitative systems in RB-NCM.(Use methods analogous to RB-FCM)

42. Study a real world model using RB-FCM and RB-NCM and compare them.

43. Obtain some characteristics of the NCR with real valued causality strength.

44. Obtain some interesting properties about NPCR. Construct a real model and compare the model by applying NPCR and FPCR.

45. Define Neutrosophic Markovian Modeling and illustrate it by a real world problem.

46. Take a particular legal issue (case) and use NCM to deal with the result or conclusions and compare it with legal statements/ judgments.

47. Construct a NCModeler (a software tool analogous to FCModeler) in creating metabolic and regulatory network models using NCMs. Hence or otherwise using a real world problem in a metabolic network prove NCMs give results closer to truth than FCMs, when the data under study involves a lot of indeterminacy.

48. Use neutrosophic mechanism for causal relations in a model and compare it with fuzzy mechanisms for causal relations.

49. Prove NCMs without circles will remain static after certain inference steps.

50. Prove NCMs with circles may have very complicated hidden patterns. Hence or otherwise prove circles enable the feedback mechanism in NCM inference process to have important consequences in many real world applications.

51. Prove some results analogous to FCMs in case of NCMs, which can be regularly divided into m-FCMs.

52. Determine an algorithm for an NCM to be divided into basic NCMs.

53. Define key vertices of an NCM and prove NCM with no key vertices has no circles i.e. they are trivial.

54. Give an example of a NCM, which is different from an FCM having several circles with one common vertex.

55. Can you prove by an illustration that NCMs with circles can give better results than FCMs with circles?

56. Construct Rule Based Neutrosophic Cognitive Maps (RB-NCM).

57. For a real world problem adopt Neutrosophic Maps and compare it with Cognitive Maps.

58. Study a real model with uncertainties using the notion of Uncertainty Neuron Neutrosophic Maps.

59. Develop a new concept of web mining procedure analysis analogous to WEMIA using NCMs instead of FCMs.

# BIBLIOGRAPHY

In this book we have collected and utilized a comprehensive list of research papers that have dealt with the field of Fuzzy Cognitive Maps. The bibliographic list given below is exhaustive and lists all those research monographs, dissertations and publications that we could get access to, and which were referred to during the course of our work on the current book. Taking into account the feasibility of referring research material over the Internet, we have also attempted to give the website addresses and links of most of the research publications which are available as PDF (Portable Document Formats).

# INDEX

The indexes given below only refer to the definition of the concepts enlisted. It would be a very cumbersome process to list each and every instance of words, and as a result we have restrained ourselves to only give the related page number where the definition of the concept first occurs.

**O**

**P**

**R**

**S**

**T**

**U**

## Other concepts

# ABOUT THE AUTHORS

**Dr. W. B. Vasantha** is an Associate Professor in the Department of Mathematics, Indian Institute of Technology Madras, Chennai, where she lives with her husband Dr. K. Kandasamy and daughters Meena and Kama. Her current interests include Smarandache algebraic structures, fuzzy theory, coding/ communication theory. In the past decade she has guided seven Ph.D. scholars in the different fields of non-associative algebras, algebraic coding theory, transportation theory, fuzzy groups, and applications of fuzzy theory of the problems faced in chemical industries and cement industries. Currently, six Ph.D. scholars are working under her guidance. She has to her credit 241 research papers of which 200 are individually authored. Apart from this, she and her students have presented around 262 papers in national and international conferences. She teaches both undergraduate and post-graduate students and has guided over 41 M.Sc. and M.Tech. projects. She has worked in collaboration projects with the Indian Space Research Organization and with the Tamil Nadu State AIDS Control Society. She has authored a Book Series, consisting of ten research books on the topic of Smarandache Algebraic Structures which were published by the American Research Press.

She can be contacted at vasantha@iitm.ac.in
You can visit her work on the web at: http://mat.iitm.ac.in/~wbv

---

**Dr. Florentine Smarandache** is an Associate Professor of Mathematics at the University of New Mexico, Gallup Campus, USA. He published over 60 books and 80 papers and notes in mathematics, philosophy, literature, rebus. In mathematics his research papers are in number theory, non-Euclidean geometry, synthetic geometry, algebraic structures, statistics, and multiple-valued logic (fuzzy logic and fuzzy set, neutrosophic logic and neutrosophic set, neutrosophic probability). He contributed with proposed problems and solutions to the Students Mathematical Competitions.
His latest interest in i*nformation fusion* were he works with Dr.Jean Dezert from ONERA (French National Establishment for Aerospace Research in Paris) in creating a new theory of plausible and paradoxical reasoning (DSmT).

In a world of chaotic alignments, traditional logic with its strict boundaries of truth and falsity has not imbued itself with the capability of reflecting the reality. Despite various attempts to reorient logic, there has remained an essential need for an alternative system that could infuse into itself a representation of the real world. Out of this need arose the system of Neutrosophy (the philosophy of neutralities, introduced by FLORENTIN SMARANDACHE), and its connected logic — Neutrosophic Logic, which is a further generalization of the theory of Fuzzy Logic.

In this book we study the concepts of Fuzzy Cognitive Maps (FCMs) and their Neutrosophic analogue, the Neutrosophic Cognitive Maps (NCMs).

Fuzzy Cognitive Maps are fuzzy structures that strongly resemble neural networks, and they have powerful and far-reaching consequences as a mathematical tool for modeling complex systems. Neutrosophic Cognitive Maps are generalizations of FCMs, and their unique feature is the ability to handle indeterminacy in relations between two concepts thereby bringing greater sensitivity into the results.

Some of the varied applications of FCMs and NCMs which has been explained by us, in this book, include: modeling of supervisory systems; design of hybrid models for complex systems; mobile robots and in intimate technology such as office plants; analysis of business performance assessment; formalism debate and legal rules; creating metabolic and regulatory network models; traffic and transportation problems; medical diagnostics; simulation of strategic planning process in intelligent systems; specific language impairment; web-mining inference application; child labor problem; industrial relations: between employer and employee, maximizing production and profit; decision support in intelligent intrusion detection system; hyper-knowledge representation in strategy formation; female infanticide; depression in terminally ill patients and finally, in the theory of community mobilization and women empowerment relative to the AIDS epidemic.

**$23.55**